\documentclass[reqno]{amsart}

\usepackage{hyperref,todonotes,comment,xcolor,xr}
\usepackage[color,matrix,arrow]{xy}
\usepackage{graphics,amssymb}
\usepackage{pagecolor,lipsum}
\usepackage{latexsym}
\usepackage{mathrsfs}
\xyoption{curve}

\definecolor{color}{HTML}{b30086}
\usepackage{hyperref}
\hypersetup{ 
colorlinks,
citecolor=color,
filecolor=color,
linkcolor=color,
urlcolor=color
}

\newtheorem{lemma}{Lemma}[section]
\newtheorem{proposition}[lemma]{Proposition}
\newtheorem{theorem}[lemma]{Theorem}
\newtheorem{corollary}[lemma]{Corollary}
\newtheorem{question}[lemma]{Question}
\newtheorem{conjecture}[lemma]{Conjecture}

\newtheorem*{theoremA}{Theorem}
\newtheorem*{corollaryA}{Corollary}
\newtheorem*{propositionA}{Proposition}

\theoremstyle{definition}
\newtheorem{example}[lemma]{Example}
\newtheorem{definition}[lemma]{Definition}
\newtheorem{remark}[lemma]{Remark}
\newtheorem{warning}[lemma]{Warning}
\newtheorem{notation}[lemma]{Notation}

\newcommand{\mfk}[1]{\mathfrak{#1}}
\newcommand{\mbb}[1]{\mathbb{#1}}
\newcommand{\mcl}[1]{\mathcal{#1}}

\newcommand{\msc}[1]{\mathscr{#1}}
\newcommand{\mbf}[1]{\mathbf{#1}}
\newcommand{\msf}[1]{\mathsf{#1}}
\newcommand{\opn}[1]{\operatorname{#1}}
\newcommand{\ot}{\otimes}
\newcommand{\Vect}{\msc{V}\!ect}
\newcommand{\uHom}{\underline{\operatorname{Hom}}}
\newcommand{\sKan}{\msc{K}\!an}
\newcommand{\sCat}{\msc{C}\!at}

\newcommand{\E}{\mathbf{E}}
\newcommand{\uE}{\underline{\mathbf{E}}}
\newcommand{\Disk}{\mathbf{Disk}}
\newcommand{\uDisk}{\underline{\mathbf{Disk}}}
\newcommand{\SM}{\mathscr{S}\!\!\mathscr{M}}
\newcommand{\dL}{\operatorname{L}}
\newcommand{\pf}{pro\text{-}fin}

\DeclareMathOperator{\Hom}{Hom}
\DeclareMathOperator{\End}{End}
\DeclareMathOperator{\Ext}{Ext}
\DeclareMathOperator{\RHom}{RHom}
\DeclareMathOperator{\Fun}{Fun}
\DeclareMathOperator{\Rep}{Rep}
\DeclareMathOperator{\Fin}{Fin}
\DeclareMathOperator{\Maps}{Maps}
\DeclareMathOperator{\Bord}{Bord}

\renewcommand{\1}{\mathbf{1}}

\renewcommand{\hat}{\widehat}

\definecolor{page_color}{HTML}{000000}
\definecolor{text_color}{HTML}{F0EAD6}

\makeindex

\title[TQFTs from derived quantum group representations]{$3$-dimensional TQFTs from derived categories of quantum group representations}
\date{August 10, 2026}

\author{Cris Negron}
\address{Department of Mathematics, University of Southern California, Los Angeles, CA 90007}
\email{cnegron@usc.edu}

\begin{document}

\maketitle

\begin{abstract}
For any finite modular tensor category $\mfk{A}$, we show that the associated derived $\infty$-category $\msc{D}(\mfk{A})$ supports a topological quantum field theory in dimension $3$.  This TQFT takes the form of a symmetric monoidal functor from an $\infty$-category of surfaces with markings by objects in $\msc{D}(\mfk{A})$, and appropriately decorated bordisms, to the $\infty$-category of dg vector spaces.  We show that the state spaces in this theory are naturally identified with linearized mapping spaces for $\msc{D}(\mfk{A})$.
\par

Though our TQFT requires markings from the derived $\infty$-category, we show that all markings can be removed after taking a homotopy truncation.  The resulting unmarked TQFT produces projective mapping class group actions on cohomology in dimension $2$, and in particular a projective $\opn{SL}_2(\mbb{Z})$-action on Hochschild cohomology.  We expect these mapping class group actions to recover those of Lentner et al.\ \cite{lentneretal23} and Schweigert-Woike \cite{schweigertwoike21}.  In dimension $3$ we obtain power-series valued knot invariants, and power-series valued invariants for modular tensor categories.  We also obtain power-series invariants for closed $3$-manifold, though these can already be calculated at the abelian level.
\par

Our derived field theories are proposed as mathematical formalizations for topological $A$-twists of certain $\mcl{N}=4$ supersymmetric QFTs, in dimension $3$.  Following physical principles, we discuss the possibility of deforming our TQFTs along local systems via an analogous (conjectural) deformation of quantum group representations along the Langlands dual group.
\end{abstract}

\tableofcontents

\section{Introduction}

This paper concerns $3$-dimensional topological field theories obtained from derived modular tensor categories. Our favorite examples to consider are categories of small quantum group representations $(\Rep_qG)_{\opn{small}}$, for a semisimple algebraic group $G$ at an even order parameter $q$ \cite{gainutdinovlentnerohrmann,negron26}, and categories of modules over rigid $C_2$-cofinite vertex operator algebras \cite{mcrae}. New exotic examples of finite modular tensor categories also appear in work of Ostrik and Utiralova \cite{ostrikutiralova}.
\par

In short, we show that the derived category $D_{fin}=D^b(\mfk{A}_{fin})$ of finite length complexes over a finite modular tensor category $\mfk{A}$ supports a $3$-dimensional topological field theory
\[
\mbb{L}_{D_{fin}}:\Bord_{3,2}^{nc}\to D(\opn{Vect})
\]
which is valued in the derived category of linear cochains. Here the superscript $nc$ indicates that our theory is \emph{non-compact}, i.e.\ is only defined on $3$-manifolds with non-trivial outgoing boundary.  For the value on a connected genus $g$ surface, we exhibit a natural identification with the derived Homs
\[
\mbb{L}_{D_{fin}}(\Sigma)\overset{\sim}\to \opn{RHom}_{\mfk{A}}(C^{\ot g},\1),
\]
where $C$ is the so-called canonical coend. Thus, for example, the value on the torus recovers Hochschild cohomology for $\mfk{A}$, and the value on the sphere recovers extensions of the unit. Furthermore, in the case where the input category $\mfk{A}$ is semisimple the theory $\mbb{L}_{D_{fin}}$ collapses to the non-compact part of classical Reshetikhin-Turaev theory, and the full Reshetikhin-Turaev theory is the unique extension of $\mbb{L}_{D_{fin}}$ to the full bordism category.

\begin{remark}
The state spaces for $\mbb{L}_{D_{fin}}$ are allowed to be infinite dimensional precisely because of the omission of closed $3$-manifolds in the non-compact bordism category.
\end{remark}

Our TQFT $\mbb{L}_{D_{fin}}$ is a projection of an analogous theory $L_{\msc{D}_{fin}}$ which is defined at the level of $\infty$-categories, and which includes markings from the derived $\infty$-category. Before continuing with an enumeration of our findings, however, let us explain the appropriate context for our study. We return to a discussion of the technical aspects of the text in Section \ref{sect:derived_intro} below.

\subsection{Mathematics of non-semisimple TQFTs}

Over the past 10 years or so there have been a number of works on ``non-semisimple" TQFTs in mathematics.  Arguably, the first results in this direction concern homology-dependent TQFTs constructed from the unrolled quantum group \cite{blanchetetal16,derenzigeermirand20}. These TQFTs are not only anomalous, but depend on the choice of a generic homology class over each surface and bordism.  
\par

Such homology dependent theories were originally constructed in work of Blanchet, Costantino, Geer, and Patureau-Mirand, specifically for $\mfk{sl}_2$, and were later generalized to arbitrary semisimple $\mfk{g}$ and also some super groups as well \cite{derenzigeermirand20,garnergeeryoung25,geeryoung25}.  Though we won't elaborate further on the nature of these TQFTs, for practical reasons, let us only say that the genericity condition on the homology class translates to a kind of projectivity constraint on representations for the unrolled quantum group, and that projectivity plays a key role in the construction of the corresponding theories.
\par

Following the homology dependent contributions from \cite{blanchetetal16,derenzigeermirand20}, it was shown in work of De Renzi, Gainutdinov, Geer, Patureau-Mirand, and Runkel that any finite (non-semisimple) modular tensor category can be used to produce a $3$-dimensional Reshetikhin-Turaev type theory
\[
L^{abelian}_\mfk{A}:\Bord^{decorated}_{3,2}\to \opn{Vect}.
\]
The theory $L^{abelian}_\mfk{A}$ depends not on the choice of a homology class, but on the presence of projectives in various strategic positions \cite{derenzietal23}.  These projectives are technically only \emph{required} in order to generate closed $3$-manifold invariants, via the use of a modified trace, but should arguably also be present in dimension $2$ \cite{derenzi21,geerpatureumirand}.
\par

It was subsequently shown in work of Costantino, Geer, Ha\"ioun, and Patureau-Mirand \cite{costantinoetalII} that any finite modular tensor category determines, additionally, an invertible $4$-dimensional theory (see also \cite{brochieretal21,beliakovaderenzi24}), and that one can recover the $3$-d theory of De Renzi et al.\ as a boundary to its $4$-d counterpart \cite{haioun}. In $4$-dimensions, projectives are required in dimensions $4$ and $3$, and even $2$ \cite{brownhaioun}, via the implementation of modified traces and the required use of ``admissible" skeins.
\par

What we'll highlight here is the apparent asymmetry in the treatment of objects, so that each non-semisimple TQFT is arguably associated to the pairing of a modular tensor category $\mfk{A}$ with its Karoubian ideal of projectives, rather than to $\mfk{A}$ itself. (This distinction is most pronounced in De Renzi's presentation of the singly extended non-semisimple Reshetikhin-Turaev theory \cite{derenzi21}. See alternatively \cite{lagiotis}.)  This suggests, to us at least, that there may be some \emph{derived} component which is being truncated when we choose to approach non-semisimple topological theories directly at the abelian level.

\begin{remark}
Though we've focused on Reshetikhin-Turaev type TQFTs, there are also non-semisimple Turaev-Viro type theories which similarly leverage projectives in an essential manner \cite{derenzigeermirand18,costantinoetal}.
\end{remark}

\subsection{Derived topological theories from twisted QFTs}

Alongside the mathematical developments discussed above, a parallel story emerged in mathematical physics. In the physical setting, the derived nature of TQFTs is seemingly commonplace.
\par

The clearest presentation of our case of interest, from a strict mathematical perspective, appears in work of Creutzig, Dimofte, Garner, and Geer \cite{creutzigetal24}.  However, the main physical principles used in \cite{creutzigetal24} are adapted from a general theory of twisting for supersymmetric QFTs which was developed in work of Witten from the early 90's \cite{witten91}. For any $\mcl{N}=4$ supersymmetric theory--whatever such a thing might be--we have two distinguished topological twists, an $A$-twist and a $B$-twist \cite{kapustinwitten06,gaiotto19}.  Both twists are interesting to consider, and are presumably connected by a type of Langlands duality.  However, in the present text we are interested in the {\it $A$-side} of this story.

Having $A$-twisted our favorite $\mcl{N}=4$ supersymmetric theory, the resulting topological theory is \emph{very} exotic from the perspectives of, say, Reshetikhin-Turaev or the cobordism hypothesis. It is not only derived in nature, but the state spaces appear as \emph{sheaves} on a stack of local systems. See for example \cite{marcus95,gaiotto18,gaiotto19}.
\par

In Creutzig et al.\ the authors posit that the $A$-twist for a particular supersymmetric $\mcl{N}=4$ theory associated to $\opn{SL}_n$, in $3$-dimensions, recovers a derived variant of the unrolled non-semisimple TQFT produced by Blanchet, Costantino, Geer, Patureau-Mirand, and De Renzi \cite{blanchetetal16,derenzigeermirand20}.  Here the presence of local systems has to do with the unrolled-ness of the quantum group, and after taking fibers at trivial local systems we reduce to a derived variant of the non-semisimple Reshetikhin-Turaev theory from De Renzi, Gainutdinov, Geer, Patureau-Mirand, and Runkel \cite{derenzietal23}.
\par

Thus, from \cite{creutzigetal24} and the surrounding physics literature, we anticipate the following: (a) There are derived topological field theories which are, at the very least, associated to derived categories of small quantum group representations. (b) These derived theories should recover pre-existing mathematical theories after a homological truncation. (c) They should also deform naturally along local systems to produce topological field theories which ``couple" to flat connections, again at least in the quantum group setting.
\par

In the present study we realize points (a) and (b), or more precisely \emph{versions} of points (a) and (b), via explicit mathematical construction. Conjectural approaches to point (c) are also discussed in Section \ref{sect:conjectures}.
\par

Let us now turn our attention to the contents of this text.

\subsection{TQFTs from derived $\infty$-categories}
\label{sect:derived_intro}

Let $\mfk{A}$ be a finite ribbon tensor category which is defined over an algebraically closed field $k$, and let $\msc{D}_{fin}=\msc{D}^b(\mfk{A}_{fin})$ be its associated derived $\infty$-category of bounded, finite length complexes.

\begin{remark}
For us, all tensor categories are presentable (cocomplete), compactly generated, and abelian. Taking the compacts, or rigid objects $\mfk{A}\rightsquigarrow \mfk{A}_{fin}$ recovers a tensor category in the usual sense of, say, \cite{bakalovkirillov01} or \cite{egno15}.  See Section \ref{sect:tensor_cats}.
\end{remark}

From such $\mfk{A}$ and $\msc{D}_{fin}$ we produce corresponding symmetric monoidal $\infty$-categories of labeled, anomalous, non-compact $3$-dimensional bordisms
\[
\Bord^{nc}_{\mfk{A}}=\Bord_{3,2}^{nc}(\mfk{A})\ \ \text{and}\ \ \Bord^{nc}_{\msc{D}_{fin}}=\Bord^{nc}_{3,2}(\msc{D}_{fin}).
\]
In the case of $\mfk{A}$, for example, the objects in this bordism category are surfaces $\Sigma_x$ with markings by objects $x_1,\dots,x_n$ in $\mfk{A}$. Morphisms $M_{\alpha}:\Sigma_x\to \Sigma'_y$ are bordism between the underlying surfaces $M:\Sigma\to \Sigma'$ which are equipped with, loosely speaking, embedded cylinders which travel from the markings on $\Sigma$ to the markings on $\Sigma'$. The terminal ends of these cylinders, on the outgoing boundary, are furthermore marked by morphisms $\alpha_j:m_j(x)\to y_j$ between specified products of the incoming labels to the outgoing labels (see Sections \ref{sect:skein_bordisms} and \ref{sect:des_bord_a}).
\par

To say things simply, the surfaces are decorated by objects in $\mfk{A}$ and bordisms are decorated by morphisms in $\mfk{A}$.  The $\infty$-ness comes in as we allow the embedded cylinders in bordisms to \emph{move} within the ambient space $M$ under homotopy. This movement institutes certain universal skein relations into the bordism category which acknowledge the ribbon structure on $\mfk{A}$. In the case of $\msc{D}_{fin}$, one has homotopical contributions both from the geometry of the bordisms and from the higher structure on $\msc{D}_{fin}$ itself. 

More generally, for any ``sufficiently ribbon-like monoidal $\infty$-category" $\msc{E}$ we construct a bordism $\infty$-category $\Bord_{\msc{E}}^{nc}$ of surfaces with markings by objects in $\msc{E}$, and $3$-dimensional cobordisms with markings by morphisms in $\msc{E}$ (Proposition \ref{prop:bord_e}). The bordism categories $\Bord_{\msc{E}}^{nc}$ are non-compact in two regards: In top dimension we disallow $3$-manifolds without boundary, and in dimension $2$ we require that each component of a surface is equipped with at least one marking from $\msc{E}$ (cf.\ \cite{derenzi21} and Section \ref{sect:mark_scheme}).

\begin{theoremA}[\ref{thm:der_lrt}, \ref{prop:der_states}]
For any finite modular tensor category $\mfk{A}$, with finite derived $\infty$-category $\msc{D}_{fin}=\msc{D}^b(\mfk{A}_{fin})$, there is a symmetric monoidal functor
\begin{equation}\label{eq:218}
L_{\msc{D}_{fin}}:\Bord_{\msc{D}_{fin}}^{nc}\to \Vect
\end{equation}
to the $\infty$-category of dg vector spaces.  Furthermore, for any connected genus $g$ surface $\Sigma$, we have a natural identification of the state spaces
\begin{equation}\label{eq:216}
L_{\msc{D}_{fin}}(\Sigma_x)\overset{\sim}\to \underline{\Maps}_{\msc{D}_{fin}}\big(C^{\ot g},m_I(x)\big).
\end{equation}
\end{theoremA}

To explain the notation in the expression \eqref{eq:216}, the object $C$ in $\mfk{A}_{fin}\subseteq \msc{D}_{fin}$ is the coadjoint object, which we define generally as the value $C=mR(\1)$ for $m:\mfk{A}\ot_k \mfk{A}\to \mfk{A}$ the product and $R:\mfk{A}\to \mfk{A}\ot_k \mfk{A}$ its right adjoint. This object is more commonly referred to as the canonical coend. Additionally, for $\msc{D}=\msc{D}(\mfk{A})$ the full unbounded derived $\infty$-category, the bifunctor
\[
\underline{\Maps}_{\msc{D}}:\msc{D}^{op}\times \msc{D}\to \Vect
\]
is the inner-Hom functor for the natural action of $\Vect$ on $\msc{D}$.  These are, more colloquially, the derived Homs for $\msc{D}$ as a linear $\infty$-category, and we recover the $\Ext$ groups via cohomology
\[
H^{\ast}\underline{\Maps}_{\msc{D}}(u,v)=\Ext^{\ast}_{\mfk{A}}(u,v).
\]
By $\underline{\Maps}_{\msc{D}_{fin}}$ we simply mean the restriction of this bifunctor the full subcategory of finite complexes.  Taking finally $I$ to be the underlying marking set for our surface $\Sigma_x$, $m_I:\msc{D}_{fin}^I\to \msc{D}_{fin}$ is the iterated product functor.

\begin{remark}
A version of Theorem \ref{thm:der_lrt} can also be proved at the $1$-categorical level. We emphasize the homotopical setting due to our interest in an expanded analysis of TQFTs for modular $\infty$-categories, in both dimensions $3$ and $4$. See Theorem \ref{thm:dis_der}, and Section \ref{sect:expansion} for further elaboration.
\end{remark}

Though our non-compactness assumption requires the appearance of non-trivial markings on surfaces in $\Bord_{\msc{D}_{fin}}^{nc}$, we show that all markings can be removed after a homotopy truncation.

\begin{theoremA}[\ref{thm:unmarked_lrt}]
Consider a finite modular tensor category $\mfk{A}$.  After adopting unit labels and applying a homotopy truncation, the theory \eqref{eq:218} specifies a unique symmetric monoidal functor $\mbb{L}_{D_{fin}}:\Bord_{3,2}^{nc}\to D(\opn{Vect})$ from the unmarked $3$-dimensional bordism category which completes a diagram
\begin{equation}\label{eq:246}
\xymatrix{
\opn{h}\Bord_{\opn{Units}(\msc{D}_{fin})}^{nc}\ar[rr]^{\opn{h}L_{\msc{D}_{fin}}|_{\opn{Units}}}\ar[dr]_{forget} & & D(\opn{Vect})\\
	& \Bord_{3,2}^{nc}\ar[ur]_{\mbb{L}_{D_{fin}}} & .
}
\end{equation}
The state spaces in the unmarked TQFT, for any connected genus $g$ surfaces $\Sigma$, are identified with the derived Homs
\[
\mbb{L}_{D_{fin}}(\Sigma)\overset{\sim}\to \RHom_{\mfk{A}}(C^{\ot g},\1).
\]
\end{theoremA}

From our description of the state spaces, we understand that the unmarked theory $\mbb{L}_{D_{fin}}$ is valued specifically in the category of coconnnective cochains $D^{\geq 0}(\opn{Vect})$.  Similarly, the restriction of the theory $L_{\msc{D}_{fin}}$ along the inclusion $\Bord_{\mfk{A}}^{nc}\to \Bord_{\msc{D}_{fin}}^{nc}$ produces a theory which takes values in the $\infty$-category $\Vect^{\geq 0}$.
\par

Since taking cohomology in the coconnective setting $H^0:\Vect^{\geq 0}\to \opn{Vect}$ is symmetric monoidal, we can reasonably speak of the cohomologies
\[
H^0\mbb{L}_{D_{fin}}\ \ \text{and}\ \ H^0L_{\msc{D}_{fin}}|_{\Bord_{\mfk{A}}^{nc}}
\]
as topological field theories. These cohomologies are shown to recover the abelian non-semisimple TQFTs from De Renzi et al.\ \cite{derenzietal23} which we discussed above.

\begin{propositionA}[\ref{prop:heart_lrt}]
Let $\mfk{A}$ be a finite modular tensor category.  The $0$-th cohomology of the theory \eqref{eq:218} is naturally identified with the marked, non-compact part of the $\opn{Vect}$-valued TQFT from De Renzi-Gainutdinov-Geer-Patureau-Mirand-Runkel
\[
H^0L_{\msc{D}_{fin}}|_{\Bord_{\mfk{A}}^{nc}}\overset{\sim}\to L_\mfk{A}^{abelian}|_{\Bord_{\mfk{A}}^{nc}}.
\]
Furthermore, the $0$-th cohomology of the unmarked theory $\mbb{L}_{D_{fin}}$ from \eqref{eq:246} recovers the unmarked, non-compact part of $L_\mfk{A}^{abelian}$,
\[
H^0\mbb{L}_{D_{fin}}\overset{\sim}\to L_{\mfk{A}}^{abelian}|_{\Bord_{3,2}^{nc}}.
\] 
\end{propositionA}

\begin{remark}
The TQFT from \cite{derenzietal23} is arguably already non-compact, as it only allows for closed $3$-manifold which come equipped with certain embedded algebraic decorations.  In particular, its unmarked part is explicitly non-compact.
\end{remark}

\begin{remark}
The ``$L$" in $L_{\msc{D}_{fin}}$ is in reference to Lyubashenko, and we refer to $L_{\msc{D}_{fin}}$ as the derived Lyubashenko-Reshetikhin-Turaev (LRT) theory.
\end{remark}

\begin{remark}
We conjecture that the unmarked theory from Theorem \ref{thm:unmarked_lrt} lifts to a TQFT at the $\infty$-categorical level $\mbb{L}_{\msc{D}_{fin}}:\Bord^{nc}_{3,2}\to \Vect$, and furthermore that the marked and unmarked theories can be integrated into a single object. The removal of markings at the $\infty$-categorical level is a highly substantive point from our perspective, and we hope to iterate on these constructions in future works. (See Sections \ref{sect:skein_bordisms}, \ref{sect:expansion}, \ref{sect:conjectures}.)
\end{remark}

As a final mathematical point for the introduction, let us highlight some of the outputs of our work in dimensions $2$ and $3$.

\subsection{Phenomena in $2$-dimensions}

To gain some kind of understanding of, say, the unmarked TQFT $\mbb{L}_{D_{fin}}$, one can consider its behaviors one surfaces and (punctured) $3$-manifolds.  In dimension $2$ we obtain a projective action of the mapping class group $\opn{MCG}(\Sigma)$ on each state space $\mbb{L}_{D_{fin}}(\Sigma)$.
\par

To elaborate, we have the subgroup $\opn{MCG}(\Sigma)$ of the non-anomalous autobordisms of $\Sigma$ which is generated by the mapping cylinders of oriented self-diffeomorphisms $\opn{Diff}(\Sigma)^{+}\to \opn{MCG}(\Sigma)\subseteq \opn{Aut}_{\text{non-anom}}(\Sigma)$.  Pulling back along the forgetful functor from the anomalous bordism category, we obtain a central extension
\[
1\to \mbb{Z}\to \tilde{\opn{MCG}}(\Sigma)\to \opn{MCG}(\Sigma)\to 1,\ \ \tilde{\opn{MCG}}(\Sigma)\subseteq \opn{Aut}_{\Bord^{nc}_{3,2}}(\Sigma).
\]
This central extension acts on $\Sigma$ in $\Bord_{3,2}^{nc}$, and hence acts on the state space $L(\Sigma)$ for any field theory $L:\Bord_{3,2}^{nc}\to S$.
\par

In our case, we recall that the state spaces for $\mbb{L}_{D_{fin}}$ are identified with cohomology for the input category $\mfk{A}$.  So we obtain the following.

\begin{corollaryA}[\ref{cor:mcg}]
Let $\mfk{A}$ be a finite modular tensor category.  For any connected genus $g$ surface $\Sigma$, the derived maps $\RHom_{\mfk{A}}(C^{\ot g},\1)$ carry the natural structure of a dg representation for the projective mapping class group $\tilde{\opn{MCG}}(\Sigma)$.
\end{corollaryA}

In genus $0$ nothing interesting happens, as the mapping class group is trivial.  In genus $1$, the mapping class group of the torus is $\opn{SL}_2(\mbb{Z})$, and there is a canonical identification of $\RHom_{\mfk{A}}(C,\1)$ with the derived lift of Hochschild cohomology $\RHom_{\End_{\opn{Vect}}(\mfk{A})}(id_{\mfk{A}},id_{\mfk{A}})$.  Hence we obtain, after taking cohomology, a projective action of $\opn{SL}_2(\mbb{Z})$ on Hochschild cohomology $HH^{\ast}(\mfk{A})$.
\par

These kinds of projective mapping class group actions on cohomology were already observed in works of Lentner, Mierach, Schweigert, and Sommerhäuser \cite{lentneretal23}, and of Schweigert and Woike \cite{schweigertwoike21}.  So, from the perspective of those texts, we are saying that the mapping class group actions from \cite{lentneretal23,schweigertwoike21} arise naturally as $2$-dimensional projection of a full $3$-dimensional derived field theory for finite non-semisimple modular tensor categories.

\begin{remark}
We have not checked explicitly that our mapping class group representations reproduce those of \cite{lentneretal23,schweigertwoike21}, though we certainly expect this to be the case.
\end{remark}

\subsection{Phenomena in $3$-dimensions}

We continue our investigation of the theory $\mbb{L}_{D_{fin}}$ in dimension $3$.  For a connected closed $3$-manifold $M$, we do not consider $M$ directly, but fix a surface $\Sigma$ and excise the interiors of two chosen embeddings
\begin{equation}\label{eq:arrangement}
\xymatrixrowsep{3mm}
\xymatrix{
\Sigma\ar[dr] & & \bar{\Sigma}\ar[dl]\\
	& M
}
\end{equation}
to obtain a bordism $\check{M}:\Sigma\to \Sigma$.  (Here we assume that the embeddings from $\Sigma$ extend to some filling $\mbf{\Sigma}$ of the surface.)  The resulting bordism is uniquely determined by $M$ when $\Sigma=S^2$ for example, though this is not the case in higher genus.
\par

The bordism $\check{M}$ obtained from such an arrangement \eqref{eq:arrangement} determines an endomorphism $\mbb{L}_{D_{fin}}(\check{M}):\mbb{L}_{D_{fin}}(\Sigma)\to \mbb{L}_{D_{fin}}(\Sigma)$ and, to wash out some cochain-level ambiguities, we take cohomology to get an endomorphism
\[
H^{\ast}\mbb{L}_{D_{fin}}(\check{M}):H^{\ast}\mbb{L}_{D_{fin}}(\Sigma)\to H^{\ast}\mbb{L}_{D_{fin}}(\Sigma).
\]
To wash out further ambiguities in the definition of $\mbb{L}_{D_{fin}}$ itself, we also take a determinant to obtain a power-series valued invariant
\begin{equation}\label{ex:342}
\opn{Inv}(\mfk{A}|\check{M}):=\sum_{n\geq 0}\det(H^n\mbb{L}_{D_{fin}}(\check{M})-X\cdot id_{H^n})\cdot t^n\ \in\ \mcl{O}(\mbb{A}^1_k)[\![t]\!],
\end{equation}
where $X$ is a formal parameter and we've adopted the notation $\mcl{O}(\mbb{A}^1_k)=k[X]$.
\par

The invariant $\opn{Inv}(\mfk{A}|\check{M})$ simply records the eigenvalues of the endomorphism $H^n\mbb{L}_{D_{fin}}(\check{M})$ in each degree.  One can use this construction to obtain $3$-manifold invariants and also invariants of framed knots.

\begin{propositionA}[{\ref{lem:3d_0}, \ref{prop:inv_knots}}]
Let $\mfk{A}$ be a modular tensor category.  In genus $0$, i.e.\ when $\Sigma=S^2$, the above construction produces power-series valued invariants of $3$-manifolds $\opn{Inv}(\mfk{A}|M):=\opn{Inv}(\mfk{A}|\check{M})$.
\par

In genus $1$, taking $M=S^3$, the above construction produces power-series valued invariants of framed knots $\nu:S^1\to \mbb{R}^3\subseteq S^3$, $\opn{Inv}(\mfk{A}|\nu):=\opn{Inv}(\mfk{A}|\check{S}^3)$.
\end{propositionA}

Though the construction $\opn{Inv}(\mfk{A}|M)$ for $3$-manifolds is clearly ``derived", and also a natural output for any $\Vect$-valued TQFT, we show that these invariants are actually determined by their constant term in $t$, $\opn{Inv}(\mfk{A}|M)=c_0(X)+O(t)$.  The constant term is furthermore a scalar which is already obtainable from the corresponding abelian theory $L^{abelian}_{\mfk{A}}$ from \cite{derenzietal23}.  So we do not truly access any derived information in this case. (See Sections \ref{sect:3d_lin} and \ref{sect:3d_0} below.)
\par

In genus $1$ however, the resulting knot invariants \emph{do} involve the higher Hochschild cohomologies $HH^{>0}(\mfk{A})$, at least in principle. We discuss the situation for $\mfk{A}=(\Rep_q\opn{SL}_2)_{\opn{small}}$ in detail in Section \ref{sect:3d_1}, where a number of interesting points arise.

\subsection{A world of homotopical TQFTs}

From a strict mathematical perspective, we have a generic interested in TQFTs for ``modular" $\infty$-categories--a notion which we'll leave somewhat ambiguous for the moment. For any modular $\infty$-category $\msc{A}$, we conjecture the existence of a $3$-dimensional TQFT from marked and unmarked bordism categories
\[
L_{\msc{A}}:\Bord^{nc}_{\msc{A}}\to \Vect\ \ \text{and}\ \ \mbb{L}_{\msc{A}}:\Bord^{nc}_{3,2}\to \Vect
\]
which we furthermore conjecture extend down to closed curves, and recovers $\msc{A}$ as the value on $S^1$. We are also interested in homotopical field theories in dimension $4$, skein modules for modular $\infty$-categories, and their connections to phenomena in dimension $3$. Examples of modular $\infty$-categories should include derived categories of finite modular $1$-categories, and presumably (quantized) derived categories of sheaves on compact symplectic varieties.
\par

If we take this broader point of view, the conclusions of the present text verify one explicit claim in this direction, or maybe a first such claim. We elaborate on the topics proposed above, and also the integration of local systems, in Sections \ref{sect:expansion} and \ref{sect:conjectures} below.

\subsection{Some metacommentary}

This paper is, to some extent, a re-articulation of non-semisimple TQFT through the lense of $\infty$-categories.  We understand, however, that this perspective lies outside of the realm of expertise for many intended readers.  For this reason we have taken a number of steps to make the text more legible for the homotopical non-expert.
\par

First, we provide specific referencing to the literature on $\infty$-categories throughout. In our ideal moments, all homotopical arguments are directly supported by a concrete simplicial analysis or an explicit reference. Second, we provide an extensive review of basic principles for simplicial categories and $\infty$-categories--specifically in the guise of weak Kan complexes--in Section \ref{sect:categories}. Cocartesian fibrations and the un/straightening equivalence are given special attention.  For the ambitious reader, we hope that the paper is more-or-less self-contained.

\subsection{Outline}
{\it Throughout the text $A^{\heartsuit}$ generally denotes a ribbon (or modular) tensor category, rather than the fraktur notation $\mfk{A}$ employed above.}  We swap over to our standard notation $A^{\heartsuit}$ from this point on.
\par
  
The theory $L_{\msc{D}_{fin}}$ is obtained from the corresponding abelian theory $L^{abelian}_{A^{\heartsuit}}$ via a somewhat subtle derivation process, the zeroth step of which already occurred in preceding work with Agustina Czenky \cite{czenkynegron}.  The first point for us is the construction of a bordism category which is sufficiently responsive to basic manipulations within the theory of higher algebra.  We then apply a ``vertical" localization to move from the discrete setting to the homotopical setting, followed by a relative Kan extension to move from bordisms over the homotopy $\infty$-category to the derived setting.  The subsequent removal of markings is a rather idiosyncratic process which has to do with further localization along a class of rather poorly behaved morphisms.  The paper proceeds as follows:

\begin{itemize}
\item {\bf Background.} We discuss the efficacy of markings in topological field theories in Section \ref{sect:skein_bordisms}, just to clarify some confusion around these issues before we begin our study in earnest. In Section \ref{sect:categories} we discuss basic categorical backgrounds, and provide a brisk overview of $\infty$-categories.  In Section \ref{sect:monoidal_infty} we recall basic notions for symmetric monoidal $\infty$-categories, then in Sections \ref{sect:fr_disk} and \ref{sect:disk_monoidal} we introduce a symmetric monoidal $\infty$-category $fr\Disk$ of framed disks which classifies ``compatible pairs" of ``homotopically balanced monoidal" $\infty$-categories. In Section \ref{sect:ribbon_framed} we show that any ribbon tensor category $A^{\heartsuit}$, in the $1$-categorical sense of the term, naturally determines $fr\Disk$-categories $(A^{\heartsuit})^{\ot}\to fr\Disk^{\ot}$ and $A^{\ot}=\opn{Ch}(A^{\heartsuit})^{\ot}\to fr\Disk^{\ot}$.
\vspace{2mm}

\item {\bf Bordism and discrete topological theories.} We construct the trivially marked, $3$-dimensional, non-compact, anomalous bordism $\infty$-category $\Bord_{\ast}^{nc}$ in Section \ref{sect:bord_star}. We subsequently construct the symmetric monoidal $\infty$-category $\Bord^{nc}_{\msc{E}}$ of $\msc{E}$-marked $3$-dimensional bordisms in Section \ref{sect:marked_bords}, where the markings come from any $fr\Disk$-monoidal $\infty$-category $\msc{E}^{\ot}\to fr\Disk^{\ot}$.  We explicitly describe the $A=\opn{Ch}(A^{\heartsuit})$-marked bordism category associated to any ribbon tensor category $A^{\heartsuit}$, then explain how results from \cite{derenzietal23,czenkynegron} provide a TQFT
\begin{equation}\label{eq:390}
L_{A}:\Bord_{A}^{nc}\to \opn{Ch}(\opn{Vect})
\end{equation}
in Section \ref{sect:dis_lrt}. We derive the TQFT $L_A$ at the $1$-categorical level, to produce a symmetric monoidal functor from the bordism category $\Bord_D^{nc}$ with labels from the discrete derived category $D(A^{\heartsuit})$, in Section \ref{sect:1_categorical}. We subsequently discuss the process, and rationalle, for lifting to the $\infty$-categorical level.
\vspace{2mm}

\item {\bf Localization and homotopical theories.} In Section \ref{sect:localization} we present a robust method of localization for monoidal $\infty$-categories, and more generally for cocartesian fibrations, along well-behaved classes of maps. Here we follow work of Hinich \cite{hinich16} in an essential way. We apply this general theory of localization to obtain from \eqref{eq:390} a symmetric monoidal functor
\[
L_{\msc{K}}:\Bord_{\msc{K}}^{nc}\to \Vect
\]
for the $3$-dimensional bordism category with labels in the homotopy $\infty$-category $\msc{K}=\msc{K}(A^{\heartsuit})$. As before, the target $\Vect$ is the derived $\infty$-category of linear cochains.\vspace{2mm}

\item {\bf Kan extension and derived theories.} We provide an analysis of Kan extensions over a base in Section \ref{sect:rel_adj_kan}. Then, in Section \ref{sect:derived_lrt}, we apply Kan extension along the natural functor $\Bord^{nc}_{\msc{K}}\to \Bord^{nc}_{\msc{D}}$ to obtain a $3$-dimensional field theory
\[
L_{\msc{D}}:\Bord^{nc}_{\msc{D}}\to \Vect
\]
for the derived $\infty$-category $\msc{D}=\msc{D}(A^{\heartsuit})$. We furthermore calculate the state spaces for $L_{\msc{D}}$ via the linear mapping spaces in $\msc{D}$, and identify the cohomology of the theory $L_{\msc{D}}$ with the abelian theory $L^{abelian}_{A^{\heartsuit}}$ from De Renzi et al.\ \cite{derenzietal23}.
\vspace{2mm}

\item {\bf Un-marking and inspection.} In Section \ref{sect:derived_w_skeins} we explain how the functor $L_{\msc{D}}$ localizes along a special class $\Theta_{\msc{D}}$ of bordisms in order to produce a TQFT from an $\infty$-category $\Bord_{\msc{D}}[\Theta^{-1}_{\msc{D}}]$ of $3$-dimensional bordisms with embedded ``universal skeins" labeled by $\msc{D}$.  After restricting to unit labels, we obtain a canonical symmetric monoidal inclusion $\Bord_{\ast}^{nc}\to \Bord_{\msc{D}}^{nc}$ and hence a symmetric monoidal functor $\Bord^{nc}_{\ast}\to \Vect$ by pulling back along this inclusion.  We then localize at the corresponding class $\Theta_{\ast}$ in $\Bord^{nc}_{\ast}$, and take a homotopy truncation to get a map of symmetric monoidal $1$-categories
\[
pre\text{-}\mbb{L}_{D_{fin}}:\opn{h}\Bord_{\ast}^{nc}[\Theta^{-1}_{\ast}]\to D(\opn{Vect}).
\]
\par

It is shown in Section \ref{sect:unmarked_lrt} that the forgetful functor $\opn{h}\Bord_{\ast}^{nc}[\Theta^{-1}_{\ast}]\to \Bord_{3,2}^{nc}$ is in fact an equivalence, so that the above map $pre\text{-}\mbb{L}_{D_{fin}}$ determines an unmarked topological field theory
\[
\mbb{L}_{D_{fin}}:\Bord_{3,2}^{nc}\to D(\opn{Vect}).
\]
We conclude Section \ref{sect:unmarked_lrt} with an analysis of the theory $\mbb{L}_{D_{fin}}$ in dimension $2$, and Section \ref{sect:dimension_3} is dedicated to a discussion of phenomena in dimension $3$.
\vspace{2mm}

\item {\bf Conjectures.} In Section \ref{sect:conjectures}, we present a number of fundamental conjectures which concern extensions of our analysis to ``modular" $\infty$-categories which are \emph{not} constructed from modular tensor categories at the abelian level, and also extensions of our theory $L_{\msc{D}}$ down to $S^1$. Following the physics literature on the subject, we further discuss couplings of derived field theories to moduli of local systems.
\vspace{2mm}

\item {\bf Appendices.} We have three appendices. The first, Appendix \ref{sect:more_fre2}, covers various minutiae for the operad of framed colored disks and the second, Appendix \ref{sect:framed_pairs}, explains how maps between balanced monoidal categories produce $fr\Disk$-monoidal $\infty$-categories. The third, Appendix \ref{sect:module_cats}, discusses module $\infty$-categories over symmetric monoidal $\infty$-categories and their inner-Homs.
\end{itemize}

\subsection{Suggested reading}
The astute reader may notice the extraordinary length of this paper. While we hope each (Sub)section serves a clear purpose, one needn't read the paper line-by-line in order to understand either the results or the arguments. We suggest a particular reading of the paper, where one might skim certain materials and read others more thoroughly, in recognition of the relative importance of the details:\vspace{1mm}
\begin{center}
Section \ref{sect:skein_bordisms} (\emph{skim}) $\to$ Section \ref{sect:categories} (\emph{reference as needed}) $\to$ Section \ref{sect:monoidal_infty} (\emph{skim}) $-\!\cdots$\vspace{2mm}\\
$\to$ Sections \ref{sect:fr_disk} \& \ref{sect:disk_monoidal} (\emph{skim}) $\to$ Section \ref{sect:ribbon_framed} (\emph{read-ish}) $-\!\cdots$\vspace{2mm}\\
$\to$ Section \ref{sect:bord_star} (\emph{read Subsection \ref{sect:bord_principles}}) $\to$ Section \ref{sect:marked_bords} (\emph{read}) $-\!\cdots$\vspace{2mm}\\
$\to$ Section \ref{sect:dis_lrt} (\emph{read}) $\to$ Section \ref{sect:1_categorical} (\emph{read}) $\to$ Section \ref{sect:localization} (\emph{skim}) $-\!\cdots$\vspace{2mm}\\
$\to$ Section \ref{sect:homtopy_lrt} (\emph{read}) $\to$ Section \ref{sect:rel_adj_kan} (\emph{read}) $\to$ Section \ref{sect:derived_lrt} (\emph{read}) $-\!\cdots$\vspace{2mm}\\
$\to$ Section \ref{sect:derived_w_skeins} (\emph{see only definitions, Lemma \ref{lem:color_change}, Theorem \ref{thm:merge_loc}, Corollary \ref{cor:der_lrt_skeins}, and Subsection \ref{sect:univ_sk}}) $-\!\cdots$\vspace{2mm}\\
$\to$ Section \ref{sect:unmarked_lrt} (\emph{read}) $\to$ Sections \ref{sect:conjectures} \& \ref{sect:dimension_3} (\emph{reader's choice}) 
\end{center}

\subsection{Acknowledgements}
Thanks to John Francis, David Jordan, and Kevin Walker for very early discussions on TQFTs in the derived setting. Thanks to Marco De Renzi, Nathan Geer, and Matt Young for many discussions on non-semisimple TQFTs. Thanks to Simon Lentner for discussions of--and sporadic writings on--crossed braided categories in the (algebro)geometric context. Thanks to Melissa Zhang and Jernej Grlj for discussions of skein lasagna modules, and possibilities for derived skein modules. Thanks to Benjamin Ha\"ioun, Patrick Kinnear, and Pavel Safronov for discussions of many things, including topological theories and local systems. Thanks to Qi You for sharing insights on Hochschild cohomology for quantum groups. Thanks to Tudor Dimofte and Nathan Benjamin for discussions concerning (T)QFT and topological twisting in physics. Thanks finally to my local colleagues at USC, Aravind Asok, Julian Chaidez, Peter Haine, and Sheel Ganatra, for answering a very large numer of questions regarding topology and homotopy theory. The author was supported by NSF CAREER Grant No.\ DMS-2239698 and Simons Collaboration Grant No.\ 999367.

\section{Remarks on marked bordisms}
\label{sect:skein_bordisms}

Before beginning, let us make some preliminary remarks on marked surfaces, decorated bordisms, and skein relations in general. This topic is not necessarily introductory, but we hope it will preempt some confusion around the issue.

\subsection{Why marked surfaces?} We'll list here three perspectives through which one might conceptualize marked surfaces. The first is through the lense of conformal field theory. Here we understand that a Reshetikhin-Turaev type theory (generally) admits a conformal boundary. The state spaces in the topological bulk are then identified with the conformal blocks for the boundary conformal theory. See for example \cite{fuchsrunkelschweigert02,runkel}.
\par

The conformal blocks in a conformal field theory naturally depend on insertions of VOA modules at points on a given surface. These insertion points are furthermore equipped with local coordinates. This is to say, we expect the state spaces in a Reshetikhin-Turaev theory to accept markings by objects in the defining category, and we expect the marking points to be equipped with such local coordinates. These local coordinates can be thought of as specifying a collection of disk embeddings, rather than embeddings of points.
\par

Alternatively, we can think of a Reshetikhin-Turaev theory itself as a boundary theory to a once extended topological theory in dimensiona $4$. This $4$-dimensional theory should be skein theoretic (cf.\ \cite{brownhaioun,costantinoetalII,haioun}), and objects in any surface category for the $4$-dimensional theory should, in particular, consist of choices of markings on the given surface (see \cite{brownhaioun}). To say that the Reshetikhin-Turaev type theory is defined for marked surface is then to say that we have well-defined values for arbitrary choices of maps from the unit $\opn{Vect}$ to the factorization homology on surfaces. (In short, markings appear in Reshetikhin-Turaev theories because markings appear in factorization homology.)

For a final perspective, we note that any extended topological field theory $Z:\mcl{B}ord^{\star}_{3,2,1}\to \msc{P}\!r_k$, valued in presentable $\infty$-categories say, determines and is determined by an associated fibration $\mcl{B}ord^{\star}_{\msc{Z}}\to \mcl{B}ord^{\star}_{3,2,1}$. For $\msc{Z}=Z(S^1)$, the $(\infty,2)$-category $\mcl{B}ord^{\star}_{\msc{Z}}$ \emph{is already} explicitly constructible as an $(\infty,2)$-category of $\msc{Z}$-marked bordisms (cf.\ Section \ref{sect:bal_description}). This $(\infty,2)$-category restricts, in various ways, to define $\infty$-categories of surfaces with embedded $\msc{Z}$-labeled disks and $3$-manifolds which are labeled--in one way or another--by morphisms in $\msc{Z}$. Hence a TQFT from a sufficiently robust bordism category for marked surfaces is closely related to the existence of an extended topological field theory.

\subsection{Extending into $3$-manifolds}

Having decided that we accept marked surfaces, one has some decisions to make regarding bordisms between such marked surfaces $M_{\alpha}:\Sigma_x\to \Sigma'_y$. Specifically, what kinds of decorations should we have in $M$ which connect the markings on the incoming boundary $\Sigma$ to the markings on the outgoing boundary $\Sigma'$?
\par

In this text we take a fairly minimalist approach to this question, and employ a minimalist solution. In particular, we construct (in Section \ref{sect:marked_bords}) a category of bordisms $\Bord^{nc}_{\msc{E}}$ with embedded cylinders which travel from the markings on the incoming boundary to markings on the outgoing boundary. We later localize to obtain a coarse approximation to bordisms with embedded skein diagrams (see Section \ref{sect:derived_w_skeins}). The minimal skein relations which appear in this process are sufficient to facilitate an \emph{un}marking procedure in Section \ref{sect:unmarked_lrt}, which leads to an unmarked TQFT $\mbb{L}_{D_{fin}}:\opn{Bord}^{nc}_{3,2}\to D(\opn{Vect})$ in Theorem \ref{thm:unmarked_lrt}.
\par

A maximalist approach to this question would seek a truly skein theoretic realization of bordisms. A ``perfect" notion of skeins should facilitate a calculation of factorization homology over surfaces for ribbon $\infty$-categories, as in the $1$-categorical setting \cite{brownhaioun}.

\subsection{Skein modules?}

The production of true skeins, or ``perfect skeins", for ribbon $\infty$-categories is, from our perspective, a deep question which warrants independent consideration. One might define a ribbon $\infty$-category efficiently as an (anomalous) $2$-dimensional topological field theory $\msc{A}$ valued in linear presentable $\infty$-categories $\msc{P}\!r_k$ \cite{mullerwoike23,steinebrunner}. Here we would work specifically with the straightened fibration over a $2$-dimensional bordism category. 
\par

While the author would state, somewhat emphatically, that they do not know how to construct skein modules for ribbon $\infty$-categories, the proposed definition of ribbon $\infty$-categories itself suggests certain constraints which one might attend to. For example, rather than considering skeins as embedded graphs in a $3$-manifold, we might think of skeins as embedded anomalous surfaces in a manifold. Here by an anomalous surface we specifically mean a surface with a chosen bounding $3$-manifold, and we should assume that the bounding $3$-manifold is ``trivial" (e.g.\ a ball in the case of $S^2$ and a doughnut in the case of the torus). One considers such an embedded surface skein, let us call it, along with its vector space obtained by evaluating under the given field theory for our modular $\infty$-category $\msc{A}$.
\par

Skein relations in this surface approach \emph{may} then be facilitated by embedded bordisms between such embedded surfaces, and so leverage some aspects of a $3$-dimensional field theory for $\msc{A}$. So, practically speaking, a perfect theory of skeins and skein relations for a modular $\infty$-category may be contingent on the pre-existence of a $3$-d field theory from a somewhat restrictive marked bordism category. In any case, this topic is completely open, and we encourage the interested reader to consider this line of inquiry, even returning to the $1$-categorical setting if necessary.

\begin{remark}
See Section \ref{sect:conjectures} for additional discussions of skeins and $4$-dimensional TQFTs, in the derived setting.
\end{remark}

\begin{remark}
In this text we model ribbon tensor categories in the $\infty$-categorical context by coupled pairs of framed $\opn{E}_2$-categories, rather than by $2$-dimensional category valued TQFTs. This choice is purely practical, or even superficial, and can be altered in accordance with what's available at the $1$-categorical level.  (See Section \ref{sect:mark_scheme}.)
\end{remark}

\section{($\infty$-)Categories}
\label{sect:categories}

To begin our study, we provide a brisk overview of tensor categories, simplicial categories, $\infty$-categories, Kan complexes, and cocartesian fibrations.  This is mostly to set notations, and also to provide a coarse reference for the non-expert.  For a detailed account we suggest Kerodon \cite{kerodon}.

\subsection{Tensor categories}
\label{sect:tensor_cats}

We follow the conventions of \cite{negron26}.  So, for us, a tensor category is a linear presentable, compactly generated abelian monoidal category for which the subcategories of compact and rigid objects agree, and for which the compacts form a tensor category in the sense of \cite{egno15}.  So, all compact objects are of finite length, have finite-dimensional Homs, and the unit is simple.  We also require that the product commutes with small colimits in each factor, and assume that the base field $k$ is algebraically closed.  Via Ind-completion, and taking the compacts, one can go back and forth between our tensor categories and those of \cite{egno15}.
\par

For a tensor category $B$ we let $B_{fin}$ denote the full subcategory of finite length = rigid objects.  We call a tensor category finite if it admits a compact projective generator, in which case $B$ is equivalent to the category of modules over a finite-dimensional algebra.  We have the standard notions of braided and ribbon tensor categories (see \cite{egno15}).

\begin{definition}
A tensor category $B$ is called ribbon if it comes equipped with a braiding $c$ and balancing transformation $\theta$ \cite[Definition 8.10.1]{egno15} for which $\theta_{x^{\ast}}=\theta_x^{\ast}$ whenever $x$ is rigid.  We call a finite ribbon tensor category modular if, in addition, any object $x$ with trivial double braiding
\[
c_{-,x}\circ c_{x,-}=id_{x\ot-}
\]
is isomorphic to an additive power of the unit $x\cong \1^{\oplus \Lambda}$.
\end{definition}

Throughout all ``modular" tensor categories are assumed to be finite.

\subsection{Simplicial sets and simplicial categories}

For $\Delta$ the category of finite linearly ordered sets and (weakly) order preserving morphisms, we let $\opn{Set}_{\Delta}$ denote the category of functors $\Fun(\Delta^{\opn{op}},\opn{Set})$. This is the category of simplicial sets.  For any finite linearly ordered set $I$ we let $\Delta^I$ denote the standard simplex $\Delta^I=\Hom_{\Delta}(-,I)$. For each $n\in \mathbb{Z}_{\geq 0}$ we take $[n]=\{0,1,\dots,n\}$ and also $\Delta^n=\Delta^{[n]}$.  For any simplicial set $K$ we let $K[n]$ denote the value $K([n])$, and we note that there is a bijection
\[
\{\text{maps }\sigma:\Delta^n\to K\text{ in }\opn{Set}_{\Delta}\}\overset{\sim}\to K[n],\ \ \sigma\mapsto \sigma(id_{[n]}),
\]
by Yoneda.
\par

The symmetric monoidal structure on $\opn{Set}$, under the cartesian product, induces a symmetric monoidal structure on the category $\opn{Set}_{\Delta}$ of simplicial sets. The category $\opn{Set}_{\Delta}$ is also complete and cocomplete, with (co)limits calculated via (co)limits in the target category $\opn{Set}$.
\par

Each object $S$ in $\opn{Set}$ becomes a simplicial set via the constant functor $S:\Delta^{\opn{op}}\to \opn{Set}$ and this gives a fully faithful, symmetric monoidal embedding
\[
\opn{Set}\hookrightarrow \opn{Set}_{\Delta}.
\]
In this way our discussions of simplicial sets, and simplicial categories, specialize to discussions of sets and plain old categories.

\begin{definition}
A simplicial category $\underline{C}$ is a category enriched in $\opn{Set}_{\Delta}$ and a simplicial functor is a functor between $\opn{Set}_{\Delta}$-enriched categories.
\end{definition}

For each $n$ and index $0\leq i\leq n$ we let $\Lambda^n_i\subseteq \Delta^n$ denote the simplicial subset whose simplices $\Delta^m\to \Lambda^n_i$ are precisely those simplices in $\Delta^n$ which factor through some face $\Delta^I\subseteq \Delta^n$ with $I\subseteq [n]$ and $[n]-\{i\}\nsubseteq I$.  We call $\Lambda^n_i$ the \emph{$i$-th horn} in $\Delta^n$.

\begin{definition}
A map of simplicial sets $\msc{X}\to \msc{S}$ is called a Kan fibration if each lifting problem
\begin{equation}\label{eq:kan_lift}
\xymatrix{
\Lambda^n_i\ar[d]\ar[r] & \msc{X}\ar[d]\\
\Delta^n\ar[r]\ar@{..>}[ur] & \msc{S},
}
\end{equation}
with $n$ arbitrary and $0\leq i\leq n$, admits a solution.  We call a simplicial set $\msc{X}$ a Kan complex if the unique map $\msc{X}\to \ast$ to the terminal simplicial set is a Kan fibration.  A morphism between Kan complexes is just a map of simplicial sets. 
\end{definition}

To clarify our language, $\msc{X}\to \msc{S}$ is a Kan fibration if and only if each commuting square as in \eqref{eq:kan_lift} admits a map $\Delta^n\to \msc{X}$ which split the diagram into two commuting triangles.  The terminal simplicial set $\ast$ is, as our notation suggests, the constant simplicial set associated to the (or rather, any) singleton.

\begin{remark}
The singular complex functor $\opn{Sing}:\opn{Top}\to \opn{Kan}$ provides a (Quillen) equivalence between compactly generated, locally Hausdorff topological spaces and Kan complexes.  The inverse is given by geometric realization.  For this reason people often refer to Kan complexes as \emph{Spaces}.
\end{remark}

\begin{definition}
A simplicial category $\underline{C}$ is called fibrant if, for each pair of objects $x$ and $y$ in $\underline{C}$, the Hom complex $\uHom_{\underline{C}}(x,y)$ is a Kan complex.
\end{definition}

We note that any simplicial category $\underline{C}$ has an underlying plain old category $C$ whose morphisms are the $0$-simplices
\[
\Hom_C(x,y):=\Hom_{\opn{Set}_{\Delta}}(\Delta^0,\uHom_{\underline{C}}(x,y)).
\]

\subsection{$\infty$-Categories}

\begin{definition}
A map of simplicial sets $F:\msc{C}\to \msc{T}$ is called an inner fibration if any lifting problem
\begin{equation}\label{eq:infty_lift}
\xymatrix{
\Lambda^n_i\ar[d]\ar[r] & \msc{C}\ar[d]^F\\
\Delta^n\ar[r]\ar@{..>}[ur] & \msc{T},
}
\end{equation}
with $n$ arbitrary and $0<i<n$, admits a solution.  We call $\msc{C}$ an $\infty$-category if the map $\msc{C}\to \ast$ to the terminal simplicial set is an inner fibration.  A functor between $\infty$-categories is just a map of simplicial sets.
\end{definition}

Note that the lifting problems here \eqref{eq:infty_lift} only involve \emph{inner} horns, i.e. the horns $\Lambda^n_i\to \Delta^n$ with $i$ between $1$ and $n-1$.  For an $\infty$-category $\msc{C}$, an \emph{object} in $\msc{C}$ is a $0$-simplex $x:\Delta^0\to \msc{C}$, a \emph{morphism} is a $1$-simplex $\alpha:\Delta^1\to \msc{C}$, and for any two maps $\alpha_{01}:x_0\to x_1$ and $\alpha_{12}:x_1\to x_2$ a composite is any third map $\alpha_{02}:x_0\to x_2$ which fits into a $2$-simplex
\[
\xymatrix{
 & x_1\ar[dr]^{\alpha_{12}} \\
x_0\ar[ur]^{\alpha_{01}}\ar[rr]_{\alpha_{02}} & & x_2.
}
\]
We note that the composite is not uniquely determined in a strict sense, but it is determined up to a contractible space of choices \cite[\href{https://kerodon.net/tag/01BW}{01BW}]{kerodon}.  So, whenever such ambiguities are acceptable, we might speak of ``the" composite $\alpha_{02}$ and write simply $\alpha_{02}=\alpha_{12}\alpha_{01}$.

\begin{definition}
Given an object $x$ in an $\infty$-category $\msc{C}$, the identity map $id_x:x\to x$ is the unique degnerate $1$-simplex $id_x:\Delta^1\to \msc{C}$ with $id_x|_{\{0\}}=id_x|_{\{1\}}=x$.  An isomorphism in $\msc{C}$ is a map $\alpha:x\to y$ which admits another map $\beta:x\to y$ for which we have composites $\alpha\beta=id_y$ and $\beta\alpha=id_x$.
\end{definition}

An ($\infty$-)subcategory in an $\infty$-category $\msc{C}$ is a simplicial subset $\msc{C}'\subseteq \msc{C}$ in which a given $n$-simplex $\sigma:\Delta^n\to \msc{C}$ lies in $\msc{C}'$ if and only if each object $\sigma|_{\{i\}}$ and each edge $\sigma|_{\Delta^{\{i,j\}}}$ lies in $\msc{C}'$.  We assume additionally that each composable pair or maps in $\msc{C}'[1]$ admits a composition in $\msc{C}'[1]$. A subcategory is called \emph{full} if each morphism $\Delta^1\to \msc{C}$ between objects in $\msc{C}'$ is also in $\msc{C}'$.

\begin{definition}\label{def:subcats}
Given a collection of objects $X_0$ in an $\infty$-category $\msc{C}$, the full subcategory spanned by $X_0$ in $\msc{C}$ is the unique full subcategory $\msc{C}'\subseteq\msc{C}$ with $\msc{C}'[0]=X_0$. Given a collection of objects $X_0$ in $\msc{C}$ and morphisms $X_1$ between objects in $X_0$ which are stable under composition, the subcategory spanned by the $X_i$ in $\msc{C}$ is the unique subcategory $\msc{C}'\subseteq\msc{C}$ with $\msc{C}'[0]=X_0$ and $\msc{C}'[1]=X_1$. 
\end{definition}

One checks directly that any subcategory $\msc{C}'$ in an $\infty$-category $\msc{C}$ is itself an $\infty$-category.

\begin{proposition}[{\cite[\href{https://kerodon.net/tag/01H1}{01H1}]{kerodon}}]
For any $\infty$-category $\msc{C}$ the subcategory $\msc{C}^{\opn{Kan}}$ spanned by all objects, but only isomorphisms between objects, is a Kan complex.
\end{proposition}

The simplicial set $\msc{C}^{\opn{Kan}}$ is usually called the core of $\msc{C}$, though we will just call it the associated Kan complex.

\subsection{Functor categories and mapping spaces}

For any $\infty$-category $\msc{C}$ and simplicial set $K$ the complex of simplicial maps
\[
\Fun(K,\msc{C}),\ \ \Fun(K,\msc{C})[n]=\Hom_{\opn{Set}_{\Delta}}(\Delta^n\times K,\msc{C}),
\]
is an $\infty$-category, and for any map of simplicial sets $L\to K$ the restriction functor
\[
\Fun(K,\msc{C})\to \Fun(L,\msc{C})
\]
is an inner fibration \cite[\href{https://kerodon.net/tag/01BT}{01BT}]{kerodon}.  For objects $x,y:\ast\to \msc{C}$ we obtain the standard mapping complex via pullback
\[
\xymatrix{
\Maps_{\msc{C}}(x,y)\ar[r]\ar[d] & \Fun(\Delta^1,\msc{C})\ar[d]^{\opn{restrict}}\\
\ast\amalg \ast\ar[r]_(.4){[x\ y]} & \Fun(\partial\Delta^1,\msc{C}).
}
\]
As the pullback of any kind of ``fibration" is again a fibration of the same type, we understand that $\Maps_{\msc{C}}(x,y)$ is an $\infty$-category at all $x$ and $y$.  Since all morphisms in $\Maps_{\msc{C}}(x,y)$ are isomorphisms \cite[\href{https://kerodon.net/tag/01DK}{01DK}]{kerodon}, this $\infty$-category is furthermore a Kan complex.

\begin{definition}
For objects in an $\infty$-category $x,y:\ast\to \msc{C}$, the associated mapping space $\Maps_{\msc{C}}(x,y)$ is the Kan complex constructed as above.
\end{definition}

We note that, in general, the mapping spaces do not admit strict composition functors.  We do however have composition functors
\[
\circ:\Maps_{\msc{C}}(x_1,x_2)\times \Maps_{\msc{C}}(x_0,x_1)\to \Maps_{\msc{C}}(x_0,x_2)
\]
which are well-defined up to homotopy \cite[\href{https://kerodon.net/tag/01PQ}{01PQ}]{kerodon}.

\subsection{Isofibrations}

\begin{definition}
An isofibration between $\infty$-categories is an inner fibration $F:\msc{C}\to \msc{D}$ which satisfies either of the following (equivalent) lifting conditions:
\begin{itemize}
\item For each isomorphism $\bar{\alpha}:\bar{x}\to \bar{y}$ in $\msc{D}$, and object $x$ in $\msc{C}$ with $F(x)=\bar{x}$, there is an isomorphism $\alpha:x\to y$ in $\msc{C}$ with $F(\alpha)=\bar{\alpha}$.\vspace{1.5mm}
\item For each isomorphism $\bar{\alpha}:\bar{x}\to \bar{y}$ in $\msc{D}$, and object $y$ in $\msc{C}$ with $F(y)=\bar{y}$, there is an isomorphism $\alpha:x\to y$ in $\msc{C}$ with $F(\alpha)=\bar{\alpha}$.
\end{itemize}
\end{definition}

Isofibrations are the ``good" maps between $\infty$-categories, especially when one wants to take limits of $\infty$-categories in a homotopically responsible manner.  See for example \cite[\href{https://kerodon.net/tag/03GF}{03GF}, \href{https://kerodon.net/tag/03H3}{03H3}, \href{lurie09p://kerodon.net/tag/033P}{033P}]{kerodon}.

\begin{remark}
Isofibrations are also referred to as categorical fibrations in the literature. Compare for example \cite{lurie09} to \cite{cisinski19,kerodon}. 
\end{remark}

\begin{remark}
We've obviously chosen to work with ``quasi-categories", or ``weak Kan complexes" as our preferred model for $\infty$-categories.  Though we are not concerned with the sinful nature of such choices, there are practical points to consider here.  In particular, if we want to embed our analysis of $3$-d TQFTs within a larger study of extended field theories, complete Segal spaces are likely preferable to quasi-categories as a model for $\infty$-categories.
\end{remark}

\subsection{Homotopy coherent nerves}

From any simplicial category $\underline{C}$ one has an associated simplicial set $\opn{N}(\underline{C})$ which is produced via a construction called the \emph{homotopy coherent nerve}.  Simplices in the nerve are given by functors
\[
\opn{N}(\msc{C})[n]=\Fun_{\opn{Cat}_{\Delta}}(\opn{Path}\Delta^n,\underline{C})
\]
from a certain path category construction \cite[\href{https://kerodon.net/tag/00KN}{00KN}, \href{https://kerodon.net/tag/00KS}{00KS}]{kerodon}.  We have $\opn{Path}\Delta^0$ the one point category and $\opn{Path}\Delta^1$ the category with two objects, trivial endomorphisms, and a single map between them.  Thus the $0$ and $1$-simplices in the nerve consist of objects and maps in the underlying plain category $C$ for $\underline{C}$.  A $2$-simplex in the nerve consists of a noncommuting diagram
\[
\xymatrix{
 & x_1\ar[dr]^g \\
x_0\ar[rr]_h\ar[ur]^{f} & & x_2 
}
\]
in the plain category $C$, along with a choice of a $1$-simplex $\zeta:\Delta^1\to \uHom_{\underline{C}}(x_0,x_2)$ for which $\zeta|_{\{0\}}=gf$ and $\zeta|_{\{1\}}=h$.
\par

The nerve construction is functorial in the sense that any simplicial functor $F:\underline{C}\to \underline{D}$ determines a map of simplicial sets $\opn{N}(F):\opn{N}(\underline{C})\to \opn{N}(\underline{D})$, and for any composite of simplicial functors we have $\opn{N}(G)\circ\opn{N}(F)=\opn{N}(G\circ F)$ 

\begin{theorem}[\cite{cordierporter86}]
If $\underline{C}$ is a fibrant simplicial category (i.e.\ is enriched in Kan complexes) then the homotopy coherent nerve $\opn{N}(\underline{C})$ is an $\infty$-category.
\end{theorem}

Recall that we have the symmetric embedding $\opn{Set}\hookrightarrow \opn{Set}_{\Delta}$ so that plain categories embed fully faithfully into simplicial categories.  In the case where $C$ is a plain category the nerve $\opn{N}(C)$ has $n$-simplices specified by tuples of maps $\{f_{ij}:0\leq i<j\leq n\}$, $f_{ij}:x_i\to x_j$, with $f_{jk}f_{ij}=f_{ik}$ whenever $i\leq j\leq k$.  The $\infty$-category $\msc{C}=\opn{N}(C)$ has the property that any inner horn $\Lambda^n_i\to \msc{C}$ admits a \emph{unique} filling to an $n$-simplex $\Delta^n\to \msc{C}$, and any $\infty$-category $\msc{C}$ with such unique horn fillings is of the form $\msc{C}=\opn{N}(C)$ for uniquely determined $C$. 

\begin{definition}
We call an $\infty$-category $\msc{C}$ discrete if, for each $n$ and $0<i<n$, the restriction functor $\Fun(\Delta^n,\msc{C})\to \Fun(\Lambda^n_i,\msc{C})$ is an isomorphism of simplicial sets.
\end{definition}

\begin{example}[Kan complexes]
For Kan complexes $\msc{X}$ and $\msc{Y}$ the functor complexes $\Fun(\msc{X},\msc{Y})$ are all Kan complexes \cite[\href{https://kerodon.net/tag/00TN}{00TN}]{kerodon}. We therefore obtain a corresponding fibrant simplicial category $\underline{\opn{Kan}}$ of Kan complexes. We take the homotopy coherent nerve to obtain the $\infty$-category of $\msc{K}\!an=\opn{N}(\underline{\opn{Kan}})$ of Kan complexes.
\end{example}

As we explained above in the generic setting, objects in $\msc{K}\!an$ are Kan complexes and morphisms are maps of Kan complexes, i.e.\ maps of simplicial sets.  A $2$-simplex
\[
\xymatrix{
	&\msc{X}_1\ar[dr]^{f_{12}} \\
\msc{X}_0\ar[rr]_{f_{02}}\ar[ur]^{f_{01}} & & \msc{X}_2
}
\]
in $\msc{K}\!an$ consists of choices of maps between Kan complexes and an additional homotopy $\zeta:\Delta^1\times \msc{X}_0\to \msc{X}_2$ from $f_{12}f_{01}$ to $f_{02}$. A map between Kan complexes $f:\msc{X}\to\msc{Y}$ is an isomorphism in $\msc{K}\!an$ if and only if it is a homotopy equivalence, in the expected sense \cite[\href{https://kerodon.net/tag/00U1}{00U1}]{kerodon}.

\begin{definition}
A Kan complex $\msc{X}$ is called contractible if the unique map $\msc{X}\to \ast$ is a homotopy equivalence.
\end{definition}

\begin{example}[$\infty$-categories]
For $\infty$-categories the functor complex $\Fun(\msc{C},\msc{D})$ is not a Kan complex in general. However, we can take the associated Kan complexes $\Fun(\msc{C},\msc{D})^{\opn{Kan}}$ to obtain a fibrant simplicial category $\underline{\opn{Cat}}_{\infty}^+$.  The nerve $\msc{C}\!at_{\infty}=\opn{N}(\underline{\opn{Cat}}^+_{\infty})$ is the $\infty$-category of $\infty$-categories.
\end{example}

Let $\msc{C}$ and $\msc{D}$ be $\infty$-categories. A transformation $\zeta:\Delta^1\times \msc{C}\to \msc{D}$ between functors $F_0$ and $F_1$ is a natural isomorphism in $\msc{C}\!at_{\infty}$, and hence lies in $\Fun(\msc{C},\msc{D})^{\opn{Kan}}$, if and only if it evaluates to an isomorphism $\zeta_x:F_0(x)\to F_1(x)$ at each $0$-simplex $x$ in $\msc{C}$ \cite[\href{https://kerodon.net/tag/01DK}{01DK}]{kerodon}.  So again, $0$, $1$, and $2$-simplices in $\sCat_{\infty}$ can be understood in concrete terms.  

\begin{example}[Plain categories]
For plain categories, the functor category of the nerves $\Fun(\opn{N}(C),\opn{N}(D))$ is just the nerve of the usual category $\opn{N}(\Fun(C,D))$ of functors and natural transformations.  The associated Kan complex $\opn{N}(\Fun(C,D))^{\opn{Kan}}$ is the nerve $\opn{N}(\Fun(C,D)^{\text{nat.isom}})$ of the groupoid of functors and natural isomorphisms.  In any case, we take the associated Kan complexes of the functor categories to get a fibrant simplicial category $\underline{\opn{Cat}}^+$ of plain categories, and take the homotopy coherent nerve to get the $\infty$-category $\msc{C}\!at=\opn{N}(\underline{\opn{Cat}})$ of plain categories. 
\end{example}

Note that $\sCat$ sits in $\sCat_{\infty}$ as the full subcategory spanned by discrete $\infty$-categories, and $\sKan$ sits in $\sCat_{\infty}$ as the full subcategory spanned by Kan complexes. Also, we can consider the underlying plain categories $\opn{Kan}$ and $\opn{Cat}_{\infty}$, of Kan complexes and $\infty$-categories in $\underline{\opn{Kan}}$ and $\underline{\opn{Cat}}_{\infty}$ respectively.  Applying the homotopy coherent nerve, we obtain inclusions of $\infty$-categories
\[
\opn{N}(\opn{Kan})\hookrightarrow \msc{K}\!an\ \ \text{and}\ \ \opn{N}(\opn{Cat}_{\infty})\hookrightarrow \msc{C}\!at_{\infty}.
\]

\begin{remark}
We have ignored all size constraints in the above examples.  Formally, we fix a sequence of universes
\[
\msf{Small}\ \subseteq\ \msf{Medium}\ \subseteq\ \msf{Large}\ \subseteq\ \msf{Huge}
\]
in order to speak of small, medium sized, large, and huge simplicial sets.  We therefore have small, medium sized, large, and huge $\infty$-categories which are all weak Kan complexes in their respective universes.  The above examples are all large $\infty$-categories. Your run-of-the-mill $\infty$-category is assumed to be of a medium size.
\end{remark}

\begin{definition}
A functor between $\infty$-categories $F:\msc{C}\to \msc{D}$ is an equivalence if it is an isomorphism in $\sCat_{\infty}$.  More directly, $F$ is an equivalence if there exists a functor $F':\msc{D}\to \msc{C}$ which admits natural isomorphisms $F'F\overset{\sim}\to id_{\msc{C}}$ and $id_{\msc{D}}\overset{\sim}\to FF'$.
\end{definition}

The following result is fundamental, and is used ad nauseam throughout the text.

\begin{lemma}\label{lem:pullback_equiv}
Given an equivalence of $\infty$-categories $p:\msc{S}\to \msc{T}$ and a pullback diagram
\[
\xymatrix{
\msc{C}\ar[r]^F\ar[d] & \msc{D}\ar[d]^q\\
\msc{S}\ar[r]_{p} & \msc{T}
}
\]
in $\opn{Cat}_{\infty}$, the map $F$ is an equivalence provided one of $p$ or $q$ is an isofibration.
\end{lemma}

We note that, for any pair of maps $p_i:\msc{C}_i\to \msc{T}$ in $\opn{Cat}_{\infty}$, the pullback of simplicial sets $\msc{C}_0\times_{\msc{T}}\msc{C}_1$ is automatically an $\infty$-category whenever one of the $p_i$ is an inner fibration. So, in almost-all situations, pullbacks in the plain category $\opn{Cat}_{\infty}$ are just pullbacks in $\opn{Set}_{\Delta}$.

\begin{proof}[Proof of Lemma \ref{lem:pullback_equiv}]
If one of $p$ of $q$ is an isofibration then the pullback diagram is a categorical pullback diagram, in the sense of \cite[\href{https://kerodon.net/tag/032Y}{032Y}]{kerodon}. Hence the result follows by \cite[\href{https://kerodon.net/tag/033M}{033M}]{kerodon}.
\end{proof}

\subsection{Language: Discrete $(\infty-)$categories}

{\it From this point on we don't distinguish between a plain category and its nerve. A functor $F:\msc{C}\to D$ from an $\infty$-category to a plain category, or vise-versa, is always understood to be a functor from $\msc{C}$ to the nerve of $D$, or vise-versa.  By a \emph{discrete} category we either mean a discrete $\infty$-category, or just a plain category, depending on the context.}

\begin{notation}
A diagram in an $\infty$-category $\msc{C}$ is a map $d:K\to \msc{C}$ from a simplicial set $K$. When $\msc{C}$ is the nerve of a fibrant simplicial category $\msc{C}=\opn{N}(\underline{C})$, and $C$ is the underling discrete/plain category for $\underline{C}$, a diagram $d:K\to \msc{C}$ is said to strictly commute if $d$ factors through the simplicial inclusion $C\to \msc{C}$,
\[
\xymatrixrowsep{3mm}
\xymatrix{
K\ar[rr]^{d}\ar@{-->}[dr]_{\exists!} & & \msc{C}\\
 & C\ar[ur] & .
}
\]
\end{notation}

When we speak of diagrams in an $\infty$-category $K\to \msc{C}$ we usually assume $K$ is a discrete category, so that a diagram in $\msc{C}$ is exactly what one expects intuitively. We most often speak of strict commutativity in the cases where $\msc{C}=\sCat_{\infty}$ or $\sKan$.

\subsection{Nerves of strict $2$-categories}

In the case of a plain category $C$, we saw that there is no interesting information in an $n$-simplex $\Delta^n\to \opn{N}(C)$ above dimension $1$.  One similarly finds that, for a $2$-category, an $n$-simplex in the nerve is nothing more than a compatible collection of $2$-simplices.  Let us elaborate.

By a strict $2$-category we mean a category enriched in plain categories.  Taking the nerves on morphisms, we can view strict $2$-categories as precisely those simplicial categories whose Hom complexes are all discrete $\infty$-categories.  A strict $2$-category $\underline{C}$ is fibrant, when considered as a simplicial category, if and only if $\underline{C}$ is enriched in groupoids.

\begin{lemma}\label{lem:2cat}
Let $\underline{C}$ be a strict $2$-category.  An $n$-simplex in the nerve $\sigma:\Delta^n\to \opn{N}(\underline{C})$ consists of the following data:
\begin{itemize}
	\item A choice of an object $x_i$ for each index $i=1,\dots,n$.\vspace{1mm}
	\item For each pair of indices $0\leq i<j\leq n$, a choice of map $\alpha_{ij}:x_i\to x_j$ in the underlying discrete category $C$.\vspace{1mm}
	\item For each triple $0\leq i<j<k\leq n$, the choice of a $2$-morphism $\zeta_{ijk}:\alpha_{jk}\alpha_{ij}\to \alpha_{ik}$.
\end{itemize}
These data are (only) required to satisfy the constraint
\[
\zeta_{ikl}(id_{\alpha_{kl}}\circ \zeta_{ijk})=\zeta_{ijl}(\zeta_{jkl}\circ id_{\alpha_{ij}}):\alpha_{kl}\alpha_{jk}\alpha_{ij}\to \alpha_{il}
\]
for each quadruple $0\leq i<j<k<l\leq n$.
\end{lemma}

\begin{proof}
The above data specify an $n$-simplex in the Duskin nerve \cite[\href{https://kerodon.net/tag/009T}{009T}]{kerodon}, and the Duskin nerve is identified with the homotopy coherent nerve in the strict setting \cite[\href{https://kerodon.net/tag/00KY}{00KY}]{kerodon}.
\end{proof}

\subsection{The homotopy category of an $\infty$-category}

Given an $\infty$-category $\msc{C}$, the homotopy category $\opn{h}\msc{C}$ is the discrete category with the same objects as $\msc{C}$ and morphisms given by equivalence classes $[\alpha]:x\to y$ of maps in $\msc{C}$.  Here $\alpha$ is equivalent to $\alpha'$ whenever there is a diagram of the form
\[
\xymatrix{
x\ar[r]^{\alpha'}\ar[d]_{id_x}\ar[dr] & y\ar[d]^{id_y}\\
x\ar[r]_{\alpha} & y
}
\] 
in $\msc{C}$.  This relation is stable under composition \cite[\href{https://kerodon.net/tag/0048}{0048}]{kerodon}, so that $\opn{h}\msc{C}$ is in fact a discrete category.
\par

Equivalently, the homotopy category $\opn{h}\msc{C}$ is defined by taking the same objects as $\msc{C}$ and morphisms given by the connected components
\[
\Hom_{\opn{h}\msc{C}}(x,y)=\pi_0\big(\Maps_{\msc{C}}(x,y)\big).
\]
The homotopy category construction is functorial, and provides a left adjoint to the inclusion $\sCat\to \sCat_{\infty}$.

\subsection{Cocartesian fibrations and straightening}
\label{sect:co_carts}

\begin{definition}
Let $q:\msc{C}\to \msc{T}$ be a functor between $\infty$-categories.  A morphism $\alpha:x\to y$ in $\msc{C}$ is called $q$-cocartesian if any lifting problem of the form
\[
\xymatrix{
	\Lambda^n_0\ar[d]\ar[r]^{\bar{\sigma}} & \msc{C}\ar[d]^q\\
	\Delta^n\ar@{..>}[ur]\ar[r] & \msc{T},
}
\]
with $n>1$ and specified initial edge $\bar{\sigma}|_{\Delta^{\{0,1\}}}=\alpha$, admits a solution.  An inner fibration $q:\msc{C}\to \msc{T}$ is called a cocartesian fibration if, for any map $\bar{\alpha}:\bar{x}\to \bar{y}$ in $\msc{T}$ and any $x$ in $\msc{C}$ with $q(x)=\bar{x}$, there is a $q$-cocartesian morphism $\alpha:x\to y$ in $\msc{C}$ with $q(\alpha)=\bar{\alpha}$.
\end{definition}

One can argue directly that cocartesian lifts of morphisms are unique up to isomorphism in $\Fun(\Delta^1,\msc{C})$. In particular, for a fixed functor $q:\msc{C}\to \msc{T}$, any two $q$-cocartesian lifts $\alpha:x\to y$ and $\alpha':x\to y'$ of a given map $\bar{\alpha}:\bar{x}\to \bar{y}$ in the base $\msc{T}$ fit into a diagram
\[
\xymatrix{
	& y\ar[dr]^{\beta} &	&\ar@{}[d]|{\text{\normalsize over}} & & \bar{y}\ar[dr]^{id} & \\
x\ar[ur]^{\alpha}\ar[rr]_{\alpha'} & & y' & & \bar{x}\ar[ur]^{\bar{\alpha}}\ar[rr]_{\bar{\alpha}} & & \bar{y}
}
\]
in which $\beta$ is an isomorphism. In fact, such lifts are unique up to a contractible space of choices \cite[\href{https://kerodon.net/tag/01VK}{01VK}]{kerodon}. 

\begin{definition}
A map of cocartesian fibrations is a strictly commuting diagram
\[
\xymatrix{
\msc{C}\ar[rr]^F\ar[dr]_{q} & & \msc{D}\ar[dl]^{p}\\
	&\msc{T}
}
\]
in which $q$ and $p$ are cocartesian, and in which $F$ sends $q$-cocartesian edges in $\msc{C}$ to $p$-cocartesian edges in $\msc{D}$.  We call such a map $F$ an equivalence of cocartesian fibrations if, after forgetting the base $\msc{T}$, it is an equivalence of $\infty$-categories.
\end{definition}

\begin{remark}
Any equivalence between cocartesian fibrations does admit an inverse over the base \cite[\href{https://kerodon.net/tag/0285}{0285}, \href{https://kerodon.net/tag/028B}{028B}]{kerodon}. Furthermore, there is an $\infty$-category of cocartesian fibrations over $\msc{T}$ in which the isomorphisms are precisely the equivalences, as defined above (see Section \ref{sect:cocart_sm}).
\end{remark}

Cocartesian fibrations serve, arguably, as the leverage point through which the theory of $\infty$-categories functions.  Essentially, every ``thing" in $\msc{C}\!at_{\infty}$ is a cocartesian (or cartesian) fibration over a specified base.  We take a moment to discuss these structures in more detail.

\begin{notation}
Given an inner fibration $q:\msc{C}\to \msc{T}$ and an object $t:\ast\to \msc{T}$, we let $\msc{C}_t$ denote the associated fiber $\msc{C}_t=\{t\}\times_{\msc{T}}\msc{C}$.
\end{notation}

We note that the fibers along any inner fibration are all $\infty$-categories, since inner fibrations are stable under pullback. Of course, at the moment we are especially interested in the case where the given map is a cocartesian fibration.

For any cocartesian fibration $q:\msc{C}\to \msc{T}$ and any simplicial set $K$, the induced functor
\[
q_\ast:\Fun(K,\msc{C})\to \Fun(K,\msc{T})
\]
is a cocartesian fibration as well.  The $q_{\ast}$-cocartesian edges in $\Fun(K,\msc{C})$ are precisely those transformations $\zeta:\Delta^1\times K\to \msc{C}$ which evaluate to a $q$-cocartesian edge $\zeta_x:\Delta^1\to \msc{C}$ at each $0$-simplex $x$ in $K$ \cite[\href{https://kerodon.net/tag/01VG}{01VG}]{kerodon}.  Hence, for any map $f:s\to t$ in $\msc{T}$ and corresponding functor $\tilde{f}=fp_1:\Delta^1\times \msc{C}_s\to \msc{T}$, there is a unique cocartesian transformation $T_f:\Delta^1\times \msc{C}_s\to \msc{C}$ which solves the lifting problem
\[
\xymatrix{
\msc{C}_s=\{0\}\times\msc{C}_s\ar[rr]^(.55){\opn{include}}\ar[d] & & \msc{C}\ar[d]^q\\
\Delta^1\times\msc{C}_s\ar[rr]_{\tilde{f}}\ar@{..>}[urr] & & \msc{T}.
}
\]
Evaluating this transformation at $1$ in $\Delta^1$, we ``travel along $f$" to obtain a \emph{transport functor} between the fibers $f_!:\msc{C}_s\to \msc{C}_t$.  One sees, via uniqueness of cocartesian lifts, that any sequence of maps $gf:s\to t\to u$ in the base produces a diagram
\[
\xymatrix{
	& \msc{C}_t\ar[dr]^{g_!}\\
\msc{C}_s\ar[ur]^{f_!}\ar[rr]_{(gf)_!} & & \msc{C}_u
}
\]
in $\sCat_{\infty}$.
\par

Taking this concept to its extreme, we find that any cocartesian fibration $q:\msc{C}\to \msc{T}$ determines an associated functor $\opn{St}(q):\msc{T}\to \sCat_{\infty}$ with $\opn{St}_q(s)=\msc{C}_s$, $\opn{St}_q(f)=f_!$, etc. It turns out that the fibration $q:\msc{C}\to \msc{T}$ is also \emph{reconstructible} from its associated functor.  Indeed, Grothendieck straightening tells us that the assignment $q\mapsto \opn{St}(q)$ extends to an equivalence
\[
\opn{St}:\msc{C}\!ocart(\msc{T})\overset{\sim}\to \Fun(\msc{T},\sCat_{\infty})
\]
from an $\infty$-category of cocartesian fibrations over a given base $\msc{T}$ (see Section \ref{sect:cocart_sm}) to the $\infty$-category of $\sCat_{\infty}$-valued functors from $\msc{T}$ \cite{lurie09,hebestreitheutsruit25,cisinskinguyen}.  The inverse
\[
\opn{Un}:\Fun(\msc{T},\sCat_{\infty})\overset{\sim}\to \msc{C}\!ocart(\msc{T})
\]
is provided by pulling back along the universal fibration $q_{\opn{univ}}:\msc{P}\sCat_{\infty}\to \sCat_{\infty}$ of pointed $\infty$-categories \cite[\href{https://kerodon.net/tag/020S}{0205}, \href{https://kerodon.net/tag/028K}{028K}]{kerodon}.
\par

Now, when we think of a ``structured" ($\infty$-)category $\msc{C}$, like a monoidal ($\infty$-)category for example, we are usually thinking about a functor $F_{\msc{C}}:T\to \sCat_{\infty}$ from some specific category $T$ which universally encodes a collection of functors and transformations which constitute the given structure.  Hence, from the perspective of Grothendieck straightening and unstraightening, we now have two options: We can think of such structured $\infty$-categories as a distinguished class of cocartesian fibrations over $T$, or as a class of functors to $\sCat_{\infty}$.
\par

Though there is some utility in working on both sides of the straightening-unstraightening equivalence, {\it it is often preferable to work with cocartesian fibrations rather than $\sCat_{\infty}$-valued functors}. The reason is, arguably, that a cocartesian fibration $\msc{C}\to \msc{T}$ contains less erroneous information than its associated functor to $\sCat_{\infty}$.

\begin{remark}
There are some size constraints to consider here, all of which we've ignored.  For example, the total space $\msc{C}$ in a cocartesian fibrations $\msc{C}\to \msc{T}$ is allowed to be large, but all of the fibers $\msc{C}_t$ are assumed to be of a medium size.
\end{remark}

\subsection{Functor categories for cocartesian fibrations}
\label{sect:fun_cc}

Consider inner fibrations $q:\msc{C}\to \msc{T}$ and $p:\msc{D}\to \msc{T}$.  The structure map on $\msc{D}$ induces an inner fibration $p_\ast:\Fun(\msc{C},\msc{D})\to \Fun(\msc{C},\msc{T})$ \cite[\href{https://kerodon.net/tag/01BV}{01BV}]{kerodon} and we take the fiber along $q$ to obtain an $\infty$-category
\[
\Fun_{\msc{T}}(\msc{C},\msc{D}):= \{q\}\times_{\Fun(\msc{C},\msc{T})}\Fun(\msc{C},\msc{D}).
\]
Objects in this $\infty$-category are functors over $\msc{T}$, morphisms are transformations over $\msc{T}$, and general $n$-simplices are strictly commuting diagrams
\[
\xymatrix{
\Delta^n\times \msc{C}\ar[rr]^{\sigma}\ar[dr]_{qp_2} & & \msc{D}\ar[dl]^p\\
	& \msc{T} & .
}
\]

\begin{definition}
For cocartesian fibrations as above, we take
\[
\Fun^{cc}_{\msc{T}}(\msc{C},\msc{D})\ \subseteq\ \Fun_{\msc{T}}(\msc{C},\msc{D})
\]
the full subcategory spanned by those functors which send $q$-cocartesian edges in $\msc{C}$ to $p$-cocartesian edges in $\msc{D}$.
\end{definition}

\section{Symmetric monoidal $\infty$-categories}
\label{sect:monoidal_infty}

We recall basic information regarding symmetric monoidal $\infty$-categories. We are especially interested in the transition from symmetric monoidal simplicial categories to symmetric monoidal $\infty$-categories. We also discuss equifibered symmetric monoidal $\infty$-categories which, to a large extent, replace $\infty$-operads in our study \cite{haugsengkock24,barkanetal25}.

\subsection{Finite pointed sets}

\begin{definition}
The category of finite pointed sets $\opn{Fin}_\ast$ is the category of finite sets, along with morphisms $\bar{f}:I\to J$ given by the choice of a subset $I_{\bar{f}}\subseteq I$ and a set map $I_{\bar{f}}\to J$.  Given such a map $\bar{f}:I\to J$ in $\Fin_{\ast}$, we let $I_{j}$ denote the preimage $\bar{f}^{-1}(j)\subseteq I$ of an element $j\in J$.
\end{definition}

By a slight abuse of notation we let $\bar{f}:I_{\bar{f}}\to J$ denote the set map underlying a map of finite pointed sets $\bar{f}:I\to J$.  The composition of maps $\bar{f}:I\to J$ and $\bar{g}:J\to K$ in $\Fin_{\ast}$ is defined by taking the subset $I_{\bar{g}\bar{f}}=\bar{f}^{-1}(J_{\bar{g}})$ along with the underlying set map
\[
\bar{g}\bar{f}|_{I_{\bar{g}\bar{f}}}:I_{\bar{g}\bar{f}}\to J_{\bar{g}}\to K.
\]

\begin{definition}
A map of finite pointed sets $\bar{f}:I\to J$ is called active if $I_{\bar{f}}=I$, and called inert if each preimage $I_j=\bar{f}^{-1}(j)$ is a singleton. For any index $i$ in $I$ we let $\bar{\rho}_i:I\to \{0\}$ denote the unique inert morphism with preimage $\bar{\rho}_i^{-1}(0)=\{i\}$.
\end{definition}

Fibrations over $\Fin_{\ast}$ play a prominent role in higher algebra. We note that any functor to $\Fin_{\ast}$ is automatically an inner fibration, since $\Fin_{\ast}$ is discrete.

\begin{definition}
Given a functor $q:\msc{C}\to \Fin_{\ast}$, a morphism $f:x\to y$ in $\msc{C}$ is called inert if it is $q$-cocartesian and the image $q(f)$ is inert in $\Fin_{\ast}$. We call $f$ active if it has active image in $\Fin_{\ast}$.
\end{definition}

\subsection{Symmetric monoidal $\infty$-categories}

\begin{definition}\label{def:sm}
A symmetric monoidal $\infty$-category is a cocartesian fibration $\msc{E}^{\ot}\to \opn{Fin}_{\ast}$ for which transport along the inert projections $\bar{\rho}_i:I\to \{0\}$ induce an equivalence
\begin{equation}\label{eq:849}
\rho_!:\msc{E}^{\ot}_I\overset{\sim}\to \msc{E}^I.
\end{equation}
A map of symmetric monoidal $\infty$-categories is a map of cocartesian fibrations
\[
\xymatrix{
\msc{E}^{\ot}\ar[rr]^F\ar[dr] & & \msc{D}^{\ot}\ar[dl]\\
	& \opn{Fin}_{\ast} & ,
}
\]
i.e.\ a functor over $\opn{Fin}_{\ast}$ which preserves cocartesian edges.
\end{definition}

In the expression \eqref{eq:849} $I$ is an arbitrary finite set, $\msc{E}$ is the fiber over the singleton $\msc{E}=\msc{E}^{\ot}_{\{0\}}$, and $\msc{E}^I$ is the exponent $\msc{E}^I=\Fun(I,\msc{E})$.

\begin{remark}
We employ no special notation for $\infty$-categories of symmetric monoidal functors, as the generic cocartesian functor category
\[
\Fun_{\opn{Fin}_{\ast}}^{cc}(\msc{E}^{\ot},\msc{D}^{\ot})
\]
from Section \ref{sect:fun_cc} already parametrizes maps of symmetric monoidal $\infty$-categories.
\end{remark}

\begin{remark}
Given any symmetric monoidal $\infty$-category $\msc{E}^{\ot}\to \Fin_{\ast}$, we note that the fiber over the empty set $\msc{E}^{\ot}_{\emptyset}$ is a contractible Kan complex.  This follows by terminality of the empty product in $\sCat_{\infty}$ and the requisite equivalence \eqref{eq:849}. Hence the transport functor along the unique active map $\bar{m}:\{1,2\}\to \{0\}$ provides a homotopically associative, and homotopically commutative product functor
\[
m_{\msc{E}}=\bar{m}_!:\msc{E}\times \msc{E}\cong \msc{E}^{\ot}_{\{1,2\}}\to \msc{E}
\]
and transport along the initial map $\bar{u}:\emptyset\to \{0\}$ provides a unit for this product
\[
u_{\msc{E}}=\bar{u}_!:\ast\cong\msc{E}^{\ot}_{\emptyset}\to \msc{E}.
\]

For any map of symmetric monoidal $\infty$-categories $F^{\ot}:\msc{E}^{\ot}\to \msc{D}^{\ot}$, compatibility between transport functors ensures the existence of homotopy commuting diagrams
\[
\xymatrix{
\msc{E}\times\msc{E}\ar[r]^{F^2}\ar[d]_{m_{\msc{E}}} & \msc{D}\times\msc{D}\ar[d]^{m_{\msc{D}}} & \ar@{}[d]|{\text{\normalsize and}} & & \ast\ar[dl]_{u_{\msc{E}}}\ar[dr]^{u_{\msc{D}}}\\
\msc{E}\ar[r]_{F} & \msc{D} & & \msc{E}\ar[rr]_F & & \msc{D}.
}
\]
The existence of higher compatibilities, up to arbitrary depth, is articulated by the corresponding transformation between the functors $F_{\msc{E}^{\ot}},F_{\msc{D}^{\ot}}:\Fin_{\ast}\to \sCat_{\infty}$ under straightening.
\end{remark}

As one expects, the homotopy truncation of any symmetric monoidal $\infty$-category remains symmetric monoidal.

\begin{lemma}\label{lem:h_symm}
For any symmetric monoidal $\infty$-category $p:\msc{E}^{\ot}\to \Fin_{\ast}$, the induced map from the homotopy category $\opn{h}p:\opn{h}\msc{E}^{\ot}\to \Fin_{\ast}$ is also a symmetric monoidal $\infty$-category. We have the fibers $(\opn{h}\msc{E}^{\ot})_I=\opn{h}(\msc{E}^{\ot}_I)$, and the transport functors $\bar{f}_!:\opn{h}\msc{E}^{\ot}_I\to \opn{h}\msc{E}^{\ot}_J$ are induced by the transport functors for $\msc{E}^{\ot}$.
\end{lemma}

\begin{proof}
It is immediate to see that every $p$-cocartesian edge in $\msc{E}^{\ot}$ is $\opn{h}p$-cocartesian in $\opn{h}\msc{E}^{\ot}$.  Hence the fibration $\opn{h}\msc{E}^{\ot}\to \Fin_{\ast}$ is cocartesian and the claim about transport follows as well.  Since transport provides an equivalence $\bar{\rho}_!:\msc{E}^{\ot}_I\overset{\sim}\to \msc{E}^I$ for each finite set $I$, by assumption, the induced map on homotopy categories $\opn{h}\msc{E}^{\ot}_I\to \opn{h}\msc{E}^I$ is also an equivalence.
\end{proof}

One obtains examples of symmetric monoidal $\infty$-categories immediately from both the discrete and simplicial settings.  We provide details below.

\subsection{Symmetric monoidal simplicial categories}

A monoidal simplicial category is a simplicial category $\underline{E}$ equipped with a simplicial functor
\[
\ot:\underline{E}\times\underline{E}\to \underline{E}
\]
along with a unit $\1$, unit isomorphisms, and associativity.  A symmetric structure on a simplicial monoidal category is a choice of a natural isomorphism $\ot\overset{\sim}\to \ot^{op}$ which is an involution and satisfies the appropriate compatibilities with the unit and associator.

\begin{definition}\label{def:e_ot}
For a symmetric monoidal simplicial category $\underline{E}$ we define a simplicial category $\underline{E}^{\ot}$ whose objects are given by pairs of a finite set $I$ \emph{equipped with a linear ordering} and an object $x$ in the product $\underline{E}^I$.  The mapping complex between objects $x$ and $y$ over finite sets $I$ and $J$, with specified orderings, is the disjoint union of the complexes
\begin{equation}\label{eq:128}
\uHom_{\underline{E}^{\ot}}(x,y)_{\bar{f}}=\prod_{j\in J}\uHom_{\underline{E}}(x_{I_j},y_j)
\end{equation}
over all $\bar{f}:I\to J$ in $\Fun_\ast$.
\end{definition}

Here each $x_{I_j}$ denotes the product
\[
x_{i_1}\ot\cdots \ot x_{i_r}:=x_{i_1}\ot(x_{i_2}\ot\cdots(x_{i_{r-1}}\ot x_{i_r})\cdots),
\]
where $I_j=\{i_1,\dots,i_r\}$ is ordered according to the ambient ordering on $I$. In the case when $I_j$ is empty we take $x_{I_j}=\1_{E}$. We note that the map of pointed sets $\bar{f}:I\to J$ underlying a map in $\underline{E}^{\ot}$ is \emph{not} required to respect the given orderings, even though the sets which parametrize objects in $\underline{E}^{\ot}$ are ordered.
\par

Composition is defined in the expected way, though it is a slightly delicate thing to write down.  Consider $n$-simplices
\[
\sigma:\Delta^n\to \uHom_{\underline{E}^{\ot}}(x,y)_{\bar{f}}\ \ \text{and}\ \ \sigma':\Delta^n\to \uHom_{\underline{E}^{\ot}}(y,z)_{\bar{g}}
\]
over maps $\bar{f}:I\to J$ and $\bar{g}:J\to K$ in $\Fin_{\ast}$. These simplices separate into tuples $\sigma=\{\sigma_j:j\in J\}$ and $\sigma'=\{\sigma'_k:k\in K\}$ relative to the product decompositions \eqref{eq:128}.  Then, for each $k$ in $K$, we have the corresponding simplices
\begin{equation}\label{eq:composition}
(\sigma'\sigma)_k:\Delta^n\to \uHom_{\underline{E}^{\ot}}(\otimes_{j\in J_k}x_{I_j},z_k)\overset{\opn{sym}^{\ast}}\longrightarrow \uHom_{\underline{E}^{\ot}}(x_{I_k},z_k)
\end{equation}
where the first map above is the composite simplex
\[
\sigma'_k(\bigotimes_{j\in J_k}\sigma_j):\Delta^n\to \uHom_{\underline{E}^{\ot}}(\otimes_{j\in J_k}x_{I_j},z_k)
\]
and $\opn{sym}^{\ast}$ is precomposition with the symmetry isomorphism $x_{I_k}\overset{\sim}\to \otimes_{j\in J_k}x_{I_j}$.  The composite simplex $\sigma'\sigma$ is, finally, the map in $\uHom_{\underline{E}^{\ot}}(x,z)_{\bar{g}\bar{f}}$ specified by the given tuple $\{(\sigma'\sigma)_k:k\in K\}$.
\par

The following is apparent.

\begin{lemma}
If $\underline{E}$ is a fibrant symmetric monoidal simplicial category, then the associated simplicial category $\underline{E}^{\ot}$ is also fibrant.
\end{lemma}

By construction the simplicial category $\underline{E}^{\ot}$ comes equipped with a forgetful functor $p:\underline{E}^{\ot}\to \Fin_{\ast}$ which sends each object $(x,I)$ to $I$, and collapses each component in the mapping complex $\uHom_{\underline{E}^{\ot}}(x,y)_{\bar{f}}$ to $\bar{f}$.

\begin{lemma}\label{lem:466}
Given a finite set $I$, any choice of a linear ordering on $I$ provides an injective equivalence $\underline{E}^I\overset{\sim}\to (\underline{E}^{\ot})_I$ to the fiber of $p:\underline{E}^{\ot}\to \Fin_{\ast}$ over $I$.
\end{lemma}

Here by an equivalence we mean a simplicial functor which is essentially surjective on the underlying plain categories, say, and which induces homotopy equivalences on the mapping complexes.

\begin{proof}
The inclusion just takes a tuple $x$ in $\underline{E}^I$ to the object $(x,I)$, where $I$ is given its specified ordering.  The induced map on Hom complexes
\[
\uHom_{\underline{E}^{I}}(x,y)=\prod_{i\in I}\uHom_{\underline{E}}(x_i,y_i)\to \uHom_{\underline{E}^{\ot}}(x,y)_{id_I} 
\]
is simply an equality.  Hence the functor $\underline{E}^I\to (\underline{E}^{\ot})_I$ is fully faithful.  The identities on both $I$ and $x$ provide isomorphisms between $(x,I)$ and $x$ paired with any alternate ordering on $I$.  So the inclusion is essentially surjective as well.
\end{proof}

We note that the fibers of the nerve are the nerve of the fibers, $\opn{N}((\underline{E}^{\ot})_I)=\opn{N}(\underline{E}^{\ot})_I$.  Hence, taking the nerve of the simplicial equivalence from Lemma \ref{lem:466} provides an equivalence of $\infty$-categories $\msc{E}^I\overset{\sim}\to (\msc{E}^{\ot})_I$.

\begin{proposition}\label{prop:symm_infty}
For a fibrant symmetric monoidal simplicial category $\underline{E}$, and $\msc{E}^{\ot}=\opn{N}(\underline{E}^{\ot})$, the functor $\opn{N}(p):\msc{E}^{\ot}\to \Fin_\ast$ is a cocartesian fibration.  Furthermore, given any map of finite pointed sets $\bar{f}:I\to J$, choice of object $x$ in $(\msc{E}^{\ot})_I$, and $x_{\bar{f}}=(x_{I_j}:j\in J)$ in $(\msc{E}^{\ot})_J$, the morphism
\[
i_{x,\bar{f}}\in \Hom_{\underline{E}^{\ot}}(x,x_{\bar{f}})_{\bar{f}}=\prod_{j\in J}\Hom_{\underline{E}_1}(x_{I_j},x_{I_j})
\]
specified by the identities $id_{x_{I_j}}$ is a $\opn{N}(p)$-cocartesian edge over $\bar{f}$.
\end{proposition}

\begin{proof}
Here $x$ is, precisely, an object in $\msc{E}^I$ along with a choice of ordering on $I$.  The products $x_{I_j}$ are all taken with respect to this ordering, exactly as above.
\par

It suffices to prove that, for any map of finite pointed sets $\bar{f}:I\to J$ and $x$ in $\msc{E}^I$, the proposed edge $i_{x,\bar{f}}:x\to x_{\bar{f}}$ is cocartesian over $\bar{f}$.  For this it suffices to show that, given any $z\in \msc{E}^{\ot}$ and $K=p(z)$, the diagram
\[
\xymatrix{
\Maps_{\msc{E}^{\ot}}(x_{\bar{f}},z)\ar[rr]^{-\circ i_{x,\bar{f}}}\ar[d] & & \Maps_{\msc{E}^{\ot}}(x,z)\ar[d]\\
\Maps_{\Fin_{\ast}}(J,K)\ar[rr]^{-\circ \bar{f}} & & \Maps_{\Fin_\ast}(I,K)
}
\]
is a pullback square in $\sKan$ \cite[\href{https://kerodon.net/tag/01TL}{01TL}]{kerodon}.  Since the spaces in the bottom row are discrete it suffices further to show that, for each $\bar{g}:J\to K$ and $\bar{h}=\bar{g}\bar{f}$, the map on the fibers
\begin{equation}\label{eq:205}
-\circ i_{x,\bar{f}}:\Maps_{\msc{E}^{\ot}}(x_{\bar{f}},z)_{\bar{g}}\to \Maps_{\msc{E}^{\ot}}(x,z)_{\bar{h}}
\end{equation}
is a homotopy equivalence.
\par

Whether or not the map \eqref{eq:205} is an equivalence can be checked at the level of the homotopy category $\opn{h}\sKan$, at which point one can replace the mapping spaces for $\msc{E}^{\ot}$ with the mapping complexes for $\underline{E}^{\ot}$ \cite[\href{https://kerodon.net/tag/02LN}{02LN}]{kerodon}. We now ask if the map
\[
-\circ i_{x,\bar{f}}:\uHom_{\underline{E}^{\ot}}(x_{\bar{f}},z)_{\bar{g}}\to \uHom_{\underline{E}^{\ot}}(x,z)_{\bar{h}}
\]
is an equivalence.  However, since the symmetry transformations on $\underline{E}$ are isomorphisms, one can check directly from the defintion of composition \eqref{eq:composition} that this map is simply an isomorphism of simplicial sets.
\end{proof}

\begin{remark}
Given $p:\underline{E}^{\ot}\to \Fin_{\ast}$ as in the statement of Proposition \ref{prop:symm_infty}, we generally let $p:\msc{E}^{\ot}\to \Fin_{\ast}$ denote the induced cocartesian fibration from the homotopy coherent nerve $\msc{E}^{\ot}=\opn{N}(\underline{E}^{\ot})$, by an abuse of notation.
\end{remark}

By the above description of cocartesian edges in $\msc{E}^{\ot}$ we see that the (nerve of the) symmetric product functor
\[
\xymatrixrowsep{3mm}
\xymatrix{
 & (\msc{E}^{\ot})_I\ar[dr]^{\opn{transport}} \\
\msc{E}^I\ar[ur]^{\sim}\ar[rr]_{\opn{product}}& & \msc{E},
}
\]
associated to any choice of ordering on $I$, is a transport functor along the unique active map $I\to \{0\}$.  The transport functors along the various inert maps $\rho_i:I\to \{0\}$ are the projections $p_i:\msc{E}^I\cong (\msc{E}^{\ot})_I\to \msc{E}$.  Hence we observe the following.

\begin{corollary}\label{cor:simp_sym_mon}
Consider a fibrant symmetric monoidal simplicial category $\underline{E}$ with associated $\infty$-category $\msc{E}$.  The cocartesian fibration $p:\msc{E}^{\ot}\to \Fin_\ast$ from Proposition \ref{prop:symm_infty} gives $\msc{E}$ the structure of symmetric monoidal $\infty$-category.
\end{corollary}

\subsection{Maps of symmetric monoidal simplicial categories}
\label{sect:underline_maps}

For any map of symmetric monoidal simplicial categories $F:\underline{E}_0\to \underline{E}_1$ we have the corresponding functor $F^{\ot}:\underline{E}_0^{\ot}\to \underline{E}_1^{\ot}$ between simplicial categories over $\Fin_\ast$.  The functor $F^{\ot}$ is defined on objects via the exponents $F^I:\underline{E}^I_0\to \underline{E}^I_1$ and on Hom complexes via the maps
\[
\begin{array}{l}
\uHom_{\underline{E}_0^{\ot}}(x,y)_{\bar{f}}=\prod_{j\in J}\uHom_{\underline{E}_0}(x_{I_j},y_j)\vspace{2mm}\\
\hspace{3cm}\longrightarrow
\prod_{j\in J}\uHom_{\underline{E}_1}(F(x)_{I_j},F(y_j))=\uHom_{\underline{E}_1^{\ot}}(F^I(x),F^J(y))_{\bar{f}}
\end{array}
\]
provided by applying $F$ then the tensor compatibilities $F(x)_{I_j}\overset{\sim}\to F(x_{I_j})$.  Taking homotopy coherent nerves produces a functor $F^{\ot}:\msc{E}_0^{\ot}\to\msc{E}_1^{\ot}$ between $\infty$-categories over $\Fin_\ast$, where we abuse notation to write simply $F^{\ot}$ for $\opn{N}(F^{\ot})$.

\begin{proposition}\label{prop:1096}
For any map of fibrant symmetric monoidal simplicial categories $F:\underline{E}_0\to \underline{E}_1$ the corresponding functor
\[
\xymatrix{
\msc{E}_0^{\ot}\ar[rr]^{F^{\ot}}\ar[dr]_{p_0} & & \msc{E}_1^{\ot}\ar[dl]^{p_1}\\
	& \Fin_\ast & 
} 
\]
is a map of symmetric monoidal $\infty$-categories.
\end{proposition}

\begin{proof}
For each morphism $\bar{f}:I\to J$ in $\opn{Fin}_{\ast}$, $x$ in $(\msc{E}_0^{\ot})_I$, and index $j\in J$ we have the diagram
\[
\xymatrix{
F(x_{i_1}\ot\dots\ot x_{i_r})\ar@{-->}[rr]^{F(id)=id} & & F(x_{I_j})\\
F(x_{i_1})\ot\dots\ot F(x_{i_r})\ar[rr]_(.65){(i_{F^{\ot}(x),\bar{f}})_j}\ar@{-->}[u]^{\cong}\ar[urr]|{F^{\ot}(i_{x,\bar{f}})_j} & & F(x)_{I_j}\ar[u]_{\cong}
}
\]
in $\msc{E}_1$.  These diagrams yields a diagram
\[
\xymatrix{
	& F(x)_{\bar{f}}\ar[dr]^{\cong}\\
F^{\ot}(x)\ar[rr]_{F^{\ot}(i_{x,\bar{f}})}\ar[ur]^{i_{F(x),\bar{f}}} & & F(x_{\bar{f}})
}
\]
in $\msc{E}^{\ot}_1$ in which the map $i_{F(x),\bar{f}}$ is $p_1$-cocartesian.  Since all isomorphisms are $p_1$-cocartesian, and cocartesian edges satisfy the $2$-of-$3$ property \cite[\href{https://kerodon.net/tag/01TS}{01TS}]{kerodon}, it follows that $F^{\ot}(i_{x,\bar{f}})$ is $p_1$-cocartesian.  Hence $F^{\ot}$ preserves cocartesian edges, and is thus a map of symmetric monoidal $\infty$-categories.
\end{proof}

\subsection{Fiber products and transfer of symmetric monoidal structures}

Throughout the text we construct symmetric monoidal $\infty$-categories by combining data from, say, a topological source and an algebraic source.  Formally such categories arise as fiber products of symmetric monoidal $\infty$-categories along pairs of symmetric monoidal functors. To address such situations we have the following.

\begin{proposition}\label{prop:fp_symm}
Given a pair of symmetric monoidal functors $F^{\ot}_i:\msc{E}^{\ot}_i\to \msc{T}^{\ot}$ in which one of the $F_i$ is an isofibration, the fiber product $\msc{E}^{\ot}_0\times_{\msc{T}^{\ot}}\msc{E}^{\ot}_1$ inherits a unique symmetric monoidal structure under which the two projections
\[
\msc{E}^{\ot}_0\times_{\msc{T}^{\ot}}\msc{E}^{\ot}_1\to \msc{E}^{\ot}_i
\]
are symmetric monoidal functors
\end{proposition}

We omit the details, but indicate some method of proof for the interested reader.

\begin{proof}[Idea of proof]
This result can be proved by hand, or can be deduced from the fact that the forgetful functor from $\SM_{\infty}$ (Defintion \ref{def:sm_infty}) to $\sCat_{\infty}$ commutes with pullbacks \cite[Corollary 3.2.2.4]{ha}.
\end{proof}

We also note that symmetric monoidality is stable under equivalence of isofibrations.

\begin{lemma}\label{lem:sym_transf}
Consider a diagram
\[
\xymatrix{
\msc{C}\ar[rr]^F\ar[dr]_p & & \msc{C}'\ar[dl]^{p'}\\
	& \Fin_{\ast}
}
\]
in which both $p$ and $p'$ are isofibrations, and $F$ is an equivalence.  If one of $p$ or $p'$ realizes $\msc{C}$ or $\msc{C}'$ as a symmetric monoidal $\infty$-category, then both $\msc{C}$ and $\msc{C}'$ are symmetric monoidal and $F$ is an equivalence of symmetric monoidal $\infty$-categories.
\end{lemma}

\begin{proof}
In this case $F$ is an equivalence of isofibrations over $\Fin_{\ast}$ \cite[\href{https://kerodon.net/tag/0285}{0295}]{kerodon}, and hence both $p$ and $p'$ are cocartesian fibrations whenever one is \cite[\href{https://kerodon.net/tag/028A}{028A}]{kerodon}.  Furthermore $F$ preserves cocartesian edges by \cite[\href{https://kerodon.net/tag/028B}{028B}]{kerodon}, and for any finite set $I$ we have a commuting diagram
\[
\xymatrix{
\msc{C}_I\ar[rr]^{F_I}\ar[d] & & \msc{C}'_I\ar[d]\\
\prod_{i\in I}\msc{C}_{\{0\}}\ar[rr]_{\prod_iF_{\{0\}}} & & \prod_{i\in I}\msc{C}'_{\{0\}},
}
\]
in $\msc{C}\!at_{\infty}$ obtained from transport along the inert maps $\bar{\rho}_i:I\to \{0\}$.  It follows that the map $\msc{C}_I\to \prod_i\msc{C}_{\{0\}}$ is an equivalence if and only if $\msc{C}'_I\to \prod_i\msc{C}'_{\{0\}}$ is an equivalence, and hence both transport functors induce equivalences by our hypothesis.  Taken together, we see that both $\msc{C}$ and $\msc{C}'$ are symmetric monoidal $\infty$-categories and that $F$ is an equivalence of symmetric monoidal $\infty$-categories.
\end{proof}

\subsection{Equifibered $\infty$-categories}

Below we consider the symmetric monoidal category of finite sets $\Fin$, under disjoint unions, and have the associated fibration $\Fin^{\ot}\to \Fin_{\ast}$ from Proposition \ref{prop:symm_infty}.

\begin{definition}[{\cite[Example 4.3.1]{barkanetal25}}]\label{def:equifib}
A symmetric monoidal $\infty$-category over $\Fin^{\ot}$, i.e.\ a symmetric monoidal $\infty$-category with a fixed symmetric monoidal functor $\pi:\msc{T}^{\ot}\to \opn{Fin}^{\ot}$, is said to be equifibered over $\Fin$ if the structure map $\pi$ is an isofibration and the diagram
\begin{equation}\label{eq:prod_pull}
\xymatrix{
\msc{T}\times \msc{T}\ar[r]\ar[d] & \msc{T}\ar[d]\\
\Fin\times \Fin\ar[r] & \Fin
}
\end{equation}
for the product functors is a pullback diagram in $\sCat_{\infty}$. A map of equifibered symmetric monoidal $\infty$-categories is a strictly commuting diagram of symmetric monoidal functors
\[
\xymatrix{
\msc{S}^{\ot}\ar[rr]\ar[dr] & & \msc{T}^{\ot}\ar[dl]\\
	& \Fin^{\ot} & .
}
\]
\end{definition}

By an \emph{equifibered $\infty$-category} over $\Fin$ we always mean a symmetric monoidal $\infty$-category which is equifibered over $\Fin$. Intuitively, a symmetric monoidal $\infty$-category $\msc{T}^{\ot}\to \Fin_{\ast}$ admits an equifibered $\infty$-category structure over $\Fin$ provided objects in the underlying $\infty$-category $\msc{T}$ are not too large, and the monoidal product on $\msc{T}$ is given by a type of disjoint union.

\begin{example}
The category $\opn{Top}^c$ of compact topological spaces with disjoint union is \emph{almost} equifibered over $\Fin$, via the connected components functor $\pi_0:\opn{Top}^c\to \Fin$. In particular, the product diagram \eqref{eq:prod_pull} is a pullback diagram, but the structure map $\opn{Top}^c\to \Fin$ fails to be an isofibration. Hence, after taking an isofibrant replacement $\opn{Top}^c\overset{\sim}\to \opn{T}^c\to \Fin$ we find that $\opn{T}^c$ is equifibered over $\Fin$. Such a replacement can be constructed by considering, for example, the category of pairs $(X,\alpha)$ where $X$ is a compact topological space and $\alpha$ is a bijection $\alpha:\pi_0(X)\overset{\sim}\to I$ from the connected components to a finite set.
\par

Similarly, the $\infty$-category $\msc{K}\!an^c$ of compact Kan complexes with disjoint union becomes an equifibered $\infty$-category over $\Fin$, via the connected components functor, after taking an isofibrant replacement.
\end{example}

\begin{example}
For any algebraic group $G$ over a field $k$, one can check that the category $\opn{Rep}(G)$ of $G$-representations admits no equifibered structure over $\Fin$.
\end{example}

\begin{definition}\label{def:T_mon}
Fix an equifibered $\infty$-category over $\Fin$, $\pi:\msc{T}^{\ot}\to \Fin^{\ot}$. A $\msc{T}$-monoidal $\infty$-category is an $\infty$-category equipped with a cocartesian fibration $q:\msc{E}^{\ot}\to\msc{T}^{\ot}$ which satisfies the following:
\begin{itemize}
\item For each object $x$ in the fiber $\msc{T}^{\ot}_K$ over a given set $K$ in $\opn{Fin}_{\ast}$, and inert projections $\rho_k:x\to x_k$ over the projections $\bar{\rho}_k:K\to \{0\}$, transport provides an equivalence
\[
\rho_!:\msc{E}^{\ot}_x\overset{\sim}\to \prod_{k\in K} \msc{E}^{\ot}_{x_k}.
\]\vspace{1mm}
\item For each object $x$ in the fiber $\msc{T}^{\ot}_K$, and cocartesian edge $f:x\to x'$ over the unique active map $\bar{f}:K\to \{0\}$, transport provides an equivalence
\[
f_!:\msc{E}^{\ot}_x\overset{\sim}\to \msc{E}^{\ot}_{x'}.
\]
\end{itemize}
A map of $\msc{T}$-monoidal $\infty$-categories is a map of cocartesian fibrations over $\msc{T}^{\ot}$.
\end{definition}

We note that any $\msc{T}$-monoidal $\infty$-category $\msc{E}^{\ot}$ comes equipped with a cocartesian fibration to $\Fin_{\ast}$, via composition. One can check, using \cite[\href{https://kerodon.net/tag/023M}{023M}]{kerodon} and the first condition above, that the induced map $\msc{E}^{\ot}\to \Fin_{\ast}$ gives $\msc{E}^{\ot}$ the structure of a symmetric monoidal $\infty$-category and, subsequently, that $q:\msc{E}^{\ot}\to \msc{T}^{\ot}$ is a symmetric monoidal functor. We compose further to observe a symmetric monoidal functor $\msc{E}^{\ot}\to \opn{Fin}^{\ot}$.

\begin{lemma}
If $\msc{T}^{\ot}\to \Fin_{\ast}$ is equifibered over $\opn{Fin}$, then any $\msc{T}$-monoidal $\infty$-category $q:\msc{E}^{\ot}\to \msc{T}^{\ot}$ is equifibered over $\opn{Fin}$ as well.
\end{lemma}

\begin{proof}
It suffices to show that the diagram
\[
\xymatrix{
\msc{E}^{\ot}_{\{0,1\}}\ar[rr]^{product}\ar[d]_{q_{\{0,1\}}} & & \msc{E}\ar[d]^{q_{\{0\}}}\\
\msc{T}^{\ot}_{\{0,1\}}\ar[rr]_{product} & & \msc{T}
}
\]
is a pullback diagram in $\sCat_{\infty}$, or equivalently that the induced map $\msc{E}^{\ot}_{\{0,1\}}\to \msc{E}\times_{\msc{T}}\msc{T}^{\ot}_{\{0,1\}}$ is an equivalence.
\par

As we explain at Lemma \ref{lem:3779} below, the product functor for $\msc{E}^{\ot}$ sends $q_{\{0,1\}}$-cocartesian edges to $q_{\{0\}}$-cocartesian edges, and hence the map $\msc{E}^{\ot}_{\{0,1\}}\to \msc{E}\times_{\msc{T}}\msc{T}^{\ot}_{\{0,1\}}$ is a map of cocartesian fibrations over $\msc{T}^{\ot}_{\{0,1\}}$. This map is therefore an equivalence if and only if its fibers over $\msc{T}^{\ot}_{\{0,1\}}$ are equivalences \cite[\href{https://kerodon.net/tag/028B}{028B}]{kerodon}. Taking the fiber over a given object $x$ in $\msc{T}^{\ot}_{\{0,1\}}$, with image $x'$ in $\msc{T}$, however, simply recovers the equivalence $\msc{E}^{\ot}_x\to \msc{E}^{\ot}_{x'}$ promised in the second point of Definition \ref{def:T_mon}.
\end{proof}

\subsection{Transport equivalences for $\msc{T}$-monoidal $\infty$-categories}
\label{sect:rho}

Fix a symmetric monoidal $\infty$-category $\msc{T}^{\ot}\to \Fin_{\ast}$ which is equifibered over $\opn{Fin}$, and a $\msc{T}$-monoidal $\infty$-category $q:\msc{E}^{\ot}\to \msc{T}^{\ot}$.

\begin{definition}\label{def:separation}
We call an object $(I_k:k\in K)$ in $\opn{Fin}^{\ot}$ separated if each set $I_k$ is a (nonempty) singleton. Given a symmetric monoidal $\infty$-category $p:\msc{T}^{\ot}\to \Fin_{\ast}$ which is equifibered over $\Fin$, via a structure map $\pi:\msc{T}^{\ot}\to \opn{Fin}^{\ot}$, an object $x$ in $\msc{T}^{\ot}$ is called separated if its image $\pi(x)$ is separated in $\opn{Fin}^{\ot}$. A map $\alpha:\tilde{x}\to x$ in $\msc{T}^{\ot}$ is called a separation of $x$ whenever $\tilde{x}$ is separated and $\alpha$ is a $p$-cocartesian edge.
\end{definition}

In $\Fin^{\ot}$ for example, a separation of an object $x=(I_k:k\in K)$ is provided by the object $\tilde{x}=(\{i\}: i\in \amalg_k I_k)$ with the apparent map $\tilde{x}\to x$. It is clear that the separation of $x$ is uniquely determined up to unique isomorphism in the overcategory $\Fin^{\ot}_{/x}$. Similarly, equifiberedness of $\msc{T}^{\ot}$ over $\Fin$ implies uniqueness of separating morphisms in $\msc{T}^{\ot}$.
\par

For a fixed object $x$ in $\msc{T}^{\ot}$ take a separation $\alpha:\tilde{x}\to x$ and a cocartesian edge $x\to x'$ over the active map $K\to \{0\}$ in $\opn{Fin}_{\ast}$. Then the composite $\tilde{x}\to x'$ is also cocartesian, and taking transport for our $\msc{T}$-monoidal $\infty$-category $\msc{E}^{\ot}$ now produces a diagram
\[
\xymatrix{
\msc{E}^{\ot}_{\tilde{x}}\ar[dr]\ar[rr]^{\alpha_!} & & \msc{E}^{\ot}_x\ar[dl]\\
	& \msc{E}_{x'}^{\ot}
}
\]
in $\sCat_{\infty}$ in which both functors to $\msc{E}_{x'}^{\ot}$ are equivalences, by the definiton of $\msc{T}$-monoidality. It follows that $\alpha_!:\msc{E}^{\ot}_{\tilde{x}}\to \msc{E}_x^{\ot}$ is an equivalence as well.
\par

The separating object $\tilde{x}$ is now over a set $I$ in $\opn{Fin}_{\ast}$, and we have the inert projections $\rho_i:\tilde{x}\to x_i$ in $\msc{T}^{\ot}$ which induce an equivalence $\msc{E}^{\ot}_{\tilde{x}}\overset{\sim}\to \prod_{i\in I}\msc{E}_{x_i}^{\ot}$. We therefore obtain a uniquely determined equivalences from $\msc{E}^{\ot}_x$,
\begin{equation}\label{eq:1206}
\xymatrix{
	& \msc{E}_{\tilde{x}}^{\ot}\ar@{..>}[dr]^{\sim}\ar@{..>}[dl]_{\sim}\\
\msc{E}^{\ot}_x\ar[rr]^{\sim}_{\exists!} & & \prod_i\msc{E}^{\ot}_{x_i}.
}
\end{equation}

\begin{definition}
Let $\msc{T}^{\ot}\to \Fin_{\ast}$ be equifibered $\Fin$. Given a $\msc{T}$-monoidal $\infty$-category $\msc{E}^{\ot}\to \msc{T}^{\ot}$ and an object $x$ in $\msc{T}^{\ot}$, the completing equivalence from \eqref{eq:1206} is referred to as the transport equivalence, and is generally denoted
\[
\rho_!:\msc{E}^{\ot}_x\overset{\sim}\to \prod_{i\in I}\msc{E}^{\ot}_{x_i}
\]
by an abuse of notation.
\end{definition}

\subsection{Language and notation for monoidal $\infty$-categories}

Given an $\infty$-category $\msc{E}$, we can speak of a symmetric monoidal \emph{structure} on $\msc{E}$. This is the choice of a symmetric monoidal $\infty$-category $\msc{E}^{\ot}\to \Fin_{\ast}$ with $\msc{E}=\msc{E}^{\ot}_{\{0\}}$, or more generally with a specified isomorphism of simplicial sets $\msc{E}\cong \msc{E}^{\ot}_{\{0\}}$. Having specified a symmetric monoidal structure on a given $\infty$-category $\msc{E}$, we may also just refer to $\msc{E}$ itself as a symmetric monoidal $\infty$-category, by an abuse of language.
\par

Given a symmetric monoidal $\infty$-category $\msc{T}^{\ot}\to \Fin_{\ast}$, we say $\msc{T}$ is equifibered over $\Fin$ if $\msc{T}^{\ot}\to \Fin_{\ast}$ is a symmetric monoidal $\infty$-category which is equifibered over $\Fin$.

Given $\infty$-categories $\msc{E}_1$ and $\msc{E}_2$ with specified symmetric monoidal structures, we can also speak of a symmetric monoidal structure on a given functor $F:\msc{E}_1\to \msc{E}_2$. This is the choice of a map $F^{\ot}:\msc{E}^{\ot}_1\to \msc{E}^{\ot}_2$ between cocartesian fibrations over $\Fin_{\ast}$ which fits into a diagram
\[
\xymatrix{
\msc{E}_1\ar[rr]^F\ar[d]_{include} & & \msc{E}_2\ar[d]^{include}\\
\msc{E}_1^{\ot}\ar[rr]_{F^{\ot}} & & \msc{E}_2^{\ot}.
}
\]

\section{$\infty$-Categories of disk embeddings}
\label{sect:fr_disk}

We construct a symmetric monoidal $\infty$-category $fr\Disk$ of $\pm$-colored framed disks, and disk embeddings. Besides our addition of $\pm$-colors, many of the constructions and results of this section are standard. See for example \cite[Remark 2.10]{ayalafrancis15}.

Following these introductory materials, we show in Section \ref{sect:disk_monoidal} that $fr\Disk$ is equifibered over $\Fin$. Throughout this study we are generically interested in $fr\Disk$-monoidal $\infty$-categories as a coarse approximation to ribbon tensor categories, in the homotopical setting. In particular, we show in Section \ref{sect:ribbon_framed} below that any ribbon tensor (1-)category admits the structure of a $fr\Disk$-monoidal $\infty$-category, and that the original ribbon structure can be recovered from this structure over $fr\Disk$.

\subsection{Colored sets}
\label{sect:col_sets}

A finite colored set is the pairing $\mcl{I}$ of a finite set $I$ with a morphism of sets $\opn{color}_{\mcl{I}}:I\to \{\pm\}$.  Maps of colored sets $f:\mcl{I}\to \mcl{J}$ are maps of the underling uncolored sets $f:I\to J$ which preserve positively colored elements.  (So, we are allowed to send a negative index $i$ to a positive index $j$, but not vise-versa.)  We let $\mbf{Fin}$ denote the category of finite colored sets.
\par

When we write $i\in \mcl{I}$ for a colored set $\mcl{I}$, we mean that $i$ is an element in the underlying uncolored set $I$. As in the uncolored setting, given a map $f:\mcl{I}\to \mcl{J}$ in $\mbf{Fin}$ we let $\mcl{I}_j$ denote the the preimage $f^{-1}(j)\subset \mcl{I}$ along with its induced coloring. We have the analogously defined category $\mbf{Fin}_{\ast}$ of finite \emph{pointed} colored sets.

\subsection{The $\infty$-category of colored disk embeddings}
\label{sect:colored_disks}

We construct a symmetric monoidal simplicial category $fr\uDisk$ of framed colored disks.  Objects are finite colored sets $\mcl{I}$, and morphisms $f:\mcl{I}\to \mcl{J}$ (in the underlying discrete category) are oriented embedding of disks
\begin{equation}\label{eq:1260}
f:D\times \mcl{I}\to D\times \mcl{J}.
\end{equation}
We assume that $f$ sends positively colored disks in $D\times \mcl{I}$ to positively colored disks in $D\times \mcl{J}$, so that $f$ recovers a map in $\mbf{Fin}$ on connected components, and we take specifically $D$ to be the standard disk $D=\{z\in \mbb{R}^2:|z|\leq 1\}$.  The mapping complexes are the subspaces
\[
\uHom_{fr\uDisk}(\mcl{I},\mcl{J})=\opn{Embed}(D\times \mcl{I},D\times \mcl{J})\ \subseteq\ \Maps(D\times \mcl{I},D\times \mcl{J})
\]
whose $n$-simplices consist of those continuous function
\[
\sigma:|\Delta^n|\times D\times \mcl{I}\to D\times \mcl{J}
\]
which evaluate at each $t\in |\Delta^n|$ to a smooth oriented embedding $\sigma_t:D\times \mcl{I}\to D\times \mcl{J}$ which recovers a map in $\mbf{Fin}$ on connected components.

\begin{remark}
One can think of negatively colored disks as copies of the standard disk in $\mbb{R}^2$ equipped with the negated orientation.
\end{remark}

The composition functions
\[
\circ:\uHom_{fr\uDisk}(\mcl{J},\mcl{K})\times\uHom_{fr\uDisk}(\mcl{I},\mcl{J})\to \uHom_{fr\uDisk}(\mcl{I},\mcl{J})
\]
are defined by taking a pair of $n$-simplices $\sigma:|\Delta^n|\times D\times \mcl{I}\to D\times \mcl{J}$ and $\sigma':|\Delta^n|\times D\times \mcl{J}\to D\times \mcl{K}$ to the usual composite
\[
|\Delta^n|\times D\times \mcl{I}\overset{\opn{diag}\times 1}\to |\Delta^n|\times |\Delta^n|\times D\times \mcl{I}\overset{1\times \sigma}\to |\Delta^n|\times D\times \mcl{J}\overset{\sigma'}\to D\times \mcl{K}.
\]

\begin{lemma}
The simplicial category $fr\uDisk$ is fibrant.
\end{lemma}

\begin{proof}
One can fill horns by noting that the geometric realization $|\Delta^n|$ retracts onto any given horn $|\Lambda^n_i|$.
\end{proof}

We have the apparent symmetric monoidal structure on $fr\uDisk$ provided by disjoint union, and taking connected components provides a symmetric monoidal functor $\pi_0:fr\uDisk\to \opn{Fin}$ to the discrete category of finite pointed sets, with disjoint union. Applying the homotopy coherent nerve therefore produces a symmetric monoidal $\infty$-category over $\Fin^{\ot}$,
\[
\pi_0^{\ot}:\opn{N}(fr\uDisk^{\ot})\to \opn{Fin}^{\ot},\ \ (\mcl{I}_k:k\in K)\mapsto (I_k:k\in K),
\]
by Proposition \ref{prop:1096}.

We give a slightly more efficient presentation of the symmetric monoidal structure on $\opn{N}(fr\Disk)$ in Section \ref{sect:snatched_disks} below, after recalling a colored variant of the framed little disk operad.

\subsection{The $\infty$-operad of colored disks}
\label{sect:fr_e2}

In addition to the symmetric monoidal $\infty$-category of framed disks we have a related disky $\infty$-operad. This is a framed and colored variant of the usual little disk operad, which we recall here.
\par

Objects in the preceding (fibrant) simplicial category $fr\underline{\E}_2$ are finite colored sets and morphisms $f:\mcl{I}\to \mcl{J}$ consist of a choice of subset $\mcl{I}_f\subseteq \mcl{I}$ along with an oriented disk embedding $f:D\times \mcl{I}_f\to D\times \mcl{J}$. For higher simplices, the mapping complexes decompose over morphisms in $\mbf{Fin}_{\ast}$,
\[
\uHom_{fr\uE_2}(\mcl{I},\mcl{J})=\coprod_{\bar{\xi}\in \Hom_{\mbf{Fin}_{\ast}}(\mcl{I},\mcl{J})}\opn{Embed}(D\times \mcl{I}_{\bar{\xi}},D\times \mcl{J})_{\bar{\xi}},
\]
where an $n$-simplex in $\opn{Embed}(D\times \mcl{I}_{\bar{\xi}},D\times\mcl{J})_{\bar{\xi}}$ is a continuous function $|\Delta^n|\times D\times \mcl{I}_{\bar{\xi}}\to D\times \mcl{J}$ which evaluates to an oriented embedding at each $t\in |\Delta^n|$, and which recovers the specified map $\bar{\xi}$ on connected components.
\par

We have the apparent forgetful functor to the discrete category of pointed colored sets $fr\underline{\E}_2\to \mbf{Fin}_{\ast}$, and can forget further to $\Fin_{\ast}$. Taking the homotopy coherent nerve, we obtain an $\infty$-category with an isofibration $fr\E_2\to \Fin_{\ast}$, and one can check that this fibration gives $fr\E_2$ the structure of an $\infty$-operad \cite[Definition 2.1.1.10]{ha}.

Now, taking disjoint unions of continuous functions
\[
\begin{array}{l}
\left(|\Delta^n|\times D\times (\amalg_{k\in K_l}\mcl{I}_k)\overset{\sigma_l}\to D\times \mcl{J}_l:l\in L\right)\mapsto\vspace{1mm}\\
\hspace{2cm}\left(\amalg_l\sigma_l:|\Delta^n|\times D\times(\amalg_{k\in K_{\bar{f}}}\mcl{I}_k)\to D\times \mcl{J}\right)
\end{array}
\]
defines a simplicial functor $fr\uDisk^{\ot}\to fr\underline{\E}_2$. We similarly define a disjoint union functor $\Fin^{\ot}\to \Fin_{\ast}$ for the symmetric monoidal category of finite sets. Via these disjoint union functors, we connect the symmetric monoidal simplicial disk category to the simplicial little disk operad via a pullback diagram
\begin{equation}\label{eq:1377}
\xymatrix{
fr\uDisk^{\ot}\ar[rr]^{\opn{disj\ union}}\ar[d] & & fr\underline{\E}_2\ar[d]\\
\Fin^{\ot}\ar[rr]_{\opn{disj\ union}} & & \Fin_{\ast}
}
\end{equation}
of simplicial categories. Taking the homotopy coherent nerve then produces a pullback diagram of $\infty$-categories.

\subsection{Snatched presentation of $fr\Disk^{\ot}$}
\label{sect:snatched_disks}

We construct a symmetric monoidal $\infty$-category $p:fr\Disk^{\ot}\to \Fin_{\ast}$ and exhibit a symmetric monoidal equivalence
\begin{equation}\label{eq:1384}
\xymatrix{
\opn{N}(fr\uDisk^{\ot})\ar[rr]^{\sim}\ar[dr] & & fr\Disk^{\ot}\ar[dl]^p\\
	& \Fin_{\ast} & .
}
\end{equation}
The point here is to provide a more user-friendly expression of the symmetric monoidal structure on the underlying $\infty$-category $fr\Disk=\opn{N}(fr\uDisk)$.

An object in $fr\Disk^{\ot}$ is a colored set $\mcl{I}$ equipped with a set map $\mu:\mcl{I}\to K$, or equivalently a continuous map $\mu:D\times \mcl{I}\to K$, where $K$ is a finite set which is taken to be completely positive if one likes. A morphism $(\mcl{I},\mu)\to (\mcl{J},\nu)$ in $fr\Disk^{\ot}$ is the pairing of a map $f:D\times \mcl{I}_f\to D\times \mcl{J}$ in $fr\E_2$ with a map $\bar{f}:K_{\bar{f}}\to L$ in $\Fin_{\ast}$ which fits into a diagram
\begin{equation}\label{eq:1419}
\xymatrix{
D\times \mcl{I}_f\ar[rr]^f\ar[d]_{\mu|_{\mcl{I}_f}} & & D\times \mcl{J}\ar[d]^{\nu}\\
K_{\bar{f}}\ar[rr]_{\bar{f}} & & L,
}
\end{equation}
and for which $I_f=\mu^{-1}K_{\bar{f}}$.
\par

In higher dimensions, an $n$-simplex $\sigma:\Delta^n\to fr\Disk^{\ot}$ is precisely the information of simplices $\sigma:\Delta^n\to fr\E_2$ and $\bar{\sigma}:\Delta^n\to \Fin_{\ast}$ which are ``linked at their vertices $\{l\}\subseteq \Delta^n$" via morphisms $\mu^l:\mcl{I}^l\to K^l$. These linking morphisms are (only) required to produce a consistency at the level of edges, in the sense that each pair of edges
\[
f_{lm}:D\times \mcl{I}_{f_{lm}}^l\to D\times \mcl{I}^m\ \ \text{and}\ \ 
\bar{f}_{lm}:K^l_{\bar{f}_{lm}}\to K^m
\]
complete diagrams as in \eqref{eq:1419}.
\par

In this alternate guise, the disjoint union functor appears as the projection
\[
\opn{disj\ union}=p_1:fr\Disk^{\ot}\to fr\E_2,\ \ \left(\sigma,\bar{\sigma},\mu\right)\to \sigma,
\]
and the other projection provides an isofibration to finite pointed sets $p=p_2:fr\Disk^{\ot}\to \Fin_{\ast}$, $(\sigma,\bar{\sigma},\mu)\mapsto \bar{\sigma}$.
\par

Similarly, for the discrete category of finite sets, we have the alternate presentation of $\opn{Fin}^{\ot}$ as the category whose objects are pairings $(I,\mu)$ of a finite set with a specified map $\mu:I\to K$, and whose morphisms consist of diagrams
\[
\xymatrix{
I\ar[r]\ar[d]_{\mu} & J\ar[d]^{\nu}\\
K\ar[r] & L.
}
\]
Under this alternate formulation, we again take connected components to obtain a functor $\pi:fr\Disk^{\ot}\to \Fin^{\ot}$ which fits into a pullback diagram
\begin{equation}\label{eq:1427}
\xymatrix{
fr\Disk^{\ot}\ar[d]_{\pi}\ar[rr]^{\opn{disj\ union}} & & fr\E_2\ar[d]\\
\Fin^{\ot}\ar[rr]_{\opn{disj\ union}} & & \Fin_{\ast}.
}
\end{equation}

\begin{lemma}\label{lem:1434}
The projection $p:fr\Disk^{\ot}\to \Fin_{\ast}$ gives $fr\Disk$ the structure of a symmetric monoidal $\infty$-category, and there is an equivalence $\opn{N}(fr\uDisk^{\ot})\overset{\sim}\to fr\Disk^{\ot}$ of symmetric monoidal $\infty$-categories over $\Fin^{\ot}$. 
\end{lemma}

\begin{proof}
Let $'\Fin^{\ot}$ denote the symmetric monoidal construction of Proposition \ref{prop:symm_infty} and $\Fin^{\ot}$ denote the alternate construction given above. We have the assignment which sends an ordered tuple of sets $(I_k:k\in K)$ to the disjoint union $\amalg_kI_k$ equipped with the unique map $\mu:\amalg_kI_k\to K$ which has $\mu^{-1}(k)=I_k$. This assignment defines a functor $'\Fin^{\ot}\to \Fin^{\ot}$ over $\Fin_{\ast}$, in both ways, which is seen to be essentially surjective and fully faithful, and thus an equivalence.

Considering the identities on $fr\E_2$ and $\Fin_{\ast}$, we now obtained a unique functor $\opn{N}(fr\uDisk^{\ot})\to fr\Disk^{\ot}$ which completes a transformation between the respective pullback diagrams $\Delta^1\times \Delta^1\to \opn{Cat}_{\infty}$ from \eqref{eq:1377} and \eqref{eq:1427}. This completing functor is an equivalence by \cite[\href{https://kerodon.net/tag/033K}{033K}]{kerodon}. The fact that the isofibration $p:fr\Disk^{\ot}\to \Fin_{\ast}$ is a cocartesian fibration which gives $fr\Disk^{\ot}$ the structure of a symmetric monoidal $\infty$-category now follows by the transfer result of Lemma \ref{lem:sym_transf}.
\end{proof}

\begin{warning}
From now on the symmetric monoidal structure on $\Fin$ is always realized via the fibration $\Fin^{\ot}\to \Fin$ of pairs $(I,\mu)$ as above, though the distinction is of little technical importance.
\end{warning}

\subsection{Linear disk arrangements}

We often consider the subcategory $Linfr\Disk^{\ot}\subseteq fr\Disk^{\ot}$ spanned by all objects, and those morphisms
\[
f:D\times\mcl{I}\to D\times \mcl{J}
\]
whose constituent disk embeddings $f_i:D_i\to D_j$ send the $x$-axis of $D_i$ into the $x$-axis of $D_j$, and have derivatives $T_0f_i:\mathbb{R}^2\to \mathbb{R}^2$ equal to the identity.  One can check that the inclusion
\[
Linfr\Disk^{\ot}\to fr\Disk^{\ot}
\]
is an equivalence of symmetric monoidal $\infty$-categories (Lemma \ref{lem:linfr_fr}).
\par

From homotopies between such linear arrangements of disks, i.e.\ from $2$-simplices in $Linfr\Disk$, one can extract some useful combinatorial data.

\subsection{Extracting framed braids}
\label{sect:braids}

From any morphism $f:(\mcl{I},\mu)\to (\mcl{J},\nu)$ in $Linfr\Disk^{\ot}$, and $j\in \mcl{J}$, the linear ordering on the $x$-axis in the $j$-th disk $D_j$ induces a linear ordering on the preimage $\mcl{I}_j\to \{j\}$. Any homotopy between two maps in $Linfr\Disk^{\ot}$, $\zeta:|\Delta^1|\times D\times \mcl{I}_f\to D\times\mcl{J}$, then traces out framed braids in the disks $D\times\mcl{J}$.  By a framed braid we mean an element in the fundamental group of the space of framed configurations of points in the disk $\opn{Br}_{fr}(K)=\pi_0(\opn{Conf}^K_{fr}(D))$, where we point the space $\opn{Conf}^K_{fr}(D)$ via a choice of an ordering on the set $K$, a map $K\to D$ which aligns the points according to the given ordering along the $x$-axis, and the points inherit their framings from the ambient framing on $D$. We note the decomposition $\opn{Br}_{fr}(K)=\opn{Br}(K)\ltimes(\mathbb{Z}^{K})$ where the integral factor tracks the windings of the framing around each embedded singleton $\{k\}\to D$.

\begin{definition}
For any map $f:(\mcl{I},\mu)\to (\mcl{J},\nu)$ in $Linfr\Disk^{\ot}$ we take
\[
\opn{Br}_{fr}(f)=\prod_{j\in J}\opn{Br}_{fr}(I_j),
\]
where $I$ and $J$ are the underlying uncolored sets for $\mcl{I}$ and $\mcl{J}$, and each $I_j$ is given its ordering induced by $f$.
\end{definition}

Expanding on what was stated above, recall that a $2$-simplex
\[
\sigma=
\xymatrix{
 & (\mcl{J},\nu)\ar[dr]^g \\
(\mcl{I},\mu)\ar[rr]_h\ar[ur]^f & & (\mcl{L},\lambda)
}
\]
in $Linfr\Disk^{\ot}$ is precisely the choice of a homotopy $\sigma:|\Delta^1|\times D\times \mcl{I}_h\to D\times \mcl{J}$ with $\sigma_0=gf|_{\mcl{I}_h}$ and $\sigma_1=h$, and with the evaluation $\sigma_t:D\times \mcl{I}_h\to D\times \mcl{J}$ at each time $t\in |\Delta^1|$ an oriented embedding.  Such a $2$-simplex specifies an element $[\sigma]=(\beta,n)\in \opn{Br}_{fr}(gf)$, where the integral factor $n\in \prod_j(\mathbb{Z}^{I_j})=\prod_{i\in I_f}\mathbb{Z}$ has each entry $n_i$ tracing the winding number of the $x$-vector in the derivative for $\sigma_t|_{D_i}:D_i\to D_j$ around $0$.  The $j$-th contribution to the braid group $\beta_j$ in $\beta=(\beta_j:j\in J)$ is given by following the centers of the disks under the given homotopy
\[
\beta_j=\left[|\Delta^1|\times I_j\overset{1\times0\times 1}\to |\Delta^1|\times D\times I_j\overset{\sigma|_{I_j}}\to D_j\right].
\]

\subsection{A remark on pairings of uncolored disks}
\label{sect:pairing}

We note that there are two symmetric embeddings $\iota_{\pm}:fr\opn{Disk}^{\ot}\to fr\Disk^{\ot}$ from the symmetric monoidal $\infty$-category of \emph{un}colored disks. These embeddings are isomorphisms onto the full subcategories spanned by totally negative, and totally positive disks respectively. (One defines the uncolored disk category in direct analogy with the colored presentation from above.)

We can consider further the full subcategory spanned by both collections of totally negative and totally positive disks. The two embeddings $\iota_{\pm}$ are then connected by a transformation $\iota:\Delta^1\times fr\opn{Disk}^{\ot}\to fr\Disk^{\ot}$ with $\iota|_{\{0\}}=\iota_-$ and $\iota|_{\{1\}}=\iota_+$, and which is furthermore an equivalence onto this mixed subcategory. For any tuple of negatively colored disks $x=(-I,\mu)$, the corresponding natural morphism $\iota_x:\Delta^1\to fr\Disk^{\ot}$ is the negating of the identity $-id_x:(-I,\mu)\to (+I,\mu)$.
\par

One can realize the map $\iota$ directly at the level of simplicial categories or by solving a cocartesian lifting problem, though the precise details are not so important for us. Our only point here is that $fr\Disk^{\ot}$ itself can be considered as a type of coupling of two copies of the uncolored disk category along with certain transition morphisms.

\section{$fr\Disk$-monoidal $\infty$-categories}
\label{sect:disk_monoidal}

We show that $fr\Disk$ is equifibered over $\Fin$ (Defintion \ref{def:equifib}). We then provide basic information regarding $fr\Disk$-subcategories, and symmetric monoidal subcategories, in $fr\Disk$-monoidal $\infty$-categories. We also discuss a relation between $fr\Disk$-monoidality and monoidality over the colored little disk $\infty$-operad $fr\E_2$.

As a fundamental example for the text, we show in Section \ref{sect:ribbon_framed} that any ribbon tensor category, in the $1$-categorical sense of the term, admits a natural $fr\Disk$-monoidal structure. Following Section \ref{sect:ribbon_framed}, we provide constructions of marked bordism categories, then begin our descent into TQFTs.

\subsection{$fr\Disk$ is equifibered over $\Fin$}

\begin{proposition}\label{prop:frdisks_efib}
The components functor $fr\Disk^{\ot}\to \Fin^{\ot}$ gives $fr\Disk$ the structure of an equifibered symmetric monoidal $\infty$-category over $\Fin$.
\end{proposition}

\begin{proof}
According to our discussion around Corollary \ref{cor:simp_sym_mon}, the product functors for $fr\Disk$ and $\Fin$ are induced by the products on the preceding simplicial categories. Hence it suffices to show that the diagram
\[
\xymatrix{
fr\uDisk\times fr\uDisk\ar[rr]^{\amalg}\ar[d] & & fr\uDisk\ar[d]\\
\Fin\times \Fin\ar[rr]^{\amalg} & & \Fin
}
\]
is a pullback diagram of simplicial categories, in the strict sense of the term. For this we need to show that the corresponding diagrams on objects and morphisms are pullbacks.
\par

The case of objects is obvious. Let us now fix finite sets $I_{\varepsilon}$ and $J_{\varepsilon}$ with colored lifts $\mcl{I}_{\varepsilon}$ and $\mcl{J}_{\varepsilon}$. Take $I=I_0\amalg I_1$ and $J=J_0\amalg J_1$, and define $\mcl{I}$ and $\mcl{J}$ similarly. We consider the diagram of simplicial sets
\begin{equation}\label{eq:1239}
\xymatrix{
\uHom_{fr\uDisk}(\mcl{I}_0,\mcl{J}_0)\times \uHom_{fr\uDisk}(\mcl{I}_1,\mcl{J}_1)\ar[rr]^(.6){\amalg}\ar[d] & & \uHom_{fr\uDisk}(\mcl{I},\mcl{J})\ar[d]\\
\Hom_{\Fin}(I_0,J_0)\times \Hom_{\Fin}(I_1,J_1)\ar[rr]^(.6){\amalg}& & \Hom_{\Fin}(I,J).
}
\end{equation}
We note that the vertical maps here are given explicitly by taking maps on connected components and forgetting the colorings
\[
(\sigma:|\Delta^n|\times D\times \mcl{L}\to D\times \mcl{M})\mapsto (\pi_0(\sigma):L\to M).
\]

Since Homs over $\Fin$ are discrete, \eqref{eq:1239} is a pullback diagram if and only if, for chosen set maps $(\bar{f}_0,\bar{f}_1)\mapsto \bar{f}$, disjoint union provides an isomorphism between the fibers
\begin{equation}\label{eq:1258}
\uHom_{fr\uDisk}(\mcl{I}_0,\mcl{J}_0)_{\bar{f}_0}\times \uHom_{fr\uDisk}(\mcl{I}_1,\mcl{J}_1)_{\bar{f}_1}\to \uHom_{fr\uDisk}(\mcl{I},\mcl{J})_{\bar{f}}.
\end{equation}
This is clear however since both spaces decompose over $J$ as the product space
\[
\prod_{j\in J}\uHom_{fr\uDisk}(\mcl{I}_j,\{\pm 0\}),\ \ \pm=\opn{color}(j),
\]
and under this decomposition the disjoint union map \eqref{eq:1258} is just the identity.
\end{proof}

The same arguments show that the uncolored disk category $fr\opn{Disk}$ is equifibered over $\Fin$ as well.

\subsection{$fr\Disk$-categories}

As a shorthand, by a $fr\Disk$-category we mean a $fr\Disk$-monoidal $\infty$-category. Similarly, by a $fr\opn{Disk}$-category we mean a $fr\opn{Disk}$-monoidal $\infty$-category.

\begin{definition}\label{def:1279}
For any $fr\Disk$-category $q:\msc{E}^{\ot}\to fr\Disk^{\ot}$ we let $\msc{E}^-$ and $\msc{E}^+$ denote the respective fibers over the negatively and positively colored singletons $\{\pm 0\}$ in $fr\Disk\subseteq fr\Disk^{\ot}$. Given a colored set $\mcl{I}$, with negative and positive indices $I_-$ and $I_+$ respectively, the colored exponent $\msc{E}^{\mcl{I}}$ for $\msc{E}^{\ot}$ is defined as the product
\[
\msc{E}^{\mcl{I}}=(\msc{E}^-)^{I_-}\times (\msc{E}^+)^{I_+}.
\]
\end{definition}

The fibers $\msc{E}^{\pm}$ are the \emph{underlying $\infty$-categories} for a given $fr\Disk$-category $\msc{E}^{\ot}\to fr\Disk^{\ot}$. According to the discussion of Section \ref{sect:pairing}, one can think of a $fr\Disk$-category as a pair of $fr\opn{Disk}$-categories coupled by a $fr\opn{Disk}$-monoidal functor.

To elaborate, we have the two embeddings $\iota_{\pm}:fr\opn{Disk}^{\ot}\to fr\Disk^{\ot}$ and associated transformation $\iota:\Delta^1\times fr\opn{Disk}^{\ot}\to fr\Disk^{\ot}$ whose image is the full subcategory spanned by consistently colored disks. Hence, any $fr\Disk$-category $\msc{E}^{\ot}\to fr\Disk^{\ot}$ pulls pack to a cocartesian fibration
\[
\iota^{\ast}\msc{E}^{\ot}\to \Delta^1\times fr\opn{Disk}^{\ot}.
\]

The fibers of the above fibration over $\{\varepsilon\}\times fr\opn{Disk}^{\ot}$ exhibits $fr\opn{Disk}$-structures $q_{\varepsilon}:(\msc{E}^{\varepsilon})^{\ot}\to fr\opn{Disk}^{\ot}$ on the underlying $\infty$-categories $\msc{E}^{\pm}$. The unique cocartesian solution to the the lifting problem
\[
\xymatrix{
\{0\}\times (\msc{E}^-)^{\ot}\ar[rr]\ar[d] & & \iota^{\ast}\msc{E}\ar[d]\\
\Delta^1\times (\msc{E}^-)^{\ot}\ar[rr]_{id\times q_-}\ar@{..>}[urr] & & \Delta^1\times fr\opn{Disk}^{\ot}
}
\]
then provides a $fr\opn{Disk}$-monoidal functor
\[
\xymatrix{
(\msc{E}^-)^{\ot}\ar[rr]^{d^{\ot}}\ar[dr]_{q_-} & & (\msc{E}^+)^{\ot}\ar[dl]^{q_+}\\
	& fr\opn{Disk}^{\ot} & .
}
\]

\begin{remark}
Straightening provides an identification between our $fr\Disk$-categories and (colored) ``framed disk algebras" in $\sCat_{\infty}$, in the sense of Ayala-Francis \cite{ayalafrancis15}.
\end{remark}

\subsection{Structure maps for $fr\Disk$-categories}

\begin{definition}\label{def:str_frdisk}
For a given $fr\Disk$-category $\msc{E}^{\ot}\to fr\Disk^{\ot}$, the unit and product functors
\[
\opn{unit}:\ast\cong \msc{E}^{\ot}_{\emptyset}\to \msc{E}^{\pm}\ \ \text{and}\ \ m:\msc{E}^{\pm}\times\msc{E}^{\pm}\cong\msc{E}^{\ot}_{\{\pm 0,\pm 1\}}\to \msc{E}^{\pm}
\]
on the underlying $\infty$-categories $\msc{E}^{\pm}$ are the transport functors along the respective like colored embeddings $\emptyset\to D$ and $D\amalg D\to D$ in $fr\Disk\subseteq fr\Disk^{\ot}$, respectively. The transition functor $d:\msc{E}^-\to \msc{E}^+$ is the transport functor along the color negating identity $id:D_-\to D_+$.
\end{definition}

We commonly understand the tuple of functors
\[
(\opn{unit}:\ast\to \msc{E}^{\pm},\ m:\msc{E}^{\pm}\times\msc{E}^{\pm}\to \msc{E}^{\pm},\ d:\msc{E}^-\to \msc{E}^+)
\]
as the underlying structure maps for such a fibration $\msc{E}^{\ot}\to fr\Disk^{\ot}$.
\par

We note that, since all disk embeddings $D\amalg D\to D$ are isotopic (Proposition \ref{prop:1456}), the product functors from Definition \ref{def:str_frdisk} are well-defined up to natural isomorphism. For the sake of specificity however, we can fix the disk embedding $D\amalg D\to D$ as the one which maps the first disks by $z\mapsto \frac{1}{3}z-\frac{1}{2}$ and the second disk by $z\mapsto \frac{1}{3}z+\frac{1}{2}$.

\begin{remark}
We commonly understand an (uncolored) $fr\opn{Disk}$-category $\msc{F}^{\ot}\to fr\opn{Disk}^{\ot}$ as a homotopical variant of a balanced monoidal category. In topological terms, the braiding transformation for the product functor is provided by the cocartesian transformation $\Delta^1\times (\msc{F}\times\msc{F})\to \msc{F}$ over the rotation of the two disks within the ambient disk $|\Delta^1|\times(D\amalg D)\to D$, and the balancing transformation $\Delta^1\times\msc{F}\to \msc{F}$ is provided by rotating a single disk within itself $|\Delta^1|\times D\to D$.
\par

Similarly, from a $fr\Disk$-category $\msc{E}^{\ot}\to fr\Disk^{\ot}$, we extract our two homotopical balanced monoidal categories $(\msc{E}^{\pm})^{\ot}\to fr\Disk^{\ot}$ and a homotopically balanced monoidal functor $d^{\ot}:(\msc{E}^-)^{\ot}\to (\msc{E}^+)^{\ot}$.
\end{remark}

\subsection{$fr\Disk$-categories are $fr\mbf{E}_2$-monoidal $\infty$-categories}
\label{sect:disks_v_disks}

We take a moment to clarify the relationship between our monoidal $\infty$-categories, over the equifibered $\infty$-category $fr\Disk$, and monoidal $\infty$-categories over the $\infty$-operad of little disks. For the purposes of this discussion, the reader need only understand that an $\infty$-operad is an $\infty$-category $\msc{O}^{\ot}$ equipped with a type of isofibration to finite pointed sets $\msc{O}^{\ot}\to \Fin_{\ast}$ \cite[Definition 2.1.1.10]{ha}. As with symmetric monoidal $\infty$-categories, the underlying $\infty$-category for an $\infty$-operad is the fiber over the singleton $\msc{O}=\msc{O}^{\ot}_{\{0\}}$.
\par

Given an $\infty$-operad $\msc{O}^{\ot}$, an $\msc{O}$-monoidal $\infty$-category is a cocartesian fibration $\msc{E}^{\odot}\to \msc{O}^{\ot}$ which has consistent fibers, as in \eqref{eq:849}. See \cite[Deﬁnitions 2.1.1.10 \& 2.1.2.13]{ha}. We note that any $\msc{O}$-monoidal $\infty$-category is itself an $\infty$-operad, and that the structure map $\msc{E}^{\odot}\to \msc{O}^{\ot}$ is a map of $\infty$-operads \cite[Proposition 2.1.2.12]{ha}.
\par

We have the disjoint union functor $\opn{Fin}^{\ot}\to \Fin_{\ast}$, $(I,\mu)\mapsto I$, and for a given operad $\msc{O}^{\ot}$ an \emph{envelope} for $\msc{O}^{\ot}$ is any $\infty$-category $\opn{Env}(\msc{O})^{\ot}$ which fits into a pullback diagram
\[
\xymatrix{
\opn{Env}(\msc{O})^{\ot}\ar[rr]\ar[d] & & \msc{O}^{\ot}\ar[d]\\
\opn{Fin}^{\ot}\ar[rr]_{\opn{disj\ union}} & & \Fin_{\ast}
}
\]
in $\sCat_{\infty}$. We assume additionally that the given diagram strictly commutes, and that the projection to $\Fin^{\ot}$ is an isofibration.

\begin{remark}
Since the structure map for any $\infty$-operad is an isofibration, the strict pullback $\opn{Fin}^{\ot}\times_{\Fin_{\ast}}\msc{O}^{\ot}$ always provides a pullback, i.e.\ limit diagram, in $\sCat_{\infty}$ \cite[\href{https://kerodon.net/tag/03G7}{03G7}]{kerodon}. So, any envelope admits a unique equivalence to the strict pullback, as a fibration over $\Fin^{\ot}$.
\end{remark}

For any envelope $\opn{Env}(\msc{O})^{\ot}$, composing with the symmetric structure map (not disjoint union map!) for $\Fin^{\ot}$ produces a functor $\opn{Env}(\msc{O})^{\ot}\to \Fin_{\ast}$ which realizes $\opn{Env}(\msc{O})^{\ot}$ as a symmetric monoidal $\infty$-category \cite[Corollary 2.2.4.5.]{ha}. The underlying $\infty$-category for the resulting fibration $\opn{Env}(\msc{O})^{\ot}\to \Fin_{\ast}$ is identified with the non-full subcategory
\[
\opn{Env}(\msc{O})\cong\opn{Act}(\msc{O})\subseteq \msc{O}^{\ot}
\]
consisting of all objects, but only active maps in $\msc{O}^{\ot}$. We have our primary example of interest.

\begin{example}
The pullback diagram from \eqref{eq:1427} realizes the symmetric monoidal $\infty$-category of framed colored disks as an envelope $fr\Disk^{\ot}=\opn{Env}(fr\E_2)$ for the $\infty$-operad of framed colored disks.
\end{example}

Now, it was shown in works of Barkan, Haugseng, Kock, and Steinebrunner that the corresponding envelope functor
\[
\opn{Env}=\Fin^{\ot}\times_{\Fin_{\ast}}-:\msc{O}\!p_{\infty}\to (\SM_{\infty})_{/\Fin^{\ot}}
\]
is fully faithful \cite{haugsengkock24}, and furthermore an equivalence onto the full subcategory of equifibered $\infty$-categories over $\Fin$ \cite{barkanetal25}. Under this equivalence $\msc{O}$-monoidal $\infty$-categories, in the sense of \cite{ha}, are identified with $\opn{Env}(\msc{O})$-monoidal $\infty$-categories, in the sense of Definition \ref{def:T_mon}. Restricting to our particular case of interest, and noting \cite[\href{https://kerodon.net/tag/01ZS}{01ZS}]{kerodon}, we obtain equivalences between $fr\E_2$-monoidal $\infty$-categories and $fr\Disk$-monoidal $\infty$-categories.

\begin{theorem}[\cite{haugsengkock24,barkanetal25}]\label{thm:hk}
The envelope functor induces equivalences of $\infty$-categories
\[
\opn{Env}:fr\E_2\text{-}\msc{C}\!ats\overset{\sim}\to fr\Disk\text{-}\msc{C}\!ats.
\]
\end{theorem}

\subsection{Symmetric subcategories in $fr\Disk$-categories}
\label{sect:subcats}

Let $\msc{E}^{\ot}\to fr\Disk^{\ot}$ be a $fr\Disk$-category and consider full subcategories $\msc{L}^{\pm}\subseteq \msc{E}^{\pm}$ which are stable under isomorphism in the underlying $\infty$-categories $\msc{E}^{\pm}$. We define the associated subcategory $\msc{L}^{\ot}\subseteq \msc{E}^{\ot}$ generated by the $\msc{L}^{\pm}$ to be the (unique) full subcategory whose fibers over $fr\Disk^{\ot}$ fit into strict pullback diagrams
\begin{equation}\label{eq:3000}
\xymatrix{
\msc{L}^{\ot}_t\ar[rr]\ar[d] & & \msc{L}^{\mcl{I}}\ar[d]\\
\msc{E}^{\ot}_t\ar[rr]_{\rho_!} & & \msc{E}^{\mcl{I}},
}
\end{equation}
where $t=(\mcl{I},\mu)$ is arbitrary. Since the inclusion $\msc{L}^{\mcl{I}}\to \msc{E}^{\mcl{I}}$ is an isofibration in this case, and since each $\rho_!$ is an equivalence, the above pullback diagram implies that the transport equivalences for $\msc{E}^{\ot}$ restrict to equivalences $\msc{L}^{\ot}_t\overset{\sim}\to \msc{L}^{\mcl{I}}$ for $\msc{L}^{\ot}$ as well (Lemma \ref{lem:pullback_equiv}).
\par

Now, according to Section \ref{sect:disks_v_disks}, we might \emph{imagine} $\msc{E}^{\ot}$ as a symmetric monoidal $\infty$-category of disks which are labeled by objects in the $\msc{E}^{\pm}$. The symmetric structure is then given by taking disjoint unions of such labeled disks. From this perspective one expects that the subcategory $\msc{L}^{\ot}$ of disks with labels in $\msc{L}^{\pm}$ forms a symmetric monoidal subcategory in $\msc{E}^{\ot}$. One can show that this is in fact the case.

\begin{lemma}\label{lem:full_to_symm}
For full subcategories $\msc{L}^{\pm}\subseteq \msc{E}^{\pm}$ which are stable under isomorphism, the associated $\infty$-category $\msc{L}^{\ot}$ is naturally symmetric monoidal and the inclusion $\msc{L}^{\ot}\to \msc{E}^{\ot}$ is a symmetric monoidal functor. Furthermore, $\msc{L}^{\ot}$ is stable under isomorphism in $\msc{E}^{\ot}$.
\end{lemma}

Note that the assertion that the inclusion $\msc{L}^{\ot}\to \msc{E}^{\ot}$ is symmetric monoidal forces the symmetric structure map $\msc{L}^{\ot}\to \opn{Fin}_{\ast}$ to be the one obtained from $\msc{E}^{\ot}$ via restriction. We delay the proof of Lemma \ref{lem:full_to_symm} to Appendix \ref{sect:full_to_symm}.

\begin{remark}
A version of Lemma \ref{lem:full_to_symm} holds for non-full subcategories $\msc{L}^{\pm}\subseteq \msc{E}^{\pm}$. The construction of the category $\msc{L}^{\ot}$ is more delicate in this case, however.
\end{remark}

\subsection{$fr\Disk$-subcategories}

\begin{proposition}\label{prop:frdisk_subcats}
Let $q:\msc{E}^{\ot}\to fr\Disk^{\ot}$ be a $fr\Disk$-category and $\msc{L}^{\pm}$ be full subcategories in the underlying $\infty$-categories $\msc{E}^{\pm}=\msc{E}^{\ot}_{\{\pm 0\}}$ which are stable under isomorphism. Suppose also that
\begin{itemize}
\item The unit maps $\opn{unit}:\ast\to \msc{E}^{\pm}$ have image in $\msc{L}^{\pm}$.\vspace{1mm}
\item The transition map $d:\msc{E}^-\to \msc{E}^+$ sends $\msc{L}^-$ into $\msc{L}^+$.\vspace{1mm}
\item The product functors $m:\msc{E}^{\pm}\times\msc{E}^{\pm}\to \msc{E}^{\pm}$ send $\msc{L}^{\pm}\times \msc{L}^{\pm}$ into $\msc{L}^{\pm}$.
\end{itemize}
Then there is a unique, full $fr\Disk$-subcategory $\msc{L}^{\ot}$ in $\msc{E}^{\ot}$ which is stable under isomorphism and has underlying $\infty$-categories $\msc{L}^{\pm}$.
\end{proposition}

The $fr\Disk$-subcategory $\msc{L}^{\ot}$ is, as one expects, precisely the symmetric subcategory constructed in Section \ref{sect:subcats}. So we are giving necessary (and sufficient) conditions for the associated symmetric subcategory $\msc{L}^{\ot}$ from Lemma \ref{lem:full_to_symm} to be a $fr\Disk$-subcategory. To assist with the proof we record a related lemma.

\begin{lemma}\label{lem:structure_transp}
For any $fr\Disk$-category $q:\msc{E}^{\ot}\to fr\Disk^{\ot}$ and morphism $f:s=(\mcl{I},\mu)\to t=(\mcl{J},\nu)$ in $fr\Disk^{\ot}$, the corresponding transport functor fits in to a diagram
\[
\xymatrix{
\msc{E}^{\ot}_s\ar[rr]^{f_!}\ar[d]_{\rho_!}^{\sim} & & \msc{E}^{\ot}_t\ar[d]_{\rho_!}^{\sim}\\
\msc{E}^{\mcl{I}}\ar[rr]_{m_f} & & \msc{E}^{\mcl{J}}
}
\]
in $\sCat_{\infty}$ in which $m_f$ decomposes as a composite of functors $D:\msc{E}^{\mcl{I}}\to \msc{E}^{\mcl{I}_0}$, $M:\msc{E}^{\mcl{I}_0}\to\msc{E}^{\mcl{J}_0}$, and $U:\msc{E}^{\mcl{J}_0}\to \msc{E}^{\mcl{J}}$ where
\begin{itemize}
\item $\mcl{I}_0$ is obtained from $\mcl{I}$ by changing some negatively colored indices positive, and $D$ applies the transition functor $d:\msc{E}^-\to \msc{E}^+$ at all such indices.\vspace{1mm}
\item $\mcl{J}_0\subseteq \mcl{J}$ is the image of $f:D\times\mcl{I}_f\to D\times \mcl{J}$, or rather $\pi_0(f)$, and $M$ is an iterated application of the product functors $m:\msc{E}^{\pm}\times \msc{E}^{\pm}\to \msc{E}^{\pm}$.\vspace{1mm}
\item $U:\msc{E}^{\mcl{J}_0}\to \msc{E}^{\mcl{J}}$ inserts the unit $\opn{unit}:\ast\to \mcl{E}^{\pm}$ at all indices in the complement $\mcl{J}\setminus \mcl{J}_0$.
\end{itemize}
\end{lemma}

The proof of Lemma \ref{lem:structure_transp} is time consuming but straightforward, and is provided in Appendix \ref{sect:str_proof}. We return to our original proposition.

\begin{proof}[Proof of Proposition \ref{prop:frdisk_subcats}]
From the description of transport provided in Lemma \ref{lem:structure_transp}, and our hypotheses, it is clear that $\msc{L}^{\ot}$ is stable under transport along arbitrary maps in $fr\Disk^{\ot}$. Since $\msc{L}^{\ot}$ is stable under isomorphism in $\msc{E}^{\ot}$, stability under transport implies that, for any object $x$ in $\msc{L}^{\ot}$, map $q(x)\to \bar{y}$ in $fr\Disk^{\ot}$, and $q$-cocartesian lift $x\to y$ in $\msc{E}^{\ot}$, the object $y$ is also in $\msc{L}^{\ot}$. Hence, by fullness, the structure map $\msc{E}^{\ot}\to fr\Disk^{\ot}$ restricts to a cocartesian fibration $\msc{L}^{\ot}\to fr\Disk^{\ot}$ and the inclusion $\msc{L}^{\ot}\to \msc{E}^{\ot}$ preserves cocartesian edges. It follows that the transport functors for $\msc{E}^{\ot}$, along maps in $fr\Disk^{\ot}$, restrict to provide transport functors for $\msc{L}^{\ot}$.
\par

We need to show that the transport functors for $\msc{L}^{\ot}$ induce the requisite decompositions of the fibers, as outlined in Definition \ref{def:T_mon}. Consider an arbitrary object $t=(\mcl{I},\mu)$ in $fr\Disk^{\ot}$ and the separation $\tilde{t}=(\mcl{I},id_{\mcl{I}})$. By assumption we have pullback diagrams
\[
\xymatrix{
\msc{L}^{\ot}_{\tilde{t}}\ar[r]\ar[d] & \msc{L}^{\mcl{I}}\ar[d]\\
\msc{E}^{\ot}_{\tilde{t}}\ar[r]^{\sim} & \msc{E}^{\mcl{I}}
}
\]
in which the vertical maps are isofibrations, which implies the restricted functor $\msc{L}^{\ot}_{\tilde{t}}\to \msc{L}^{\mcl{I}}$ is an equivalence. The transport equivalence $\msc{E}^{\ot}_{t}\overset{\sim}\to \msc{E}^{\mcl{I}}$ similarly restricts to an equivalence $\msc{L}^{\ot}_{t}\overset{\sim}\to \msc{L}^{\mcl{I}}$. From the diagram
\[
\xymatrix{
\msc{L}^{\ot}_{\tilde{t}}\ar[rr]\ar[dr]_{\sim} & & \msc{L}^{\ot}_t\ar[dl]^{\sim}\\
 & \msc{L}^{\mcl{I}}
}
\]
we conclude now that the transport functor $\msc{L}^{\ot}_{\tilde{t}}\to \msc{L}^{\ot}_t$ is an equivalence.
\par

Consider now the maps $t=(\mcl{I},\mu)\to t'=\mcl{I}$ just given by the identity on $\mcl{I}$. From the diagram
\[
\xymatrix{
\msc{L}^{\ot}_{\tilde{t}}\ar[rr]^{\sim}\ar[dr]_{\sim} & & \msc{L}^{\ot}_{t}\ar[dl]\\
	& \msc{L}^{\ot}_{t'}
}
\]
we conclude that transport provides an equivalence $\msc{L}^{\ot}_t\overset{\sim}\to \msc{L}^{\ot}_{t'}$ as well. Writing out $\mu:\mcl{I}\to K$ and taking $t_k=\mcl{I}_k$ for $k\in K$, we similarly observer a diagram
\[
\xymatrix{
\msc{L}^{\ot}_{\tilde{t}}\ar[dr]\ar[rr]^{\sim}\ar@{-->}[ddr]_{\sim} & & \msc{L}^{\ot}_t\ar[dl]\ar@{-->}[ddl]^{\sim}\\
	& \prod_{k\in K}\msc{L}^{\ot}_{t_k}\ar@{-->}[d]|{\sim}\\
	& \msc{L}^{\mcl{I}}
}
\]
to find that transport provides a final equivalence $\msc{L}^{\ot}_t\overset{\sim}\to \prod_k \msc{L}^{\ot}_{t_k}$. So we see that the cocartesian fibration $\msc{L}^{\ot}\to fr\Disk^{\ot}$ satisfies the factorization constraints from Definition \ref{def:equifib}, and verify the $fr\Disk$-monoidal structure on $\msc{L}^{\ot}$.
\end{proof}

\section{Disk categories from ribbon tensor categories}
\label{sect:ribbon_framed}

We show that any discrete ribbon tensor category $A^{\heartsuit}$ produces a gaggle of $fr\Disk$-monoidal $\infty$-categories $q_{\star}:A^{\ot}_{\star}\to fr\Disk^{\ot}$ (Proposition \ref{prop:fr_cochains} and Corollary \ref{cor:1714}). Here $A$ is the category of cochains $A=\opn{Ch}(A^{\heartsuit})$ and $\star$ is a modifier. Our constructions follow from a general result which asserts that any functor between balanced monoidal categories $d:E^-\to E^+$ determines a $fr\Disk$-monoidal $\infty$-category $q:E^{\ot}\to fr\Disk^{\ot}$.
\par

In order to support our analysis of unbounded complexes over cocomplete categories, like $A^{\heartsuit}$, we recall basic information concerning pro-finite completion.  Pro-finite objects are essential tools for us, as a pro-finite category is always on the ``other side" of a finitely generated linear abelian category under a duality.

\subsection{$fr\Disk$-categories from discrete balanced categories}
\label{sect:framed_ribbon}

Recall that a balanced monoidal category $E$ is a (discrete) braided monoidal category which is equipped with a natural isomorphism $\theta:id_E\to id_E$ which satisfies
\[
\theta_{x\ot y}=(\theta_x\ot \theta_y)c_{yx}c_{xy},
\]
for any pair of objects $x$ and $y$ in $E$, where $c$ is the braiding. A balanced functor between balanced monoidal categories $E$ and $D$, with associated transformations $\theta^E$ and $\theta^D$, is a braided monoidal functor $F:E\to D$ which satisfies $F\theta^E=\theta^DF$.

A balanced tensor category $E$ is a ribbon tensor category if there is an equality $\theta_x^{\ast}=\theta_{x^{\ast}}$ whenever $x$ is rigid in $E$.

We consider a pair of balanced monoidal categories $E^-$ and $E^+$ equipped with a balanced functor $d:E^-\to E^+$.  For such a pairing, and colored set $\mcl{I}$ we consider the colored exponent $E^{\mcl{I}}$, as in Definition \ref{def:1279}, and for any map $f:(\mcl{I},\mu)\to (\mcl{J},\nu)$ in $Linfr\Disk^{\ot}$, with underlying disk arrangement $D\times \mcl{I}_f\to D\times \mcl{J}$, we let
\begin{equation}\label{eq:m_f}
m_f:E^{\mcl{I}}\to E^{\mcl{J}}
\end{equation}
denote the associated product with constituent factors
\[
m_f^j:E^{\mcl{I}_j}\to E^{\opn{color}(j)},\ \ x\mapsto d^{\varepsilon_{i_1}}(x_{i_1})\ot\dots\ot d^{\varepsilon_{i_r}}(x_r).
\]
Here each $\mcl{I}_j=\{i_1,\dots, i_r\}$ is given the ordering specified by the map $f$ (see Section \ref{sect:braids}) and
\[
d^{\varepsilon_i}=\left\{
\begin{array}{ll}
id & \text{if }\opn{color}(i)=\opn{color}(j)\\
d & \text{if }\opn{color}(i)=-\opn{color}(j).
\end{array}\right.
\]
We refer the the map $m_f$ from \eqref{eq:m_f} as the \emph{colored product functor} for the triple $\{E^-,E^+,d\}$ over $f$.
\par

To be clear, $m_f$ sends the excess factors $E^{\mcl{I}\setminus \mcl{I}_{f}}$ to a point, i.e.\ $m_f$ is a composite of the projection $E^{\mcl{I}}\to E^{\mcl{I}_{f}}$ with the aforementioned multiplications, and for any index $j\in \mcl{J}$ with empty preimage the map \eqref{eq:m_f} just chooses the unit object $\1_{\opn{color}(j)}:\ast\to E^{\opn{color}(j)}$. For any composite $gf:\mcl{I}\to \mcl{L}$ in $Linfr\E_2$ we also have the natural isomorphism
\[
m_gm_f\overset{\cong}\to m_{gf}
\]
provided by the associators and the monoidal structure on $d$.

\begin{proposition}\label{prop:a_pm}
Given any triple $E=(E^-,E^+,d)$ consisting of balanced monoidal categories and a balanced monoidal functor $d:E^-\to E^+$, there is a canonically defined $fr\Disk$-category $q:E^{\ot}\to fr\Disk^{\ot}$ whose fiber over each colored set $(\mcl{I},\mu)$ is the colored exponent $E^{\ot}_{(\mcl{I},\mu)}=E^{\mcl{I}}$ and whose transport functor along any map $f:(\mcl{I},\mu)\to (\mcl{J},\nu)$ in $Linfr\Disk^{\ot}$ is the colored product functor
\[
f_!=m_f:E^{\mcl{I}}\to E^{\mcl{J}}.
\]
\end{proposition}

The proof is covered in Appendix \ref{sect:framed_pairs}.

\subsection{Describing the fibrations $q:E^{\ot}\to fr\Disk^{\ot}$}
\label{sect:A_ot_des}

The fibration $q:E^{\ot}\to fr\Disk^{\ot}$ is defined by taking a full subcategory $E^{\ot}\subseteq \msc{E}_d^{\ot}$ in a certain cocartesian fibration $\tilde{q}:\msc{E}_d^{\ot}\to fr\Disk^{\ot}$ associated to the map $d$.  We take in particular $\msc{E}^{\ot}_d=E^{\ot}_d\times_{fr\opn{Disk}^{\ot}}fr\Disk^{\ot}$, for $E_d^{\ot}\to fr\opn{Disk}^{\ot}$ the $fr\opn{Disk}$-category defined via the balanced monoidal category $E_d$ from Section \ref{sect:e_ot_const}. The fibration $E_d^{\ot}\to fr\opn{Disk}$ is constructed explicitly in Section \ref{sect:bal_description}, and we thus obtain an explicit description of the fibration $q:E^{\ot}\to fr\Disk^{\ot}$.  We record the relevant low-dimensional information for $q$.

Before providing this description however, we recall some preliminary topological information. For a $2$-simplex
\[
\xymatrix{
	& (\mcl{J},\nu)\ar[dr]^{g}\\
(\mcl{I},\mu)\ar[rr]_{h}\ar[ur]^{f} & & (\mcl{L},\lambda)
}
\]
in $Linfr\Disk^{\ot}$ with specifying homotopy $\sigma:|\Delta^1|\times D\times\mcl{I}_{h}\to D\times \mcl{L}$ from $gf$ to $h$ the associated framed braid $[\sigma]=(\beta,n)\in \opn{Br}_{fr}(gf)$ from Section \ref{sect:braids} has underlying symmetry in $\prod_jS_{I_l}$ exchanging the ordering on $\mcl{I}_l$ induced by $gf$ with the ordering induced by $h$.  Hence $\beta$ provides, via the braidings on $E^{\pm}$, a natural isomorphism $\beta:m_{gf}\to m_{h}$. The framed braid itself $[\sigma]=(\beta,n)$ then induces a transformation
\begin{equation}\label{eq:zeta_beta}
\zeta_{\beta,n}=\beta\circ m_{gf}(\theta^n)\circ \opn{assoc}:m_gm_f\cong m_{gf}\to m_h
\end{equation}
in which the integral parameter $n\in \prod_j\mathbb{Z}^{I_j}$ acts via the balancing transformaions on the various factors.
\par

Returning to the fibration $E^{\ot}\to fr\Disk^{\ot}$, the objects in $E^{\ot}$ consist of a choice of finite colored set $\mcl{I}$ with a map to a finite set $\mu:\mcl{I}\to K$ and an object $x$ in $E^{\mcl{I}}$. For $f:(\mcl{I},\mu)\to (\mcl{J},\nu)$ in the specified subcategory $Linfr\Disk^{\ot}$, a map $\alpha:(x,\mcl{I},\mu)\to (y,\mcl{J},\nu)$ in $E^{\ot}$ over $f$ is the choice of a map $\alpha:m_f(x)\to y$.  A $2$-simplex
\[
\xymatrix{
	& y\ar[dr]^{\kappa} &	&\ar@{}[d]|{\text{\normalsize over}} & & (\mcl{J},\nu)\ar[dr]^g\ar@{}[d]|(.6){\sigma} & & \ar@{}[d]|{\text{\normalsize in $Linfr\Disk^{\ot}$}}\\
x\ar[ur]^{\alpha}\ar[rr]_{\gamma} & & z & & (\mcl{I},\mu)\ar[ur]^f\ar[rr]_h & & (\mcl{L},\lambda) &
}
\]
is the choice of maps $\alpha:m_f(x)\to y$, $\kappa:m_g(y)\to z$, and $\gamma:m_h(x)\to z$ which produce a diagram
\[
\xymatrix{
m_gm_f(x)\ar[d]_{\zeta_{\beta,n}}\ar[rr]^{m_g(\alpha)} & & m_g(y)\ar[d]^{\kappa}\\
m_h(x)\ar[rr]_{\gamma} & & z
}
\]
in $E^{\mcl{L}}$, where $(\beta,n)$ is the framed braid specified by $\sigma$.
\par

One should recall at this point that $E^{\ot}$ is a symmetric monoidal $\infty$-category and that the fibration $q:E^{\ot}\to fr\Disk^{\ot}$ is a symmetric monoidal functor, as is the restriction to linear disk arrangements $q:E^{\ot}|_{Linfr\Disk^{\ot}}\to Linfr\Disk^{\ot}$. From the above description we see that the underlying $\infty$-category for $E^{\ot}|_{Linfr\Disk^{\ot}}$ is the $\infty$-category of disks labeled by objects in $E^{\pm}$ and (linear) disk embeddings. The symmetric structure is given by taking disjoint unions of such labeled disks.
\par

Since the inclusion $Linfr\Disk^{\ot}\to fr\Disk^{\ot}$ is an equivalence (Lemma \ref{lem:linfr_fr}), and the structure map $q:E^{\ot}\to fr\Disk^{\ot}$ is an isofibration, we have that the inclusion $E^{\ot}|_{Linfr\E_2}\to E^{\ot}$ is an equivalence as well (Lemma \ref{lem:pullback_equiv}).  From this one obtains a precise description of the homotopy category for $E^{\ot}$, as a discrete symmetric monoidal category.

\begin{lemma}
For $(E^+,E^-,d)$ as in Proposition \ref{prop:a_pm}, the homotopy category $\opn{h}E^{\ot}_{\{0\}}$ of the underlying $\infty$-category for the symmetric monoidal structure $E^{\ot}\to \opn{Fin}_{\ast}$ has objects given by objects $x$ in $E^{\mcl{I}}$, across varying $\mcl{I}$, and morphisms $[\alpha_f]:x\to y$ given by the choice of a linear arrangement of disks $f:D\times \mcl{I}\to D\times\mcl{J}$ coupled with a morphism $\alpha:m_f(x)\to y$. Two such maps $[\alpha_f]:x\to y$ and $[\alpha'_{f'}]:x\to y$ are equal in $\opn{h}E^{\ot}_{\{0\}}$ if and only if they are related by the action of the framed braid group, i.e.\ if there is a framed braid $(\beta,n)$ from $f$ to $f'$ for which we have an equality
\[
\alpha'\circ \zeta_{\beta,n}|_x=\alpha.
\]
The symmetric monoidal structure is given by taking disjoint unions of disks and (cartesian) products of morphisms in $E^{\pm}$.
\end{lemma}

Here, again, by a framed braid $(\beta,n)$ \emph{from $f$ to $f'$} we mean a framed braid for which the associated elements in the symmetric groups $\bar{\beta}_j\in S_{I_j}$ are the unique permutations which exchange the ordering on $I_j$ defined by $f$ with the ordering defined by $f'$.

\subsection{Pro-finite objects}
\label{sect:pf_obj}

We show below that the category of unbounded cochains $A=\opn{Ch}(A^{\heartsuit})$ over a ribbon tensor category produces a $fr\Disk$-category $A^{\ot}\to fr\Disk^{\ot}$.  Here the positive fiber is a pro-finite variant of $\opn{Ch}(A^{\heartsuit})$ and the negative fiber is $\opn{Ch}(A^{\heartsuit})^{op}$.  Both unboundedness and pro-finiteness cause some technical issues, which we resolve over the next two subsections.
\par

Fix a tensor category $A^{\heartsuit}$, or more generally a finitely generated linear abelian category (Definition \ref{def:lin_ab}). First, the category $A^{\heartsuit}_{\pf}$ of pro-finite objects in $A^{\heartsuit}$ consists of objects $x$ in $A^{\heartsuit}$ which come equipped with a filtered collection $\mcl{B}_x$ of inclusions $i^{\lambda}:x^{\lambda}\to x$, over a filtered category $\Lambda$, which satisfy the following:
\begin{itemize}
\item An inclusion $i':x'\to x$ belongs to $\mcl{B}_x$ if and only if there is a map $i^{\lambda}:x^{\lambda}\to x$ in the collection which factors through $i'$.\vspace{1mm}
\item Each quotient $\pi_{\lambda}:x\to x_{\lambda}=x/x^{\lambda}$ has $x_{\lambda}$ of finite length.\vspace{1mm}
\item The natural map
\begin{equation}\label{eq:can_lim}
x\to \varprojlim_\lambda x_{\lambda}
\end{equation}
is an isomorphism.\vspace{1mm}
\end{itemize}
Morphisms between objects in $A^{\heartsuit}_{\pf}$ are maps $f:x\to y$ for which $x\times_y(y^{\mu}\to y)$ lies in $\mcl{B}_x$ whenever $y^{\mu}\to y$ lies in $\mcl{B}_y$.  Equivalently, we require each composite $x\to y\to y_{\mu}$ to factor through some projection $x\to x_{\lambda}\to y_\mu$.  We call the collection $\mcl{B}_x$ for pro-finite $x$ the \emph{base} for $x$ around $0$, and a general map $f:x\to y$ in $A^{\heartsuit}$ between pro-finite objects is called \emph{continuous} if it is a map in $A^{\heartsuit}_{\pf}$.

\begin{remark}
The indexing category $\Lambda$ here is not small but, after identifying $A^{\heartsuit}$ with a corepresentation category \cite[Theorem 5.1]{takeuchi77}, one sees that $\Lambda$ is is fact essentially small.  So the limit $\varprojlim_{\lambda}x_{\lambda}$ always exists in $A^{\heartsuit}$.
\end{remark}

\begin{remark}
Filteredness, and the first point above, tells us that the intersection $i^{\lambda\mu}:x^{\lambda}\times_xx^{\mu}\to x$ is in $\mcl{B}_x$ whenever the constituent maps $i^{\lambda}:x^{\lambda}\to x$ and $i^{\mu}:x^{\mu}\to x$ are in $\mcl{B}_x$.
\end{remark}

We note that $A^{\heartsuit}$ and $A^{\heartsuit}_{\pf}$ share the subcategory of finite length objects.  In particular, each finite length object is naturally pro-finite and each map between finite length objects is continuous.  So we have fully faithful inclusions
\[
A^{\heartsuit}\leftarrow A^{\heartsuit}_{fin}\to A^{\heartsuit}_{\pf}.
\]
The objects in $A^{\heartsuit}_{fin}$ are precisely the discrete objects in $A^{\heartsuit}_{\pf}$, in the sense that $0\to x$ is in the base around $0$ for pro-finite $x$ if and only if $x$ is of finite length.
\par

When $A^{\heartsuit}$ is a tensor category, we endow $A^{\heartsuit}_{\pf}$ with the monoidal structure induced by the completed tensor product
\[
x\ot y:=x\hat{\ot}y=\varprojlim_{\lambda,\mu}x_{\lambda}\ot y_{\mu}.
\]
The object $x\ot y$ is naturally pro-finite, with canonical base given by all subobjects $w\to x\ot y$ which contain the kernel of some structural projection $x\ot y\to x_{\lambda}\ot y_{\mu}$.  When $A^{\heartsuit}$ is braided (resp.\ balanced) the braiding (resp.\ balancing structure) induces a unique braiding (resp.\ balancing structure) on $A^{\heartsuit}_{\pf}$ under which the inclusion $A^{\heartsuit}_{fin}\to A^{\heartsuit}_{\pf}$ is a map of braided (resp.\ balanced) monoidal categories.
\par

In the following statement, the opposite category $(A^{\heartsuit})^{op}$ is taken to be monoidal under the opposite product $x\ot^{op}y=y\ot x$. When $A^{\heartsuit}$ is braided the opposite category is naturally braided via the same operations $c^{op}_{xy}=c_{xy}$ and any subsequent twist for $A^{\heartsuit}$ is a twist for $(A^{\heartsuit})^{op}$ as well.

\begin{lemma}\label{lem:1572}
For any tensor category $A^{\heartsuit}$, the duality functor provides an equivalence of monoidal categories
\begin{equation}\label{eq:1564}
-^{\ast}:(A^{\heartsuit})^{op}\overset{\sim}\to A^{\heartsuit}_{\pf}.
\end{equation}
When $A^{\heartsuit}$ is furthermore ribbon, this equivalence is an equivalence of balanced monoidal categories.
\end{lemma}

We note that the expression for the dual $x^{\ast}:=\underline{\Hom}(x,\1)$ is valid even when $x$ is of infinite length, where $\underline{\Hom}$ is the inner-Hom functor for the tensor action of $A^{\heartsuit}$ on itself. In the above expression \eqref{eq:1564} our duality functor is, more specifically, this inner-Hom functor on $A^{\heartsuit}$.

\begin{proof}
First recall that duality is an equivalence on the full subcategory of finite length objects, $-^{\ast}:(A^{\heartsuit}_{fin})^{op}\overset{\sim}\to A^{\heartsuit}_{fin}$. We understand that any object $x$ in $A^{\heartsuit}$ is canonically expressed as the filtered colimit $x=\varinjlim_{\lambda}x_\lambda$ of its finite length subobjects $x_\lambda\to x$, and that all maps in $A^{\heartsuit}$ send finite subobjects to finite subobjects.  Since the duality functor exchanges colimits with limits, we now have the canonical expression for the dual $x^{\ast}=\varprojlim x^{\ast}_{\lambda}$ which endows $x^{\ast}$ with the base provided by the kernels of the projections $x^{\ast}\to x^{\ast}_{\lambda}$.  So we see that the duality functor $-^{\ast}:(A^{\heartsuit})^{op}\to A^{\heartsuit}$ naturally has image in the non-full subcategory $A^{\heartsuit}_{\pf}$ of pro-finite objects.  The inverse to $-^{\ast}$ is given by the continuous (right) dual ${^{\ast}-}:A^{\heartsuit}_{\pf}\to (A^{\heartsuit})^{op}$, ${^{\ast}z}=\varinjlim_{\tau}{^{\ast}z}_\tau$.
\par

As for the tensor structure, we have the natural isomorphisms $(x_{\lambda}\ot y_{\mu})^{\ast}\overset{\sim}\to y^{\ast}_{\mu}\ot x^{\ast}_{\lambda}$ for finite-length objects which induce unique natural isomorphisms $(x\ot y)^{\ast}\to y^{\ast}\ot x^{\ast}$ for arbitrary objects in $A^{\heartsuit}$.  This gives $-^{\ast}$ its monoidal structure.  Compatibility with the braiding follows from the identification $(c_{xy})^{\ast}= c_{x^{\ast}y^{\ast}}$ whenever $x$ and $y$ are of finite length \cite[Lemma 2.1.11]{bakalovkirillov01}, which forces such an equality at all objects, and the identification for the twist $\theta^{\ast}_x=\theta_{x^{\ast}}$ is precisely the ribbon condition. 
\end{proof}

We note that the opposite category $(A^{\heartsuit})^{op}$ itself is both complete and cocomplete, since $A^{\heartsuit}$ is both complete and cocomplete.  Hence Lemma \ref{lem:1572} tells us that the category of pro-finite objects $A^{\heartsuit}_{\pf}$ is both complete and cocomplete as well, and also abelian.  When $A^{\heartsuit}$ is finite, and thus has enough projectives, we see that $A^{\heartsuit}_{\pf}$ has enough injectives.

\begin{lemma}\label{lem:pf_limits}
Let $A^{\heartsuit}$ be a tensor category.
\begin{enumerate}
\item The forgetful functor $A^{\heartsuit}_{\pf}\to A^{\heartsuit}$ preserves all limits.\vspace{1mm}
\item Given any essentially small cofiltered diagram $K^{op}\to A^{\heartsuit}_{fin}$ in which all of the transition maps $x_{\alpha}\to x_{\beta}$ are surjective, and any pro-finite $x$, a map $f:x\to \varprojlim_{\alpha}x_{\alpha}$ in $A^{\heartsuit}$ is a (continuous) isomorphism in $A^{\heartsuit}_{\pf}$ if and only if the following hold:\vspace{1mm}
\begin{enumerate}
\item Each composite $f_{\alpha}:x\to x_{\alpha}$ factors through a structural projection $\pi_{\mu}:x\to x_{\mu}$, and the resulting map $x_{\mu}\to x_{\alpha}$ is surjective.\vspace{1mm}
\item Each structural projection $\pi_{\lambda}:x\to x_{\lambda}$ factors through some composite $f_{\beta}:x\to x_{\beta}$.
\end{enumerate}
\end{enumerate}
\end{lemma}

\begin{proof}
(1) Follows from the fact that the duality functor $-^{\ast}=\underline{\Hom}(-,\1):(A^{\heartsuit})^{op}\to A^{\heartsuit}$ is left exact and restricts to an equivalence onto $A^{\heartsuit}_{\pf}$. (2) Continuity requires that each composite $f_{\alpha}$ factors through some $\pi_{\lambda}$.  The fact that the (now continuous) map $f:x\to \varprojlim_{\alpha} x_{\alpha}$ is an isomorphism if and only if each projection $\pi_{\lambda}$ factor through some $f_{\beta}$ follows from the dual statement in $A^{\heartsuit}$.  Namely, one sees that a map from a filtered colimit in which all transition functions are injective $\varinjlim_{\alpha}y_{\alpha}\to y$ is an isomorphism if and only if each constituent map $y_{\alpha}\to y$ factors through some inclusion $y_{\mu}\to y$ of a finite length subobject, and the resulting map $y_{\alpha}\to y_{\mu}$ is injective, and each inclusion $y_{\lambda}\to y$ from a finite length subobject factors through some $y_{\beta}\to y$.  The first condition tells us that the map $\varinjlim_{\alpha}y_{\alpha}\to y$ is injective and the latter forces surjectivity.
\end{proof}

Applying Proposition \ref{prop:a_pm} to Lemma \ref{lem:1572}, we observe the following.

\begin{proposition}\label{prop:fr_heartsuit}
For any ribbon tensor category $A^{\heartsuit}$, the triple
\[
(A^{\heartsuit})^-=(A^{\heartsuit})^{op},\ \ (A^{\heartsuit})^+=A^{\heartsuit}_{\pf},\ \ -^{\ast}:(A^{\heartsuit})^{op}\to A^{\heartsuit}_{\pf}
\]
defines a $fr\Disk$-monoidal $\infty$-category $q^{\heartsuit}:(A^{\heartsuit})^{\ot}\to fr\Disk^{\ot}$ as in Proposition \ref{prop:a_pm}.
\end{proposition}

\subsection{Pro-finite complexes}
\label{sect:pf_complexes}

In this work we have a fundamental interest in \emph{cochains} over ribbon tensor categories. We consider pro-finite completion in this specific setting.

For complexes $\opn{Ch}(A^{\heartsuit})$ over a tensor category $A^{\heartsuit}$, degree-wise duality provides an equivalence
\begin{equation}\label{eq:1608}
-^{\ast}:\opn{Ch}(A^{\heartsuit})^{op}\overset{\sim}\to \opn{Ch}(A^{\heartsuit}_{\pf}),\ \ (x^{\ast})^n=(x^{-n})^{\ast}.
\end{equation}
This degree-wise duality functor is, however, exactly the inner-Hom functor for the tensor action of $\opn{Ch}(A^{\heartsuit})$ on itself.  We note also that $\opn{Ch}(A^{\heartsuit})$ is itself a tensor category, with finite objects precisely the subcategory
\[
\opn{Ch}(A^{\heartsuit})_{fin}=\opn{Ch}^b(A^{\heartsuit}_{fin})
\]
of bounded complexes of finite length objects.

\begin{lemma}
For any tensor category $A^{\heartsuit}$, the forgetful functor
\[
\opn{Ch}(A^{\heartsuit}_{\pf})\to \opn{Ch}(A^{\heartsuit})
\]
identifies complexes over $A^{\heartsuit}_{\pf}$ with the non-full subcategory $\opn{Ch}(A^{\heartsuit})_{\pf}$ of pro-finite cochains.
\end{lemma}

We claim, in fact, that these two pro-finite subcategories in $\opn{Ch}(A^{\heartsuit})$ are \emph{equal}.

\begin{proof}
The above discussion at \eqref{eq:1608} tells us that the inclusion $\opn{Ch}(A^{\heartsuit}_{\pf})\to \opn{Ch}(A^{\heartsuit})_{\pf}$ is fully faithful and essentially surjective. To see that every object in $\opn{Ch}(A^{\heartsuit})_{\pf}$ actually lies in the image we note first that limits in $\opn{Ch}(A^{\heartsuit})$ are calculated degree-wise.  Second, we note that for any basic inclusion $x^{\lambda}\to x$ into a pro-finite complex, and intermediate subobject $(x^{\lambda})^n\to (x')^n\to x^n$, there is an intermediate subcomplex $x^{\lambda}\to x'\to x$ with $x'$ equal to $(x')^n$ in degree $n$. Thus $x'$ lies in the base as well.
\par

We now have that the base around $0$ for $x$ produces, in each degree, a base around $0$ for $x^n$.  So each $x^n$ lies in $A^{\heartsuit}_{\pf}$.  The fact that the differential is continuous follows from the fact that $x$ is the limit of its finite quotients $x=\varprojlim_{\lambda}x_{\lambda}$.  So we see that $x$ lies in $\opn{Ch}(A^{\heartsuit}_{\pf})$.
\end{proof}

By the materials from Section \ref{sect:pf_obj} we understand that the category $\opn{Ch}(A^{\heartsuit}_{\pf})$ carries a completed product so that the duality equivalence \eqref{eq:1608} is an equivalence of monoidal categories, and of balanced categories when $A^{\heartsuit}$ is ribbon.  Since the product $v\ot w$ of complexes in $\opn{Ch}(A^{\heartsuit})$ has degree $n$ component $(v\ot w)^n=\oplus_{n_1+n_2=n}v^{n_1}\ot w^{n_2}$, it is easy to see that the product on $\opn{Ch}(A^{\heartsuit}_{\pf})$ is given explicit by the formula
\[
x\ot y=\cdots \to (x\hat{\ot} y)^{n-1}\to (x\hat{\ot}y)^n\to \cdots,\ \ \text{with}\ \ (x\hat{\ot}y)^n=\prod_{n_1+n_2=n}x^{n_1}\hat{\ot}y^{n_2}.
\]
In particular, the inclusions from the full subcategories of partially bounded, locally finite complexes
\[
\xymatrixrowsep{3mm}
\xymatrix{
\opn{Ch}^-(A^{\heartsuit}_{fin})\ar[dr] & & \opn{Ch}^+(A^{\heartsuit}_{fin})\ar[dl]\\
	& \opn{Ch}(A^{\heartsuit}_{\pf})
	}
\]
are both monoidal.  We record the following for the sake of referencing.

\begin{corollary}\label{cor:ch_prepair}
Given a ribbon tensor category $A^{\heartsuit}$, duality provides an equivalence of balanced monoidal categories
\[
-^{\ast}:\opn{Ch}(A^{\heartsuit})^{op}\overset{\sim}\to \opn{Ch}(A^{\heartsuit}_{\pf}).
\]
\end{corollary}

\begin{proof}
In this case $\opn{Ch}(A^{\heartsuit})$ is ribbon, under the signed braiding $c_{xy}=\oplus_{n,m}(-1)^{nm}c^{\heartsuit}_{x^ny^m}$ and balancing structure $\theta_x=\oplus_n \theta^{\heartsuit}_{x^n}$.  So this is just a particular instance of Lemma \ref{lem:1572}.
\end{proof}

The following bit of homological algebra is fundamental.

\begin{proposition}\label{prop:tensor_qiso}
For any tensor category $A^{\heartsuit}$, the monoidal products on $\opn{Ch}(A^{\heartsuit})$ and $\opn{Ch}(A^{\heartsuit}_{\pf})$ preserve quasi-isomorphisms in each factor.
\end{proposition}

\begin{proof}
Via duality, it suffices to prove the result for the product on $\opn{Ch}(A^{\heartsuit})$. Since the tensor functor $x\ot -$ for generic $x$ in $\opn{Ch}(A^{\heartsuit})$ commutes with mapping cones, it suffices to show that $x\ot y$ is acyclic whenever $y$ is acyclic.  For this we write $y$ as a directed colimit $y=\varinjlim_n y(n)$ for the bounded subcomplexes
\[
y(n)=\cdots \to 0\to y^{-n}\to y^{-n+1}\to \cdots\to y^{n-1}\to Z^n(y)\to 0\to \cdots.
\]
Since the product $x\ot -$ commutes with arbitrary colimits we have $x\ot y=\varinjlim_nx\ot y(n)$.
\par

For each constituent complex $x\ot y(n)$ the inclusion $x\ot Z^{-n}(y)\to x\ot y(n)$ is a quasi-isomorphism.  This is clear by a simple spectral sequence argument using the filtration on $x\ot y(n)$ provided by the degree on $x$ \cite[Theorem 5.5.11]{weibel95}.
\par

We note that the degree $m$ cocycles in $x\ot Z^{-n}(y)$ are precisely the subspace $Z^{m+n}(x)\ot Z^{-n}(y)$, and that in $x\ot y(n+1)$ these cocycles are all bounded by the degree $m-1$ subspace $Z^{m+n}(x)\ot y^{-n-1}\subseteq (x\ot y)^{m-1}$. Hence the sequence
$x\ot Z^{-n}(y)\to x\ot y(n)\to x\ot y(n+1)$ induces the $0$ map on cohomology
\[
\xymatrix{
	& H^{\ast}\big(x\ot y(n)\big)\ar[dr]\\
H^{\ast}\big(x\ot Z^{-n}(y)\big)\ar[rr]_0\ar[ur]^{\cong} & & H^{\ast}\big(x\ot y(n+1)\big),
}
\]
and it follows that each inclusion $x\ot y(n)\to x\ot y(n+1)$ induces the $0$ map on cohomology. Thus we conclude that the colimit $\varinjlim_nH^{\ast}\big(x\ot y(n)\big)$ vanishes.
\par

Now, any finitely generated abelian category (Definition \ref{def:lin_ab}) is Grothendieck abelian, so that directed colimits are exact.  Hence, from the expressions
\[
B^m\big(x\ot y\big)=\varinjlim_n B^m\big(x\ot y(n)\big)\ \ \text{and}\ \ Z^m\big(x\ot y)=\varinjlim_n Z^m\big(x\ot y(n)\big)
\]
and the exact sequences
\[
0\to B^m\big(x\ot y(n)\big)\to Z^m\big(x\ot y(n)\big)\to H^m\big(x\ot y(n)\big)\to 0,
\]
we conclude in this case that, for all $m\in \mbb{Z}$, $H^m(x\ot y)=\varinjlim_nH^{\ast}\big(x\ot y(n)\big)=0$.
\end{proof}

\subsection{Disk categories from ribbon cochains}
\label{sect:fr_cochains}

\begin{proposition}\label{prop:fr_cochains}
For any ribbon tensor category $A^{\heartsuit}$, with associated ribbon tensor category $A=\opn{Ch}(A^{\heartsuit})$ of unbounded cochains, the pair of balanced monoidal categories $A^-=A^{op}$ and $A^+=A_{\pf}$, along with the balanced duality equivalence
\[
-^{\ast}:A^{op}\overset{\sim}\to A_{\pf},
\]
determines a $fr\Disk$-monoidal $\infty$-category $q:A^{\ot}\to fr\Disk^{\ot}$.  The fibers of $A^{\ot}$ over colored sets are the colored exponents $A^{\ot}_{(\mcl{I},\mu)}=A^{\mcl{I}}$ and the transport functor along any map $f:(\mcl{I},\mu)\to (\mcl{J},\nu)$ in $Linfr\Disk^{\ot}$ is the associated colored product functor $m_f:A^{\mcl{I}}\to A^{\mcl{J}}$ from \eqref{eq:m_f}.
\end{proposition}

The $fr\Disk$-category of unbounded cochains contains all other homologically oriented $fr\Disk$-categories of interest.

\begin{corollary}\label{cor:1714}
For the following full subcategories $A^{\ot}_{\star}\subseteq A^{\ot}$, the structure map $q:A^{\ot}\to fr\Disk^{\ot}$ restricts to a cocartesian fibration $q_{\star}:A^{\ot}_{\star}\to fr\Disk^{\ot}$ which gives $A^{\ot}_{\star}$ the structure of a $fr\Disk$-monoidal $\infty$-category.
\begin{itemize}
\item The full subcategory $A^{\ot}_{fin}$ spanned by the finite length complexes $A^-_{fin}=\opn{Ch}^b(A^{\heartsuit}_{fin})^{op}$ and $A^+_{fin}=\opn{Ch}^b(A^{\heartsuit}_{fin})$, and their higher powers $A^{\mcl{I}}_{fin}\subseteq A^{\mcl{I}}$.\vspace{1mm}
\item The full subcategory $A^{\ot}_{l\text{-}fin}$ spanned by partially bounded locally finite complexes $A^-_{l\text{-}fin}=\opn{Ch}^-(A^{\heartsuit}_{fin})^{op}$ and $A^+_{l\text{-}fin}=\opn{Ch}^+(A^{\heartsuit}_{fin})$, and their higher powers $A^{\mcl{I}}_{l\text{-}fin}\subseteq A^{\mcl{I}}$.\vspace{1mm}
\item The full subcategory $A^{\ot}_{conn}$ spanned by (co)connective complexes $A^-_{conn}=\opn{Ch}^{\leq 0}(A^{\heartsuit})^{op}$ and $A^+_{conn}=\opn{Ch}^{\geq 0}(A^{\heartsuit}_{\pf})$, and their higher powers $A^{\mcl{I}}_{conn}\subseteq A^{\mcl{I}}$.
\end{itemize}
Furthermore, in each case the inclusion $A^{\ot}_{\star}\to A^{\ot}$ is a $fr\Disk$-monoidal functor.
\end{corollary}

\begin{proof}
Apply Proposition \ref{prop:frdisk_subcats} to the $fr\Disk$-category $A^{\ot}\to fr\Disk^{\ot}$.
\end{proof}

By the explicit construction of the $fr\Disk$-category associated to a triple as in Proposition \ref{prop:a_pm}, we also have the apparent $fr\Disk$-monoidal inclusion $(A^{\heartsuit})^{\ot}\to A^{\ot}$ which identifies $(A^{\heartsuit})^{\ot}$ with the full subcategory in $A^{\ot}$ spanned by complexes concentrated in degree $0$.

\begin{lemma}
The inclusions $(A^{\heartsuit})^{\mcl{I}}\to A^{\mcl{I}}$ extend to a $fr\Disk$-monoidal functor
\[
\xymatrix{
(A^{\heartsuit})^{\ot}\ar[rr]\ar[dr] & & A^{\ot}\ar[dl]\\
	& fr\E_2 & 
}
\]
which identifies $(A^{\heartsuit})^{\ot}$ with the full subcategory in $A^{\ot}$ spanned by complexes which are concentrated in degree $0$.
\end{lemma}

\begin{notation}
Given a ribbon tensor category $A^{\heartsuit}$, $q:A^{\ot}\to fr\Disk^{\ot}$ always denotes the $fr\Disk$-monoidal $\infty$-category of unbounded cochains, as in Proposition \ref{prop:fr_cochains}.
\end{notation}

\section{Anomalous bordisms with topology}
\label{sect:bord_star}

We construct a symmetric monoidal $\infty$-category $\Bord_{\ast}^{nc}$ of non-compact, anomalous, trivially marked bordisms in $3$-dimensions.  In Section \ref{sect:marked_bords} we use this $\infty$-category, and a canonical symmetric monoidal functor $\partial_{out}:\Bord_{\ast}^{nc}\to fr\Disk$, to construct a symmetric monoidal $\infty$-category $\Bord^{nc}_{\msc{E}}$ of $3$-dimensional $\msc{E}$-marked bordisms which we associate to any $fr\Disk$-monoidal $\infty$-category $q:\msc{E}^{\ot}\to fr\Disk^{\ot}$.
\par

The majority of this section is spent dealing with various technical minutiae and general tomfoolery.  So we actually suggest reading the synopsis provided in Section \ref{sect:bord_principles}, then briskly skimming the remaining contents. The conscientious reader can return to the section at a later point to digest the details, as the need arises.

\subsection{Marked bordisms, in principle}
\label{sect:bord_principles}

At a basic level, our trivially marked bordism $\infty$-category $\Bord_{\ast}^{nc}$ appears as follows: Objects are surfaces marked by colored disks $\zeta:D\times \mcl{I}\to \Sigma$.  Maps between such marked surfaces $M_f:\Sigma_{\zeta}\to \Sigma'_{\eta}$ are $3$-dimensional bordisms $M:\Sigma\to \Sigma'$ which are equipped with an additional oriented embedding of cylinders
\[
f:\opn{Cyl}\times \mcl{I}\to M.
\]
This embedding preserves boundaries, and hence defines two disk embeddings at the boundaries $f_{\Sigma}:D\times \mcl{I}\to\Sigma$ and $f_{\Sigma'}:D\times \mcl{J}\to \Sigma'$.  We require that $f_{\Sigma}$ recovers the given markings $\zeta$ on $\Sigma$, and that $f_{\Sigma'}$ maps \emph{into} the image of the markings $\eta$ on $\Sigma'$. So, at the outgoing boundary we recover an embedding of colored disks $\partial_{out}(f):D\times\mcl{I}\to D\times \mcl{J}$, i.e.\ a map in $fr\Disk$.
\par

For the higher structure, we allow the cylinders $f:\opn{Cyl}\times \mcl{I}\to M$ to move within $M$ under homotopy.  In this way, an $n$-simplex in $\Bord_{\ast}^{nc}$ traces out an $n$-simplex in $fr\Disk$ at the outgoing boundary, and we observe a map of $\infty$-categories $\partial_{out}:\Bord_{\ast}^{nc}\to fr\Disk$.  Furthermore, $\Bord_{\ast}^{nc}$ admits the expected symmetric monoidal structure under disjoint union, and the map $\partial_{out}$ is the linear term of a corresponding symmetric monoidal functor to the $\infty$-category of framed colored disks
\[
\partial_{out}^{\ot}:(\Bord_{\ast}^{nc})^{\ot}\to fr\Disk^{\ot}.
\]
\par

Below we take special care to accomodate the anomaly, deal appropriately with collars around the boundary, mod out by diffeomorphism in dimension $3$, etc. However, for almost all of what is done in this paper, the above description of the bordism category will suffice. So, upon first engagement, the reader might only skim the contents of this section, then proceed directly to subsequent materials.

\subsection{Marked surfaces and bordisms}
\label{sect:bordisms}

By a marked surface $\Sigma_{\zeta}$ we mean a compact oriented surface $\Sigma$ with a choice of an oriented embedding
\[
\zeta:D\times\mcl{I}\to \Sigma.
\]
A bordism $M_{f}:\Sigma_{\zeta}\to \Sigma'_{\eta}$ between marked surfaces is an oriented $3$-dimensional bordism between the underlying surfaces
\[
M,\ \phi:\Sigma\amalg \bar{\Sigma}'\overset{\sim}\to \partial M,
\]
along with a choice of collar around the boundary
\[
\tilde{\phi}:\Sigma\times[0,\delta)\amalg \bar{\Sigma}'\times(1-\delta,1]\hookrightarrow M
\]
and a choice of embedding $f:\opn{Cyl}_a\times\mcl{I}\to M$ from cylinders $\opn{Cyl}_a=D\times[0,a]$, with $a>0$, whose intersections at the boundaries
\[
p_2:(\opn{Cyl}_a\times\mcl{I})\times_M\Sigma=D\times\mcl{I}\to \Sigma\ \ \text{and}\ \ p_2':(\opn{Cyl}_a\times\mcl{I})\times_M\Sigma'=D\times\mcl{I}\to \Sigma'
\]
recover the markings on the incoming boundary, and send $D\times\mcl{I}$ into the image of the markings on the outgoing boundary.
\par

We consider the collar on $M$ only up to equivalence, where two collars are equivalent if they are equalized by a third, i.e.\ if there is a diagram
\[
\xymatrixcolsep{-1mm}
\xymatrix{
	& \Sigma\times[0,\varepsilon)\amalg \bar{\Sigma}'\times(1-\varepsilon,1]\ar[dl]_{incl}\ar[dr]^{incl}\\
\Sigma\times[0,\delta)\amalg \bar{\Sigma}'\times(1-\delta,1]\ar[dr]_{\tilde{\phi}} & & \Sigma\times[0,\delta')\amalg \bar{\Sigma}'\times(1-\delta',1]\ar[dl]^{\tilde{\phi}'}\\
	& M
}
\]
at some $\varepsilon\leq \delta,\delta'$.  We also assume that $f:\opn{Cyl}_a\times\mcl{I}\to M$ is constant towards the boundary of $M$ in the sense that
\[
f|_{D\times[0,\varepsilon)\times\mcl{I}}:(z,t,i)\mapsto (\zeta(z,i),t)\ \ \in\ \Sigma\times[0,\varepsilon)\ \subseteq M
\]
\[
\text{and}\ \ f|_{D\times(a-\varepsilon,a]\times\mcl{I}}:(z,t,i)\mapsto (\zeta(z,i),1-a+t)\ \ \in\ \bar{\Sigma}'\times(1-\varepsilon,1]\ \subseteq M
\]
for some small $\varepsilon>0$.

\subsection{Mapping complexes at a fixed bordism}
\label{sect:maps1}

Given marked surfaces $\zeta:D\times\mcl{I}\to \Sigma$ and $\eta:D\times\mcl{J}\to\Sigma'$, and an oriented bordism $M:\Sigma\to \Sigma'$ (with fixed collar), we take
\[
\uHom(\Sigma_{\zeta},\Sigma'_{\eta};M)_a
\]
the subcomplex in $\Maps(\opn{Cyl}_a\times\mcl{I},M)$ whose $n$-simplices consist of continuous maps $\sigma:|\Delta^n|\times \opn{Cyl}_a\times \mcl{I}\to M$ which evaluate at each $t\in |\Delta^n|$ to a smooth embedding $\sigma_t:\opn{Cyl}_a\times\mcl{I}\to M$ whose behaviors around the boundary are as prescribed in Section \ref{sect:bordisms}.
\par

For each smooth oriented isomorphism $w:[0,a]\overset{\sim}\to [0,b]$ with $\partial w/\partial t=1$ around the boundaries we have the induced oriented isomorphism $\opn{Cyl}_a\overset{\sim}\to \opn{Cyl}_b$ which induces a further isomorphism of simplicial complexes
\[
w^{\ast}:\uHom(\Sigma_{\zeta},\Sigma_{\eta};M)_b\overset{\cong}\to \uHom(\Sigma_{\zeta},\Sigma_{\eta};M)_a.
\]
For $\opn{Param}$ the opposite category of positive real numbers and oriented isomorphisms $w:[0,a]\overset{\sim}\to [0,b]$ as above we then have the functor
\[
\uHom(\Sigma_{\zeta},\Sigma_{\eta};M)_?:\opn{Param}\to \opn{Kan}.
\]

\begin{definition}
Given marked surfaces $\zeta:D\times\mcl{I}\to \Sigma$ and $\eta:D\times\mcl{J}\to\Sigma'$, and an oriented bordism $M:\Sigma\to \Sigma'$, we take
\[
\uHom(\Sigma_{\zeta},\Sigma'_{\eta};M)=\varinjlim_{\opn{Param}}\uHom(\Sigma_{\zeta},\Sigma'_{\eta};M)_a.
\]
\end{definition}

Note that $n$-simplices in the Hom complex are, explicitly, equivalence classes of $n$-simplices $\sigma:|\Delta^n|\times \opn{Cyl}_a\times \mcl{I}\to M$, where two such simplices are equivalent if they are related by such an oriented reparametrization in the linear factor
\[
\xymatrix{
\opn{Cyl}_a\times\mcl{I}\ar[rr]^{\cong}\ar[dr] & & \opn{Cyl}_b\times \mcl{I}\ar[dl]\\
	& M & .
}
\]
Due to this freedom of reparametrization, we generally consider maps $\opn{Cyl}\times\mcl{I}\to M$ only from the standard cylinder $\opn{Cyl}=\opn{Cyl}_1=D\times[0,1]$.

\subsection{Isomorphisms between bordisms}

Fix surfaces $\Sigma$ and $\Sigma'$.  We define the groupoid $\Bord(\Sigma,\Sigma')$ whose objects are oriented bordisms $M:\Sigma\to \Sigma'$, with fixed equivalence class of a collar, and whose morphisms are oriented diffeomorphisms $\chi:M\to M'$ which are constant towards the boundary.  More precisely, for choices of collars $\opn{Collar}\to M$ and $\opn{Collar}'\to M'$, we require the existence of a subcollar
\[
\begin{tikzpicture}
\draw (0,.8) node {$\opn{Collar}_{\varepsilon}$};
\draw (0,.45) node {$=$};
\draw (-3.5,0) node {$\opn{Collar}$};
\draw (0,0) node {$\Sigma\times[0,\varepsilon)\amalg\bar{\Sigma}'\times(1-\varepsilon,1]$};
\draw (3.5,0) node {$\opn{Collar}'$};
\draw [thick, <-] (-2.9,0) -- (-2.1,0);
\draw [thick, ->] (2.1,0) -- (2.9,0);
\end{tikzpicture}
\]
along which $\chi|_{\opn{Collar}_{\varepsilon}}=id_{\opn{Collar}_{\varepsilon}}$.

For example, for any bordism $M:\Sigma\to \Sigma'$, the resulting bordisms
\[
(\Sigma\times[-1,0])\coprod_{\Sigma}M\ \ \text{and}\ \ M\coprod_{\Sigma'}(\Sigma'\times [1,2])
\]
obtained from gluing are both isomorphic to $M$ in $\Bord(\Sigma,\Sigma')$.

We note that each isomorphism $\chi:M\to M'$ in $\Bord(\Sigma,\Sigma')$ determines an isomorphism between the fixed mapping complexes
\[
\chi_{\ast}:\uHom(\Sigma_{\zeta},\Sigma'_{\eta};M)\overset{\sim}\to \uHom(\Sigma_{\zeta},\Sigma'_{\eta};M').
\]

\begin{definition}\label{def:hom_bord}
For marked surfaces $\Sigma_{\zeta}$ and $\Sigma_{\eta}$, we define the Hom complex
\[
\uHom_{\underline{\mfk{Bord}}_{\ast}}(\Sigma_{\zeta},\Sigma'_{\eta})=\varinjlim_{\opn{Bord}(\Sigma,\Sigma')}\uHom(\Sigma_{\zeta},\Sigma'_{\eta};M)
\]
and the composition operations
\[
\circ:\uHom_{\underline{\mfk{Bord}}_{\ast}}(\Sigma'_{\eta},\Sigma''_{\xi})\times\uHom_{\underline{\mfk{Bord}}_{\ast}}(\Sigma_{\zeta},\Sigma'_{\eta})\to \uHom_{\underline{\mfk{Bord}}_{\ast}}(\Sigma_{\zeta},\Sigma''_{\xi})
\]
are defined by gluing along the boundary
\begin{equation}\label{eq:1285}
M'_{g}\circ M_{f}=M_{f}\coprod_{\Sigma'_{\eta}} M'_{g}.
\end{equation}
\end{definition}

To clarify what we mean by \emph{gluing along the boundary}, let us consider $n$-simplices
\[
\tau:|\Delta^n|\times \opn{Cyl}_a\times\mcl{I}\to M\ \ \text{and}\ \ \omega:|\Delta^n|\times\opn{Cyl}_b\times\mcl{J}\to M'
\]
in $\uHom(\Sigma_{\zeta},\Sigma'_{\eta})$ and $\uHom(\Sigma'_{\eta},\Sigma''_{\xi})$ respectively. The induced map of disks $\tau_{\Sigma'}:|\Delta^n|\times D\times \mcl{I}\to D\times \mcl{J}$ obtained by intersecting at the outgoing boundary $\Sigma'$ produces an inclusion
\[
\tilde{\tau}_{\Sigma'}:|\Delta^n|\times\opn{Cyl}_b\times\mcl{I}\to |\Delta^n|\times \opn{Cyl}_b\times\mcl{J}
\]
and we restrict along this inclusion to obtain an embedding $\omega|_{\tau}:|\Delta^n|\times \opn{Cyl}_b\times\mcl{I}\to M'$. We have the apparent decomposition
\[
|\Delta^n|\times \opn{Cyl}_{a+b}\times \mcl{I}=(|\Delta^n|\times \opn{Cyl}_a\times \mcl{I})\coprod_{(|\Delta^n|\times D\times \mcl{I})}(|\Delta^n|\times \opn{Cyl}_b\times\mcl{I})
\]
and the two maps
\[
\xymatrixrowsep{3mm}
\xymatrix{
|\Delta^n|\times \opn{Cyl}_a\times\mcl{I}\ar[dr]_{\tau} & & |\Delta^n|\times \opn{Cyl}_b\times\mcl{I}\ar[dl]^{\omega|_{\tau}}\\
	& M\amalg_{\Sigma'}M'
}
\]
now induce a map from the pushout
\[
\omega\circ\tau:|\Delta^n|\times \opn{Cyl}_{a+b}\times \mcl{I}\to M\amalg_{\Sigma'}M'.
\]
The glued bordism $M\amalg_{\Sigma'}M'$ along with this $n$-simplex $\omega\circ\tau$ is the composite $M'_{\omega}\circ M_{\tau}$ of \eqref{eq:1285}.

\subsection{The trivially marked bordism category}

The $\ast$-marked, non-anomalous, simplicial bordism category $\underline{\mfk{Bord}}_{\ast}=\underline{\mfk{Bord}}_{3,2}(\ast)$ is the category of oriented surfaces with colored markings $\zeta:D\times \mcl{I}\to \Sigma$, and mapping spaces as in Definition \ref{def:hom_bord}.  This simplicial category is symmetric monoidal under disjoint union.  Note that we have the symmetric monoidal forgetful functor
\[
\underline{\mfk{Bord}}_{\ast}\to \mfk{Bord}_{3,2}
\]
to the usual (discrete) bordism category of unmarked oriented surfaces and unmarked bordisms. This non-anomalous bordism category is not precisely what we want to consider however, as we want to introduce the anomaly.
\par

First, we have the usual anomalous bordism category $\Bord_{3,2}$ which consists of oriented surfaces $\Sigma$ equipped with a Lagrangian $\lambda\subseteq H_1(\Sigma,\mathbb{R})$ and bordisms $M:\Sigma\to \Sigma'$ equipped with an integer $n_M$.  Composition is defined by gluing bordisms and translating the integer parameters via the Maslov index.  See for example \cite[Section 4.3]{derenzietal23}.
\par

The anomalous bordism category is again symmetric monoidal under disjoint union, and we have the symmetric monoidal forgetful functor $\Bord_{3,2}\to \mfk{Bord}_{3,2}$.  We then construct the marked, anomalous, simplicial bordism category via pullback
\begin{equation}\label{eq:1333}
\xymatrix{
\underline{\Bord}_{\ast}:=\Bord_{3,2}\times_{\mfk{Bord}_{3,2}}\underline{\mfk{Bord}}_{\ast}\ar[rr]\ar[d] & & \underline{\mfk{Bord}}_{\ast}\ar[d]\\
\Bord_{3,2}\ar[rr] & & \mfk{Bord}_{3,2}.
}
\end{equation}
Here the objects are marked surfaces $\Sigma_{\zeta,\lambda}$ with a Lagrangian, and the mapping complexes are the fiber products
\[
\Hom_{\Bord}(\Sigma_{\lambda},\Sigma'_{\lambda'})\times_{\Hom_{\mfk{Bord}}(\Sigma,\Sigma')}\uHom_{\underline{\mfk{Bord}}_{\ast}}(\Sigma_{\zeta},\Sigma_{\eta}).
\]
This is again a symmetric monoidal simplicial category under disjoint union.

\begin{definition}
The anomalous $\ast$-marked bordism $\infty$-category $\Bord_{\ast}=\Bord_{3,2}(\ast)$ is the homotopy coherent nerve of the symmetric monoidal simplicial category
\[
\Bord_{\ast}:=\opn{N}(\underline{\Bord}_{\ast}),
\]
were $\underline{\Bord}_{\ast}$ is as in \eqref{eq:1333}
\end{definition}

\subsection{Anomalous notations!}

{\it Throughout this document we work with anomalous bordisms $\Bord_{\ast}$, but we always ignore the anomaly data both in our notations and in our arguments.} We only work in the anomalous setting because our theories are, quite simply, anomalous.

\subsection{Non-compact bordisms}
\label{sect:nc_bord}

The non-compact bordism category $\Bord^{nc}_{3,2}$ in $\Bord_{3,2}$ is the symmetric monoidal subcategory consisting of all surfaces but only those bordisms $M:\Sigma\to \Sigma'$ in which each component $M_c$ in $M$ has non-vanishing outgoing boundary.

In the marked setting, we have the full symmetric monoidal simplicial subcategory $\underline{\mfk{Bord}}_{\ast}^{surf.nc}$ of surfaces with nontrivial markings.  Specifically, we consider marked surfaces $D\times\mcl{I}\to \Sigma$ in which each component $\Sigma_c$ in $\Sigma$ is marked by at least one disk.  We take the fiber product with the anomalous, unmarked, non-compact bordisms to get the symmetric monoidal simplicial category
\[
\underline{\Bord}_{\ast}^{nc}=\Bord^{nc}_{3,2}\times_{\mfk{Bord}_{3,2}}\underline{\mfk{Bord}}_{\ast}^{surf.nc}.
\]

\begin{definition}\label{def:bord_nc}
The non-compact $\ast$-marked bordism $\infty$-category $\Bord_{\ast}^{nc}$ is the homotopy coherent nerve
\[
\Bord_{\ast}^{nc}:=\opn{N}(\underline{\Bord}_{\ast}^{nc}).
\]
\end{definition}

The $\infty$-category $\Bord_{\ast}^{nc}$ admits the apparent description as the subcategory in $\Bord_{\ast}$ spanned by nontrivially marked surfaces and those bordisms $M_f:\Sigma_{\zeta}\to \Sigma'_{\eta}$ for which each component in the underlying $3$-manifold $M_c\subseteq M$ has non-vanishing outgoing boundary.

\begin{remark}
The term ``non-compact" is more about the possible targets for a symmetric monoidal functor $L:\Bord^{nc}_{\star}\to \msc{V}$.  The point is that coevaluation bordisms $\opn{coev}:\bar{\Sigma}\amalg \Sigma\to \emptyset$ are omitted when we move from $\Bord_{\star}$ to $\Bord^{nc}_{\star}$, so that the value of such a field theory $L$ on a surface is allowed to be non-compact, or more specifically non-rigid, in $\msc{V}$ in the non-compact setting.
\end{remark}

\begin{remark}
When working in the non-compact setting there is an apparent choice of orientation.  Namely, by applying the orientation reversal equivalence $\Bord^{op}_{3,2}\overset{\sim}\to \Bord_{3,2}$ we exchange outgoing and incoming boundaries, and hence exchange a non-compact category of bordisms with non-vanishing outgoing boundary with an alternate non-compact category of bordisms with non-vanishing incoming boundary.  One should think of this incoming vs outgoing choice as largely an issue of preference.  For example, we prefer to allow for vanishing incoming boundary so that $S^2$ is an algebra in $\Bord^{nc}_{3,2}$, rather than a coalgebra, and so that all surfaces in $\Bord^{nc}_{3,2}$ are $S^2$-modules in $\Bord^{nc}_{3,2}$.
\end{remark}

\subsection{Disk arrangements at the boundary}

We have a symmetric monoidal simplicial functor $\underline{d}_{out}:\underline{\mfk{Bord}}_{\ast}\to fr\uDisk$ given by intersecting at the outgoing boundary.  In particular, we send each marked surface $\zeta:D\times\mcl{I}\to \Sigma$ to the colored marking set $\mcl{I}$ in $fr\uDisk$ and for an $n$-simplex
\[
\sigma:|\Delta^n|\times \opn{Cyl}\times \mcl{I}\to M,
\]
with $M$ a bordism $M:\Sigma\to \Sigma'$ between marked surfaces $\Sigma_{\zeta}$ and $\Sigma_{\eta}$, we intersect with the outgoing disks to obtain a well-defined $n$-simplex
\[
\xymatrix{
|\Delta^n|\times D\times \mcl{I}\ar@{-->}[rr]^{\exists!\ \underline{d}_{out}(\sigma)}\ar[d] & & D\times\mcl{J}\ar[d]^{\eta}\\
|\Delta^n|\times \opn{Cyl}\times \mcl{I}\ar[rr]_{\sigma}	& & M
}
\]
We apply the construction from Section \ref{sect:underline_maps}, and the homotopy coherent nerve, or alternatively a construction as in Section \ref{sect:snatched_disks}, to obtain a map of symmetric monoidal $\infty$-categories $d_{out}^{\ot}:\mfk{Bord}^{\ot}_{\ast}\to fr\Disk^{\ot}$. We then restrict along the symmetric monoidal functor $\Bord^{\ot}_{\ast}\to \mfk{Bord}_{\ast}^{\ot}$ to obtain a map from the anomalous bordism category.

\begin{definition}\label{def:partial_out}
The symmetric monoidal functor
\[
\partial_{out}^{\ot}:\Bord_{\ast}^{\ot}\to fr\Disk^{\ot}
\]
is the composite of the projection $\Bord_{\ast}^{\ot}\to \mfk{Bord}_{\ast}^{\ot}$ with the map $d_{out}^{\ot}$ from above. We let
\[
\partial_{out}:\Bord_{\ast}\to fr\Disk
\]
denote the linear term of the functor $\partial^{\ot}_{out}$, i.e.\ the fiber of $\partial^{\ot}_{out}$ over $\{0\}$ in $\opn{Fin}_{\ast}$.
\end{definition}

We usually restrict further to the symmetric monoidal subcategory $(\Bord^{nc}_{\ast})^{\ot}$ of non-compact bordisms, though we still employ the same notation $\partial^{\ot}_{out}$.

\section{An $\infty$-category of marked bordisms}
\label{sect:marked_bords}

Given a $fr\Disk$-monoidal $\infty$-category $q:\msc{E}^{\ot}\to fr\Disk^{\ot}$, we define an associated symmetric monoidal $\infty$-category $\Bord_{\msc{E}}^{nc}$ of $\msc{E}$-marked bordisms.  In the case of cochains $A=\opn{Ch}(A^{\heartsuit})$ over a ribbon tensor category $A^{\heartsuit}$, we describe the bordism $\infty$-category $\Bord_A^{nc}$ in explicit terms. In Section \ref{sect:dis_lrt} we show, subsequently, that results from \cite{derenzietal23,czenkynegron} imply the existence of a symmetric monoidal functor $\opn{L}_{A_{fin}}:\Bord^{nc}_{A_{fin}}\to \opn{Ch}(\opn{Vect})$ whenever the input category $A^{\heartsuit}$ is finite and modular. 

\subsection{Bordisms with $fr\Disk$-markings}

Recall our symmetric monoidal functor $\partial_{out}^{\ot}:\Bord_{\ast}^{\ot}\to fr\Disk^{\ot}$ from the $3$-dimensional bordism $\infty$-category of Section \ref{sect:nc_bord}.

\begin{definition}
Given a $fr\Disk$-monoidal $\infty$-category $q:\msc{E}^{\ot}\to fr\Disk^{\ot}$, the $\msc{E}$-marked bordism category is the fiber product
\[
\Bord_{\msc{E}}^{\ot}:=\Bord_{\msc{E}}^{\ot}\times_{fr\Disk^{\ot}}\msc{E}^{\ot}.
\]
The non-compact $\msc{E}$-marked bordism category is the subcategory
\[
(\Bord_{\msc{E}}^{nc})^{\ot}:=(\Bord^{nc}_{\ast})^{\ot}\times_{\Bord_{\ast}^{\ot}}\Bord_{\msc{E}}^{\ot}.
\]
\end{definition}

\begin{remark}
In our notation $\Bord_{\msc{E}}$, one can take ``$\msc{E}$" to be the underlying triple $(\msc{E}^-,\msc{E}^+,d)$. In our case of special interest $A=\opn{Ch}(A^{\heartsuit})$, for ribbon $A^{\heartsuit}$, we simply confuse the distinction between the triple $(A^{op},A_{\pf},d)$ and $A$ itself, and recognize $\Bord_A^{\ot}$ as the symmetric monoidal $\infty$-category of surfaces and bordisms with labels from $A$ and $A_{\pf}$.
\end{remark}

We note that the forgetful functors $\Bord_{\msc{E}}^{\ot}\to \Bord^{\ot}_{\ast}$ and $(\Bord_{\msc{E}}^{nc})^{\ot}\to (\Bord_{\ast}^{nc})^{\ot}$ are cocartesian fibrations, and hence the corresponding maps to $\Fin_{\ast}$ are cocartesian fibrations.

\begin{proposition}\label{prop:bord_e}
The functor $\Bord_{\msc{E}}^{\ot}\to \opn{Fin}_{\ast}$ endows $\Bord_{\msc{E}}$ with a symmetric monoidal structure under which the cocartesian fibration $\Bord_{\msc{E}}^{\ot}\to \Bord_{\ast}^{\ot}$ is a map of symmetric monoidal $\infty$-categories.  Similarly, $\Bord_{\msc{E}}^{nc}$ is naturally symmetric monoidal and the inclusion $(\Bord_{\msc{E}}^{nc})^{\ot}\to \Bord_{\msc{E}}^{\ot}$ is a symmetric monoidal functor.
\end{proposition}

\begin{proof}
Since the symmetric monoidal functor $\msc{E}^{\ot}\to fr\Disk^{\ot}$ is an isofibration, this result is an immediate consequence of Proposition \ref{prop:fp_symm}.
\end{proof}

\begin{notation}\label{not:labeled_unlabeled}
Given an object $\Sigma_x$, or morphism $M_{\alpha}$ in $\Bord_{\msc{E}}$ we let $\Sigma_{\opn{top}(x)}$ and $M_{\opn{top}(\alpha)}$ denote the corresponding images in $\Bord_{\ast}$ under the forgetful functor/projection $\Bord_{\msc{E}}\to \Bord_{\ast}$. For a morphism $M_f:\Sigma_{\zeta}\to \Sigma'_{\eta}$ in $\Bord_{\ast}$ we let
\begin{equation}\label{eq:2109}
m_f:(\Bord_{\msc{E}})_{\Sigma_{\zeta}}=\msc{E}^{\mcl{I}}\to (\Bord_{\msc{E}})_{\Sigma'_{\eta}}=\msc{E}^{\mcl{J}}
\end{equation}
denote the corresponding transport functor $m_f=(M_f)_!$.
\end{notation}

In the expression \eqref{eq:2109}, $\mcl{I}$ and $\mcl{J}$ are the marking sets for $\Sigma_{\zeta}$ and $\Sigma'_{\eta}$ respectively.  We also note that the transport functor along $M_f$ is identified with transport for $\msc{E}^{\ot}$ along the corresponding map $\partial_{out}(M_f)$ in $fr\Disk$. So $m_f$ is always a type of colored product for $\msc{E}^{\mcl{I}}$, in the language of Section \ref{sect:framed_ribbon}.

\subsection{Describing the $\infty$-category $\Bord_A$}
\label{sect:des_bord_a}

We consider a ribbon tensor category $A^{\heartsuit}$ and its associated $fr\Disk$-category $q:A^{\ot}\to fr\E_2$ of unbounded cochains, as in Proposition \ref{prop:fr_cochains}. This fibration was described explicitly in Section \ref{sect:A_ot_des}.
\par

From our description of $A^{\ot}$ over $fr\Disk^{\ot}$ we deduce an explicit description of the bordism category $\Bord_A$ over $\Bord_{\ast}$.  Specifically, we can describe the equivalent subcategory
\[
Linfr\Disk\times_{fr\Disk}\Bord_{A}=(Linfr\Disk\times_{fr\Disk}\Bord_{\ast})\times_{fr\Disk}A^{\ot}_{\{0\}}
\]
spanned by all objects and all morphisms whose disk arrangements at the outgoing boundary lie in $Linfr\Disk$ (Corollary \ref{cor:whatever}).

Objects in $\Bord_A$ are surfaces with embedded, colored disks which are marked by objects in $A^{\pm}$,
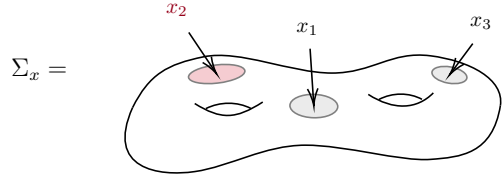
\[
\scalebox{.8}{
\tikzset{every picture/.style={line width=0.75pt}} 
\begin{tikzpicture}[x=0.75pt,y=0.75pt,yscale=-1,xscale=1]
\draw    (128.57,84.03) .. controls (139.09,63) and (150.35,56.82) .. (173.11,52.17) .. controls (195.88,47.52) and (227.83,62.32) .. (253.53,62.36) .. controls (279.22,62.41) and (290.66,48.39) .. (315.39,51.53) .. controls (340.12,54.67) and (361.78,66.18) .. (359.31,82.75) ;
\draw    (128.57,84.03) .. controls (111.25,122.89) and (145.28,127.35) .. (167.55,124.17) .. controls (189.81,120.98) and (211.13,107.65) .. (236.83,107.6) .. controls (262.52,107.56) and (296.45,117.93) .. (315.39,116.52) .. controls (334.32,115.12) and (354.36,113.34) .. (359.31,82.75) ;
\draw    (167.55,81.48) .. controls (184.11,91.69) and (193.09,91.52) .. (210.23,80.54) ;
\draw    (174.55,84.66) .. controls (183.38,80.2) and (194.39,80.2) .. (202.79,84.57) ;
\draw    (276.42,75.74) .. controls (292.26,85.95) and (306.89,85.3) .. (317.24,74.8) ;
\draw    (283.12,78.93) .. controls (291.56,74.47) and (302.1,74.47) .. (310.13,78.84) ;
\draw  [color={rgb, 255:red, 0; green, 0; blue, 0 }  ,draw opacity=0.5 ][fill={rgb, 255:red, 208; green, 2; blue, 27 }  ,fill opacity=0.2 ] (163.67,65.01) .. controls (163.33,61.88) and (170.97,58.44) .. (180.73,57.31) .. controls (190.49,56.19) and (198.68,57.82) .. (199.01,60.95) .. controls (199.35,64.08) and (191.71,67.52) .. (181.95,68.64) .. controls (172.19,69.76) and (164,68.13) .. (163.67,65.01) -- cycle ;
\draw  [color={rgb, 255:red, 0; green, 0; blue, 0 }  ,draw opacity=0.5 ][fill={rgb, 255:red, 0; green, 0; blue, 0 }  ,fill opacity=0.08 ] (227.28,82.67) .. controls (227.38,78.62) and (234.18,75.52) .. (242.46,75.74) .. controls (250.75,75.96) and (257.38,79.42) .. (257.28,83.46) .. controls (257.17,87.51) and (250.38,90.61) .. (242.1,90.39) .. controls (233.81,90.17) and (227.18,86.71) .. (227.28,82.67) -- cycle ;
\draw  [color={rgb, 255:red, 0; green, 0; blue, 0 }  ,draw opacity=0.5 ][fill={rgb, 255:red, 0; green, 0; blue, 0 }  ,fill opacity=0.08 ] (316.41,61.67) .. controls (316.87,58.89) and (322.14,57.48) .. (328.18,58.54) .. controls (334.21,59.6) and (338.74,62.71) .. (338.27,65.5) .. controls (337.81,68.29) and (332.55,69.69) .. (326.51,68.64) .. controls (320.47,67.58) and (315.95,64.46) .. (316.41,61.67) -- cycle ;
\draw    (163.83,36.87) -- (180.23,61.32) ;
\draw [shift={(181.34,62.98)}, rotate = 236.15] [color={rgb, 255:red, 0; green, 0; blue, 0 }  ][line width=0.75]    (10.93,-3.29) .. controls (6.95,-1.4) and (3.31,-0.3) .. (0,0) .. controls (3.31,0.3) and (6.95,1.4) .. (10.93,3.29)   ;
\draw    (239.92,47.07) -- (242.15,81.07) ;
\draw [shift={(242.28,83.07)}, rotate = 266.25] [color={rgb, 255:red, 0; green, 0; blue, 0 }  ][line width=0.75]    (10.93,-3.29) .. controls (6.95,-1.4) and (3.31,-0.3) .. (0,0) .. controls (3.31,0.3) and (6.95,1.4) .. (10.93,3.29)   ;
\draw    (341.99,43.88) -- (328.53,61.98) ;
\draw [shift={(327.34,63.59)}, rotate = 306.62] [color={rgb, 255:red, 0; green, 0; blue, 0 }  ][line width=0.75]    (10.93,-3.29) .. controls (6.95,-1.4) and (3.31,-0.3) .. (0,0) .. controls (3.31,0.3) and (6.95,1.4) .. (10.93,3.29)   ;

\draw (148.12,17.57) node [anchor=north west][inner sep=0.75pt]  [color={rgb, 255:red, 150; green, 1; blue, 26 }  ,opacity=1 ] [align=left] {$\displaystyle x_{2}$};
\draw (229.99,29.15) node [anchor=north west][inner sep=0.75pt]   [align=left] {$\displaystyle x_{1}$};
\draw (338.91,25.78) node [anchor=north west][inner sep=0.75pt]   [align=left] {$\displaystyle x_{3}$};
\draw (70,60) node {{\Large $\Sigma_x =$}};
\draw (380,110) node {{\Large .}};
\end{tikzpicture}}
\]
The positive colored disks are labeled by objects in $\opn{Ch}(A^{\heartsuit}_{\pf})=A_{\pf}=A^+$ and the negatively colored disks are labeled by objects in $\opn{Ch}(A^{\heartsuit})^{op}=A^{op}=A^-$.  We denote such an object by $\Sigma_x$, where $\Sigma$ is the underlying oriented surface and $x=(x_A,\opn{top}(x))$ is the pairing of the tuple of objects $x_A$ in $A^{\mcl{I}}$ with the underlying configuration of disks $\opn{top}(x):D\times\mcl{I}\to \Sigma$.

Morphisms $M_{\alpha}:\Sigma_x\to \Sigma'_y$ in the subcategory $Linfr\Disk\times_{fr\Disk}\Bord_A$ consist of a bordisms $M_{\opn{top}(\alpha)}:\Sigma_{\opn{top}(x)}\to \Sigma_{\opn{top}(y)}$ between the underlying $\ast$-marked surfaces along with morphisms $\alpha_j:x_{I_j}\to y_j$ in $A^{\opn{color}(j)}$ over each disk in the outgoing boundary
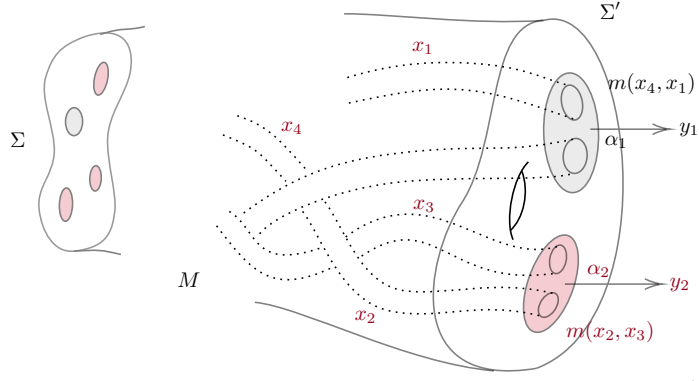
\[
\scalebox{.8}{
\tikzset{every picture/.style={line width=0.75pt}} 
\begin{tikzpicture}[x=0.75pt,y=0.75pt,yscale=-1,xscale=1]
\draw [color={rgb, 255:red, 0; green, 0; blue, 0 }  ,draw opacity=0.5 ]   (411.57,251.46) .. controls (328.39,265.12) and (340.82,184.56) .. (365.05,147.22) .. controls (389.28,109.88) and (380.32,34.87) .. (416.25,33.16) ;
\draw  [color={rgb, 255:red, 0; green, 0; blue, 0 }  ,draw opacity=0.5 ][fill={rgb, 255:red, 208; green, 2; blue, 27 }  ,fill opacity=0.2 ] (410.97,228.72) .. controls (418.79,230.18) and (429.5,217.76) .. (434.9,200.96) .. controls (440.3,184.17) and (438.34,169.37) .. (430.52,167.91) .. controls (422.7,166.45) and (411.99,178.87) .. (406.59,195.67) .. controls (401.19,212.46) and (403.15,227.26) .. (410.97,228.72) -- cycle ;
\draw  [color={rgb, 255:red, 0; green, 0; blue, 0 }  ,draw opacity=0.5 ][fill={rgb, 255:red, 0; green, 0; blue, 0 }  ,fill opacity=0.08 ] (434.41,141.19) .. controls (443.9,141.09) and (451.16,124.2) .. (450.61,103.47) .. controls (450.07,82.75) and (441.93,66.03) .. (432.44,66.13) .. controls (422.95,66.23) and (415.69,83.12) .. (416.23,103.85) .. controls (416.78,124.57) and (424.91,141.29) .. (434.41,141.19) -- cycle ;
\draw    (406.6,121.94) .. controls (396.99,131.22) and (391.54,158.48) .. (397.03,170.74) ;
\draw    (402.47,127.67) .. controls (407.5,140.66) and (403.46,157.83) .. (395.27,164.92) ;
\draw  [dash pattern={on 0.84pt off 2.51pt}]  (293,69) .. controls (313,63) and (337,58) .. (370,60) .. controls (403,62) and (411,68) .. (431.45,73.91) ;
\draw  [dash pattern={on 0.84pt off 2.51pt}]  (294,85) .. controls (320.34,79.43) and (344.79,74.33) .. (369,77) .. controls (393.21,79.67) and (411,87) .. (436.46,95.49) ;
\draw  [color={rgb, 255:red, 0; green, 0; blue, 0 }  ,draw opacity=0.5 ] (431.45,73.91) .. controls (434.88,72.87) and (438.79,76.87) .. (440.18,82.83) .. controls (441.56,88.79) and (439.9,94.46) .. (436.46,95.49) .. controls (433.03,96.53) and (429.12,92.53) .. (427.74,86.57) .. controls (426.35,80.61) and (428.01,74.94) .. (431.45,73.91) -- cycle ;
\draw  [dash pattern={on 0.84pt off 2.51pt}]  (238,167) .. controls (248.08,155.09) and (283.13,139.88) .. (322,135) .. controls (360.87,130.12) and (418.88,134.57) .. (434.9,129.1) ;
\draw  [dash pattern={on 0.84pt off 2.51pt}]  (228,154) .. controls (263,128) and (292,119) .. (339,115) .. controls (386,111) and (391,117) .. (437.03,107.24) ;
\draw  [color={rgb, 255:red, 0; green, 0; blue, 0 }  ,draw opacity=0.5 ] (437.03,107.24) .. controls (441.14,107.76) and (443.99,113.07) .. (443.4,119.11) .. controls (442.81,125.14) and (439.01,129.61) .. (434.9,129.1) .. controls (430.79,128.58) and (427.94,123.26) .. (428.53,117.23) .. controls (429.12,111.19) and (432.92,106.72) .. (437.03,107.24) -- cycle ;
\draw  [color={rgb, 255:red, 0; green, 0; blue, 0 }  ,draw opacity=0.5 ] (427.19,174.24) .. controls (429.93,175.04) and (431.27,179.63) .. (430.17,184.5) .. controls (429.08,189.38) and (425.96,192.68) .. (423.22,191.88) .. controls (420.48,191.08) and (419.14,186.48) .. (420.24,181.61) .. controls (421.33,176.74) and (424.45,173.44) .. (427.19,174.24) -- cycle ;
\draw  [color={rgb, 255:red, 0; green, 0; blue, 0 }  ,draw opacity=0.5 ] (424.11,204.85) .. controls (426.18,206.8) and (425.62,211.45) .. (422.87,215.23) .. controls (420.11,219.02) and (416.2,220.5) .. (414.13,218.55) .. controls (412.07,216.6) and (412.62,211.95) .. (415.38,208.17) .. controls (418.14,204.38) and (422.05,202.9) .. (424.11,204.85) -- cycle ;
\draw  [dash pattern={on 0.84pt off 2.51pt}]  (265,153) .. controls (277.12,169) and (289.29,209.46) .. (320.24,215.56) .. controls (351.19,221.66) and (374.73,208.54) .. (414.13,218.55) ;
\draw  [dash pattern={on 0.84pt off 2.51pt}]  (279,147) .. controls (291.12,163) and (301.32,199.02) .. (328.46,202.69) .. controls (355.61,206.36) and (384.71,194.84) .. (424.11,204.85) ;
\draw  [dash pattern={on 0.84pt off 2.51pt}]  (224,93) .. controls (242,99) and (248,105) .. (266.58,126.06) ;
\draw  [dash pattern={on 0.84pt off 2.51pt}]  (215,105) .. controls (237,114) and (235,113) .. (254.57,134.13) ;
\draw [color={rgb, 255:red, 0; green, 0; blue, 0 }  ,draw opacity=0.5 ]   (445.47,101.41) -- (493.4,101.41) ;
\draw [shift={(495.4,101.41)}, rotate = 180] [color={rgb, 255:red, 0; green, 0; blue, 0 }  ,draw opacity=0.5 ][line width=0.75]    (10.93,-3.29) .. controls (6.95,-1.4) and (3.31,-0.3) .. (0,0) .. controls (3.31,0.3) and (6.95,1.4) .. (10.93,3.29)   ;
\draw [color={rgb, 255:red, 0; green, 0; blue, 0 }  ,draw opacity=0.5 ]   (429.5,198.58) -- (488.61,198.58) ;
\draw [shift={(490.61,198.58)}, rotate = 180] [color={rgb, 255:red, 0; green, 0; blue, 0 }  ,draw opacity=0.5 ][line width=0.75]    (10.93,-3.29) .. controls (6.95,-1.4) and (3.31,-0.3) .. (0,0) .. controls (3.31,0.3) and (6.95,1.4) .. (10.93,3.29)   ;
\draw  [dash pattern={on 0.84pt off 2.51pt}]  (292,169) .. controls (347,128) and (382.77,190.03) .. (427.19,174.24) ;
\draw  [dash pattern={on 0.84pt off 2.51pt}]  (302,181) .. controls (329.96,162.34) and (343,173) .. (365.99,183.81) .. controls (388.99,194.63) and (406.73,193.92) .. (423.22,191.88) ;
\draw  [dash pattern={on 0.84pt off 2.51pt}]  (219,153) .. controls (237,164) and (232.65,190.64) .. (276,176) ;
\draw  [dash pattern={on 0.84pt off 2.51pt}]  (211,162) .. controls (233,191) and (246,206) .. (283.54,188.14) ;
\draw [color={rgb, 255:red, 0; green, 0; blue, 0 }  ,draw opacity=0.5 ]   (291,29) .. controls (322,29) and (371.75,36.16) .. (416.25,33.16) ;
\draw [color={rgb, 255:red, 0; green, 0; blue, 0 }  ,draw opacity=0.5 ]   (235,210) .. controls (265,213) and (329,248) .. (385,253) ;
\draw [color={rgb, 255:red, 0; green, 0; blue, 0 }  ,draw opacity=0.4 ]   (128.26,43.39) .. controls (112.63,49.09) and (107.83,55.52) .. (103.84,68.7) .. controls (99.85,81.89) and (109.83,100.83) .. (109.19,115.82) .. controls (108.56,130.8) and (98.03,137.18) .. (99.67,151.67) .. controls (101.32,166.16) and (109.17,179.03) .. (121.33,177.93) ;
\draw [color={rgb, 255:red, 0; green, 0; blue, 0 }  ,draw opacity=0.4 ]   (128.26,43.39) .. controls (157.09,34.09) and (159.46,54.03) .. (156.56,66.95) .. controls (153.65,79.87) and (143.36,92.03) .. (142.66,107.02) .. controls (141.96,122) and (148.65,142) .. (147.13,153.02) .. controls (145.61,164.04) and (143.79,175.68) .. (121.33,177.93) ;
\draw  [color={rgb, 255:red, 0; green, 0; blue, 0 }  ,draw opacity=0.5 ][fill={rgb, 255:red, 208; green, 2; blue, 27 }  ,fill opacity=0.2 ] (140.46,59.46) .. controls (138.19,59.2) and (135.47,63.58) .. (134.4,69.25) .. controls (133.33,74.92) and (134.31,79.73) .. (136.58,79.99) .. controls (138.86,80.26) and (141.57,75.87) .. (142.64,70.2) .. controls (143.71,64.53) and (142.74,59.72) .. (140.46,59.46) -- cycle ;
\draw  [color={rgb, 255:red, 0; green, 0; blue, 0 }  ,draw opacity=0.5 ][fill={rgb, 255:red, 0; green, 0; blue, 0 }  ,fill opacity=0.08 ] (121.7,87.93) .. controls (118.75,87.91) and (116.31,91.81) .. (116.25,96.64) .. controls (116.2,101.48) and (118.55,105.42) .. (121.51,105.44) .. controls (124.46,105.47) and (126.9,101.57) .. (126.96,96.73) .. controls (127.01,91.9) and (124.66,87.96) .. (121.7,87.93) -- cycle ;
\draw  [color={rgb, 255:red, 0; green, 0; blue, 0 }  ,draw opacity=0.5 ][fill={rgb, 255:red, 208; green, 2; blue, 27 }  ,fill opacity=0.2 ] (116.78,138.28) .. controls (114.49,138.39) and (112.52,143.16) .. (112.38,148.92) .. controls (112.24,154.69) and (113.98,159.28) .. (116.27,159.17) .. controls (118.55,159.06) and (120.52,154.3) .. (120.66,148.53) .. controls (120.8,142.76) and (119.06,138.17) .. (116.78,138.28) -- cycle ;
\draw  [color={rgb, 255:red, 0; green, 0; blue, 0 }  ,draw opacity=0.5 ][fill={rgb, 255:red, 208; green, 2; blue, 27 }  ,fill opacity=0.2 ] (135.26,124.31) .. controls (133.26,124.4) and (131.55,128.03) .. (131.44,132.41) .. controls (131.33,136.79) and (132.87,140.27) .. (134.87,140.17) .. controls (136.88,140.08) and (138.59,136.45) .. (138.69,132.07) .. controls (138.8,127.69) and (137.26,124.21) .. (135.26,124.31) -- cycle ;
\draw [color={rgb, 255:red, 128; green, 128; blue, 128 }  ,draw opacity=1 ]   (416.25,33.16) .. controls (471,18) and (495,220) .. (411.57,251.46) ;
\draw [color={rgb, 255:red, 128; green, 128; blue, 128 }  ,draw opacity=1 ]   (138,42) .. controls (153,37) and (154,40) .. (166,37) ;
\draw [color={rgb, 255:red, 128; green, 128; blue, 128 }  ,draw opacity=1 ]   (121.33,177.93) .. controls (138,177) and (138,177) .. (151,180) ;

\draw (500,96) node [anchor=north west][inner sep=0.75pt]  [rotate=-0.19] [align=left] {$\displaystyle y_{1}$};
\draw (494,193) node [anchor=north west][inner sep=0.75pt]  [color={rgb, 255:red, 150; green, 1; blue, 26 }  ,opacity=1 ,rotate=-0.45] [align=left] {$y_{2}$};
\draw (455,65) node [anchor=north west][inner sep=0.75pt]   [align=left] {$\displaystyle m(x_4 ,x_1)$};
\draw (332.13,45) node [anchor=north west][inner sep=0.75pt]  [color={rgb, 255:red, 150; green, 1; blue, 26 }  ,opacity=1 ]  [align=left] {$\displaystyle x_{1}$};
\draw (250,95) node [anchor=north west][inner sep=0.75pt]  [color={rgb, 255:red, 150; green, 1; blue, 26 }  ,opacity=1 ] [align=left] {$x_{4}$};
\draw (453,104) node [anchor=north west][inner sep=0.75pt]  [rotate=-0.19] [align=left] {$\alpha_{1}$};
\draw (442,185) node [anchor=north west][inner sep=0.75pt]  [color={rgb, 255:red, 150; green, 1; blue, 26 }  ,opacity=1 ,rotate=-0.45] [align=left] {$\alpha _{2}$};
\draw (429,220) node [anchor=north west][inner sep=0.75pt]  [color={rgb, 255:red, 150; green, 1; blue, 26 }  ,opacity=1 ,rotate=-0.45] [align=left] {$m(x_{2} ,x_{3})$};
\draw (296.18,213.62) node [anchor=north west][inner sep=0.75pt]  [color={rgb, 255:red, 150; green, 1; blue, 26 }  ,opacity=1 ] [align=left] {$\displaystyle x_{2}$};
\draw (332.93,144) node [anchor=north west][inner sep=0.75pt]  [color={rgb, 255:red, 150; green, 1; blue, 26 }  ,opacity=1 ] [align=left] {$\displaystyle x_{3}$};
\draw (185,190) node [anchor=north west][inner sep=0.75pt]   [align=left] {$\displaystyle M$};
\draw (450,20) node [anchor=north west][inner sep=0.75pt]   [align=left] {$\displaystyle \Sigma'$};
\draw (80,102) node [anchor=north west][inner sep=0.75pt]   [align=left] {$\displaystyle \Sigma$};
\draw (510,260) node {.};
\end{tikzpicture}}
\]
Here we note that the outgoing arrangement of disks aligns each incoming tuple $D\times \mcl{I}_j$ linearly along the $x$-axis of $D_j\subseteq \Sigma'$, so that the corresponding marking objects $\{x_i:i\in \mcl{I}_j\}$ are naturally ordered via the topology of $M_{\opn{top}(\alpha)}$.  So it makes sense to label each outgoing disk $D_j$ by the chosen morphism $\alpha_j$ from the \emph{ordered} product of the incoming objects to the outgoing objects,
\[
\alpha_j:x_{I_j}=m(x_{i_1},\dots,x_{i_r})\to y_j.
\]

\begin{remark}
This ordered product includes an application of duality whenever we encounters a color transition $\opn{color}(i_a)=-\opn{color}(j)$, as in Section \ref{sect:framed_ribbon}.
\end{remark}

A $2$-simplex
\[
\xymatrix{
	& \Sigma'_y\ar[dr]^{M'_{\delta}}\ar@{}[d]|(.6){\sigma} & &\ar@{}[d]|{\text{\normalsize over}} & & \Sigma'_{\eta}\ar[dr]^{M'_{\opn{top}(\delta)}}\ar@{}[d]|(.6){h=\opn{top}(\sigma)} \\
\Sigma_x\ar[ur]^{M_{\alpha}}\ar[rr]_{M''_{\gamma}} & & \Sigma''_z & & \Sigma_{\zeta}\ar[ur]^{M_{\opn{top}(\alpha)}}\ar[rr]_{M''_{\opn{top}(\gamma)}} & & \Sigma''_{\xi}
}
\]
in $Linfr\Disk\times_{fr\Disk}\Bord_{A}$ consists of an underlying $2$-simplex in $Linfr\Disk\times_{fr\Disk}\Bord_{\ast}$, which is explicitly represented by a choice of a diffeomorphism $\phi:M'\circ M \to M''$ and a homotopy
\[
h:|\Delta^1|\times \opn{Cyl}\times \mcl{I}\to M''
\]
from the cylindrical diagram $\phi(\opn{top}(\delta)\circ \opn{top}(\alpha))$ to the cylindrical diagram $\opn{top}(\gamma)$, along with an algebraic compatibility at the terminal boundary.  For the algebraic compatibility, the homotopy $h$ provides a homotopy $|\Delta^1|\times D\times \mcl{I}_l\to D_l$ between linear disk arrangements at the outgoing boundary $\xi:D\times\mcl{L}\to \Sigma''$, and this boundary homotopy traces out a collection of framed braids $(\beta_l,n_l)$ in the outgoing disks $D_l$.  We require that the associated actions of the framed braid transformation $\zeta_{\beta_l,n_l}$ from \eqref{eq:zeta_beta} completes a diagram
\[
\xymatrixcolsep{12mm}
\xymatrix{
 & y_{J_l}\ar[dr]^{\delta_l} &\\
\bigotimes_{j\in J_l}(x_{I_j})\ar[ur]^{\otimes_j\alpha_{j}}\ar[r]_(.6){\zeta_{\beta_l,n_l}} & x_{I_l}\ar[r]_{\xi_l} & z_l
}
\]
for each index $l$ in $\mcl{L}$.

The symmetric monoidal structure on $\Bord_A$ is given by disjoint union.

\begin{remark}
One should view the allowance of homotopies as the institution of a universal skein relation, where the twistings and rotations at the outgoing boundary are compensated for algebraically via the appropriate application of the ribbon structure on $A^-$ or $A^+$.
\end{remark}

\subsection{Including bounding constraints}

We recall that we have fully faithful inclusions
\[
\xymatrix{
A^{\ot}_{fin}\ar[r]\ar[dr] & A^{\ot}\ar[d] & A^{\ot}_{l\text{-}fin}\ar[dl]\ar[l]\\
	& fr\Disk
}
\]
for the $fr\Disk$-categories of finite and locally finite complexes (Proposition \ref{prop:fr_cochains}).  Pulling back, we obtain full, symmetric monoidal subcategories
\[
\Bord_{A_{fin}}\to \Bord_A\leftarrow \Bord_{A_{l\text{-}fin}}
\]
spanned by surfaces labeled by finite and locally finite cochains, respectively.  In particular, the above explicit description of $\Bord_A$ restricts to provide explicit descriptions of both $\Bord_{A_{fin}}$ and $\Bord_{A_{l\text{-}fin}}$.

\subsection{Comments on general marking schemes}
\label{sect:mark_scheme}

Given a symmetric monoidal $\infty$-category $\msc{T}^{\ot}$, and a corresponding ``$\msc{T}$-flavored" bordism $\infty$-category, by which we mean a symmetric monoidal $\infty$-category of bordisms with a symmetric monoidal functor
\[
\Bord_{triv}^{\ot}\to \msc{T}^{\ot}\times \Bord_{3,2},
\]
we obtain a corresponding class of marked bordisms just as above. Namely, for any symmetric monoidal cocartesian fibration $\msc{A}^{\ot}\to \msc{T}^{\ot}$ we get an $\infty$-category of $\msc{A}$-marked bordisms via pullback
\[
\Bord_{\msc{A}}=\Bord_{triv}\times_{\msc{T}^{\ot}}\msc{A}^{\ot}.
\]
This situation is what we would call a \emph{$3$-dimensional $\msc{T}$-flavored marking scheme}.
\par

We've provided above a $3$-dimensional $fr\Disk$-flavored marking scheme, but certainly there are other natural schemes one might consider. We have some interest, for example, in the case where $\msc{T}$ is the $2$-dimensional bordism category (see \cite{costello,brochierwoike}). Though we operate within the scheme of framed disks, we note that many of our arguments through the remainder of the text apply equally well under any marking scheme.
\par

The main point of contention, as far as our arguments are concerned, is the existence (or non-existence) of an abelian theory for cochains under the given scheme.

\section{Discrete Lyubashenko-Reshetikhin-Turaev theories}
\label{sect:dis_lrt}

For any modular tensor category $A^{\heartsuit}$--which we recall is a special type of ribbon tensor category--we construct a symmetric monoidal functor $\opn{L}_A:\Bord^{nc}_A\to \opn{Ch}(\opn{Vect})_{\pf}$ from the $\infty$-category of $A$-labeled bordisms.  We describe the values $\opn{L}_A(\Sigma_x)$ for a connected marked surface $\Sigma$ via mapping complexes in $A$.
\par

The field theory $\opn{L}_A$ is constructed by completing pre-existing TQFTs from De Renzi-Gainutdinov-Geer-Patureau-Mirand-Runkel \cite{derenzietal23} and Czenky-Negron \cite{czenkynegron}.

\subsection{Discrete LRT in the finite setting}

\begin{theorem}\label{thm:fin_dis_lrt}
Let $A^{\heartsuit}$ be a finite modular tensor category and $A_{fin}^{\ot}\to fr\Disk$ be the associated $fr\Disk$-monoidal $\infty$-category of finite cochains.  The field theories from \cite{derenzietal23,czenkynegron} determine a corresponding symmetric monoidal functor
\[
\dL_{A_{fin}}:\Bord^{nc}_{A_{fin}}\to \opn{Ch}^+(\opn{Vect}_{fin}).
\]
\end{theorem}

\begin{proof}[Proof/Construction]
The inclusion
\[
Linfr\Disk\times_{fr\Disk}\Bord^{nc}_{A_{fin}}\to \Bord^{nc}_{A_{fin}}
\]
is an equivalence (Corollary \ref{cor:whatever}) which is bijective on objects, so that the induced map on homotopy categories is an isomorphism.  Hence the homotopy category $\opn{h}\Bord^{nc}_{A_{fin}}$ is explicitly identified as the discrete category of $A_{fin}$-labeled surfaces $\Sigma_x$ and equivalence classes of bordisms $M_{\alpha}:\Sigma_x\to \Sigma_y$ in which $\partial_{out}(M_{\opn{top}(\alpha)})$ consists of a linear arrangement of disks $D\times \mcl{I}\to D\times \mcl{J}$.
\par

Consider maps $M_{\alpha},M_{\gamma}:\Sigma_x\to \Sigma'_y$ in $Linfr\Disk\times_{fr\Disk}\Bord^{nc}_{A_{fin}}$ with corresponding disk arrangements $f,g:D\times\mcl{I}\to D\times\mcl{J}$ in $Linfr\Disk$. These maps are identified in the homotopy category if and only if there is a smooth homotopy between the constituent cylinder embeddings $\opn{Cyl}\times \mcl{I}\to M$, as in Section \ref{sect:maps1}, for which the resulting homotopy $\sigma:|\Delta^1|\times D\times\mcl{I}\to D\times\mcl{J}$ at the outgoing boundary, and associated framed braid $[\sigma]=(\beta,n)$, produce a diagram
\[
\xymatrix{
 & m_g(x)\ar[dr]^{\gamma_g} \\
m_f(x)\ar[rr]_{\alpha_f}\ar[ur]^{\zeta_{\beta,n}} & & y
}
\]
in $A_{fin}^{\mcl{J}}$.  Here $\zeta_{\beta,n}$ is the transformation provided by applications of the braiding and twist on $\opn{Ch}^b(A^{\heartsuit}_{fin})$, as in \eqref{eq:zeta_beta}.  The symmetric monoidal structure on $\opn{h}\Bord^{nc}_{A_{fin}}$ is given by disjoint union, as expected.
\par

Now for $\mathsf{Bord}_{A_{fin}}^{nc}$ the ``reduced" bordism category from \cite[Definition 6.1]{czenkynegron}, we have a symmetric monoidal functor
\[
p:\opn{h}\Bord^{nc}_{A_{fin}}\to \msf{Bord}_{A_{fin}}^{nc}
\]
which sends each marked surface $\Sigma_x$ to the same surface $\Sigma$ marked by the centers of the corresponding disks $\mcl{I}\to D\times \mcl{I}\to \Sigma$, and with the $i$-th marking labeled by the object $x_i$ in $x=(x_1,\dots,x_r)$.  We send a bordism $M_{\alpha}:\Sigma_x\to \Sigma_y$ to the ribbon bordism $p(M_{\alpha})$ obtained as follows: 
\begin{enumerate}
\item[(a)] We replace the cylindrical diagram $\opn{Cyl}\times \mcl{I}\to M$ with the ribbon diagram $[0,1]\times \mcl{I}\to \opn{Cyl}\times \mcl{I}\to M$ obtained by following the centers of the cylinders, and endow each constituent path with the framing inherited by the framing on the ambient cylinders.\vspace{1mm}
\item[(b)] In a small collar around the outgoing boundary $(1-\delta,1]\times \Sigma'\to M$, and at each marking disk $D_j\to \Sigma'\to M$, we take the ribbons mapping to the $x$-axis in $D_j$
\[
(1-\delta,1]\times \mcl{I}_j\to (1-\delta,1]\times D_j\to M,
\]
and merge them along the $x$-axis into a single string emerging from the center of $D_j$,
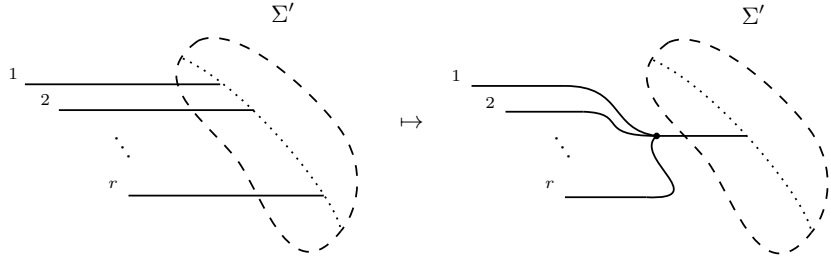
\[
\scalebox{.9}{
\tikzset{every picture/.style={line width=0.75pt}} 
\begin{tikzpicture}[x=0.75pt,y=0.75pt,yscale=-1,xscale=1]
\draw  [dash pattern={on 4.5pt off 4.5pt}]  (192.96,63.39) .. controls (214.74,53.49) and (255.54,90.39) .. (271.56,111.96) .. controls (287.59,133.52) and (281.98,158.73) .. (272.51,169.52) ;
\draw  [dash pattern={on 4.5pt off 4.5pt}]  (192.96,63.39) .. controls (159.81,86.77) and (203.38,111.96) .. (216.64,132.64) .. controls (229.89,153.33) and (245.99,203.7) .. (272.51,169.52) ;
\draw  [dash pattern={on 0.84pt off 2.51pt}]  (184.44,73.28) .. controls (222.32,90.37) and (272.51,155.13) .. (272.51,169.52) ;
\draw    (96.36,87.67) -- (205.27,87.67) ;
\draw    (115.3,102.06) -- (224.21,102.06) ;
\draw    (154.13,149.73) -- (263.04,149.73) ;
\draw  [dash pattern={on 4.5pt off 4.5pt}]  (455.28,64.29) .. controls (477.07,54.39) and (517.86,91.29) .. (533.89,112.86) .. controls (549.91,134.42) and (544.3,159.63) .. (534.83,170.42) ;
\draw  [dash pattern={on 4.5pt off 4.5pt}]  (455.28,64.29) .. controls (422.14,87.67) and (465.7,112.86) .. (478.96,133.54) .. controls (492.22,154.23) and (508.32,204.6) .. (534.83,170.42) ;
\draw  [dash pattern={on 0.84pt off 2.51pt}]  (446.76,74.18) .. controls (484.64,91.27) and (534.83,156.03) .. (534.83,170.42) ;
\draw    (345.43,88.57) -- (397.52,88.57) ;
\draw    (364.37,102.96) -- (407.93,102.96) ;
\draw    (397.52,150.63) -- (442.97,150.63) ;
\draw    (448.66,116.45) -- (498.85,116.45) ;
\draw    (407.93,102.96) .. controls (433.5,103.86) and (413.62,116.45) .. (448.66,116.45) ;
\draw    (397.52,88.57) .. controls (430.66,88.57) and (422.14,112.86) .. (448.66,116.45) ;
\draw    (442.97,150.63) .. controls (480.85,152.43) and (434.45,128.15) .. (448.66,116.45) ;

\draw (232.26,40.89) node [anchor=north west][inner sep=0.75pt]   [align=left] {$\displaystyle \Sigma '$};
\draw (494.58,41.79) node [anchor=north west][inner sep=0.75pt]   [align=left] {$\displaystyle \Sigma '$};
\draw (149,113) node [anchor=north west][inner sep=0.75pt]  [rotate=-55.49] [align=left] {$\displaystyle \cdots $};
\draw (394.35,113) node [anchor=north west][inner sep=0.75pt]  [rotate=-55.49] [align=left] {$\displaystyle \cdots $};
\draw (85.66,77.22) node [anchor=north west][inner sep=0.75pt]  [font=\scriptsize] [align=left] {$\displaystyle 1$};
\draw (103.65,91.62) node [anchor=north west][inner sep=0.75pt]  [font=\scriptsize] [align=left] {$\displaystyle 2$};
\draw (141.53,140.19) node [anchor=north west][inner sep=0.75pt]  [font=\scriptsize] [align=left] {$\displaystyle r$};
\draw (384.91,140.19) node [anchor=north west][inner sep=0.75pt]  [font=\scriptsize] [align=left] {$\displaystyle r$};
\draw (332.83,78.12) node [anchor=north west][inner sep=0.75pt]  [font=\scriptsize] [align=left] {$\displaystyle 1$};
\draw (351.77,93.41) node [anchor=north west][inner sep=0.75pt]  [font=\scriptsize] [align=left] {$\displaystyle 2$};
\draw (303.94,105.65) node [anchor=north west][inner sep=0.75pt]   [align=left] {$\displaystyle \mapsto $};
\filldraw (448.66,116.45) circle (1pt);
\end{tikzpicture}}
\]
\item[(c)] We label the unique vertex generated in this process by the underlying map $\alpha_f:m_f(x)\to y$ for the bordism $M_{\alpha}$, where $f=\partial_{out}(M_{\opn{top}(\alpha)})$.
\end{enumerate}
The fact that $p$ is well-defined on equivalence classes of morphisms follows by the relation \cite[Section 5.2, U3]{czenkynegron} on maps in $\msf{Bord}^{nc}_{A_{fin}}$, and $p$ respects composition by the relation \cite[Section 5.2, U1]{czenkynegron}. 
\par

Having constructed such a symmetric monoidal functor $p:\opn{h}\Bord^{nc}_{A_{fin}}\to \msf{Bord}^{nc}_{A_{fin}}$ we obtain the TQFT $\dL_{A_{fin}}$ by pulling back the TQFT $Z^{\ast}$ from \cite[Theorem 9.3]{czenkynegron} along $p$,
\[
\dL_{A_{fin}}:=\left\{\ \Bord^{nc}_{A_{fin}}\to \opn{h}\Bord^{nc}_{A_{fin}}\overset{p}\to \msf{Bord}^{nc}_{A_{fin}}\overset{Z^{\ast}}\to \opn{Ch}^+(\opn{Vect}_{fin}).\right.
\]
\end{proof}

\subsection{Merging bordisms}
\label{sect:merge}

Let $\zeta:D\times \mcl{I}\to \Sigma$ be a marked connected surface.  A basic merging bordism for $\Sigma_{\zeta}$ is a bordism $M_f:\Sigma_{\zeta}\to \Sigma_{\eta}$ to the same surface $\Sigma$ equipped with a single positive marking $\eta:D\to \Sigma$, and with $M_f$ satisfying the following:
\begin{enumerate}
\item[(0)] $\Sigma$ admits an auxiliary embedding $\tau:D\to \Sigma$ which contains the images of both $\zeta$ and $\eta$ in $\tau(D)$.\vspace{1.5mm}
\item  The underlying $3$-manifold for $M_f$ is the product $\Sigma\times [0,1]$.\vspace{1.5mm}
\item The mappings from the cylinders $f:\opn{Cyl}\times \mcl{I}\to \Sigma\times [0,1]$ factor through the image of the auxiliary disk embedding $\tau:D\to \Sigma$.\vspace{1.5mm}
\item  The corresponding map
\[
f\times_{(\Sigma\times [0,1])}(D\times [0,1]):\opn{Cyl}\times\mcl{I}=D\times [0,1]\times \mcl{I}\to D\times [0,1]
\]
is induced by a smooth homotopy $h:|\Delta^1|\times D\times \mcl{I}\to D$, in the sense that $f(z,t,i)=(h(t,z,i),t)$.
\end{enumerate}
For example, a basic merging bordism might appear as
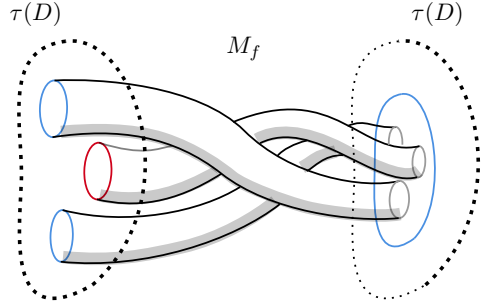
\[
\scalebox{.9}{
\tikzset{every picture/.style={line width=0.75pt}} 
\begin{tikzpicture}[x=0.75pt,y=0.75pt,yscale=-1,xscale=1]
\draw [line width=1.5]  [dash pattern={on 1.69pt off 2.76pt}]  (254.92,113.63) .. controls (218,118.06) and (231.48,161.99) .. (230.16,192.63) .. controls (228.83,223.27) and (222.34,276.79) .. (269.9,249.47) ;
\draw [line width=1.5]  [dash pattern={on 1.69pt off 2.76pt}]  (269.9,249.47) .. controls (323.97,202.22) and (294.22,108.46) .. (254.92,113.63) ;
\draw  [color={rgb, 255:red, 74; green, 144; blue, 226 }  ,draw opacity=1 ][fill={rgb, 255:red, 255; green, 255; blue, 255 }  ,fill opacity=0.18 ][line width=0.75]  (248.7,135.5) .. controls (244.52,135.27) and (240.87,142.17) .. (240.56,150.89) .. controls (240.25,159.62) and (243.38,166.88) .. (247.56,167.11) .. controls (251.74,167.34) and (255.38,160.45) .. (255.69,151.72) .. controls (256.01,142.99) and (252.87,135.73) .. (248.7,135.5) -- cycle ;
\draw  [color={rgb, 255:red, 208; green, 2; blue, 27 }  ,draw opacity=1 ][fill={rgb, 255:red, 255; green, 255; blue, 255 }  ,fill opacity=0.75 ][line width=0.75]  (273.45,170.2) .. controls (269.27,169.97) and (265.63,176.86) .. (265.32,185.59) .. controls (265,194.32) and (268.14,201.58) .. (272.32,201.81) .. controls (276.49,202.04) and (280.14,195.15) .. (280.45,186.42) .. controls (280.76,177.69) and (277.63,170.43) .. (273.45,170.2) -- cycle ;
\draw  [color={rgb, 255:red, 74; green, 144; blue, 226 }  ,draw opacity=1 ][fill={rgb, 255:red, 255; green, 255; blue, 255 }  ,fill opacity=1 ][line width=0.75]  (253.21,207.78) .. controls (249.77,207.59) and (246.76,213.89) .. (246.47,221.84) .. controls (246.19,229.79) and (248.74,236.39) .. (252.17,236.58) .. controls (255.6,236.77) and (258.61,230.47) .. (258.9,222.52) .. controls (259.19,214.57) and (256.64,207.97) .. (253.21,207.78) -- cycle ;
\draw [line width=0.75]  [dash pattern={on 0.84pt off 2.51pt}]  (440.59,113.63) .. controls (403.67,118.06) and (417.16,161.99) .. (415.84,192.63) .. controls (414.51,223.27) and (408.02,276.79) .. (455.58,249.47) ;
\draw [line width=1.5]  [dash pattern={on 1.69pt off 2.76pt}]  (455.58,249.47) .. controls (509.65,202.22) and (479.9,108.46) .. (440.59,113.63) ;
\draw [color={rgb, 255:red, 74; green, 144; blue, 226 }  ,draw opacity=1 ]   (442.28,143.08) .. controls (425.09,145.7) and (427.74,171.9) .. (427.13,190.06) .. controls (426.51,208.22) and (427.11,239.77) .. (449.26,223.58) ;
\draw [color={rgb, 255:red, 74; green, 144; blue, 226 }  ,draw opacity=1 ][line width=0.75]    (449.26,223.58) .. controls (469.47,202.18) and (460.58,140.01) .. (442.28,143.08) ;
\draw [color={rgb, 255:red, 0; green, 0; blue, 0 }  ,draw opacity=1 ]   (248.69,135.5) .. controls (274.49,134.16) and (328.75,137.99) .. (354.47,163.49) .. controls (380.19,188.99) and (420.61,196.32) .. (441.34,191.49) ;
\draw  [color={rgb, 255:red, 155; green, 155; blue, 155 }  ,draw opacity=1 ] (441.34,191.49) .. controls (438.35,191.33) and (435.76,195.85) .. (435.55,201.58) .. controls (435.34,207.32) and (437.6,212.1) .. (440.6,212.26) .. controls (443.59,212.43) and (446.18,207.91) .. (446.39,202.17) .. controls (446.6,196.44) and (444.34,191.66) .. (441.34,191.49) -- cycle ;
\draw    (247.56,167.11) .. controls (284.45,161.62) and (314.9,166.31) .. (341.78,187.46) .. controls (368.66,208.61) and (414.75,214.77) .. (440.6,212.26) ;
\draw  [color={rgb, 255:red, 155; green, 155; blue, 155 }  ,draw opacity=1 ] (452.2,172.15) .. controls (450.23,172.04) and (448.49,175.83) .. (448.32,180.62) .. controls (448.15,185.4) and (449.6,189.37) .. (451.57,189.48) .. controls (453.54,189.59) and (455.28,185.79) .. (455.45,181.01) .. controls (455.62,176.22) and (454.16,172.26) .. (452.2,172.15) -- cycle ;
\draw    (272.32,201.81) .. controls (293.57,206.65) and (315.07,204.44) .. (341.78,187.46) ;
\draw [color={rgb, 255:red, 0; green, 0; blue, 0 }  ,draw opacity=0.5 ]   (273.45,170.2) .. controls (298.61,178) and (305.95,173.43) .. (314.42,171.95) ;
\draw    (363.28,170.48) .. controls (393.25,155.71) and (427.13,195.58) .. (450.87,189.44) ;
\draw [color={rgb, 255:red, 0; green, 0; blue, 0 }  ,draw opacity=1 ]   (348.3,158.67) .. controls (406.28,135.04) and (419.26,179.3) .. (452.19,172.15) ;
\draw    (252.17,236.58) .. controls (299.43,241.35) and (348.95,224.37) .. (374.36,204.44) ;
\draw [color={rgb, 255:red, 0; green, 0; blue, 0 }  ,draw opacity=1 ]   (253.2,207.78) .. controls (289.52,216.92) and (330.05,208.87) .. (352.21,194.84) ;
\draw    (389.34,185.98) .. controls (393.9,183.77) and (398.46,177.12) .. (409.54,177.12) ;
\draw [color={rgb, 255:red, 0; green, 0; blue, 0 }  ,draw opacity=1 ]   (373.08,177.51) .. controls (375.66,174.17) and (382.83,170.48) .. (389.34,169) ;
\draw [color={rgb, 255:red, 0; green, 0; blue, 0 }  ,draw opacity=1 ]   (412.14,160.88) .. controls (422.57,156.45) and (425.82,162.36) .. (439.27,161.5) ;
\draw  [draw opacity=0] (435.15,171.4) .. controls (435.11,170.95) and (435.09,170.5) .. (435.09,170.03) .. controls (435.05,165.25) and (436.71,161.36) .. (438.81,161.33) .. controls (440.9,161.31) and (442.64,165.16) .. (442.68,169.94) .. controls (442.69,170.41) and (442.68,170.86) .. (442.65,171.31) -- (438.89,169.99) -- cycle ; \draw  [color={rgb, 255:red, 155; green, 155; blue, 155 }  ,draw opacity=1 ] (435.15,171.4) .. controls (435.11,170.95) and (435.09,170.5) .. (435.09,170.03) .. controls (435.05,165.25) and (436.71,161.36) .. (438.81,161.33) .. controls (440.9,161.31) and (442.64,165.16) .. (442.68,169.94) .. controls (442.69,170.41) and (442.68,170.86) .. (442.65,171.31) ;  
\draw [color={rgb, 255:red, 0; green, 0; blue, 0 }  ,draw opacity=0.21 ][line width=5.25]    (252.26,164.21) .. controls (270.25,162.67) and (291.52,162.2) .. (309.44,167.98) .. controls (327.37,173.76) and (341.47,184.02) .. (347.76,189.21) .. controls (354.05,194.39) and (377.17,203.12) .. (389.5,206.18) .. controls (401.83,209.24) and (429.53,212.32) .. (441.54,209.96) ;
\draw [color={rgb, 255:red, 0; green, 0; blue, 0 }  ,draw opacity=0.21 ][line width=5.25]    (253.4,234.01) .. controls (281.42,236.84) and (299.72,233.07) .. (320.31,228.35) .. controls (340.89,223.63) and (360.91,211.84) .. (371.77,202.88) ;
\draw [color={rgb, 255:red, 0; green, 0; blue, 0 }  ,draw opacity=0.21 ][line width=5.25]    (273.99,199.58) .. controls (299.72,204.77) and (320.31,197.22) .. (339.18,186.38) ;
\draw [color={rgb, 255:red, 0; green, 0; blue, 0 }  ,draw opacity=0.21 ][line width=5.25]    (361.48,168.45) .. controls (371,162.98) and (387.21,162.32) .. (409.52,174.11) .. controls (431.82,185.9) and (445.88,191.05) .. (451.83,185.9) ;
\draw [color={rgb, 255:red, 0; green, 0; blue, 0 }  ,draw opacity=0.21 ][line width=5.25]    (386.64,184.02) .. controls (396.16,178.55) and (392.36,178.83) .. (405.51,175.06) ;

\draw (222.31,90) node [anchor=north west][inner sep=0.75pt]   [align=left] {$\displaystyle \tau ( D)$};
\draw (446.34,90) node [anchor=north west][inner sep=0.75pt]   [align=left] {$\displaystyle \tau ( D)$};
\draw (343.07,109.82) node [anchor=north west][inner sep=0.75pt]   [align=left] {$\displaystyle M_{f}$};
\draw (500,240) node [anchor=north west][inner sep=0.75pt]   [align=left] {{\large .}};
\end{tikzpicture}}
\]

\begin{definition}
Given two markings $\Sigma_{\zeta}$ and $\Sigma_{\eta}$ of the same underlying surface $\Sigma$, a merging bordism $M_f:\Sigma_{\zeta}\to \Sigma_{\eta}$ is a bordism which decomposes over the components as
\[
M_f=\coprod_{c\in \pi_0(\Sigma)}M_{f^c}^c:\Sigma_{\zeta}=\amalg_c\Sigma^c_{\zeta^c}\to \Sigma_{\eta}=\amalg_c\Sigma^c_{\eta^c}
\]
with each $M^c_{f^c}$ either an isomorphism, or isomorphic to a basic merging bordism in $\Maps_{\Bord_{\ast}}(\Sigma_{\zeta},\Sigma_{\eta})$, defined as above.
\end{definition}

\begin{lemma}
For any marked surface $\Sigma_{\zeta}$ in $\Bord_{\ast}$, there is a choice of single positive markings $\eta:D\times \pi_0(\Sigma)\to \Sigma$ on the components of $\Sigma$ which admits a merging bordism $M_f:\Sigma_{\zeta}\to \Sigma_{\eta}$.  We may assume further that the outgoing diagram $\partial_{out}(M_f)$ lies in $Linfr\Disk$.
\end{lemma}

\begin{proof}
It suffices to deal with the case where $\Sigma$ is connected with specified markings $\zeta:D\times\mcl{I}\to \Sigma$. We consider a path $[0,1]\to \Sigma$ which passes through the centers of the $\zeta(D_i)$, and by taking a tubular neighborhood we obtain an embedding from a disk $\eta:D\to \Sigma$ which contains the centers of the $\zeta(D_i)$ in its interior.  We can take a smooth homotopy $\kappa:|\Delta^1|\times D\times \mcl{I}\to \Sigma$ which shrinks the image of each $D_i$ so that $\zeta'=\kappa_1:D\times \mcl{I}\to \Sigma$ has image in $\eta(D)$.
\par

We consider now $\zeta'$ as a disk embedding $\zeta':D\times \mcl{I}\to D$ by pulling back along $\eta$, and we are free to assume that the images of the disks $D_i$ are arbitrarily small in $D$.  In particular, we assume that each $D_i$ has image in the interior $D^o$.  We consider, for each index $i$, the differential $A_i^{-1}=T_0(\zeta')_i\in \opn{GL}_2(\mbb{R})$. Since $A_i$ has positive determinant we can connect the identity matrix $I_2$ to $A_i$ via a smooth path $A_i(t):|\Delta^1|\to \opn{GL}_2(\mbb{R})$.
\par

We take now $A(t):|\Delta^1|\times D\times \mcl{I}\to D$ with $A(t)|_{D_i}=s_i^{-1}A_i(t)s_i$ where $s_i:\mbb{R}^2\to \mbb{R}^2$ is the translation which sends the center of $\zeta'(D_i)$ to $0$.  Since the disk is compact, we can assume the images $\zeta'(D_i)$ are sufficiently small so that $A(t)$ is a diffeomorphism at each time $t$ which does in fact have image in the unit disk $D\subseteq \mbb{R}^2$.  Then $\zeta''=A(1):D\times \mcl{I}\to D$ is a disk embedding which sends the center of each $D_i$ to the $x$-axis in $D$ and has differentials $T_0(\zeta'')_i=id_{\mbb{R}^2}$.  We are free to assume that $\kappa(t)$ and $A(t)$ are constant in neighborhoods of both $0$ and $1$, so that we obtain smooth homotopies
\[
h^A,h^{\kappa}:|\Delta^1|\times D\times \mcl{I}\to D\to \Sigma
\]
which define respective marked bordisms
\[
f^A,f^{\kappa}:\opn{Cyl}\times \mcl{I}=D\times |\Delta^1|\times\mcl{I}\to \Sigma\times|\Delta^1|,\ \ f^{\star}(z,t,i)\mapsto (h^{\star}(z,t,i),t).
\]
The map $M_{f^{\kappa}}$ is an isomorphism, and $M_{f^A}$ satisfies (0)--(4), so that the composite $M_f=M_{f^A}M_{f^{\kappa}}$ is a merging bordism.
\end{proof}

\subsection{State spaces and merging bordisms}

For each marked surface $\zeta:D\times\mcl{I}\to \Sigma$ in $\Bord_{\ast}^{nc}$ we recall that the projection
\begin{equation}\label{eq:1214}
(\Bord_{A_{fin}}^{nc})_{\Sigma_{\zeta}}\to (A_{fin}^{\ot})_{\mcl{I}}=A_{fin}^{\mcl{I}}
\end{equation}
is an isomorphism of simplicial sets, since $\Bord_{A_{fin}}^{nc}$ is constructed from the fibration $A_{fin}^{\ot}\to fr\Disk$ via pullback.  Now for a bordism $f:\opn{Cyl}\times \mcl{I}\to M$ between between $\ast$-marked surfaces $M_f:\Sigma_{\zeta}\to \Sigma'_{\eta}$, an edge $M_{\alpha}:\Sigma_x\to \Sigma_y$ over $M_f$ is cocartesian in $\Bord_{A_{fin}}^{nc}$ if and only if its projection  to $A_{fin}^{\ot}$ is cocartesian over $fr\E_2$.  From this we conclude that the isomorphisms from \eqref{eq:1214} fit into a diagram
\[
\xymatrix{
\Delta\times (\Bord_{A_{fin}}^{nc})_{\Sigma_{\zeta}}\ar[rr]^{T_{M_f}}\ar[d]_{\cong} & & \Bord_{A_{fin}}^{nc}\ar[d]^{project}\\
\Delta\times A_{fin}^{\mcl{I}}\ar[rr]_{T_{\partial_{out}(M_f)}} & & A_{fin}^{\ot},
}
\]
where each $T_{\star}$ is the cocartesian transformation over the given map. Evaluating at $1:\ast\to \Delta^1$ now produces a diagram for the transport functors
\[
\xymatrix{
(\Bord_{A_{fin}}^{nc})_{\Sigma_{\zeta}}\ar[rr]^{(M_f)_!}\ar[d]_{\cong} & & (\Bord_{A_{fin}}^{nc})_{\Sigma'_{\eta}}\ar[d]^{\cong}\\
A_{fin}^{\mcl{I}}\ar[rr]_{\partial_{out}(M_f)_!} & & A_{fin}^{\mcl{J}}
}
\]
where we recall that the map $\partial_{out}(M_f)_!$ is the colored product functor $m_{\partial_{out}(M_f)}$.  The above diagram can furthermore be chosen to strictly commute.

After we identify the fibers in $\Bord_{A_{fin}}^{nc}$ with their respective colored exponents for $A_{fin}$, we denote the transport functor $(M_f)_!$ simply by $m_f:A_{fin}^{\mcl{I}}\to A_{fin}^{\mcl{J}}$, as in Notation \ref{not:labeled_unlabeled}.

\begin{lemma}\label{lem:1288}
Let $\Sigma_{\zeta}$ be a marked surface in $\Bord_{\ast}^{nc}$.  For any merging bordism $M_{f}:\Sigma_{\zeta}\to \Sigma_{\eta}$ the corresponding transport functor provides a natural isomorphism
\[
\dL_{A_{fin}}\circ T_{M_f}:\Delta^1\times (\Bord_{A_{fin}}^{nc})_{\Sigma_{\zeta}}\to \opn{Ch}^+(\opn{Vect}_{fin}).
\]
Equivalently, evaluating at each object $\Sigma_x$ over $\Sigma_{\zeta}$ produces a natural isomorphism between the state spaces
\[
\dL_{A_{fin}}(\Sigma_x)\overset{\sim}\to\dL_{A_{fin}}(\Sigma_{m_f(x)}).
\]
\end{lemma}

\begin{proof}
We need only demonstrate a single cocartesian edge $M_{\alpha}:\Sigma_x\to \Sigma_{m_f(x)}$ over $M_f$ which is sent to an isomorphism in $\opn{Ch}^+(\opn{Vect}_{fin})$ under $\dL_{A_{fin}}$.  In the case where $\Sigma$ is connected this is covered in (the proof of) \cite[Lemma 8.4]{czenkynegron}.  In the case where $\Sigma$ is not connected the result follows by the connected case and the monoidal structure on $\dL_{A_{fin}}$.
\end{proof}

\begin{proposition}[{\cite{derenzietal23,czenkynegron}}]\label{prop:pre_states}
Consider a connected genus $g$ surface $\Sigma$ with arbitrary markings $\zeta:D\times \mcl{I}\to\Sigma$.  Any merging bordism $M_{f}:\Sigma_{\zeta}\to \Sigma_{\eta}$ specifies a natural isomorphism
\[
\dL_{A_{fin}}(\Sigma_-)\overset{\sim}\to \Hom^{\ast}_{A_{fin}}\left(C^{\ot g},m_f-\right),
\]
of functors from $(\Bord^{nc}_{A_{fin}})_{\Sigma_{\zeta}}=A_{fin}^{\mcl{I}}$ to $\opn{Ch}^+(\opn{Vect}_{fin})$.
\end{proposition}

In the above expression $C$ is the \emph{coend} for $A^{\heartsuit}$, i.e.\ the value $C=mR(\1)$ where $R:A^{\heartsuit}\to A^{\heartsuit}\ot_{\opn{Vect}} A^{\heartsuit}$ is right adjoint to the product functor $m:A^{\heartsuit}\ot_{\opn{Vect}}A^{\heartsuit}\to A^{\heartsuit}$.  In the case of representations of a Hopf algebra $A^{\heartsuit}=\opn{Rep}(H)$, for example, $C$ is the coadjoint representation $C=H^{\ast}$.

\begin{proof}[Proof of Proposition \ref{prop:pre_states}]
By \cite[Theorem 8.5]{czenkynegron}--which follows \cite[Proposition 4.17]{derenzietal23}--there is such a natural isomorphism
\[
\dL_{A_{fin}}|_{(\Bord_{A_{fin}}^{nc})_{\Sigma_{\eta}}}\overset{\sim}\to \Hom^{\ast}_{A_{fin}}(C^{\ot g},-)
\]
when $\eta:D\to \Sigma$ is a single positive marking.  So the result follows by Lemma \ref{lem:1288}.
\end{proof}

\subsection{Discrete LRT theory for unbounded cochains}

Below we consider the symmetric monoidal category $\opn{Ch}(\opn{Vect})_{\pf}=\opn{Ch}(\opn{Vect}_{\pf})$ of pro-finite linear cochains, as in Section \ref{sect:pf_complexes}.

\begin{proposition}\label{prop:dis_lrt}
Let $A^{\heartsuit}$ be a finite modular tensor category, $A=\opn{Ch}(A^{\heartsuit})$, and $q:A^{\ot}\to fr\Disk$ be the associated $fr\Disk$-category of unbounded cochains (see Proposition \ref{prop:fr_cochains}).  There is a unique symmetric monoidal functor
\[
\dL_{A}:\Bord_{A}^{nc}\to \opn{Ch}(\opn{Vect})_{\pf}
\]
which fits into a commuting diagram
\begin{equation}\label{eq:129}
\xymatrixrowsep{4mm}
\xymatrix{
\Bord^{nc}_{A_{fin}}\ar[rr]^(.45){\dL_{A_{fin}}}\ar[d] & & \opn{Ch}^+(\opn{Vect}_{fin})\ar[d]\\
\Bord^{nc}_{A}\ar[rr]_(.45){\dL_{A}} & & \opn{Ch}(\opn{Vect})_{\pf},
}
\end{equation}
and for which the natural map $\dL_{A}(\Sigma_x)\to \varprojlim_{\lambda}\dL_{A_{fin}}(\Sigma_{x_{\lambda}})$ is an isomorphism whenever $x=\varprojlim_{\lambda}x_{\lambda}$ is an expression of $x$ as a filtered limit of finite cochains as in Lemma \ref{lem:pf_limits}.
\end{proposition}

\begin{proof}
Each object $x$ in a colored product $A^{\mcl{I}}$ admits a structural expression as a filtered limit $x=\varprojlim_{\lambda}x_{\lambda}$ over finite length (bounded) complexes.  We define the state spaces in our proposed theory directly as the limits
\[
\dL_{A}(\Sigma_x)=\varprojlim_{\lambda}\dL_{A_{fin}}(\Sigma_{x_{\lambda}}).
\]
Via the mutual refinement property, one sees that the above limit is always identified with the limit over any system as in Lemma \ref{lem:pf_limits}. Consider any bordism $M_{\alpha}:\Sigma_x\to \Sigma'_y$ with $y$ consisting of finite, bounded complexes.  Such $M_{\alpha}$ consists of a bordism $M_f:\Sigma_{\opn{top}(x)}\to \Sigma_{\opn{top}(y)}$ and a map $\alpha_f:m_f(x)\to y$ in $A^{\mcl{J}}$.
\par

Since the product on $\opn{Ch}(A^{\heartsuit})$ commutes with filtered colimits and preserves injections, the product the dual category $\opn{Ch}(A^{\heartsuit})_{\pf}$ commutes with filtered limits and preserves surjections.  Hence $m_f(x)=\varprojlim_{\lambda} m_f(x_{\lambda})$ with each map $m_f(x)\to m_f(x_{\lambda})$ surjective.  In particular the kernels of these maps provide a base for the topology on $m_f(x)$ around $0$, and continuity of $\alpha_f:m_f(x)\to y$ implies a (unique) factorization through $m_f(x_{\lambda})\to y$ at all sufficiently large $\lambda$,
\[
\xymatrixrowsep{4.5mm}
\xymatrix{
m_f(x)\ar[dr]\ar[rr]^{\alpha_f} & & y\\
	& m_f(x_{\lambda})\ar[ur]_{\alpha_f^{\lambda}} & .
}
\]
Hence we obtain a collection of compatible bordisms $M^{\lambda}_{\alpha}:\Sigma_{x_{\lambda}}\to \Sigma'_y$ for large $\lambda$, and so obtain a well-defined morphism
\[
\dL_A(M_{\alpha})=\dL_{A_{fin}}(M^{\lambda}_{\alpha})\circ \pi_{\lambda}
\]
where $\pi_{\lambda}$ is the structure map $\Sigma_x\to \Sigma_{x_{\lambda}}$.
\par

In the case of a bordism $M_{\alpha}:\Sigma_x\to \Sigma'_y$ with general marking object $y$ we have $\dL_A(\Sigma'_y)=\varprojlim_{\mu}\dL_{A^{\heartsuit}_{fin}}(\Sigma'_{y_{\mu}})$ and so define
\[
\dL_A(M_{\alpha})=\varprojlim_{\mu} \dL_A(M_{\alpha,\mu}):\dL_A(\Sigma_x)\to \dL_A(\Sigma'_y)
\]
where each $M_{\alpha,\mu}$ is the composite of $M_{\alpha}$ with the projection $\Sigma'_y\to \Sigma'_{y_{\mu}}$.  For a composite of bordisms $N_{\beta}M_{\alpha}:\Sigma_x\to \Sigma'_y\to \Sigma''_z$ with $z$ bounded and finite, the diagram(s)
\[
\xymatrix{
\Sigma_x\ar[r]^{M_{\alpha}}\ar[d] & \Sigma'_y\ar[r]^{N_{\beta}}\ar[d] & \Sigma''_z\\
\Sigma_{x_{\lambda}}\ar[r]_{M_{\alpha,\mu}^{\lambda}} & \Sigma'_{y_{\mu}}\ar[ur]_{N_{\beta}^{\mu}}
}
\]
give $\dL_A(N_{\beta}M_{\alpha})=\dL_A(N_{\beta})\dL_A(M_{\beta})$.  For general $z$ the equality $\dL_A(N_{\beta}M_{\alpha})=\dL_A(N_{\beta})\dL_A(M_{\beta})$ follows from the finite case.
\par

For the symmetric monoidal structure we have
\[
\dL_A(\Sigma_x)\ot\dL_A(\Sigma'_y)=\varprojlim_{\lambda,\mu}\dL_{A_{fin}}(\Sigma_{x_{\lambda}})\ot\dL_{A_{fin}}(\Sigma'_{y_{\mu}})
\]
and
\[
\dL_A(\Sigma_x\amalg\Sigma'_y)=\varprojlim_{\lambda,\mu}\dL_{A_{fin}}(\Sigma_{x_{\lambda}}\amalg\Sigma'_{y_{\mu}})
\]
in $\opn{Ch}(\opn{Vect})_{\pf}$. Hence we have a unique natural isomorphism $\dL_A(-)\ot \dL_A(-)\overset{\sim}\to \dL(-\amalg-)$ which fits into diagrams
\[
\xymatrix{
\dL_A(\Sigma_x)\ot\dL_A(\Sigma'_y)\ar[rr]^{\sim}\ar[d] & & \dL_A(\Sigma_x\amalg\Sigma'_y)\ar[d]\\
\dL_{A_{fin}}(\Sigma_{x_{\lambda}})\ot\dL_{A_{fin}}(\Sigma'_{y_{\mu}})\ar[rr]^{\sim} & & \dL_{A_{fin}}(\Sigma_{x_{\lambda}}\amalg\Sigma'_{y_{\mu}})
}
\]
at each $\lambda$ and $\mu$.
\end{proof}

\begin{definition}
For a finite modular tensor category $A^{\heartsuit}$, and $A=\opn{Ch}(A^{\heartsuit})$, the discrete LRT theory is the symmetric monoidal functor
\[
\opn{L}_A:\Bord^{nc}_A\to \opn{Ch}(\opn{Vect})_{\pf}
\]
constructed in Theorem \ref{thm:fin_dis_lrt} and Proposition \ref{prop:dis_lrt}.
\end{definition}

\begin{remark}
The letters LRT reference Lyubashenko, Reshetikhin, and Turaev. It would also be reasonable to reference De Renzi, Gainutdinov, Geer, Patureau-Mirans, and Runkel here, though the label DGGPR is a bit too rambunctious for our liking.
\end{remark}

\subsection{State spaces for discrete LRT}

\begin{lemma}\label{lem:ext_states}
Let $A^{\heartsuit}$ be a finite modular tensor category, and $A=\opn{Ch}(A^{\heartsuit})$.  For any connected genus $g$ surface $\Sigma$ with fixed markings $\zeta:D\times\mcl{I}\to \Sigma$, and any map $f:D\times \mcl{I}\to D$ to a positively colored disk in $fr\Disk$, there is a natural isomorphism
\[
\dL_{A}(\Sigma_{-})\overset{\sim}\to \Hom^\ast_{A}(C^{\ot g},m_f-)
\]
of functors from $(\Bord_{A}^{nc})_{\Sigma_{\zeta}}=A^{\mcl{I}}$ to $\opn{Ch}(\opn{Vect})_{\pf}$.
\end{lemma}

\begin{proof}
We note that all embeddings $f:D\times \mcl{I}\to D$ are smoothly homotopic (Proposition \ref{prop:1456}) so that all of the associated product functors $m_f:A^{\mcl{I}}\to A^+=A_{\pf}$, i.e.\ all transport functors along the maps $f$, are naturally isomorphic.  It follows that the functor $\Hom^{\ast}_A(C^{\ot g},m_f-)$ is independent of the choice of $f$, up to natural isomorphism.
\par

Now, for any marking object $x$ in $A^{\mcl{I}}$ with structural expression $x=\varprojlim_{\lambda} x_{\lambda}$, continuity of $m_f$ gives the product $m_f(x)$ as the limit $m_f(x)\overset{\sim} =\varprojlim_{\lambda} m_f(x_{\lambda})$.  Hence we have a natural isomorphism
\[
\Hom_{A}^\ast(C^{\ot g},m_f(x))\overset{\sim}\to \varprojlim_{\lambda}\Hom_{A_{fin}}^\ast(C^{\ot g},m_f(x_{\lambda}))
\]
between functors from $A^{\mcl{I}}$ to $\opn{Ch}(\opn{Vect})_{\pf}$.  By construction we have a natural isomorphism
\[
\dL_{A}(\Sigma_x)\overset{\sim}\to \varprojlim_{\lambda}\dL_{A_{fin}}(\Sigma_{x_{\lambda}})
\]
as well. Hence the natural isomorphisms $\dL_{A_{fin}}(\Sigma_{x_{\lambda}})\overset{\sim}\to \Hom_{A}^\ast(C^{\ot g},m_f(x_{\lambda}))$ from Proposition \ref{prop:pre_states} imply the claimed identification
\[
\dL_{A}(\Sigma_x)\overset{\sim}\to \Hom_{A}^\ast(C^{\ot g},m_f(x)).
\]
\end{proof}

We recall that any product functor $m_f:A^{\mcl{I}}\to A^{+}=A_{\pf}$ sends homotopy equivalences in each factor to homotopy equivalences in $A_{\pf}$, and hence sends any product of homotopy equivalences $\xi:x\to x'$ to a homotopy equivalence in $A_{\pf}$. This follows by the duality of Corollary \ref{cor:ch_prepair} and the fact that the product on $\opn{Ch}(A^{\heartsuit})$ preserves homotopy equivalence in each factor. Furthermore, the functor
\[
\Hom_{A}^\ast(C^{\ot g},-):A_{\pf}\to \opn{Ch}(\opn{Vect})_{\pf}
\]
preserves homotopy equivalences as well.  It follows from the identification of Lemma \ref{lem:ext_states}, and symmetric monoidality, that the LRT theory preserves homotopy equivalences.

\begin{corollary}\label{cor:dis_htop}
Consider a marked surface $\zeta:D\times\mcl{I}\to \Sigma$. The symmetric monoidal functor $\dL_{A}:\Bord^{nc}_{A}\to \opn{Ch}(\opn{Vect})_{\pf}$ sends homotopy equivalences in the fiber $(\Bord^{nc}_{A})_{\Sigma_{\zeta}}=A^{\mcl{I}}$ to homotopy equivalences in $\opn{Ch}(\opn{Vect})_{\pf}$.
\end{corollary}

\section{Deriving LRT at the $1$-categorical level}
\label{sect:1_categorical}

Consider the discrete derived category $D=D(A^{\heartsuit})$ of a modular tensor category $A^{\heartsuit}$. We have the induced ribbon structure on $D$ and corresponding duality equivalence $-^{\ast}:D(A^{\heartsuit})^{op}\overset{\sim}\to D(A^{\heartsuit}_{\pf})$. This equivalence gives $D(A^{\heartsuit})$ a $fr\Disk$-structure $D^{\ot}\to fr\Disk^{\ot}$ via Proposition \ref{prop:a_pm}, and we produce the subsequent marked bordism category $\Bord^{nc}_D$.

Low dimensional simplices in $\Bord^{nc}_D$ are described precisely as in Section \ref{sect:des_bord_a}, from which one deduces an explicit description of the homotopy truncation $\opn{h}\Bord^{nc}_D$. Vaguely, $\opn{h}\Bord^{nc}_D$ is (equivalent to) a symmetric monoidal category of surfaces with markings from the discrete categories $D(A^{\heartsuit})$ and $D(A^{\heartsuit}_{\pf})$, and ribbon bordisms between such surfaces. The underlying graphs for such ribbon bordisms are required to be trees in which each component has a unique internal vertex, and a unique external vertex which is attached to the outgoing surface.

\begin{theorem}\label{thm:dis_der}
Let $A^{\heartsuit}$ be a finite modular tensor category with corresponding discrete derived category $D=D(A^{\heartsuit})$. There is a symmetric monoidal functor
\[
L_D:\opn{h}\Bord^{nc}_D\to D(\opn{Vect}_{\pf})
\]
which is obtained by ``deriving" the theory from Proposition \ref{prop:dis_lrt}.
\end{theorem}

We note that $L_D$ is a map between $1$-categories, and can therefore be approached in completely $1$-categorical terms. We outline a proof of Theorem \ref{thm:dis_der} below in two parts. Following the proof, we discuss various means of ``expanding" the theory from Theorem \ref{thm:dis_der}, and requisite liftings to the $\infty$-categorical setting.

\subsection{Step 1: Introducing the homotopy category}

We consider the homotopy category $K=K(A^{\heartsuit})$ with its $fr\Disk$-structure induced by the duality equivalence $-^{\ast}:K(A^{\heartsuit})^{op}\overset{\sim}\to K(A^{\heartsuit}_{\pf})$. We show that discrete LRT induces a unique symmetric monoidal functor from $\opn{h}\Bord^{nc}_K$ which fits into a diagram
\begin{equation}\label{eq:3059}
\xymatrix{
\opn{h}\Bord^{nc}_A\ar[rr]^{\opn{L}_A}\ar[d] & & \opn{Ch}(\opn{Vect}_{\pf})\ar[d]^{quotient}\\
\opn{h}\Bord^{nc}_K\ar[rr]^{\opn{L}_K} & & D(\opn{Vect}_{\pf})
}
\end{equation}

\begin{proof}[Proof for Theorem \ref{thm:dis_der}, Step 1]
Consider the equivalence relation $\sim_{htop}$ on maps in $\opn{h}\Bord^{nc}_A$, where two morphisms $M_{\alpha},M_{\beta}:\Sigma_x\to \Sigma_y$ are equivalent whenever they have the same underlying topology
\[
M_{\opn{top}(\alpha)}=M_f=M_{\opn{top}(\beta)},
\]
and the outgoing morphisms $\alpha,\beta:m_f(x)\to y$ are related by a tuple of cochain homopies in $A^{\pm}$. We have the symmetric monoidal functor $\opn{h}\Bord^{nc}_A\to \opn{h}\Bord^{nc}_K$ provided by the balanced projection $A\to K$, and this induces a symmetric equivalence
\[
\opn{h}\Bord^{nc}_A/\sim_{htop}\overset{\sim}\to \opn{h}\Bord^{nc}_K.
\]

We note that each colored product functor $m:A^{\mcl{I}}\to A^+$ along a linear disk arrangement $D\times \mcl{I}\to D$ sends cochain homotopies in the incoming factors to cochain homotopies in $A^+$, and that the Hom functors $\Hom_A^{\ast}(C^{\ot g},-)$ preserve cochain homotopies as well. Thus, from the description of the state spaces provided in Lemma \ref{lem:ext_states}, and the fact that the functor
\[
\opn{Ch}(\opn{Vect}_{\pf})\to D(\opn{Vect}_{\pf})
\]
is a quotient along cochain homotopy, we understand that there is a unique symmetric monoidal functor $\opn{L}_K$ which completes a diagram as in \eqref{eq:3059}.
\end{proof}

\subsection{Step 2: Deriving LRT}

\begin{proof}[Proof of Theorem \ref{thm:dis_der}, Step 2]
The localization functor $K(A^{\heartsuit}_{\pf})\to D(A^{\heartsuit}_{\pf})$ admits a fully faithful right adjoint $R:D(A^{\heartsuit}_{\pf})\to K(A^{\heartsuit}_{\pf})$ which is, by nonsense, lax monoidal and compatible with the balanced structure. The existence of such an adjoint follows from the fact that $K(A^{\heartsuit}_{\pf})$ admits enough $K$-injectives (see Section \ref{sect:kan_htop_der}). Hence $R$ induces a symmetric monoidal functor $\opn{h}\Bord_R^{nc}:\opn{h}\Bord_D^{nc}\to \opn{h}\Bord_K^{nc}$, and we obtain a symmetric monoidal functor to $D(\opn{Vect}_{\pf})$ via composition
\[
\opn{L}_D=\opn{L}_K\circ \opn{h}\Bord_R^{nc}:\opn{h}\Bord^{nc}_D\to D(\opn{Vect}_{\pf}).
\]
This completes the proof of Theorem \ref{thm:dis_der}.
\end{proof}

One can show that the functor $\opn{L}_D:\opn{h}\Bord_D^{nc}\to D(\opn{Vect}_{\pf})$ is a left Kan extension for the functor $\opn{L}_K$, so that $\opn{L}_D$ is a ``derivation" of discrete LRT. One can also check that the state spaces for $\opn{L}_D$ over a genus $g$ marked surface $\Sigma_{\zeta}$ are identified with derived Homs
\[
\opn{L}_D(\Sigma_-)\overset{\sim}\to \opn{RHom}_{A^{\heartsuit}}(C^{\ot g},m_f-),
\]
in direct analogy with Lemma \ref{lem:ext_states} (cf.\ Proposition \ref{prop:der_states}). From this description of the state spaces one sees that $\opn{L}_D$ restricts to a $D(\opn{Vect})$ valued theory from the full subcategory of surfaces labeled by objects in the bounded derived category $D^b(A^{\heartsuit}_{fin})$ of finite length complexes, i.e.\ that we can forgo pro-finite completion in the finite bounded setting,
\[
\opn{L}_{D_{fin}}:\opn{h}\Bord^{nc}_{D_{fin}}\to D(\opn{Vect}).
\]

\subsection{Unmarking the discrete derived theory}
\label{sect:pre_unmark}

We recall that all nonempty surfaces in $\opn{h}\Bord_D^{nc}$ are equipped with markings from the discrete derived category. In a skein theoretic setting, we could remove all such markings by employing unit labels on surfaces and noting that the particular positions of unit markings do not matter, modulo skein relations.
\par

In our setting we can still label all surfaces by the identity object $\1_D$ in $D^b(A^{\heartsuit}_{fin})$ to obtain a symmetric inclusion
\[
\opn{h}\Bord^{nc}_{\ast}\to \opn{h}\Bord_{D_{fin}}^{nc}
\]
from the trivially marked bordism category. We then institute certain coarse skein relations on $\opn{h}\Bord^{nc}_{\ast}$ via a localization process which we outline in Section \ref{sect:derived_w_skeins}. While we defer the details, and direct the interested reader to Sections \ref{sect:derived_w_skeins} and \ref{sect:unmarked_lrt} below, let us record our ultimate conclusion.

\begin{theoremA}[\ref{thm:unmarked_lrt}]
For any finite modular tensor category $A^{\heartsuit}$, the discrete derived theory from Theorem \ref{thm:dis_der} induces an unmarked, anomalous topological field theory in $3$ dimensions
\[
\mbb{L}_{D_{fin}}:\Bord^{nc}_{3,2}\to D(\opn{Vect}).
\]
For each genus $g$ surface $\Sigma$, the corresponding state space admits a natural identification with derived Homs $\mbb{L}_{D_{fin}}(\Sigma)\cong \opn{RHom}_{A^{\heartsuit}}(C^{\ot g},\1)$.
\end{theoremA}

As we explain in Section \ref{sect:mcg}, the existence of such a field theory implies a natural, projective action of each mapping class group $\opn{MCG}(\Sigma)$ on the associated derived Homs $\opn{RHom}_{A^{\heartsuit}}(C^{\ot g},\1)$ (cf. \cite{lentneretal23,schweigertwoike21}). We also discuss implications in dimension $3$ in Section \ref{sect:dimension_3}.

\subsection{Expanding from Theorem \ref{thm:dis_der}?}
\label{sect:expansion}

At the center of our project, even beyond the present paper, is the following question:\vspace{1.5mm}

{\it -- To what extent do the symmetric monoidal functors $\opn{L}_D$ and $\mbb{L}_{D_{fin}}$ \emph{expand} into a system of (T)QFTs, in analogy with the fusion or 1-categorical settings? --}\vspace{1.5mm}

To list a few explicit points:\vspace{1.5mm}

\noindent(A) {\it For a given finite modular tensor category $A^{\heartsuit}$, is there a once extended, anomalous $3$-d topological theory $\mfk{L}'$ whose circle value $\mfk{L}'(S^1)$ recovers the derived category for $A^{\heartsuit}$? Does such an extended theory imply a natural marking scheme which supersedes the one employed in our construction of $\opn{L}_D$? (See Conjectures \ref{conj:0} and \ref{conj:3} below, as well as \cite{bartlettetal,derenzi21,lagiotis} for the $1$-categorical setting.)}\vspace{1.5mm}

\noindent(B) {\it Is there an extended \emph{$4$-dimensional} topological theory associated to the derived category of a modular tensor category, and can we recover our $3$-dimensional theories from such a $4$-d theory via a boundary construction? (See Question \ref{quest:4} below and compare with \cite{haioun}.)}\vspace{1.5mm}

\noindent(C) {\it Is there an effective notion of skeins, and skein relations, for the derived category which allows one to calculate values in the $4$-dimensional theory from (b)? In particular, can one calculate factorization homology using skeins?  Thinking further about (a), can this skein theory be used to implement skein relations into a sufficiently receptive $3$-dimensional marking scheme? (See again Sections \ref{sect:skein_bordisms} and \ref{sect:mark_scheme}, Question \ref{quest:4.5} below, and compare with \cite{costantinoetalII,brownhaioun}.)}\vspace{1.5mm}

Questions (A)--(C) all ask, in the non-semisimple/derived setting, if there is an ecosystem of interconnected topological field theories which parallel those one finds in the fusion setting, for say Crane-Yetter vs.\ Chern-Simons/Reshetikhin-Turaev theories. We record, finally, a principle concerning flat connections which is unique to the non-semisimple world.\vspace{1.5mm}

\noindent(D) {\it Does the theory $\opn{L}_D$, or some enhanced version of it, couple naturally to flat connections? (See Conjectures \ref{conj:5} and \ref{conj:6} below, and also \cite{kinnear}.)}\vspace{1.5mm}

We defer further discussions of point (D) to Section \ref{sect:conjectures}, though it is certainly a deep and interesting story, both in mathematics and in physics.

Thinking seriously about any of the questions above requires one to think about the derived category in some locale where one has an effective \emph{theory} of categories. It is apparent that the triangulated setting will not suffice for this task, which then prompts our transition to the world of $\infty$-categories.

\subsection{Derivating LRT in $\infty$-categories}

The remainder of the text--or more precisely, the next five sections--is dedicated to the construction of a lift for the marked theory from Theorem \ref{thm:dis_der} to the level of $\infty$-categories. We suggest this lift as a projection of a more expansive system of TQFTs, as proposed above, and so as an indication that such a system should in fact exist. We return to the fundamental questions (A--D), and provide a more thorough discussion of these issues, in Section \ref{sect:conjectures} at the conclusion of the text.
\par

In lifting the theory $\opn{L}_D$ to the homotopical setting, we follows essentially the same processes outlined above. We begin in Sections \ref{sect:localization} and \ref{sect:homtopy_lrt} with the movement from the abelian category of cochains to the homotopy $\infty$-category, then Kan extend to produce a theory for the derived $\infty$-category. In comparing with Step 1 above, we note that the homotopy $\infty$-category $\msc{K}(A^{\heartsuit})$ is no longer obtained from chains over $A^{\heartsuit}$ via a quotient construction, but via localization.

\begin{remark}
The methods employed below are not sensitive to our specific choice of marking scheme (see Section \ref{sect:mark_scheme}). What's required is the existence of a preceding abelian theory at the level of cochains, for the given marking scheme, as in Proposition \ref{prop:dis_lrt}.
\end{remark}

\section{Localization for monoidal $\infty$-categories}
\label{sect:localization}

We claim that the discrete LRT theory $\opn{L}_A:\Bord^{nc}_A\to \opn{Ch}(\opn{Vect})_{\pf}$ from Proposition \ref{prop:dis_lrt} \emph{localizes} to provide a symmetric monoidal functor, i.e.\ field theory, from a bordism $\infty$-category $\Bord^{nc}_{\msc{K}}$ with labels in the homotopy $\infty$-category. While this is more-or-less an immediate consequence of Corollary \ref{cor:dis_htop}, one needs to understand the behaviors of localization when applied to various types of fibrations in order to effectively proceed in this direction.
\par

In this section we present an overview of localization for cocartesian fibrations, and monoidal $\infty$-categories, which essentially follows work of Hinich on the topic \cite{hinich16}. (See also \cite[Propositions 2.2.1.9 \& 4.1.7.4]{ha} and \cite[Appendix A]{nikolausscholze18}.)

\subsection{Vertical localization}
\label{sect:vert_loc}

Let $W$ be a collection of morphisms in an $\infty$-category $\msc{C}$. Recall that a localization for $\msc{C}$ along $W$ is the choice of a functor $F:\msc{C}\to \msc{C}'$ for which restriction
\[
F^{\ast}:\Fun(\msc{C}',\msc{D})\to \Fun(\msc{C},\msc{D}),
\]
at arbitrary $\msc{D}$, provides an equivalence onto the full subcategory $\Fun(\msc{C},\msc{D})_W\subseteq \Fun(\msc{C},\msc{D})$ spanned by those functors which send all morphisms in $W$ to isomorphisms in $\msc{D}$.  Via a pushout construction \cite[\href{https://kerodon.net/tag/05ZL}{05ZL}]{kerodon}, or fibrant replacement \cite[Proposition 3.1.3.7]{lurie09}, we understand that localizations always exist.

\begin{proposition}[{\cite[Proposition 2.1.4]{hinich16}}]\label{prop:fib_localize}
Let $q:\msc{C}\to \msc{T}$ be a cocartesian fibration and $W\subseteq \msc{C}[1]$ be a collection of edges.  Suppose that every edge in $W$ lies in the fiber $\msc{C}_t$ over some object $t$ in $\msc{T}$, that each fiber $W_t=W\cap \msc{C}_t[1]$ is stable under isomorphism in $\Fun(\Delta^1,\msc{C}_t)$, and that the transport functor $f_!:\msc{C}_t\to \msc{C}_{u}$ over each map $f:t\to u$ in $\msc{T}$ sends $W_t$ into $W_u$.  Then there is a pairing of a localization $F:\msc{C}\to \msc{C}[W^{-1}]$ with a map $q':\msc{C}[W^{-1}]\to \msc{T}$ which fts into a strictly commuting diagram
\begin{equation}\label{eq:3037}
\xymatrix{
\msc{C}\ar[rr]^F\ar[dr]_q & & \msc{C}[W^{-1}]\ar[dl]^{q'}\\
	& \msc{T} & ,
}
\end{equation}
and for which the following hold:
\begin{itemize}
\item[(i)] $q':\msc{C}[W^{-1}]\to \msc{T}$ is an isofibration.\vspace{1mm}
\item[(ii)] For any isofibration $\msc{D}\to\msc{T}$, restriction provides an equivalence
\[
F^{\ast}:\Fun_{\msc{T}}(\msc{C}[W^{-1}],\msc{D})\overset{\sim}\to \Fun_{\msc{T}}(\msc{C},\msc{D})_W.
\]
\item[(iii)] $q':\msc{C}[W^{-1}]\to \msc{T}$ is a cocartesian fibration, and $F$ is a map of cocartesian fibrations over $\msc{T}$.\vspace{1mm}
\item[(iv)] For each $t$ in $\msc{T}$ the induced map on the fibers $F_t:\msc{C}_t\to (\msc{C}[W^{-1}])_t$ realizes the target category as a localization for $\msc{C}_t$ along the class $W_t$.\vspace{1mm}
\item[(v)] For any cocartesian fibration $\msc{D}\to \msc{T}$, restriction provides an equivalence
\[
F^{\ast}:\Fun_{\msc{T}}^{cc}(\msc{C}[W^{-1}],\msc{D})\overset{\sim}\to \Fun^{cc}_{\msc{T}}(\msc{C},\msc{D})_W.
\] 
\end{itemize}
Furthermore, any such pairing $(F,q')$ in which $q'$ is an isofibration automatically satisfies (ii)--(v).
\end{proposition}

As in Section \ref{sect:fun_cc}, $\Fun_{\msc{T}}(\msc{C},\msc{D})_W$ denotes the fiber
\[
\Fun_{\msc{T}}(\msc{C},\msc{D})_W=\{q\}\times_{\Fun(\msc{C},\msc{T})}\Fun(\msc{C},\msc{D})_W
\]
and $\Fun_{\msc{T}}^{cc}(\msc{C},\msc{D})_W$ is the full subcategory spanned by all functors which preserve cocartesian edges, i.e.\ maps of cocartesian fibrations.

\begin{proof}
From either \cite[Proposition 3.1.3.7]{lurie09} or \cite[\href{https://kerodon.net/tag/05ZL}{05ZL}]{kerodon}, one sees that there is a construction of the localization $F':\msc{C}\to \msc{C}''$ for which $F'$ is injective, as a map of simplicial sets.  In this case the restriction equivalence
\[
\Fun(\msc{C}'',\msc{D})\to \Fun(\msc{C},\msc{D})_W
\]
is furthermore an isofibration \cite[\href{https://kerodon.net/tag/01F3}{01F3}]{kerodon}.  Hence, considering the case $\msc{D}=\msc{T}$, we see that the structure map $q$ lifts uniquely to a map $q'':\msc{C}''\to \msc{T}$ which completes the proposed diagram.

We can factor $q''$ into a sequence $\msc{C}''\overset{w}\to \msc{C}'\overset{q'}\to \msc{T}$ in which $w$ is an equivalence and $q'$ is an isofibration \cite[\href{https://kerodon.net/tag/02VM}{02VM}]{kerodon}.  Then the composition $F=wF':\msc{C}\to \msc{C}'$ is also a localization of $\msc{C}$ along $W$, and we claim that this choice of localization, $\msc{C}[W^{-1}]=\msc{C}'$ coupled with $q':\msc{C}'\to \msc{T}$, satisfies (i)--(iv). Indeed, we claim that \emph{any} choice of a localization $F:\msc{C}\to \msc{C}'$ with an isofibration $q':\msc{C}'\to \msc{T}$ which completes a diagram as in \eqref{eq:3037} satisfies (ii)--(v).

(i) Holds by construction.  (ii) Given an isofibration $\msc{D}\to \msc{T}$, we have a diagram
\[
\xymatrixrowsep{4mm}
\xymatrixcolsep{4mm}
\xymatrix{
\Fun_{\msc{T}}(\msc{C}',\msc{D})\ar[rr]\ar[dr]\ar[dd] & & \Fun(\msc{C}',\msc{D})\ar[dd]\ar[dr]^{\sim}\\
 & \Fun_{\msc{T}}(\msc{C},\msc{D})_W\ar[rr]\ar[dd] & & \Fun(\msc{C},\msc{D})_W\ar[dd] \\
\ast\ar[dr]_{=}\ar[rr]_(.6){q'} & & \Fun(\msc{C}',\msc{T})\ar[dr]^{\sim}\\
 & \ast\ar[rr]_q & & \Fun(\msc{C},\msc{T})_W.
}
\]
in which the vertical maps are isofibrations \cite[\href{https://kerodon.net/tag/01H7}{01H7}]{kerodon} and the front and back faces are pullback diagrams.  It follows that the induced map on the fibers $\Fun_{\msc{T}}(\msc{C}',\msc{D})\to \Fun_{\msc{T}}(\msc{C},\msc{D})_W$ is an equivalence as well, by Lemma \ref{lem:pullback_equiv}.

(iii) In \cite[Section 2.2.2]{hinich16} Hinich constructs a localization $\msc{E}$ for $\msc{C}$ by localizing the functor $F_{\msc{C}}^{\opn{marked}}:\msc{T}\to \sCat_{\infty}^{\opn{marked}}$ associated to $\msc{C}$.  The associated localized functor $F_{\msc{E}}:\msc{T}\to \sCat_{\infty}$ comes equipped with a transformation $F_{\msc{C}}\to F_{\msc{E}}$ so that we have an induced map of cocartesian fibrations
\[
\xymatrix{
\msc{C}\ar[rr]^l\ar[dr]_q & & \msc{E}\ar[dl]\\
	& \msc{T} & .
}
\]
By (ii) this triangle extends to a diagram
\begin{equation}\label{eq:1859}
\xymatrix{
	& \msc{C}'\ar[dr]\ar[dd]|{\hole}^(.3){q'}\\
\msc{C}\ar[rr]|(.4)l\ar[dr]_{q}\ar[ur]^{F} & & \msc{E}\ar[dl]\\
	& \msc{T} & .
}
\end{equation}
in which the sub-diagram $\Lambda^3_3\to \sCat_{\infty}$ strictly commutes and the map of $\infty$-categories $\msc{C}'\to \msc{E}$ is an equivalence, by uniqueness of localization.  Since $q'$ is an isofibration it follows, via the above diagram and the fact that $l$ is a map of cocartesian fibrations, that $q'$ is in fact a cocartesian fibration \cite[\href{https://kerodon.net/tag/028A}{028A}]{kerodon} and that $F$ is a map of cocartesian fibrations \cite[\href{https://kerodon.net/tag/028B}{028B}]{kerodon}.
\par

(iv) By \cite[Proposition 2.1.4]{hinich16} the localization $\msc{E}$ considered above has the property that all of its fibers $\msc{E}_t$ over $\msc{T}$ localize the respective fibers $\msc{C}_t$.  (This property is actually clear from its construction via the localization of $F_{\msc{C}}$.)  Via the equivalence $\msc{C}'\to \msc{E}$ over $\msc{T}$, and the above diagram \eqref{eq:1859}, it follows that the fibers $\msc{C}'_t$ similarly localize the fibers $\msc{C}_t$ \cite[\href{https://kerodon.net/tag/023M}{023M}]{kerodon}.
\par

(v) Suppose that $\msc{D}\to \msc{T}$ is a cocartesian fibration.  Given the equivalence of (ii), it suffices to show that a functor $G:\msc{C}\to \msc{D}$ over $\msc{T}$ preserves cocartesian edges if and only if its localized map $G':\msc{C}'\to \msc{D}$ preserves cocartesian edges.  However, since $F$ is essentially surjective, and since cocartesian lifts of edges in the base are unique, this just follows from the fact that $F$ itself preserves cocartesian edges.
\end{proof}

\begin{definition}\label{def:vertical}
Given a cocartesian fibration $q:\msc{C}\to \msc{T}$, a collection of edges $W$ in $\msc{C}$ is said to be a vertical class over $\msc{T}$ if it satisfies the hypotheses of Proposition \ref{prop:fib_localize}.  Given such such a vertical class $W$, a vertical localization over $\msc{T}$ is a choice of strictly commuting diagram
\[
\xymatrix{
\msc{C}\ar[rr]^(.45){F}\ar[dr]_{q} & & \msc{C}[W^{-1}]\ar[dl]^{q'}\\
	& \msc{T}
}
\]
in which $q'$ is a cocartesian fibration and $F$ is a localization functor which preserves cocartesian edges.
\end{definition}

Proposition \ref{prop:fib_localize} says that such vertical localizations for cocartesian fibrations always exist. We also observe that vertical localizations are stable under pullback.

\begin{proposition}\label{prop:vert_pullback}
Let $q:\msc{C}\to \msc{T}$ be a cocartesian fibration, $W\subseteq\msc{C}[1]$ be a vertical class over $\msc{T}$, and $\msc{C}[W^{-1}]\to \msc{T}$ be the associated vertical localization over $\msc{T}$.  Consider any functor $\msc{R}\to \msc{T}$, and let $W_{\msc{R}}$ be the collection of all maps in the fiber product $\msc{C}_{\msc{R}}=\msc{R}\times_{\msc{T}}\msc{C}$ which are sent to degenerate edges in $\msc{R}$ and to $W$ in $\msc{C}$.  Then $W_{\msc{R}}$ is a vertical class in $\msc{C}_{\msc{R}}$ and the natural map
\begin{equation}\label{eq:1893}
\msc{C}_{\msc{R}}[W_{\msc{R}}^{-1}]\to \msc{R}\times_{\msc{T}}(\msc{C}[W^{-1}])
\end{equation}
is an equivalence of cocartesian fibrations over $\msc{R}$.
\end{proposition}

\begin{proof}
We note that the fibers of $\msc{C}_{\msc{R}}$ over $\msc{R}$ are identified with the fibers of $\msc{C}$ over $\msc{T}$, and under these identifications the transport functors for $\msc{C}_{\msc{R}}$ are identified with the transport functors for $\msc{C}$.  So it is clear that $W_{\msc{R}}$ is a vertical class.  We consider the claim that the map \eqref{eq:1893} is a equivalence.
\par
 
By Proposition \ref{prop:fib_localize} (v), the map of cocartesian fibrations $\msc{C}_{\msc{R}}\to \msc{R}\times_{\msc{T}}\msc{C}[W^{-1}]$ induced by the map of cocartesian fibrations $\msc{C}\to \msc{C}[W^{-1}]$ localizes to provide a map of cocartesian fibrations
\[
\opn{nat}:\msc{C}_{\msc{R}}[W_{\msc{R}}^{-1}]\to \msc{R}\times_{\msc{T}}(\msc{C}[W^{-1}])
\]
over $\msc{R}$.  This is our functor \eqref{eq:1893}.  Note that for each object $r:\ast\to \msc{R}$, with image $t$ in $\msc{T}$, the induced map on fibers $(\msc{C}_{\msc{R}})_r\to \msc{C}_t$ is an isomorphism which sends the class $(W_{\msc{R}})_r$ bijectively to $W_t$.  Hence, from the diagram
\[
\xymatrix{
(\msc{C}_{\msc{R}})_r\ar[rr]^(.4){\cong}\ar[d] & & \msc{C}_t\ar[d]\\
(\msc{C}_{\msc{R}}[W_{\msc{R}}^{-1}])_r\ar[rr]_(.35){\opn{nat}_r} & & \left(\msc{R}\times_{\msc{T}}(\msc{C}[W^{-1}])\right)_r=(\msc{C}[W^{-1}])_t
}
\]
and Proposition \ref{prop:fib_localize} (iv), we find that each fiber $\opn{nat}_r$ over $\msc{R}$ is an equivalence.  It follows that $\opn{nat}$ itself is an equivalence \cite[\href{https://kerodon.net/tag/023M}{023M}]{kerodon}.
\end{proof}

\subsection{Aside: Generation of classes}

\begin{lemma}[Morphism generation]\label{lem:gen_mor}
Let $\msc{E}^{\ot}\to fr\Disk^{\ot}$ be a $fr\Disk$-category and $W^{\pm}$ be classes of maps in the underlying $\infty$-categories $\msc{E}^{\pm}$ which are stable under isomorphism in $\Fun(\Delta^1,\msc{E}^{\pm})$ and contain all identities. Suppose also that
\begin{itemize}
\item The transition map $d:\msc{E}^-\to \msc{E}^+$ sends $W^-$ into $W^+$.\vspace{1mm}
\item The product functors $m:\msc{E}^{\pm}\times\msc{E}^{\pm}\to \msc{E}^{\pm}$ send $(W^{\pm})^2$ into $W^{\pm}$.
\end{itemize}
There is a (unique) vertical class $W$ in $\msc{E}^{\ot}$ whose fibers $W_{(\mcl{I},\mu)}$ are precisely the pullbacks of the product classes $W^{\mcl{I}}$ along the transport equivalences $\rho_!:\msc{E}^{\ot}_{(\mcl{I},\mu)}\to \msc{E}^{\mcl{I}}$.
\end{lemma}

To recall, the transition functor $d:\msc{E}^-\to \msc{E}^+$ is the transport functor along the color negating identity $id:D_-\to D_+$, and the product functors $m:\msc{E}^{\pm}\times \msc{E}^{\pm}$ are the transport functors along any like colored disk embedding $D\amalg D\to D$.

\begin{proof}
For the union $W=\cup_{(\mcl{I},\mu)}W_{(\mcl{I},\mu)}$, stability under transport follows by the description of transport along maps in $fr\Disk^{\ot}$ provide in Lemma \ref{lem:structure_transp}.
\end{proof}

\begin{definition}
For a $fr\Disk$-category $q:\msc{E}^{\ot}\to fr\Disk^{\ot}$ and classes $W^{\pm}\subseteq \msc{E}^{\pm}[1]$ as in Lemma \ref{lem:gen_mor}, we call the corresponding collection $W\subseteq \msc{E}^{\ot}[1]$ from Lemma \ref{lem:gen_mor} the class generated by $W^{\pm}$.
\end{definition}

\subsection{Localization for monoidal $\infty$-categories}

In general, one can use Proposition \ref{prop:fib_localize} to effectively localize $\msc{O}$-monoidal $\infty$-categories along classes of maps which are stable under any collection of ``generating morphisms" in $\msc{O}^{\ot}$, for a given $\infty$-operad or equifibered symmetric monoidal $\infty$-category $\msc{O}^{\ot}$.  We cover two specific instances which are relevant for us.

\begin{proposition}\label{prop:fr_loc}
Let $q:\msc{E}^{\ot}\to fr\Disk^{\ot}$ be a $fr\Disk$-category, and $W^{\pm}$ be collections of edges in the underlying $\infty$-categories $\msc{E}^{\pm}$ which are stable under isomorphism in $\Fun(\Delta^1,\msc{E}^{\pm})$ and contain all identities.  Suppose that the $W^{\pm}$ are stable under the product functors on $\msc{E}^{\pm}$ as well as the transition functor $d:\msc{E}^{-}\to \msc{E}^+$. Then the localization $q':\msc{E}^{\ot}[W^{-1}]\to fr\Disk^{\ot}$ along the corresponding vertical class (see Lemma \ref{lem:gen_mor}) is a $fr\Disk$-monoidal $\infty$-category.
\end{proposition}

\begin{proof}
Take $\msc{G}^{\ot}=\msc{E}^{\ot}[W^{-1}]$. By Proposition \ref{prop:fib_localize}, $q':\msc{G}^{\ot}\to fr\Disk$ is a cocartesian fibration and the fibers of the localization functor $\msc{E}^{\ot}\to \msc{G}^{\ot}$ induce equivalences $\msc{E}^{\ot}_t[W^{-1}_t]\overset{\sim}\to \msc{G}_t^{\ot}$, over each object $t$ in $fr\Disk^{\ot}$.
\par

Let $f:t=(\mcl{I},\mu)\to t'$ be a cocartesian edge in $fr\Disk^{\ot}$ over the active surjection $K\to \{0\}$ in $\opn{Fin}_{\ast}$. We have the diagram
\[
\xymatrix{
\msc{E}^{\ot}_t\ar[rr]^{f_!}\ar[dr]_{\rho_!} & & \msc{E}^{\ot}_{t'}\ar[dl]^{\rho_!}\\
	& \msc{E}^{\mcl{I}}
}
\]
in which $f_!$ is an equivalence, by the definition of $fr\Disk$-monoidality. So we see, from the definition of $W$, that $f_!^{-1}(W_{t'})=W_{t}$. Similarly, for the inert projections $\rho_k:t\to t_k$ the transport equivalence $\msc{E}^{\ot}_t\overset{\sim}\to \prod_{k\in K}\msc{E}^{\ot}_{t_k}$ recovers $W_t$ as the preimage of the product $\prod_{k\in K} W_{t_k}$. Hence, from the diagrams
\[
\xymatrix{
\msc{E}_t^{\ot}[W_t^{-1}]\ar[d]_{\sim}\ar[r]^{\sim} & \msc{E}_{t'}^{\ot}[W_{t'}^{-1}]\ar[d]^{\sim} &  
\msc{E}_t^{\ot}[W_t^{-1}]\ar[r]^(.45){\sim}\ar[d]_{\sim} & \prod_k\msc{E}^{\ot}_{t_k}[W_{t_k}^{-1}]\ar[d]^{\sim}\\
\msc{G}^{\ot}_t\ar[r]_{transp} & \msc{G}^{\ot}_{t'} & 
\msc{G}^{\ot}_t\ar[r]_(.45){transp} & \prod_k\msc{G}^{\ot}_{t_k}
}
\]
we see that transport for $\msc{G}^{\ot}$ induces equivalences $\msc{G}^{\ot}_t\overset{\sim}\to \msc{G}^{\ot}_{t'}$ and  $\msc{G}^{\ot}_t\overset{\sim}\to\prod_k\msc{G}^{\ot}_{t_k}$. This verifies $fr\Disk$-monoidality of the localization.
\end{proof}

The analogous statement for symmetric monoidal $\infty$-categories is well-known.

\begin{proposition}[{\cite[Proposition 3.2.2]{hinich16}}]\label{prop:sym_loc}
Let $q:\msc{F}^{\ot}\to \Fin_{\ast}$ be a symmetric monoidal $\infty$-category, and $W_0\subseteq \msc{F}[1]$ be a class of maps which is stable under isomorphism in $\Fun(\Delta^1,\msc{F})$ and contains all identities.  Suppose also that $W_0$ is stable under the product functor $m:\msc{F}\times \msc{F}\to \msc{F}$, and define $W\subseteq \msc{F}^{\ot}[1]$ to be the collection of edges in the varied fibers $(\msc{F}^{\ot})_I$ which are sent to $W^I$ under the transport equivalences $(\msc{F}^{\ot})_I\overset{\sim}\to \msc{F}^I$.
\par

The collection $W$ is a vertical class in $\msc{F}^{\ot}$ and the corresponding cocartesian fibration $q':\msc{F}^{\ot}[W^{-1}]\to \Fin_{\ast}$ realizes $\msc{F}^{\ot}[W^{-1}]$ as a symmetric monoidal $\infty$-category.
\end{proposition}

One can show that in fact all ``symmetric monoidal localizations" are of the form considered in Proposition \ref{prop:sym_loc}.  We record a variation on this point.

\begin{lemma}\label{lem:adv_loc}
Let $q:\msc{F}^{\ot}\to \Fin_{\ast}$ be a symmetric monoidal $\infty$-category and $W\subseteq \msc{F}^{\ot}[1]$ be a class of maps for which $q(\alpha)$ is a degenerate whenever $\alpha\in W$.  Suppose that there is a localization $l:\msc{F}^{\ot}\to \msc{F}^{\ot}[W^{-1}]$ which fits into a strictly commuting diagram
\[
\xymatrix{
\msc{F}^{\ot}\ar[dr]_q\ar[rr]^l & & \msc{F}^{\ot}[W^{-1}]\ar[dl]^{q'}\\
	& \Fin_{\ast} & ,
}
\]
and that the following hold:
\begin{enumerate}
\item[(a)] $q'$ is a cocartesian fibration which gives $\msc{F}^{\ot}[W^{-1}]$ the structure of a symmetric monoidal $\infty$-category.\vspace{1mm}
\item[(b)] $l$ is a map of symmetric monoidal $\infty$-categories.\vspace{1mm}
\item[(c)] The natural map $\msc{F}[W_{\{0\}}^{-1}]\to \msc{F}^{\ot}[W^{-1}]_{\{0\}}$ is an equivalence.
\end{enumerate}
Then a map of symmetric monoidal $\infty$-categories $F:\msc{F}^{\ot}\to \msc{D}^{\ot}$ factors to provide a (uniquely associated) map of symmetric monoidal $\infty$-categories $F':\msc{F}^{\ot}[W^{-1}]\to \msc{D}^{\ot}$ if and only if $F_{\{0\}}:\msc{F}\to \msc{D}$ sends $W_{\{0\}}$ into the class of isomorphisms in $\msc{D}$.
\end{lemma}

\begin{proof}
Throughout we let $\opn{Isom}(\msc{C})\subseteq \msc{C}[1]$ denote the class of isomorphisms in an $\infty$-category $\msc{C}$. Take $\msc{G}^{\ot}=\msc{F}^{\ot}[W^{-1}]$ and consider the class $U$ in $\msc{F}^{\ot}[1]$ obtained as the union of all of the maps
\[
U_I=l^{-1}_I\left(\opn{Isom}(\msc{G}^{\ot}_I)\right)\ \subseteq\ \msc{F}^{\ot}_I[1],
\]
across all finite sets $I$.  Via the diagrams
\[
\xymatrix{
\msc{G}^{\ot}_I\ar[rr]^{\sim}_{\bar{\rho}_!} & & \msc{G}^I\\
\msc{F}^{\ot}_I\ar[u]^{l_I}\ar[rr]^{\sim}_{\bar{\rho}_!} & & \msc{F}^I\ar[u]_{l^I_{\{0\}}}
}
\]
in $\sCat_{\infty}$ we understand that each fiber $U_I$ is the preimage of $U_{\{0\}}^I$ under the transport equivalence $\msc{F}^{\ot}_I\overset{\sim}\to \msc{F}^I$.  Furthermore from the general compatibility with transport
\[
\xymatrix{
\msc{G}^{\ot}_I\ar[rr]^{\bar{f}_!} & & \msc{G}^{\ot}_J\\
\msc{F}^{\ot}_I\ar[rr]_{\bar{f}_!}\ar[u]^{l_I} & & \msc{F}^{\ot}_J\ar[u]_{l_J}
}
\]
we understand that $\bar{f}_!(U_I)\subseteq U_J$ for all $\bar{f}:I\to J$ in $\Fin_{\ast}$.  Hence $U$ is the vertical class generated by $U_{\{0\}}$ over $\Fin_{\ast}$, as in Proposition \ref{prop:sym_loc}, and via Proposition \ref{prop:fib_localize} we observe an equivalence of symmetric monoidal $\infty$-categories
\[
\xymatrix{
\msc{F}^{\ot}[U^{-1}]\ar[rr]^{\sim}\ar[dr] & & \msc{G}^{\ot}\ar[dl]^{q'}\\
	& \Fin_{\ast} & .
}
\]
See also \cite[\href{https://kerodon.net/tag/0285}{0285}]{kerodon}.
\par

Consider now any symmetric monoidal functor $F:\msc{F}^{\ot}\to \msc{D}^{\ot}$ and suppose that $F(U_{\{0\}})\subseteq \opn{Isom}(\msc{D})$.  Then by compatibility of $F$ with transport we have $F(U_I)\subseteq \opn{Isom}(\msc{D}^{\ot}_I)$ for all $I$, and hence $F(U)\subseteq \opn{Isom}(\msc{D}^{\ot})$.  So we obtain an induced map of symmetric monoidal $\infty$-categories $F:\msc{G}^{\ot}\to \msc{D}^{\ot}$ by Proposition \ref{prop:fib_localize} (v).
\par

We claim finally that $F$ sends $U_{\{0\}}$ into $\opn{Isom}(\msc{D})$ if and only if $F$ sends $W_{\{0\}}$ into $\opn{Isom}(\msc{D})$.  Clearly the former property implies the latter, and if we suppose that $F$ sends $W_{\{0\}}$ into isomorphisms then $F_{\{0\}}:\msc{F}\to \msc{D}$ admits a factorization
\[
\xymatrix{
	& \msc{G}\ar[dr]^{\exists !}\\
\msc{F}\ar[rr]_{F_{\{0\}}}\ar[ur]^{l_{\{0\}}} & & \msc{D}
}
\]
by condition (c).  But this diagram implies $F$ sends $U_{\{0\}}$ into $\opn{Isom}(\msc{D})$ as well, by the definition of $U$.
\end{proof}

\subsection{Derived $\infty$-categories as $fr\Disk$-categories}
\label{sect:fr_htop_der}

Let $A^{\heartsuit}$ be a (discrete) ribbon tensor category and $q:A^{\ot}\to fr\Disk^{\ot}$ be the corresponding $fr\Disk$-monoidal $\infty$-category of unbounded cochains.  Recall the values
\[
A^-=A^{op}=\opn{Ch}(A^{\heartsuit})^{op}\ \ \text{and}\ \ A^+=A_{\pf}=\opn{Ch}(A^{\heartsuit}_{\pf}).
\]
We have the classes $\opn{Htop}^{\pm}\subseteq A^{\pm}[1]$ of homotopy equivalences.  These classes are stable under duality $-^{\ast}:A^-\to A^+$ and under the tensor products $A^{\pm}\times A^{\pm}\to A^{\pm}$.  Stability under the tensor product, in particular, follows from the fact that the tensor product preserves homotopy equivalences between morphisms in each factor.  We then have the class $\opn{Htop}$ in $A^{\ot}$ generated by $\{\opn{Htop}^-,\opn{Htop}^+\}$, and localization produces a $fr\Disk$-category
\[
q_{\msc{K}}:\msc{K}^{\ot}\to fr\Disk^{\ot}
\]
whose underlying $\infty$-categories are the homotopy $\infty$-categories
\[
\msc{K}^-=\msc{K}^{op}=\msc{K}(A^{\heartsuit})^{op}\ \ \text{and}\ \ \msc{K}^+=\msc{K}_{\pf}=\msc{K}(A^{\heartsuit}_{\pf}),
\]
by the calculations $\msc{K}^{\pm}=A^{\pm}[(\opn{Htop}^{\pm})^{-1}]$ \cite[Proposition 1.3.4.5]{ha}.  We note that the localization $q_{\msc{K}}:\msc{K}^{\ot}\to fr\Disk^{\ot}$ is the unique $fr\Disk$-category which admits a monoidal functor
\[
\xymatrix{
A^{\ot}\ar[dr]\ar[rr] & & \msc{K}^{\ot}\ar[dl]\\
	& fr\Disk^{\ot}
	}
\]
with the fibers $A^{\pm}\to \msc{K}^{\pm}$ over the $\pm$-colored disks inverting homotopy equivalences and inducing equivalences $A^{\pm}[(\opn{Htop}^{\pm})^{-1}]\overset{\sim}\to \msc{K}^{\pm}$.  Such uniqueness is deducible from Proposition \ref{prop:fib_localize} and \cite[\href{https://kerodon.net/tag/028B}{028B}]{kerodon}.
\par

We can also consider the collections $\opn{Qiso}^{\pm}\subseteq A^{\pm}[1]$ of quasi-isomorphisms.  These classes are again stable under duality and tensoring, by Proposition \ref{prop:tensor_qiso}, and hence generate a vertical class $\opn{Qiso}$ in $A^{\ot}$. We therefore localize to obtain a $fr\Disk$-category
\[
q_{\msc{D}}:\msc{D}^{\ot}\to fr\Disk^{\ot}
\]
with underlying $\infty$-categories
\[
\msc{D}^-=\msc{D}^{op}=\msc{D}(A^{\heartsuit})^{op}\ \ \text{and}\ \ \msc{D}^+=\msc{D}_{\pf}=\msc{D}(A^{\heartsuit}_{\pf})
\]
\cite[Proposition 1.3.5.15]{ha}.  (See also Section \ref{sect:kan_htop_der}.) As in the homotopy situation, the fibration $q_{\msc{D}}:\msc{D}^{\ot}\to fr\Disk^{\ot}$ is the unique $fr\Disk$-structure on the derived $\infty$-category for which the localization functors $A^{\pm}\to \msc{D}^{\pm}$ lift to a $fr\Disk$-monoidal functor.

\section{Homotopy LRT theories}
\label{sect:homtopy_lrt}

We construct a homotopical variant $L_{\msc{K}}:\Bord^{nc}_{\msc{K}}\to \Vect_{\pf}$ of the discrete Lyubashenko-Reshetikhin-Turaev theory from Section \ref{sect:dis_lrt}, and calculate the state spaces in this theory.  The state spaces are, in particular, identified with inner mapping spaces for the action of $\Vect_{\pf}$ on $\msc{K}(A^{\heartsuit}_{\pf})$.
\par

In order to tame the topology in the target category, we show that $L_{\msc{K}}$ also restricts to a $\Vect$-valued TQFT $L_{\msc{K}_{l\text{-}fin}}:\Bord^{nc}_{\msc{K}_{l\text{-}fin}}\to \Vect$ from the $\infty$-category of bordisms with labels in the homotopy $\infty$-category of locally finite cochains.

\subsection{Homotopical vector spaces}

\begin{definition}
We take
\[
\Vect:=\opn{Ch}(\opn{Vect})[\opn{Qiso}^{-1}]\ \ \text{and}\ \ \Vect_{\pf}:=\opn{Ch}(\opn{Vect})_{\pf}[\opn{Qiso}^{-1}],
\]
and refer to these localizations as the $\infty$-categories of homotopical vector spaces and homotopical pro-finite vector spaces, respectively.
\end{definition}

Since quasi-isomorphisms are preserved under the tensor products on $\opn{Ch}(\opn{Vect})$ and $\opn{Ch}(\opn{Vect})_{\pf}$, by Proposition \ref{prop:tensor_qiso}, we see that these localizations inherit unique symmetric monoidal structures under which the localization functors $\opn{Ch}(\opn{Vect})_{\star}\to \Vect_{\star}$ are symmetric monoidal functors. This follows by Proposition \ref{prop:sym_loc}.

\begin{remark}
As with any abelian category, one can explicitly construct the $\infty$-categories $\Vect$ and $\Vect_{\pf}$ by taking the dg nerves of the dg categories of unbounded cochains $\opn{Ch}_{\opn{dg}}(\opn{Vect})$ and $\opn{Ch}_{\opn{dg}}(\opn{Vect}_{\pf})$, respectively \cite[Proposition 1.3.4.5]{ha}.  (We note that in this case the classes of homotopy equivalences and quasi-isomorphisms agree, so that there are no distinctions between the homotopy and derived $\infty$-categories.)
\par

However, one can also construct $\Vect$ as the stabilization of the $\infty$-category $\msc{K}\!an_k$ of linear Kan complexes--aka simplicial vector spaces.
\end{remark}

We observe that the pair of fully faithful, symmetric monoidal functors from $\opn{Ch}^{\pm}(\opn{Vect}_{fin})$ to $\opn{Ch}(\opn{Vect})$ and $\opn{Ch}(\opn{Vect})_{\pf}$ (see Section \ref{sect:pf_complexes}) localizes to produce fully faithful, symmetric monoidal embeddings
\[
\xymatrixrowsep{3mm}
\xymatrix{
	& \Vect_{l\text{-}fin}^{\pm}\ar[dr]\ar[dl] \\
\Vect & & \Vect_{\pf}.
}
\]
As expected, we've taken $\Vect_{l\text{-}fin}^{\pm}=\opn{Ch}^{\pm}(\opn{Vect}_{fin})[\opn{Qiso}^{-1}]$ in the above expression.

\subsection{Discrete LRT with homotopical target}

Let $A^{\heartsuit}$ be a finite modular tensor category with corresponding $fr\Disk$-category $A^{\ot}\to fr\Disk^{\ot}$ of unbounded cochains.  We have the symmetric monoidal functor $\opn{L}_A^{\ot}:(\Bord_A^{nc})^{\ot}\to \opn{Ch}(\opn{Vect})_{\pf}^{\ot}$ from Proposition \ref{prop:dis_lrt} and compose with the localization functor $\opn{Ch}(\opn{Vect})_{\pf}^{\ot}\to \Vect_{\pf}^{\ot}$ to obtain a map of symmetric monoidal $\infty$-categories
\[
L_A^{\ot}:(\Bord_A^{nc})^{\ot}\to \Vect_{\pf}^{\ot}.
\]
This is just a version of the discrete LRT theory with a homotopical target.

\subsection{Homotopy LRT theories}

Below we speak of a diagram in the $\infty$-category $\SM_{\infty}$ of symmetric monoidal $\infty$-categories.  We recall the construction of this $\infty$-category in Section \ref{sect:cocart_sm} of the appendix, but it suffices to understand that a $2$-simplex in $\SM_{\infty}$ is a choice of composable functors between symmetric monoidal $\infty$-categories
\[
\xymatrix{
\msc{F}_0^{\ot}\ar[rr]^{F_{02}}\ar[dr]|{F_{01}}\ar[ddr]_{q_0} & & \msc{F}^{\ot}_2\ar[ddl]^{q_2}\\
	& \msc{F}^{\ot}_1\ar[ur]|{F_{12}}\ar[d]|{q_1}\\
	& \Fin_{\ast}
}
\]
and a choice of natural isomorphism $\sigma:\Delta^1\times\msc{F}^{\ot}_0\to \msc{F}^{\ot}_2$ in $\Fun_{\Fin_{\ast}}(\msc{F}_0,\msc{F}_2)$ with $\sigma_0=F_{12}F_{01}$, $\sigma_1=F_{02}$, and $q_2\sigma=id_{q_0}$.

\begin{proposition}\label{prop:htop_lrt}
Consider a finite modular tensor category $A^{\heartsuit}$ and let $\msc{K}^{\ot}=\msc{K}(A^{\heartsuit})^{\ot}\to fr\Disk^{\ot}$ be the associated $fr\Disk$-category for the homotopy $\infty$-category of unbounded cochains over $A^{\heartsuit}$.  There is a unique symmetric monoidal functor $L_{\msc{K}}:\Bord_{\msc{K}}^{nc}\to \Vect_{\pf}$ from the $\infty$-category of $\msc{K}$-labeled bordisms which completes a diagram
\[
\xymatrix{
\Bord_A^{nc}\ar[rr]^{L_A}\ar[dr]_{loc} & & \Vect_{\pf}\\
	& \Bord_{\msc{K}}^{nc}\ar[ur]_{L_{\msc{K}}}
}
\]
in $\SM_{\infty}$.
\end{proposition}

\begin{proof}
We obtain the fibration $\msc{K}^{\ot}\to fr\Disk^{\ot}$ via vertical localization of the fibration $A^{\ot}\to fr\Disk^{\ot}$ along homotopy equivalence, so that the map
\[
(\Bord_A^{nc})^{\ot}\to (\Bord_{\msc{K}}^{nc})^{\ot}
\]
realizes $(\Bord_{\msc{K}}^{nc})^{\ot}$ as the symmetric monoidal localization of $(\Bord_{A}^{nc})^{\ot}$ along the pullback class $\opn{Htop}_{(\Bord_{\ast}^{nc})^{\ot}}$ (Proposition \ref{prop:vert_pullback}). Furthermore, taking the fiber at $\{0\}$ in $\Fin$, which is the same as taking the fiber along the underlying $\infty$-category $\Bord_{\ast}^{nc}\to (\Bord_{\ast}^{nc})^{\ot}$, produces an equivalence
\[
\Bord^{nc}_A[\opn{Htop}^{-1}_{\Bord^{nc}_{\ast}}]\overset{\sim}\to \Bord_{\msc{K}}^{nc} 
\]
of cocartesian fibrations over $\Bord_{\ast}^{nc}$, again by Proposition \ref{prop:vert_pullback}.  Thus Lemma \ref{lem:adv_loc} tells us that the functor
\[
(\Bord^{nc}_{A})^{\ot}\to (\Bord^{nc}_{\msc{K}})^{\ot}
\]
is universal amongst symmetric monoidal functors which invert the class $\opn{Htop}_{\Bord_{\ast}}$ in the underlying $\infty$-category.
\par

Now, according to the construction from Proposition \ref{prop:vert_pullback}, we can describe the class $\opn{Htop}_{\Bord_{\ast}}$ explicitly.  These are simply bordisms $M_{\alpha}:\Sigma_x\to \Sigma_y$ whose underlying topology $M_{\opn{top}(\alpha)}$ is the identity map, i.e.\ the $3$-manifold $M=\Sigma\times [0,1]$ with embedded cylinders $\opn{Cyl}\times \mcl{I}\to \Sigma\times [0,1]$ traveling in straight lines from $\Sigma_{\opn{top}(x)}$ to $\Sigma_{\opn{top}(y)}=\Sigma_{\opn{top}(x)}$, and whose labeling morphisms $\alpha_i:x_i\to y_i$ are all homotopy equivalences.  Hence, by Corollary \ref{cor:dis_htop}, the LRT theory
\[
L_A:\Bord_A^{nc}\to \Vect_{\pf}
\]
sends all maps in $\opn{Htop}_{\Bord_{\ast}}$ to isomorphisms in $\Vect_{\pf}$. We therefore obtain the claimed factorization
\begin{equation}\label{eq:2666}
\xymatrix{
(\Bord_A^{nc})^{\ot}\ar[rr]^{L_A^{\ot}}\ar[dr]_{loc^{\ot}} & & \Vect_{\pf}^{\ot}\\
	& (\Bord_{\msc{K}}^{nc})^{\ot}\ar@{-->}[ur]_{\exists\ L_{\msc{K}}^{\ot}} & 
}
\end{equation}
in $\SM_{\infty}$.

For uniqueness, restriction along the localization functor provides an equivalence
\[
\Fun^{cc}_{\Fin_{\ast}}((\Bord_{\msc{K}}^{nc})^{\ot},\Vect_{\pf}^{\ot})\overset{\sim}\to \Fun^{cc}_{\Fin_{\ast}}((\Bord_A^{nc})^{\ot},\Vect_{\pf}^{\ot})_{\opn{Htop}},
\]
by Proposition \ref{prop:fib_localize} (v).  Hence the homotopy fiber
\[
\Fun^{cc}_{\Fin_{\ast}}((\Bord_{\msc{K}}^{nc})^{\ot},\Vect_{\pf}^{\ot})\times^{\opn{htop}}_{\Fun^{cc}_{\Fin_{\ast}}((\Bord_A^{nc})^{\ot},\Vect_{\pf}^{\ot})_{\opn{Htop}}}\{L_A\}
\]
is contractible \cite[\href{https://kerodon.net/tag/033M}{033M}]{kerodon}.  By defintion of the homotopy pullback \cite[\href{https://kerodon.net/tag/032Z}{032Z}]{kerodon}, this contractible space parametrizes those symmetric monoidal functors which complete such a diagram \eqref{eq:2666} in $\SM_{\infty}$.
\end{proof}

\begin{definition}
For any finite modular tensor category $A^{\heartsuit}$, with corresponding homotopy $\infty$-category $\msc{K}=\msc{K}(A^{\heartsuit})$, we refer to the symmetric monoidal functor
\[
L_{\msc{K}}:\Bord^{nc}_{\msc{K}}\to \Vect_{\pf}
\]
from Proposition \ref{prop:htop_lrt} as the homotopy LRT theory.
\end{definition}

\subsection{Calculating the state spaces}

We recall that the state spaces in discrete LRT are identified with the Hom complexes
\[
\opn{L}_A(\Sigma_x)\overset{\sim}\to\Hom_{A}^{\ast}(C^{\ot g},m_f(x)),
\]
whenever $\Sigma$ is connected of genus $g$ and $m_f:A^{\mcl{I}}\to A^+=A_{\pf}$ is the product associated to any embedding $D\times \mcl{I}\to D$ to a positively colored disk.  We can obtain the Hom complexes in a manner internal to $A_{\pf}$ by taking inner-Homs for the action of $\opn{Ch}(\opn{Vect})_{\pf}$ on $A_{\pf}$ (Lemma \ref{lem:a_innhom}).  We recall also that the object $C$ here is the canonical coend in $A^{\heartsuit}$.
\par

For $\msc{K}_{\pf}=\msc{K}(A^{\heartsuit}_{\pf})$, the action of $\opn{Ch}(\opn{Vect}_{\pf})$ on $\opn{Ch}(A^{\heartsuit}_{\pf})$ localizes to endow $\msc{K}_{\pf}$ with a natural module structure over $\Vect_{\pf}$ (Corollary \ref{cor:htop_lin}), and the $\Vect_{\pf}$-module category $\msc{K}_{\pf}$ admits inner-Homs.  In particular, for any bounded complex of finite length objects $x$ in $\msc{K}_{\pf}$ we have the associated internal mapping spaces
\begin{equation}
\underline{\Maps}_{\msc{K}_{\pf}}(x,-):\msc{K}_{\pf}\to \Vect_{\pf}
\end{equation}
(Corollary \ref{cor:K_innerhoms}). The following is an analog to Lemma \ref{lem:ext_states} for homotopy LRT.

\begin{proposition}\label{prop:htop_states}
Consider a finite modular tensor category $A^{\heartsuit}$ and the corresponding $fr\Disk$-category $\msc{K}^{\ot}\to fr\Disk^{\ot}$ for the homotopy $\infty$-category. For a connected genus $g$ surface $\Sigma$ with markings $\zeta:D\times\mcl{I}\to \Sigma$, and any choice of map $f:D\times \mcl{I}\to D$ to a positive disk in $fr\Disk$, there is a natural isomorphism
\[
L_{\msc{K}}(\Sigma_-)\overset{\sim}\to \underline{\Maps}_{\msc{K}_{\pf}}(C^{\ot g},m_f-)
\]
of functors from $\msc{K}^{\mcl{I}}\cong(\Bord_{\msc{K}}^{nc})_{\Sigma_{\zeta}}$ to $\Vect_{\pf}$.
\end{proposition}

\begin{proof}
Given any $fr\Disk$-category $\msc{E}^{\ot}\to fr\Disk^{\ot}$ we have the canonical identification
\[
(\Bord_{\msc{E}}^{nc})_{\Sigma_{\zeta}}=\{\Sigma_{\zeta}\}\times_{\Bord_{\ast}^{nc}}\Bord_{\ast}^{nc}\times_{fr\Disk^{\ot}}\msc{E}^{\ot}=\{\mcl{I}\}\times_{fr\Disk^{\ot}}\msc{E}^{\ot}\overset{\sim}\to \msc{E}^{\mcl{I}}.
\]
We consider the fibration $\Bord_{\msc{K}}^{nc}\to \Bord_{\ast}^{nc}$ and, taking the fiber at $\Sigma_{\zeta}$, the diagram from Proposition \ref{prop:htop_lrt} produces a diagram
\[
\xymatrix{
A^{\mcl{I}}\ar[rr]^{L_A|_{\Sigma_{\zeta}}}\ar[dr]_{loc} & & \Vect_{\pf}\\
	& \msc{K}^{\mcl{I}}\ar[ur]_{L_{\msc{K}}|_{\Sigma_{\zeta}}} & .
}
\]
The fact that the map $A^{\mcl{I}}\to \msc{K}^{\mcl{I}}$ is the localization functor along homotopy equivalences follows by Proposition \ref{prop:vert_pullback}.
\par

Via the isomorphism $L_A\overset{\sim}\to loc\circ \Hom^{\ast}_{A}(C^{\ot g},m_f-)$ we obtain an isomorphic diagram
\[
\xymatrix{
A^{\mcl{I}}\ar[rr]^(.45){loc\Hom^{\ast}_A(C^{\ot},m_f-)}\ar[dr]_{loc} & & \Vect_{\pf}\\
	& \msc{K}^{\mcl{I}}\ar[ur]_{L_{\msc{K}}|_{\Sigma_{\zeta}}} & ,
}
\]
and hence conclude that $L_{\msc{K}}|_{\Sigma_{\zeta}}$ is the unique functor obtained from $\Hom^{\ast}_A(C^{\ot g},m_f-)$ via localization.  Rather, $L_{\msc{K}}|_{\Sigma_{\zeta}}$ is isomorphic to any other functor which completes such a diagram.
\par

By Corollary \ref{cor:K_innerhoms} and the fact that the localization functor $A^{\ot}\to \msc{K}^{\ot}$ is a map of $fr\Disk$-categories, we have a diagram
\[
\xymatrix{
A^{\mcl{I}}\ar[r]^{m_f}\ar[d]_{loc} & A^+\ar[d]_{loc}\ar[rr]^(.4){\Hom^{\ast}(C^{\ot g},-)} & & \opn{Ch}(\opn{Vect})_{\pf}\ar[d]^{loc}\\
\msc{K}^{\mcl{I}}\ar[r]_{m_f} & \msc{K}^+\ar[rr]_(.4){\underline{\Maps}(C^{\ot g},-)} & & \Vect_{\pf}.
}
\]
in $\sCat_{\infty}$. This diagram reduces, via a sequence of horn fillings, to provide a diagram
\[
\xymatrix{
A^{\mcl{I}}\ar[rr]^(.45){loc\Hom^{\ast}_A(C^{\ot},m_f-)}\ar[dr]_{loc} & & \Vect_{\pf}\\
	& \msc{K}^{\mcl{I}}\ar[ur]_{\underline{\Maps}(C^{\ot g},m_f-)} & 
}
\]
in $\sCat_{\infty}$ and a subsequent isomorphism $L_{\msc{K}}|_{\Sigma_{\zeta}}\overset{\sim}\to \underline{\Maps}_{\msc{K}_{\pf}}(C^{\ot g},m_f-)$.
\end{proof}

\begin{remark}
As is clear from the proof, or more directly from Corollary \ref{cor:K_innerhoms} in the appendix, each mapping space at fixed $x$ in $\msc{K}^+$ is identified with the Hom complex
\[
\underline{\Maps}_{\msc{K}_{\pf}}(C^{\ot g},x)\cong loc\circ \Hom^{\ast}_{A}(C^{\ot g}, x).
\]
So, as an object, this complex is not mysterious.  Its behaviors on higher simplices in $\msc{K}$ represent the only non-trivial point.  However, even this higher structure can be understood concretely if one adopts any explicit construction of the homotopy $\infty$-category.
\end{remark}

\subsection{Restricting to locally finite complexes}
\label{sect:kfin_lrt}

We are interested in producing a version of the field theory $L_{\msc{K}}$ which takes values in $\Vect$, rather than in the pro-finite completion $\Vect_{\pf}$. It is relatively clear that the restriction of $L_{\msc{K}}$ to the subcategory of bordisms with locally finite labels has the desired property. We argue the point explicitly below.
\par

Take $\msc{K}=\msc{K}(A^{\heartsuit})$ for a ribbon tensor category $A^{\heartsuit}$ and consider the corresponding $fr\Disk$-category $\msc{K}^{\ot}\to fr\Disk^{\ot}$ from Section \ref{sect:fr_htop_der}. Let $\msc{K}_{l\text{-}fin}^{\pm}\subseteq \msc{K}^{\pm}$ be the full subcategories spanned by complexes which are isomorphic (homotopic) to complexes in the full subcategories $\msc{K}^{\pm}(A^{\heartsuit}_{fin})$ in $\msc{K}^-(A^{\heartsuit})$ and $\msc{K}^+(A^{\heartsuit}_{\pf})$, respectively.  Define now $\msc{K}_{l\text{-}fin}^{\ot}\subseteq \msc{K}^{\ot}$ as the full subcategory whose fibers along various $t=(\mcl{I},\mu)$ in $fr\Disk^{\ot}$ are the pullbacks
\begin{equation}\label{eq:2770}
\xymatrix{
(\msc{K}_{l\text{-}fin}^{\ot})_t\ar[rr]\ar[d] & & \msc{K}_{l\text{-}fin}^{\mcl{I}}\ar[d]\\
\msc{K}^{\ot}_{\mcl{I}}\ar[rr]^{\sim}_{\rho_!} & & \msc{K}^{\mcl{I}},
}
\end{equation}
as in Section \ref{sect:subcats}.
\par

Since locally finite complexes are stable under tensoring and duality, Proposition \ref{prop:frdisk_subcats} tells us that $\msc{K}^{\ot}_{l\text{-}fin}$ is a $fr\Disk$-category, in the expected way.

\begin{lemma}\label{lem:k_lfin}
The map $\msc{K}_{l\text{-}fin}^{\ot}\to fr\Disk^{\ot}$ inherited from $\msc{K}^{\ot}$ via restriction is a cocartesian fibration which gives $\msc{K}_{l\text{-}fin}^{\ot}$ the structure of $fr\Disk^{\ot}$-category, and the inclusion $\msc{K}_{l\text{-}fin}^{\ot}\to \msc{K}^{\ot}$ is a $fr\Disk$-monoidal functor.
\end{lemma}

For any complex $x$ in $\msc{K}^+(A^{\heartsuit}_{fin})$ the linear mapping space $\underline{\Maps}_{\msc{K}_{\pf}}(C^{\ot g},x)$ lies in the symmetric subcategory $\Vect_{l\text{-}fin}^+$ in $\Vect_{\pf}$.  Since the inclusion $\Vect_{l\text{-}fin}^+\to \Vect_{\pf}$ is fully faithful, the corresponding map of symmetric monoidal $\infty$-categories
\[
(\Vect_{l\text{-}fin}^+)^{\ot}\to \Vect_{\pf}^{\ot}
\]
is also fully faithful \cite[\href{https://kerodon.net/tag/01VB}{01VB}]{kerodon}.  So, since all objects in $\msc{K}_{l\text{-}fin}^{\pm}$ are isomorphic to an object in $\msc{K}^{\pm}(A^{\heartsuit}_{fin})$, the entire restricted TQFT
\[
L_{\msc{K}}^{\ot}|_{(\Bord_{\msc{K}_{l\text{-}fin}}^{nc})^{\ot}}
\]
admits a unique symmetric monoidal lift to the $\infty$-category of locally finite, bounded below, homotopical vector spaces
\[
\xymatrix{
	& (\Vect_{l\text{-}fin}^+)^{\ot}\ar[dr]\\
(\Bord_{\msc{K}_{l\text{-}fin}}^{nc})^{\ot}\ar[rr]_{L^{\ot}_{\msc{K}}}\ar@{-->}[ur]^{\exists!} & & \Vect_{\pf}^{\ot}.
}
\]
Via a consideration of the symmetric monoidal embedding $\Vect_{l\text{-}fin}^+\to \Vect$, we observe the following.

\begin{proposition}\label{prop:kfin_lrt}
Let $A^{\heartsuit}$ be a finite modular tensor category and $\msc{K}_{l\text{-}fin}^{\ot}\to fr\Disk^{\ot}$ be the associated $fr\Disk$-category of partially bounded, locally finite cochains.  There is a symmetric monoidal functor
\[
L_{\msc{K}_{l\text{-}fin}}:\Bord_{\msc{K}_{l\text{-}fin}}^{nc}\to \Vect
\]
which is obtained from $L_{\msc{K}}$ via restriction. Furthermore, for any connected genus $g$ surface $\Sigma_{\zeta}$ with marking set $\mcl{I}$, and embedding $f:D\times\mcl{I}\to D$ to a positively colored disk, there is a natural isomorphism
\[
L_{\msc{K}_{l\text{-}fin}}(\Sigma_{-})\overset{\sim}\to \underline{\Maps}_{\msc{K}}(C^{\ot g},m_f-)
\]
of functors from $\msc{K}^{\mcl{I}}_{l\text{-}fin}\cong(\Bord_{\msc{K}_{l\text{-}fin}}^{nc})_{\Sigma_{\zeta}}$ to $\Vect$.
\end{proposition}

\begin{proof}
To ease notation we take $\msc{K}_{l\text{-}fin}=\msc{K}_{l\text{-}fin}^+$ and $\Vect_{l\text{-}fin}=\Vect^+_{l\text{-}fin}$. The only issue is the identification of the state spaces.  First, via the action of $\opn{Ch}^+(\opn{Vect}_{fin})$ on $\opn{Ch}^+(A^{\heartsuit}_{fin})$ we obtain a natural action of $\Vect_{l\text{-}fin}$ on the homotopy $\infty$-category $\msc{K}_{l\text{-}fin}$ via localization (Proposition \ref{prop:mod_loc}).  From the inner-Homs for $\opn{Ch}^+(A^{\heartsuit}_{fin})$, which are just given by the Hom-complexes, we also obtain inner-Homs for $\msc{K}_{l\text{-}fin}$ via localization (Proposition \ref{prop:loc_innerhoms}).  In particular, the inner-Homs for $\msc{K}_{l\text{-}fin}$ agree with those of the ambient category $\msc{K}_{\pf}$ of pro-finite complexes,
\[
\underline{\Maps}_{\msc{K}_{l\text{-}fin}}= \underline{\Maps}_{\msc{K}_{\pf}}|_{\msc{K}_{fin}^{op}\times \msc{K}_{l\text{-}fin}}.
\]
Similarly, we have the inclusion $\msc{K}_{l\text{-}fin}\to \msc{K}=\msc{K}(A^{\heartsuit})$, inner-Homs $\underline{\Maps}_{\msc{K}}$ for the action of $\Vect$ on $\msc{K}$, and an identification of inner-Homs
\[
\underline{\Maps}_{\msc{K}_{l\text{-}fin}}= \underline{\Maps}_{\msc{K}}|_{\msc{K}_{fin}^{op}\times \msc{K}_{l\text{-}fin}},
\]
so that we may write $\underline{\Maps}_{\msc{K}_{l\text{-}fin}}(C^{\ot g},m_f-)=\underline{\Maps}_{\msc{K}}(C^{\ot},m_f-)|_{(\msc{K}_{l\text{-}fin}^{\ot})_{\mcl{I}}}$.
\end{proof}

\begin{remark}
One can (easily) show that the theory $L_{\msc{K}_{l\text{-}fin}}$ is the unique one obtained from $L_{A_{l\text{-}fin}}=L_A|_{\Bord_{A_{l\text{-}fin}}^{nc}}$ via localization.
\end{remark}

\section{Adjoints and Kan extension over a base}
\label{sect:rel_adj_kan}

From our theory $L_{\msc{K}}:\Bord_{\msc{K}}^{nc}\to \Vect_{\pf}$ for the homotopy $\infty$-category we can produce a theory $L_{\msc{D}}:\Bord_{\msc{D}}^{nc}\to \Vect_{\pf}$ for the derived $\infty$-category as well.  We proceed, in particular, via a Kan extension along the localization functor $\Bord_{\msc{K}}^{nc}\to \Bord_{\msc{D}}^{nc}$.  In order to control this extension to a sufficient degree--for example, in order to calculate the state spaces in $L_{\msc{D}}$--we want to understand what such a Kan extension \emph{actually} looks like.
\par

In this section we give a general presentation of relative Kan extension, and relative adjoints for fibrations over a fixed base.  We employ relative constructions in order to account for behaviors of monoidal structures under Kan extension, for example.

\begin{remark}
There is a more complicated notion of operadic Kan extension \cite[Definition 3.1.2.2]{ha}. For our purposes, we expect that one can employ operadic Kan extension instead of relative Kan extension, and obtain the same output. So the distinction is not a deep one for us, and the reader might consult \cite[Section 3.1.2]{ha} for an alternate point of view.
\end{remark}

\subsection{Relative adjoints}

\begin{definition}
We say a functor between isofibrations
\[
\xymatrix{
\msc{C}\ar[rr]^F\ar[dr]_q & & \msc{E}\ar[dl]^p\\
	& \msc{T}
}
\]
admits a right adjoint relative to $\msc{T}$ if there is another functor
\[
\xymatrix{
\msc{E}\ar[rr]^G\ar[dr]_p & & \msc{C}\ar[dl]^q\\
	& \msc{T}
}
\]
and a transformation $\varepsilon:GF\to id_{\msc{C}}$ in $\Fun_{\msc{T}}(\msc{C},\msc{C})$ which exhibits $G$ as right adjoint to $F$, in the non-relative sense of the term \cite[\href{https://kerodon.net/tag/02EJ}{02EJ}]{kerodon}.
\end{definition}

We note that this notion is symmetric with respect to units and counits.

\begin{lemma}
If $F:\msc{C}\to \msc{E}$ is a functor between isofibrations over a given base $\msc{T}$, and $G:\msc{E}\to \msc{C}$ is a right adjoint relative to $\msc{T}$ with corresponding counit $\varepsilon:GF\to id_{\msc{C}}$, then there is a unit transformation $u:id_{\msc{E}}\to FG$ which lies in $\Fun_{\msc{T}}(\msc{E},\msc{E})$.
\end{lemma}

\begin{proof}
Let $q:\msc{C}\to \msc{T}$ and $p:\msc{E}\to \msc{T}$ be the structure maps and let $u':FG\to id_{\msc{E}}$ be any choice of unit.  We have the diagram
\[
\xymatrix{
	& GFG\ar[dr]^{\varepsilon G}\\
G\ar[ur]^{G(u')}\ar[rr]_{id_G} & & G
}
\]
in $\Fun(\msc{E},\msc{C})$ which provides a diagram
\[
\xymatrix{
	& p\ar[dr]^{id}\\
p\ar[ur]^{p(u')}\ar[rr]_{id} & & p
}
\]
after applying $q$.  So we see $p(u')$ is isomorphic to the identity on $p$ in the functor category $\Fun(\Delta^1\times \msc{E},\msc{T})$.  Since the map
\begin{equation}\label{eq:3421}
p_\ast:\Fun(\Delta^1\times \msc{E},\msc{E})\to \Fun(\partial\Delta^1\times \msc{E},\msc{E})\times_{\Fun(\partial\Delta^1\times \msc{E},\msc{T})}\Fun(\Delta^1\times \msc{E},\msc{T})
\end{equation}
is an isofibration \cite[\href{https://kerodon.net/tag/01F1}{01F1}]{kerodon} we can find a diagram
\begin{equation}\label{eq:2877}
\xymatrix{
FG\ar[r]^{id} & FG\\
id_{\msc{E}}\ar[u]^{u'}\ar[ur]\ar[r]_{id} & id_{\msc{E}}\ar[u]_{u}
}
\end{equation}
in $\Fun(\msc{E},\msc{E})$ which maps to the paired diagrams
\[
\xymatrix{
FG\ar[r]^{id} & FG &\ar@{}[d]|{\text{\normalsize and}} & &  p\ar[r]^{id} & p\\
id_{\msc{E}}\ar[r]_{id} & id_{\msc{E}} & & & p\ar[u]^{p(u')}\ar[r]_{id}\ar[ur]|{id} & p\ar[u]_{id}
}
\]
under \eqref{eq:3421}.  The transformation $u:id_{\msc{E}}\to FG$ in \eqref{eq:2877} provides the requisite unit transformation.
\end{proof}

Since the unit and counit transformations have constant projections to the base, we can pull back along any map $\msc{R}\to \msc{T}$ to obtain an adjunction over $\msc{R}$ from any adjunction over $\msc{T}$.

\begin{proposition}[{\cite[Proposition 7.3.2.5]{ha}}]\label{prop:reladj_pullback}
Suppose a map between isofibrations
\[
\xymatrix{
\msc{C}\ar[rr]^F\ar[dr] & & \msc{E}\ar[dl]\\
	& \msc{T}
}
\]
admits a right adjoint $G:\msc{E}\to \msc{C}$ relative to $\msc{T}$.  Then for any functor $\msc{R}\to \msc{T}$ the pullback
\[
G_{\msc{R}}:\msc{R}\times_{\msc{T}}\msc{E}\to \msc{R}\times_{\msc{T}}\msc{C}
\]
is a right adjoint to $F_{\msc{R}}$ relative to $\msc{R}$.
\end{proposition}

As for the existence of such adjunctions, we have a convenient reduction to the fibers in the cocartesian case.

\begin{proposition}[{\cite[Proposition 7.3.2.6]{ha}}]\label{prop:reladj_exists}
A map between cocartesian fibrations
\[
\xymatrix{
\msc{C}\ar[rr]^F\ar[dr] & & \msc{E}\ar[dl]\\
	& \msc{T}
}
\]
admits a right adjoint relative to $\msc{T}$ if and only if, for each $t$ in $\msc{T}$, the fiber $F_t:\msc{C}_t\to \msc{E}_t$ admits a right adjoint.
\end{proposition}

\subsection{Adjoints and (relative) Kan extension}

Let $l:\msc{C}\to \msc{C}'$ be a functor between $\infty$-categories.  We recall that a left Kan extension for a functor $F:\msc{C}\to \msc{E}$, along $l$, is the choice of a functor $F':\msc{C}'\to \msc{E}$ and transformation $\alpha:F\to F'l$ which realizes, at each $x$ in $\msc{C}'$, the value $F'(x)$ as a colimit of the diagram
\[
\xymatrixrowsep{3mm}
\xymatrix{
	& \msc{C}\ar@{..>}[dr]^{F}\\
\msc{C}_{/x}\ar[rr]\ar@{..>}[ur]^{forget} & & \msc{E}.
}
\]
(See \cite[\href{https://kerodon.net/tag/02Y7}{02Y7}]{kerodon}.)  Here $\msc{C}_{/x}$ is the relative overcategory
\[
\msc{C}_{/x}:=\msc{C}\times_{\msc{C}'}\msc{C}'_{/x},
\]
i.e.\ the $\infty$-category which parametrizes objects $\tilde{x}$ in $\msc{C}$ equipped with a specified map $l(\tilde{x})\to x$.

The following is standard in $1$-category theory, and should be well-known in the $\infty$-categorical setting as well.

\begin{lemma}\label{lem:2916}
Let $l:\msc{C}\to \msc{C}'$ be any functor which admitis a right adjoint $R:\msc{C}'\to \msc{C}$. Then any functor $F:\msc{C}\to \msc{E}$ admits a left Kan extension along $l$, and this Kan extension can be calculated as the composite $F'=FR:\msc{C}'\to \msc{E}$ along with the shifted unit transformation $F(u):F\to F'l$.
\end{lemma}

\begin{proof}
At each $x$ in $\msc{C}'$ we have the transformation
\begin{equation}\label{eq:2922}
F|_{\msc{C}_{/x}}\overset{F(u)}\to F'l|_{\msc{C}_{/x}}\overset{F'(\opn{can})}\to \underline{F'(x)}
\end{equation}
in $\Fun(\msc{C}_{/x},\msc{E})$, where $\underline{F'(x)}$ is the constant diagram of value $F'(x)$.  The transformation $F'l|_{\msc{C}_{/x}}\to \underline{F'(x)}$ is obtained as follows: We have the coslice equivalence
\[
\msc{C}_{/x}\overset{\sim}\to \msc{C}'\times_{\Fun(\{0\},\msc{C}')}\Fun(\Delta^1,\msc{C}')\times_{\Fun(\{1\},\msc{C}')}\{x\}
\]
\cite[\href{https://kerodon.net/tag/02GH}{02GH}]{kerodon} which specifies, and is specified by the choice of a canonical transformation
\[
\opn{can}:\ast\to \Fun(\Delta^1\times\msc{C}_{/x},\msc{C}')
\]
with $\opn{can}_0=l|_{\msc{C}_{/x}}$ and $\opn{can}_1=\underline{x}$.  Rather, the coslice equivalence chooses a transformation $\opn{can}:l|_{\msc{C}_{/x}}\to \underline{x}$, and we apply $F'$ to obtain the transformation $F'(\opn{can}):F'l|_{\msc{C}_{/x}}\to F'(\underline{x})$.  All we need to know about this transformation is that it evaluates at each object $t_{\tilde{x}}=(\tilde{x},t:l(\tilde{x})\to x)$ in $\msc{C}_{x/}$ to the chosen map $\opn{can}_{t_{\tilde{x}}}=t:l(x)\to x$.

To say that $F'$ is a left Kan extension for $F$ along $l$ is, by definition, to say that the transformation \eqref{eq:2922} realizes $F'(x)$ as a colimit for the diagram $F|_{\msc{C}_{/x}}$, for each $x$ in $\msc{C}'$.  Let us denote the transformation \eqref{eq:2922} by $\opn{can}^F_u$ and argue this point.
\par

Since the functor $l$ admits a right adjoint $R:\msc{C}'\to \msc{C}$, the counit map $\varepsilon_x:lR(x)\to x$ is a terminal object $(R(x),\varepsilon_x)$ in the relative overcategory $\msc{C}_{x/}$ \cite[\href{https://kerodon.net/tag/02J9}{02J9}, \href{https://kerodon.net/tag/038S}{038S}]{kerodon}.  Hence the transformation $\opn{can}^F_u$ realizes $F'(x)$ as a colimit of the given diagram if and only if its restriction along the inclusion $(R(x),\varepsilon_x):\ast\to \msc{C}_{/x}$ realizes $F'(x)$ as a colimit for the trivial diagram $FR(x):\ast\to \msc{E}$ \cite[\href{https://kerodon.net/tag/02XW}{02XW}]{kerodon}. But now we simply observe
\[
(\opn{can}^F_u)_{(Rx,\varepsilon_x)}=\left(FR(x)\overset{F(u_{R(x)})}\to FRlR(x)\overset{FR(\varepsilon_x)}\to FR(x)\right).
\]
By the definition of a unit and counit for an adjunction \cite[\href{https://kerodon.net/tag/02EJ}{02EJ}]{kerodon} this composite is isomorphic to the identity, and in particular an isomorphism. Hence $\opn{can}^F_u|_{(Rx,\varepsilon_x)}$ is a colimit diagram, $\opn{can}^F_u$ realizes all $F'(x)$ as a colimit for the corresponding diagram $F|_{\msc{C}_{/x}}$, and $F'$ is realized as a left Kan extension for $F$.
\end{proof}

As with adjoints, we are interested in a relative notion of Kan extension.

\begin{definition}
Consider a partial diagram of isofibrations
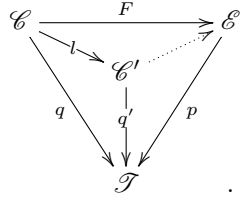
\[
\xymatrixrowsep{2mm}
\xymatrix{
\msc{C}\ar[ddddr]_q\ar[rr]^{F}\ar[dr]|l & & \msc{E}\ar[ddddl]^p\\
	& \msc{C}'\ar[ddd]|(.4){q'}\ar@{..>}[ur]\\\\\\
	& \msc{T} & .
	}
\]
A left Kan extension $F':\msc{C}'\to \msc{E}$ of $F$ along $l$ is said to be a relative left Kan extension over $\msc{T}$ if $F'$ is a map over $\msc{T}$ and if the given transformation $\alpha:F\to F'l$ satisfies $p(\alpha)=id_q$.
\end{definition}

We obtain immediately a relative version of Lemma \ref{lem:2916}.

\begin{lemma}\label{lem:adj_kan}
Consider isofibrations $q:\msc{C}\to \msc{T}$ and $q':\msc{C}'\to \msc{T}$.  If $l:\msc{C}\to \msc{C}'$ admits a right adjoint over $\msc{T}$ then any functor $F:\msc{C}\to \msc{E}$ over $\msc{T}$ admits a left Kan extension over $\msc{T}$.  In particular, for such a right adjoint $\{R,\varepsilon\}$, any choice $u:id_{\msc{C}}\to Rl$ of a unit over $\msc{T}$ provides a left Kan extension $\{FR,F(u)\}$ over $\msc{T}$.
\end{lemma}

\begin{remark}
We expect that one can establish the existence of relative Kan extensions, in the cocartesian setting, by checking at the fibers (cf.\ \cite[Proposition 4.3.3.10]{lurie09}).  However, we have no need for such a result. 
\end{remark}

\subsection{Homotopy vs.\ derived $\infty$-categories}
\label{sect:kan_htop_der}

We consider a standard example concerning the construction of the derived $\infty$-category from the homotopy $\infty$-category, via reflective localization.
\par

Consider a finitely generated linear abelian category $A^{\heartsuit}$, in the precise sense of Definition \ref{def:lin_ab}, and its homotopy $\infty$-category $\msc{K}=\msc{K}(A^{\heartsuit})$.  Suppose also that $A^{\heartsuit}$ has enough projectives, so that $\msc{K}$ has enough $K$-projective complexes.  Explicitly, we assume that each complex $x$ admits a quasi-isomorphism $P\to x$ from a complexes $P$ for which the functor $\Hom^{\ast}_A(P,-)$ preserves quasi-isomorphisms, and hence for which the inner-Hom functor $\underline{\Maps}_{\msc{K}}(P,-)$ preserves quasi-isomorphisms. See for example \cite[Theorems 2.3, 2.6]{krizmay95} or \cite[Lemma 13.3]{drinfeld04}.
\par

Dually, the $\infty$-category $\msc{K}_{\pf}=\msc{K}(A^{\heartsuit}_{\pf})$ admits enough $K$-injectives. We take
\[
\opn{Inj}\msc{K}^-=(\opn{Proj}\msc{K})^{op}\ \ \text{and}\ \ \opn{Inj}\msc{K}^+=\opn{Inj}\msc{K}_{\pf}
\]
the full subcategories spanned by $K$-injectives in $\msc{K}^{-}=\msc{K}(A^{\heartsuit})^{op}$ and $\msc{K}^+=\msc{K}(A^{\heartsuit})_{\pf}$, respectively.  We note that isomorphisms in these subcategories, which are homotopy equivalences, agree precisely with quasi-isomorphisms, and the statement that we have enough $K$-projectives is precisely the statement that the full subcategories
\[
\opn{Inj}\msc{K}^{\pm}\subseteq \msc{K}^{\pm}
\]
are reflective \cite[\href{https://kerodon.net/tag/02F5}{01F5}]{kerodon}.  We therefore obtain the following.

\begin{lemma}
Let $A^{\heartsuit}$ be a finitely generated linear abelian category with enough projectives, and take $\msc{K}^{\pm}$ one of $\msc{K}(A^{\heartsuit})^{op}$ or $\msc{K}(A^{\heartsuit}_{\pf})$.  The fully faithful inclusion $R^{\pm}:\opn{Inj}\msc{K}^{\pm}\to \msc{K}^{\pm}$ admits a left adjoint
\[
l:\msc{K}^{\pm}\to \opn{Inj}\msc{K}^{\pm}.
\]
Furthermore, this left adjoint realizes $\opn{Inj}\msc{K}^{\pm}$ as a localization of $\msc{K}^{\pm}$ along the class of quasi-isomorphisms.
\end{lemma}

\begin{proof}
Given the above analysis, this is an immediate consequence of reflective localization. See in particular \cite[\href{https://kerodon.net/tag/02FA}{02FA}, \href{https://kerodon.net/tag/02FE}{02FE}, \href{https://kerodon.net/tag/04JL}{04JL}]{kerodon}.
\end{proof}

Uniqueness for localizations now gives the following.

\begin{corollary}\label{cor:2953}
Let $A^{\heartsuit}$ be a finitely generated linear abelian category with enough projectives, and let $\msc{D}^{\pm}=\msc{K}^{\pm}[\opn{Qiso}^{-1}]$ denote one of the derived $\infty$-categories $\msc{D}^-=\msc{D}(A^{\heartsuit})^{op}$ or $\msc{D}^+=\msc{D}(A^{\heartsuit}_{\pf})$.  The following hold:
\begin{enumerate}
\item The localization functor $loc:\msc{K}^{\pm}\to \msc{D}^{\pm}$ admits a fully faithful right adjoint $R^{\pm}:\msc{D}^{\pm}\to \msc{K}^{\pm}$.\vspace{1mm}
\item For each $x$ in $\msc{K}^{\pm}$, the unit transformation $u_x:x\to R^{\pm}\circ loc(x)$ is a quasi-isomorphism, and for each $y$ in $\msc{D}^{\pm}$ the counit transformation $\varepsilon_y:loc\circ R^{\pm}(y)\to y$ is an isomorphism.\vspace{1mm}
\item The right adjoint $R^{\pm}$ is an equivalence onto the full subcategory $\opn{Inj}\msc{K}^{\pm}$ spanned by $K$-injectives in $\msc{K}^{\pm}$.\vspace{1mm}
\item The localization functor restricts to an equivalence $\opn{Inj}\msc{K}^{\pm}\overset{\sim}\to \msc{D}^{\pm}$.
\end{enumerate}
\end{corollary}

Let us suppose now that $A^{\heartsuit}$ is a finite ribbon tensor category and consider the corresponding cocartesian fibration $q:\msc{K}^{\ot}\to fr\Disk^{\ot}$ for the homotopy $\infty$-category (see Section \ref{sect:fr_htop_der}).  We consider the full subcategories $\opn{Inj}\msc{K}^{\pm}$ and, as in Section \ref{sect:subcats}, and the associated full symmetric monoidal subcategory $\opn{Inj}\msc{K}^{\ot}$ in $\msc{K}^{\ot}$ whose fibers, along objects $t=(\mcl{I},\mu)$ in $fr\Disk^{\ot}$, fit into pullback squares
\[
\xymatrix{
\opn{Inj}\msc{K}^{\ot}_t\ar[d]\ar[rr] & & \opn{Inj}\msc{K}^{\mcl{I}}\ar[d]\\
\msc{K}^{\ot}_t\ar[rr]_{\rho_!} & & \msc{K}^{\mcl{I}}.
}
\]
We have the following particular instance of Lemma \ref{lem:full_to_symm}.

\begin{lemma}\label{lem:inj_k}
Let $A^{\heartsuit}$ be a finite ribbon tensor category. The subcategory $\opn{Inj}\msc{K}^{\ot}$ of $K$-injective complexes is a symmetric monoidal $\infty$-category, and the inclusion $\opn{Inj}\msc{K}^{\ot}\to \msc{K}^{\ot}$ is a symmetric monoidal functor. Furthermore, $\opn{Inj}\msc{K}^{\ot}$ is stable under isomorphism in $\msc{K}^{\ot}$.
\end{lemma}

As a consequence of isoclosure we see that the structure map $\msc{K}^{\ot}\to fr\Disk^{\ot}$ restricts to an isofibration $\opn{Inj}\msc{K}^{\ot}\to fr\Disk^{\ot}$.

\begin{proposition}\label{prop:3009}
Suppose $A^{\heartsuit}$ is a finite ribbon tensor category.  For $q:\msc{K}^{\ot}\to fr\Disk^{\ot}$ and $q':\msc{D}^{\ot}\to fr\Disk^{\ot}$ the corresponding cocartesian fibration for the homotopy and derived $\infty$-categories, the following hold:
\begin{enumerate}
\item The (vertical) localization functor $loc:\msc{K}^{\ot}\to \msc{D}^{\ot}$ admits a relative right adjoint $R:\msc{D}^{\ot}\to \msc{K}^{\ot}$ over $fr\Disk^{\ot}$.\vspace{1mm}
\item The right adjoint $R$ is fully faithful, and an equivalence onto the full subcategory $\opn{Inj}\msc{K}^{\ot}\subseteq \msc{K}^{\ot}$ spanned by $K$-injectives.\vspace{1mm}
\item The localization functor $loc:\msc{K}^{\ot}\to \msc{D}^{\ot}$ restricts to an equivalence of isofibrations $\opn{Inj}\msc{K}^{\ot}\overset{\sim}\to \msc{D}^{\ot}$.
\end{enumerate}
\end{proposition}

\begin{proof}
(1) Over any object $(\mcl{I},\mu)$ the fiber of localization reproduces the localization functor
\[
\opn{loc}^{\mcl{I}}:\msc{K}^{\mcl{I}}\to \msc{D}^{\mcl{I}},
\]
by Proposition \ref{prop:fib_localize}, and each fiber $\opn{loc}^{\mcl{I}}$ admits a fully faithful right adjoint $R^{\mcl{I}}:\msc{D}^{\mcl{I}}\to \msc{K}^{\mcl{I}}$ by Corollary \ref{cor:2953}.  Hence $loc$ admits a right adjoint $R:\msc{D}^{\ot}\to \msc{K}^{\ot}$ over $fr\Disk^{\ot}$, by Proposition \ref{prop:reladj_exists}.
\par

(2) \& (3) By Proposition \ref{prop:reladj_pullback} the fiber of $R$ over any $(\mcl{I},\mu)$ recovers the right adjoint $R^{\mcl{I}}$ to the localization functor $\opn{loc}^{\mcl{I}}$, and we know that this functor has image in the full subcategory $\opn{Inj}\msc{K}^{\mcl{I}}$ by Corollary \ref{cor:2953} (3).  Furthermore, again after reducing to the fibers, the same result tells us that the unit and counit transformations provide natural isomorphisms
\[
x\overset{\sim}\to R\circ loc(x)\ \ \text{and}\ \ loc\circ R(y)\overset{\sim}\to y
\]
whenever $x$ is in $\opn{Inj}\msc{K}^{\ot}$ and $y$ is in $\msc{D}^{\ot}$.  However, the unit and counits for the fibers $R^{\mcl{I}}$ and $loc^{\mcl{I}}$ are just obtained from the unit and counit transformations for $R$ and $loc$.  In particular, these transformations realize $R$ and $loc|_{\opn{Inj}\msc{K}^{\ot}}$ as mutually inverse equivalences over $fr\Disk^{\ot}$.
\end{proof}

\begin{corollary}\label{cor:3075}
Suppose we are in the situation of Proposition \ref{prop:3009}. The structure map $\msc{K}^{\ot}\to fr\Disk^{\ot}$ restricts to a cocartesian fibration $\opn{Inj}\msc{K}^{\ot}\to fr\Disk^{\ot}$. Furthermore, the right adjoint $R:\msc{D}^{\ot}\to \msc{K}^{\ot}$ to localization is a fully faithful symmetric monoidal functor which restricts to an equivalence of cocartesian fibration $\msc{D}^{\ot}\to \opn{Inj}\msc{K}^{\ot}$.
\end{corollary}

\begin{proof}
The first and third claims follow by \cite[\href{https://kerodon.net/tag/028A}{028A}]{kerodon}, and the second claim follows, subsequently, from the fact that the right adjoint factors as a composite
\[
\msc{D}^{\ot}\overset{\sim}\to \opn{Inj}\msc{K}^{\ot}\hookrightarrow \msc{K}^{\ot}
\]
of an equivalence of cocartesian fibrations composed with an inclusion of symmetric monoidal $\infty$-categories.
\end{proof}

\begin{remark}
Since $K$-injective complexes are stable under tensoring, in the complete setting, the only thing that is stopping the inclusion $\opn{Inj}\msc{K}^{\ot}\to \msc{K}^{\ot}$ from being a map of cocartesian fibrations over $fr\Disk^{\ot}$ is the unit.  Specifically, the unit in the category of $K$-injectives is provided by an injective resolution $I$ of the unit in $\msc{K}^{\pm}$, so that the inclusion does not preserve unit object. See Lemma \ref{lem:3386} below.
\end{remark}

\begin{corollary}\label{cor:der_relkan}
Suppose $A^{\heartsuit}$ is a finite ribbon tensor category. Then for any functor $\msc{T}\to fr\Disk^{\ot}$, any $\infty$-category $\msc{E}$ over $\msc{T}$, and any map $F:\msc{T}\times_{fr\Disk^{\ot}}\msc{K}^{\ot}\to \msc{E}$ of $\infty$-categories over $\msc{T}$, a relative left Kan extension
\[
F':\msc{T}\times_{fr\Disk^{\ot}}\msc{D}^{\ot}\to \msc{E}
\]
along $loc_{\msc{T}}:\msc{T}\times_{fr\Disk^{\ot}}\msc{K}^{\ot}\to \msc{T}\times_{fr\Disk^{\ot}}\msc{D}^{\ot}$ exists.  This extension is explicitly given by the composite $F'=FR_{\msc{T}}$.
\end{corollary}

\begin{proof}
Apply Proposition \ref{prop:3009} and Lemma \ref{lem:adj_kan}.
\end{proof}

\subsection{Uniqueness of Kan extension over a discrete base}

We consider the following particular expression of the universal property for Kan extensions.

\begin{proposition}\label{prop:kan_univ}
Consider functors $F:\msc{C}\to \msc{E}$ and $l:\msc{C}\to \msc{C}'$.  For any left Kan extension $F':\msc{C}'\to \msc{E}$, the structural transformation $\alpha:F\to F'l$ realizes $F'$ as an initial object in the relative undercategory
\[
\Fun(\msc{C}',\msc{E})_{F/}=\Fun(\msc{C}',\msc{E})\times_{\Fun(\msc{C},\msc{E})}\Fun(\msc{C},\msc{E})_{F/}.
\]
\end{proposition}

\begin{proof}
For any choice of functor $G:\msc{C}'\to \msc{E}$, the usual universal property for Kan extensions tells us that the sequence
\[
H_{\Fun(\msc{C}',\msc{E})}(F',G)\overset{l^{\ast}}\to H_{\Fun(\msc{C},\msc{E})}(F'l,Gl)\overset{-\circ \alpha}\to H_{\Fun(\msc{C},\msc{E})}(F,Gl)
\]
is a homotopy equivalence \cite[\href{https://kerodon.net/tag/0309}{0309}]{kerodon}.  (Here the $H_{\star}$ are arbitrary choices of $\Hom$ functors \cite[\href{https://kerodon.net/tag/03MR}{03MR}]{kerodon}.)  When we choose $G=F'$, this functor sends the identity $id_F$ to $\alpha$.  We therefore have a natural isomorphism
\[
H_{\Fun(\msc{C}',\msc{E})}(F',-)\overset{\sim}\to H_{\Fun(\msc{C},\msc{E})}(F,-l)
\]
of functors from $\Fun(\msc{C}',\msc{E})$ to $\sKan$ and subsequent equivalence between the corresponding cocartesian fibrations
\[
\xymatrix{
\Fun(\msc{C}',\msc{E})_{F'/}\ar[rr]^(.4){\sim}\ar[dr]& & \Fun(\msc{C}',\msc{E})\times_{\Fun(\msc{C},\msc{E})}\Fun(\msc{C},\msc{E})_{F/}\ar[dl]\\
	& \Fun(\msc{C}',\msc{E})
	}
\]
\cite[\href{https://kerodon.net/tag/02GT}{02GT}, \href{https://kerodon.net/tag/01QC}{01QC}]{kerodon} which sends the object $id_{F'}$ to $(F',\alpha)$.  Since $id_{F'}$ is initial in $\Fun(\msc{C}',\msc{E})_{F'/}$ \cite[\href{https://kerodon.net/tag/02J2}{02J2}]{kerodon}, $\alpha$ is initial in the given fiber product.
\end{proof}

This universality statement directly implies uniqueness of Kan extensions over a discrete base.

\begin{proposition}\label{prop:kan_unique}
Let $Z$ be a discrete category, and consider a partial diagram
\begin{equation}\label{eq:2925}
\xymatrixrowsep{2mm}
\xymatrix{
\msc{C}\ar[ddddr]_q\ar[rr]^{F}\ar[dr]|l & & \msc{E}\ar[ddddl]^p\\
	& \msc{C}'\ar[ddd]|(.4){q'}\ar@{..>}[ur]\\\\\\
	& Z
	}
\end{equation}
of maps between cocartesian fibrations over $Z$.  Suppose also that one of the following conditions hold:
\begin{enumerate}
\item[i)] For each $x'$ in $\msc{C}'$ the relative overcategory $\msc{C}_{/x'}$ is nonempty.
\item[ii)] $\msc{E}$ has no initial object.
\item[iii)] $\msc{E}$ admits an initial object $o$ for which $\opn{Aut}_Z(p(o))=1$.
\end{enumerate}
If a left Kan extension of $F$ over $Z$ exists, then this extension is uniquely determined up to a contractible spaces of choices.  In particular, the full subcategory
\[
\msc{L\!K}\!an(F;Z)\ \subseteq\ \Fun_Z(\msc{C}',\msc{E})\times_{\Fun_Z(\msc{C},\msc{E})}\Fun_Z(\msc{C},\msc{E})_{F/}
\]
spanned by left Kan extensions over $Z$ is a contractible Kan complex.
\end{proposition}

\begin{proof}
Consider two $\infty$-categories over $Z$, $q_i:\msc{B}_i\to Z$. Since $Z$ is discrete, the subcomplex
\[
\Fun_{Z}(\msc{B}_0,\msc{B}_1)\ \subseteq\ \Fun(\msc{B}_0,\msc{B}_1)
\]
is the subcategory spanned by all functors over $Z$ and all transformations $\zeta:F\to F'$ between such functors for which $q_1(\zeta):q_0\to q_0$ is the identity.  Furthermore, transformations in $\Fun_{Z}(\msc{B}_0,\msc{B}_1)$ satisfy the $2$-of-$3$ property in that for any diagram
\[
\xymatrix{
	& F'\ar[dr]^{\zeta_1}\\
F\ar[ur]^{\zeta_0}\ar[rr]_{\zeta_2} & & F''
}
\]
in which two of the $\zeta_j$ are in $\Fun_Z(\msc{B}_0,\msc{B}_1)$, the third transformation is also in $\Fun_Z(\msc{B}_0,\msc{B}_1)$.
\par

Considering now the situation \eqref{eq:2925}. If we have two left Kan extensions over $Z$, say $F_0,F_1:\msc{C}'\to \msc{E}$ with transformations $\alpha_i:F\to F_il$, then the universal property of Kan extensions promises the existence of a natural isomorphism $\zeta:F_0\to F_1$ which completes a diagram
\[
\xymatrix{
	& F_0l\ar[dr]^{\zeta l}\\
F\ar[ur]^{\alpha_0}\ar[rr]_{\alpha_1} & & F_1 l.
}
\]
Since the $\alpha_i$ are transformations over $Z$ we conclude that $\zeta l$ is a transformation over $Z$ as well, and hence that $\zeta$ itself is a transformation over $Z$.
\par

To see this clearly, we consider an object $x'$ in $\msc{C}'$ and the corresponding map $p(\zeta_{x'}):q'(x')\to q'(x')$.  If $x'$ admits a map $t:l(x)\to x'$ from some $x$ in $\msc{C}$ then, via cocartesianness of $q$ and $q'$ and the fact that $l$ is a map of cocartesian fibrations, we can assume that $x$ and $t$ lie in a single fiber over $Z$, i.e.\ that $q(x)=q'(x')$ and $q'(t)=id_{q'(x')}$.  In this case, we have $p(\zeta_{l(x)})=id_{q(x)}$ and the diagram
\[
\xymatrix{
q'(l(x))\ar[rr]^{p(\zeta_{l(x)})=id}\ar[d]_{q'(t)=id} & & q'(l(x))\ar[d]^{q'(t)=id}\\
q'(x')\ar[rr]_{p(\zeta_{x'})} & & q'(x')
}
\]
in $Z$ implies $p(\zeta_{x'})=id$.
\par

On the other hand, if $x'$ admits no maps from an object of the form $l(x)$ then $F_0(x)$ and $F_1(x)$ are both colimits over the empty diagram $\msc{C}_{/x}\to \msc{E}$, and are hence initial objects in $\msc{E}$ which lie in the same fiber over $Z$.  In this case assumption (iii) must hold, which forces $p(\zeta_{x'})=id_{q'(x)}$.
\par

Let us take now
\[
\Pi_Z=\Fun_Z(\msc{C}',\msc{E})\times_{\Fun_Z(\msc{C},\msc{E})}\Fun_Z(\msc{C},\msc{E})_{F/}
\]
and
\[
\Pi=\Fun(\msc{C}',\msc{E})\times_{\Fun(\msc{C},\msc{E})}\Fun(\msc{C},\msc{E})_{F/}.
\]
We note that $\Pi_Z$ includes into $\Pi$ as the subcategory spanned by functors and transformations over $Z$.  The above arguments tell us that any map $\zeta:\{F_0,\alpha_0\}\to \{F_1,\alpha_1\}$ between left Kan extensions over $Z$ in $\Pi$ is a map in $\Pi_Z$.  So we see that the full subcategory
\[
\msc{L\!K}\!an(F;Z)\ \subseteq\ \Pi
\]
spanned by left Kan extensions over $Z$, all of which are initial in $\Pi$ by Proposition \ref{prop:kan_univ}, already lies in $\Pi_Z$.  Since the full subcategory spanned by initial objects in any $\infty$-category is contractible--provided it's not empty--we conclude that $\msc{L\!K}\!an(F;Z)$ is a contractible full subcategory in $\Pi_Z$.
\end{proof}

\begin{remark}
Na\"ive versions of Proposition \ref{prop:kan_unique} do not hold when the base $Z$ is allowed to be an $\infty$-category. Also, one can check in examples that \emph{some} hypotheses are needed in order to establish uniqueness for relative Kan extensions.
\end{remark}

As an example we can consider the case of a discrete operad. When $Z$ is an $\infty$-operad and $p:\msc{E}\to Z$ is a map of $\infty$-operads, then one of the conditions (ii) or (iii) must hold.  In particular, if $\msc{E}$ has an initial object that object lies in the contractible fiber $\msc{E}_{\emptyset}$ over the empty set in $\Fin_{\ast}$.  Such an object maps to an object in $Z_{\emptyset}$, and contractibility of this discrete category implies that all automorphism groups of objects in $Z_{\emptyset}$ are trivial.

\begin{corollary}\label{cor:kan_unique_operad}
Suppose that we have a partial diagram of cocartesian fibrations as in \eqref{eq:2925}, with $Z$ discrete.  Suppose also that $Z$ is an $\infty$-operad and that $p:\msc{E}\to Z$ is a map of $\infty$-operads.  If a Kan extension $F':\msc{C}'\to \msc{E}$ over $Z$ exists, then this extension is unique up to a contractible space of choices.
\end{corollary}

Our main example of interest here is $Z=\Fin_{\ast}$, though Corollary \ref{cor:kan_unique_operad} applies to other operads as well, including those that control associative monoidal $\infty$-categories and module $\infty$-categories over (symmetric) monoidal $\infty$-categories.
\par

Considering the case where the functor $F:\msc{C}\to \msc{E}$ is the identity on $\msc{C}$, in the statement of Proposition \ref{prop:kan_unique}, we obtain the following.

\begin{corollary}\label{cor:unique_kan_ra}
If a map of cocartesian fibrations over a discrete base $Z$ admits a relative right adjoint over $Z$, then this adjoint is unique up to a contractible space of choices.
\end{corollary}

\begin{proof}
Given such a relative right adjoint $G$ to $l$ we only note that condition (i) holds, as we always have the object $\{G(x'),\varepsilon_{x'}:lG(x')\to x'\}$ in $\msc{C}_{/x'}$.
\end{proof}

\begin{definition}
Given a situation as in the statement of Proposition \ref{prop:kan_unique}, we denote the relative left Kan extension of $F:\msc{C}\to \msc{E}$ along $l:\msc{C}\to \msc{C}'$, if it exists, by $\opn{Lan}_l(F;Z):\msc{C}'\to \msc{E}$.
\end{definition}

When no confusion will arise we drop the subscript $l$, and write simply $\opn{Lan}(F;Z)$ for $\opn{Lan}_l(F;Z)$.

\section{Deriving LRT at the $\infty$-categorical level}
\label{sect:derived_lrt}

For a finite modular tensor category $A^{\heartsuit}$ with derived $\infty$-category $\msc{D}=\msc{D}(A^{\heartsuit})$, we construct a topological field theory
\[
L_{\msc{D}}:\Bord^{nc}_{\msc{D}}\to \Vect_{\pf}
\]
from the $3$-dimensional anomalous bordism category with labels in $\msc{D}$. This theory is constructed by taking a Kan extension of the homotopical theory $L_{\msc{K}}$ along the localization functor $\Bord^{nc}_{\msc{K}}\to \Bord^{nc}_{\msc{D}}$. We then apply results from Section \ref{sect:rel_adj_kan} and Appendix \ref{sect:module_cats} to calculate the state spaces for $L_{\msc{D}}$ via maps in $\msc{D}$.
\par

We show that the restriction $L_{\msc{D}_{fin}}$ of $L_{\msc{D}}$ to the full subcategory $\Bord^{nc}_{\msc{D}_{fin}}$ of bordisms with rigid labels defines a field theory which is valued in $\Vect$, rather than $\Vect_{\pf}$. Finally we show that the cohomology of this theory recovers the original abelian theory from De Renzi et al.\ \cite{derenzietal23},
\[
H^0(L_{\msc{D}_{fin}}|_{\Bord^{nc}_{A^{\heartsuit}_{fin}}})\overset{\sim}\to \text{Original Lyubashenko theory }L_{A^{\heartsuit}_{fin}}.
\]
The theory $L_{\msc{D}_{fin}}$ can therefore be viewed as a \emph{derivation} of the original abelian Lyubashenko theory from \cite{derenzietal23}, or conversely, we can view the abelian theory as a homological truncation of a fundamentally derived theory $L_{\msc{D}_{fin}}$.

\subsection{Derived LRT theories}

\begin{theorem}\label{thm:der_lrt}
Consider a finite modular tensor category $A^{\heartsuit}$, and let $\msc{K}^{\ot}\to fr\Disk^{\ot}$ and $\msc{D}^{\ot}\to fr\Disk^{\ot}$ be the corresponding $fr\Disk$-categories for the homotopy and derived $\infty$-categories of unbounded cochains (see Section \ref{sect:fr_htop_der}).  The symmetric monoidal functor $L_{\msc{K}}:\Bord_{\msc{K}}^{nc}\to \Vect_{\pf}$ from Proposition \ref{prop:htop_lrt} admits a left Kan extension
\[
L_{\msc{D}}:\Bord_{\msc{D}}^{nc}\to \Vect_{\pf}
\]
and this Kan extension admits a canonical symmetric monoidal structure.
\end{theorem}

Let us expand on the above statement in order to clarify the point about symmetric monoidality.  We have the symmetric monoidal functor $loc^{\ot}:(\Bord_{\msc{K}}^{nc})^{\ot}\to (\Bord_{\msc{D}}^{nc})^{\ot}$ induced by the localization functors $\msc{K}^{\pm}\to\msc{D}^{\pm}$ on the marking categories (see Section \ref{sect:fr_htop_der}, Proposition \ref{prop:fr_loc}).  As these localization functors from $\msc{K}^{\pm}$ admit right adjoints, the corresponding localization functor for bordisms admits a relative right adjoint
\[
R^{\ot}:(\Bord_{\msc{D}}^{nc})^{\ot}\to (\Bord_{\msc{K}}^{nc})^{\ot}
\]
over $(\Bord_{\ast}^{nc})^{\ot}$, and this right adjoint recovers the right adjoints on the marking categories $\msc{K}^{\mcl{I}}\cong(\Bord_{\msc{K}}^{\ot})_{\Sigma_{\zeta}}$ and $\msc{D}^{\mcl{I}}\cong(\Bord_{\msc{D}}^{\ot})_{\Sigma_{\zeta}}$ over each marked surface $\Sigma_{\zeta}$ in $\Bord_{\ast}\subseteq (\Bord_{\ast})^{\ot}$ (Propositions \ref{prop:reladj_exists}, \ref{prop:reladj_pullback}).
\par

Now, given that we have such a right adjoint, relative left Kan extensions along arbitrary functors $F:(\Bord_{\msc{D}}^{nc})^{\ot}\to \msc{E}$ over $(\Bord_{\ast})^{\ot}$ exist, by (Corollary \ref{cor:der_relkan}). Furthermore, all such extensions are obtained via precomposition with $R^{\ot}$.  In particular, we have the relative left Kan extension
\begin{equation}
L_{\msc{D}}^{\ot}:=\opn{Lan}_{loc^{\ot}}\left(L^{\ot}_{\msc{K}};(\Bord_{\ast}^{nc})^{\ot}\right)=L^{\ot}_{\msc{K}}\circ R^{\ot}:
\end{equation}
\[
(\Bord_{\msc{D}}^{nc})^{\ot}\to \Vect_{\pf}^{\ot}.
\]
We also understand that the fiber of $L_{\msc{D}}^{\ot}$ along the inclusion $\Bord^{nc}_{\ast}\to (\Bord^{nc}_{\ast})^{\ot}$ recovers the left Kan extension $L_{\msc{D}}$ of the underlying functor $L_{\msc{K}}:\Bord_{\msc{K}}^{nc}\to \Vect_{\pf}$ along localization.  We are claiming, in the statement of Theorem \ref{thm:der_lrt}, that the map $L_{\msc{D}}^{\ot}$ is a map of cocartesian fibrations over $\Fin_{\ast}$, i.e.\ that it is a symmetric monoidal functor.

\begin{proof}[Proof of Theorem \ref{thm:der_lrt}]
Since the map $L^{\ot}_{\msc{K}}:(\Bord_{\msc{K}})^{\ot}\to \Vect^{\ot}_{\pf}$ is symmetric monoidal it suffices to establish symmetric monoidality of the right adjoint $R^{\ot}:(\Bord_{\msc{D}}^{nc})^{\ot}\to (\Bord_{\msc{K}}^{nc})^{\ot}$. However, this point is apparent, though we take a moment to account for the details.

In the notation of Proposition \ref{prop:reladj_pullback}, $R^{\ot}=R_{(\Bord_{\ast})^{\ot}}$ is the unique functor which completes a strictly commuting diagram
\[
\xymatrixrowsep{4mm}
\xymatrixcolsep{4mm}
\xymatrix{
(\Bord_{\msc{D}}^{nc})^{\ot}\ar[rr]\ar@{-->}[dr]\ar[dddr] & & \msc{D}^{\ot}\ar[dr]^{R}\ar[dddr]|(.39){\hole}\\
 & (\Bord_{\msc{K}}^{nc})^{\ot}\ar[rr]\ar[dd] & & \msc{K}^{\ot}\ar[dd] \\\\
 & (\Bord_{\ast}^{nc})^{\ot}\ar[rr] & & fr\Disk^{\ot}.
}
\]
By Corollary \ref{cor:3075} the functor $R$ is symmetric monoidal, and the two projections
\[
(\Bord_{\ast}^{nc})^{\ot}\leftarrow (\Bord_{\msc{D}}^{nc})^{\ot}\to \msc{D}^{\ot}
\]
are symmetric monoidal functors as well. Furthermore, by Proposition \ref{prop:fp_symm}, the cocartesian edges in $(\Bord_{\msc{K}}^{nc})^{\ot}$, as a fibration over $\Fin_{\ast}$, are precisely those edges which have cocartesian image in both $\msc{K}^{\ot}$ and $(\Bord^{nc}_{\ast})^{\ot}$. From these two points it follows immediately that the completing map $R^{\ot}:(\Bord_{\msc{D}}^{nc})^{\ot}\to (\Bord_{\msc{K}}^{nc})^{\ot}$ is a map of cocartesian fibrations over $\Fin_{\ast}$, i.e.\ a symmetric monoidal functor.
\end{proof}

\begin{remark}
As far as uniqueness, $L_{\msc{D}}$ is uniquely determined as a map in $\sCat_{\infty}$. The author cannot determine if $L_{\msc{D}}$ is uniquely determined as a functor over $\Bord_{\ast}^{nc}$.  In any case, having specified $L_{\msc{D}}$, the symmetric monoidal structure \emph{is} uniquely determined up to a contractible space of choices by Corollary \ref{cor:unique_kan_ra}.
\end{remark}

\begin{definition}
Given a finite modular tensor category $A^{\heartsuit}$, we refer to the symmetric monoidal functor
\[
L_{\msc{D}}:\Bord^{nc}_{\msc{D}}\to \Vect_{\pf}
\]
from Theorem \ref{thm:der_lrt} as the derived LRT theory for $\msc{D}=\msc{D}(A^{\heartsuit})$.
\end{definition}

\subsection{Restriction to the finite derived category}

For a modular tensor category $A^{\heartsuit}$ we consider the full subcategory $\msc{D}_{fin}^{\ot}$ spanned by those complexes with bounded, finite length cohomology in $\msc{D}^{\ot}$.  More precisely, as in Section \ref{sect:kfin_lrt}, the fibers in this $\infty$-category over objects $t=(\mcl{I},\mu)$ in $fr\Disk^{\ot}$ are the pullbacks
\[
\xymatrix{
(\msc{D}_{fin}^{\ot})_t\ar[rr]^{\sim}\ar[d] & & \msc{D}_{fin}^{\mcl{I}}\ar[d]\\
\msc{D}^{\ot}_t\ar[rr]_{\rho_!}^{\sim} & & \msc{D}^{\mcl{I}}.
}
\]
By Proposition \ref{prop:frdisk_subcats}, $\msc{D}_{fin}^{\ot}$ is a $fr\Disk$-subcategory in $\msc{D}^{\ot}$.

\begin{theorem}\label{thm:der_fin_lrt}
Let $A^{\heartsuit}$ be a finite modular tensor category, and $\msc{D}_{fin}^{\ot}\to fr\Disk^{\ot}$ be the associated $fr\Disk$-category of complexes with finite cohomology. There is a $\Vect$-valued topological field theory
\[
L_{\msc{D}_{fin}}:\Bord_{\msc{D}_{fin}}^{nc}\to \Vect
\]
which is obtained from $L_{\msc{D}}$ via restriction.
\end{theorem}

\begin{proof}
We have the symmetric monoidal functor $L_{\msc{K}_{l\text{-}fin}}:\Bord_{\msc{K}_{l\text{-}fin}}^{nc}\to \Vect$ of Proposition \ref{prop:kfin_lrt}, which is obtained from $L_{\msc{K}}$ via restriction and a change of the target category according to the pair of fully faithful embeddings
\[
\Vect_{\pf}\leftarrow \Vect_{l\text{-}fin}^+\to \Vect.
\]

Since the right adjoint $R:\msc{D}^{\ot}\to \msc{K}^{\ot}$ restricted to $\msc{D}_{fin}^{\ot}$ has image in $\msc{K}_{l\text{-}fin}^{\ot}$, the explicit expression $L_{\msc{D}}=L_{\msc{K}}R$ gives an identification of $\Vect_{l\text{-}fin}^+$-valued functors
\[
L_{\msc{D}}|_{\Bord_{\msc{D}_{fin}}^{nc}}=L_{\msc{K}}R|_{\Bord_{\msc{D}_{fin}}}\cong L_{\msc{K}_{l\text{-}fin}}R|_{\Bord_{\msc{D}_{fin}}}.
\]
After embedding into $\Vect$, we obtain the proposed symmetric monoidal functor
\[
L_{\msc{D}_{fin}}=L_{\msc{K}_{l\text{-}fin}}R|_{\Bord_{\msc{D}_{fin}}}:\Bord_{\msc{D}_{fin}}^{nc}\to \Vect.
\]
\end{proof}

\subsection{Calculation of the state spaces}

\begin{proposition}\label{prop:der_states}
Let $\Sigma_{\zeta}$ be a connected genus $g$ surface with markings $\zeta:D\times \mcl{I}\to \Sigma$.  For any embedding $f:D\times \mcl{I}\to D$ into a positively colored disk there is a natural isomorphism
\[
L_{\msc{D}}(\Sigma_-)\overset{\sim}\to \underline{\Maps}_{\msc{D}}(C^{\ot g},m_f-)
\]
of functors from the fiber $\msc{D}^{\mcl{I}}\cong (\Bord^{nc}_{\msc{D}})_{\Sigma_{\zeta}}$ to $\Vect_{\pf}$.
\end{proposition}

Before proving Proposition \ref{prop:der_states} we establish a helpful lemma.

\begin{lemma}\label{lem:3386}
Let $A^{\heartsuit}$ be a finite ribbon tensor category, and consider the subcategory $fr\Disk^{\ot}_{surj}$ in $fr\Disk^{\ot}$ spanned by all objects and all morphisms $f:\mcl{I}\to \mcl{J}$ whose corresponding disk embedding $f:D\times\mcl{I}_f\to D\times\mcl{J}$ is surjective on connected components. The inclusion $\opn{Inj}\msc{K}^{\ot}\to \msc{K}^{\ot}$ pulls back to a map of cocartesian fibrations
\begin{equation}\label{eq:3388}
\xymatrixcolsep{2mm}
\xymatrix{
	fr\Disk^{\ot}_{surj}\times_{fr\Disk^{\ot}}\opn{Inj}\msc{K}^{\ot}\ar[rr]\ar[dr] & & fr\Disk^{\ot}_{surj}\times_{fr\Disk^{\ot}}\msc{K}^{\ot}\ar[dl]\\
		& fr\Disk^{\ot}_{surj} & .
}
\end{equation}
\end{lemma}

\begin{proof}
For any map $f:(\mcl{I},\mu)\to (\mcl{J},\nu)$ in $fr\Disk^{\ot}_{surj}$, Lemma \ref{lem:structure_transp} tells us that the corresponding transport functor
\[
f_!:\opn{Inj}\msc{K}^{\ot}_{(\mcl{I},\mu)}\to \opn{Inj}\msc{K}^{\ot}_{(\mcl{J},\nu)}
\]
is expressible as a composite of duality functors and product functors. Since the subcategories $\opn{Inj}\msc{K}^{\pm}$ are identified under duality, by construction, we see that the subcategory $\opn{Inj}\msc{K}^{\ot}\subseteq \msc{K}^{\ot}$ is stable under transport if and only if objects in $\opn{Inj}\msc{K}^{\pm}$ are stable under the respective product functors on $\msc{K}^{\pm}$. Since $\opn{Inj}\msc{K}^{\ot}$ is full in $\msc{K}^{\ot}$ and stable under isomorphism, stability under transport is equivalent to the claim that any cocartesian edge $\alpha:x\to y$ from an object in $\opn{Inj}\msc{K}^{\ot}$ also has $y$ in $\opn{Inj}\msc{K}^{\ot}$.
\par

Considering $\opn{Inj}\msc{K}^{-}=(\opn{Proj}\msc{K})^{op}$, and the duality equivalence with $\opn{Inj}\msc{K}_{\pf}$, it suffices to show that a product $P\ot P'$ of two $K$-projectives remains $K$-projective. For this we observe the formula
\[
\Hom^{\ast}_A(P\ot P',y)=\Hom^{\ast}_A(P,\uHom(P',y))
\]
via inner-Homs for the action of $A=\opn{Ch}(A^{\heartsuit})$ on itself, and note that the functor $\uHom(P',-)$ sends acyclic complexes to acyclic complexes. For this latter point one can consider a finite projective generator $G$ in $A^{\heartsuit}$ and note the formula
\[
\Hom_A^{\ast}(G,\uHom(P',-))=\Hom_A^{\ast}(G\ot P',-)=\Hom_A^{\ast}(P',G^{\ast}\ot-).
\]
The left hand functor preserves acyclics if and only if $\uHom(P',-)$ preserved acyclics, and the right hand functor preserves acyclics via exactness of the action functor $G^{\ast}\ot-$ and $K$-projectivity of $P'$.
\end{proof}

We return to our analysis of the state spaces in derived LRT theory.

\begin{proof}[Proof of Proposition \ref{prop:der_states}]
By our construction via the right adjoint we have
\[
L_{\msc{D}}(\Sigma_-)=L_{\msc{K}}(\Sigma_{R-})\overset{\sim}\to \underline{\Maps}_{\msc{K}}\left(C^{\ot g},m_f(R-)\right),
\]
where the equivalence with inner-Homs comes from Proposition \ref{prop:htop_states}.  We also note that the map $f:D\times \mcl{I}\to D$ in this case is in the subcategory $fr\Disk_{surj}\subseteq fr\Disk$.
\par

By Corollary \ref{cor:3075} and Lemma \ref{lem:3386} the symmetric monoidal functor $R:\msc{D}^{\ot}\to \msc{K}^{\ot}$ restricts to a map of cocartesian fibrations
\[
fr\Disk^{\ot}_{surj}\times_{fr\Disk^{\ot}}\msc{D}^{\ot}\to fr\Disk^{\ot}_{surj}\times_{fr\Disk^{\ot}}\msc{K}^{\ot}
\]
over $fr\Disk^{\ot}_{surj}$, so that there is a natural isomorphism $m_fR\overset{\sim}\to Rm_f$.  This gives an equivalence
\[
L_{\msc{D}}(\Sigma_-)\overset{\sim}\to \underline{\Maps}_{\msc{K}}\left(C^{\ot g},Rm_f-\right).
\]
By the calculation of inner-Homs for the derived $\infty$-category provided in Proposition \ref{prop:D_innerhom_formula}, we have $\underline{\Maps}_{\msc{D}}(loc-,-)\cong \underline{\Maps}_{\msc{K}}(-,R-)$, so that finally
\[
L_{\msc{D}}(\Sigma_-)\overset{\sim}\to \underline{\Maps}_{\msc{D}}(C^{\ot g},m_f-).
\]
\end{proof}

\subsection{Recovering discrete LRT at the heart}

Consider the full subcategory $(\msc{D}^{\heartsuit}_{fin})^{\ot}$ in $\msc{D}_{fin}^{\ot}$ spanned by complexes with cohomology concentrated in degree $0$.  This is a $fr\Disk$-category by Proposition \ref{prop:frdisk_subcats}, the inclusion $(\msc{D}^{\heartsuit}_{fin})^{\ot}\to \msc{D}^{\ot}$ is a $fr\Disk$-monoidal functor, and the sequence
\[
(A^{\heartsuit})^{\ot}\to A^{\ot}\to \msc{D}^{\ot}
\]
induces an equivalence $(A^{\heartsuit})^{\ot}\overset{\sim}\to (\msc{D}_{fin}^{\heartsuit})^{\ot}$ \cite[\href{https://kerodon.net/tag/028B}{0288}]{kerodon}. Hence we obtain a symmetric monoidal equivalence $\Bord_{A^{\heartsuit}}^{nc}\overset{\sim}\to \Bord_{\msc{D}_{fin}^{\heartsuit}}^{nc}$. We therefore have two field theories
\[
\opn{L}_{A^{\heartsuit}}=\opn{L}_A|_{\Bord_{A^{\heartsuit}}^{nc}}\ \ \text{and}\ \ L_{\msc{D}_{fin}^{\heartsuit}}=L_{\msc{D}_{fin}}|_{\Bord_{\msc{D}_{fin}^{\heartsuit}}^{nc}},
\]
the first of which is valued in the full symmetric monoidal subcategory $\opn{Vect}\subseteq \opn{Ch}(\opn{Vect})$ and, via the identification of the state spaces in Proposition \ref{prop:der_states}, the second field theory is (up to natural isomorphism) valued in the full subcategory $\Vect^{\geq 0}$ of coconnective cochains.
\par

Now, the $0$-th cohomology functor $H^0:\Vect^{\geq 0}\to \opn{Vect}$ admits a natural symmetric monoidal structure, which one can realizes via localization for example.  We therefore have the $\opn{Vect}$-valued symmetric monoidal functor
\[
H^0L_{\msc{D}_{fin}^{\heartsuit}}:\Bord_{\msc{D}_{fin}^{\heartsuit}}^{nc}\to \opn{Vect}.
\]
We compare this cohomological LRT theory to the original Lyubashenko theory from \cite{derenzietal23}.  As a preliminary point, however, we record the following.

\begin{lemma}\label{lem:sym_transform}
A transformation
\[
\xymatrix{
\Delta^1\times\msc{G}_0^{\ot}\ar[rr]^{\zeta}\ar[dr]_{q_0p_2} & & \msc{G}_1^{\ot}\ar[dl]^{q_1}\\
	& \Fin_{\ast}
}
\]
between symmetric monoidal functors $F_0$ and $F_1$ is a natural isomorphism if and only if, for each $x$ in the underlying $\infty$-category $\msc{G}_0$, $\zeta$ evaluates to an isomorphism $\zeta_x:F_0(x)\to F_1(x)$ in $\msc{G}_1$.
\end{lemma}

\begin{proof}
As in the statement, we take $F_0=\zeta|_{\{0\}}$ and $F_1=\zeta|_{\{1\}}$.  For an object $x$ in a general fiber $(\msc{G}_0^{\ot})_I$ we have that $\zeta_x$ is an isomorphism if and only if the induced map
\[
(\zeta_x)_{\ast}:\pi_0\Maps_{\msc{G}_1^{\ot}}(y,F_0x)\to \pi_0\Maps_{\msc{G}^{\ot}_1}(y,F_1x)
\]
is an isomorphism at all $y$.  For the inert projections $\bar{\rho}_i:I\to \{0\}$ we have the cocartesian lifts $\rho_i:x\to x_i$ in $\msc{G}_0^{\ot}$ which map to cocartesian lifts $F_{\varepsilon}(\rho_i):F_{\varepsilon}x\to F_{\varepsilon}x_i$, since the $F_{\varepsilon}$ are symmetric monoidal functors. For each $i$ in $I$, we also have the diagrams
\[
\xymatrix{
F_{0}x\ar[r]^{\zeta_x}\ar[d]_{F_0(\rho_i)} & F_{1}x\ar[d]^{F_1(\rho_i)} \\
F_{0}x_i\ar[r]_{\zeta_{x_i}} & F_1x_i.
}
\]
\par

By the factorization property on morphisms for $\infty$-operads \cite[Definition 2.1.1.10]{ha} we see that the map $\zeta_x$ is an isomorphism if and only if all $\zeta_{x_i}$ are isomorphisms.  In particular, if we assume that $\zeta_{\{0\}}$ is a natural isomorphism then $\zeta_x$ is an isomorphism at all $x$ in $\msc{G}_0^{\ot}$.  This is sufficient to show that $\zeta$ itself is an isomorphism in the functor category $\Fun(\msc{G}_0^{\ot},\msc{G}_1^{\ot})$ \cite[\href{https://kerodon.net/tag/01DK}{01DK}]{kerodon}. We now seek an inverse to $\zeta$ in the relative functor category $\Fun_{\Fin_{\ast}}(\msc{G}_0^{\ot},\msc{G}_1^{\ot})$.
\par

The inverse to $\zeta$ in the non-relative functor category is uniquely determined, up to isomorphism, as the transformation $\zeta':F_1\to F_0$ which completes a diagram
\begin{equation}\label{eq:3485}
\xymatrix{
	& F_1\ar@{-->}[dr]^{\zeta'}\\
F_0\ar[rr]_{id}\ar[ur]^{\zeta} & & F_0.	
}
\end{equation}
Since the induced map $(q_1)_{\ast}:\Fun(\msc{G}_0^{\ot},\msc{G}_1^{\ot})\to \Fun(\msc{G}_0^{\ot},\Fin_{\ast})$ is an inner fibration \cite[\href{https://kerodon.net/tag/01H0}{01H0}]{kerodon}, and any isomorphism is therefore both a $(q_1)_{\ast}$-cocartesian and cartesian edge, we can produce a candidate inverse $\zeta'$ to $\zeta$ in $\Fun_{\Fin_{\ast}}(\msc{G}^{\ot}_0,\msc{G}^{\ot}_1)$ by solving the lifting diagram
\begin{equation}\label{eq:3492}
\xymatrix{
\Lambda^2_0\ar[r]\ar[d] & \Fun(\msc{G}_0^{\ot},\msc{G}_1^{\ot})\ar[d]\\
\Delta^2\ar[r]_(.3){q_0}\ar@{..>}[ur] & \Fun(\msc{G}_0^{\ot},\Fin_{\ast})
}
\end{equation}
in which the bottom map is constant of value $q_0$. This solution provides a $2$-simplex in the subcategory $\Fun_{\Fin_{\ast}}^{cc}(\msc{G}_{0}^{\ot},\msc{G}_1^{\ot})$ which realizes $\zeta'$ as a left inverse to $\zeta$, as a symmetric monoidal transformation, and one subsequently completes a diagram of the form
\[
\xymatrixrowsep{5mm}
\xymatrixcolsep{4mm}
\xymatrix{
F_0\ar[rr]^{\zeta}\ar[ddr]_{id}\ar[drrr]|(.57){\hole} & & F_1\ar[ddl]_{\zeta'}\ar[dr]^{id} \\
& & & F_1\\
 & F_0\ar[urr]_{\zeta}
}
\]
by solving a lifting problem as in \eqref{eq:3492} for the horn inclusion $\Lambda^3_0\to \Delta^3$ to see that $\zeta'$ is a $2$-sided inverse to $\zeta$ in $\Fun_{\Fin_{\ast}}^{cc}(\msc{G}_0^{\ot},\msc{G}_1^{\ot})$.
\end{proof}

\begin{proposition}\label{prop:heart_lrt}
For any finite modular tensor category $A^{\heartsuit}$, there is a diagram
\[
\xymatrix{
\Bord_{A^{\heartsuit}}^{nc}\ar[rr]^{L_{A^{\heartsuit}}}\ar[dr]_{\sim} & & \opn{Vect}\\
	& \Bord_{\msc{D}^{\heartsuit}_{fin}}\ar[ur]_{H^0L_{\msc{D}_{fin}^{\heartsuit}}}
}
\]
in $\SM_{\infty}$.  In particular, there is a natural isomorphism of symmetric monoidal functors $L_{A^{\heartsuit}}\overset{\sim}\to H^0L_{\msc{D}_{fin}^{\heartsuit}}|_{\Bord_{A^{\heartsuit}}}$.
\end{proposition}

\begin{proof}
We have the symmetric monoidal transformation
\[
\alpha^{\ot}:L_{A^{\heartsuit}}^{\ot}\cong L_{\msc{K}_{fin}}^{\ot}|_{(\Bord_{A^{\heartsuit}}^{nc})^{\ot}}\to L_{\msc{D}_{fin}^{\heartsuit}}^{\ot}|_{(\Bord_{A^{\heartsuit}}^{nc})^{\ot}}
\]
provided by our construction of $L_{\msc{D}}^{\ot}$ via Kan extension.  Taking cohomology, and noting that the sequence $\opn{Vect}^{\ot}\to \Vect^{\ot}\overset{H^0}\to \opn{Vect}^{\ot}$ is the identity, we obtain another natural transformation
\[
H^0\alpha^{\ot}:L_{A^{\heartsuit}}^{\ot}\to H^0L_{\msc{D}_{fin}^{\heartsuit}}^{\ot}|_{(\Bord_{A^{\heartsuit}}^{nc})^{\ot}}.
\]
We claim that $H^0\alpha^{\ot}$ is a natural isomorphism.
\par

To test if $H^0\alpha^{\ot}$ is a natural isomorphism, it suffices to check that the underlying transformation
\[
H^0\alpha:L_{A^{\heartsuit}}\to H^0L_{\msc{D}_{fin}^{\heartsuit}}|_{\Bord_{A^{\heartsuit}}^{nc}}
\]
evaluates to an isomorphism at each object $\Sigma_x$ in the underlying $\infty$-category $\Bord_{A^{\heartsuit}}^{nc}$, by Lemma \ref{lem:sym_transform}.  Furthermore, since $H^0\alpha$ is symmetric monoidal we reduce to the case where $\Sigma$ is connected and, since merging bordisms induce isomorphisms on state spaces, we reduce further to the case of a surface with a single positive marking.  In this case, for such singly marked $\Sigma_{\zeta}$ in $\Bord_{\ast}^{nc}$, we have the localization functor
\[
l:\msc{K}_{\pf}=(\Bord^{nc}_{\msc{K}})_{\Sigma_{\zeta}}\to \msc{D}_{\pf}=(\Bord^{nc}_{\msc{D}})_{\Sigma_{\zeta}}
\]
with right adjoint $R$, and the transformation
\[
\alpha:L_{A^{\heartsuit}}(\Sigma_x)=L_{\msc{K}}(\Sigma_x)\to L_{\msc{K}}(\Sigma_{Rl(x)})=L_{\msc{D}}(\Sigma_{l(x)})
\]
is induced by the unit map $u_x:x\to Rl(x)$ for the adjunction between $R$ and $l$.  In particular, $\alpha$ is obtained by applying the functor $L_{\msc{K}}$ to the bordism $\Sigma_{u_x}:\Sigma_x\to \Sigma_{Rl(x)}$ over $id_{\Sigma_{\zeta}}$ with the outgoing disk labeled by the morphism $u_x$.
\par

As relayed in Corollary \ref{cor:2953}, this object $Rl(x)$ is $K$-injective and $u_x:x\to Rl(x)$ is a quasi-isomorphism.  So, after applying the identification of state spaces from Proposition \ref{prop:der_states}, we have
\[
\alpha_{\Sigma_x}=-\circ u_x:\Hom_{A^{\heartsuit}}(C^{\ot g},x)\to \Hom^{\ast}_{A^{\heartsuit}}(C^{\ot g},Rl(x))=\underline{\Maps}_{\msc{K}}(C^{\ot g},Rl(x)).
\]
(Here we adopt the explicit construction of the homotopy $\infty$-category as the dg nerve of the dg category of linear cochains.) Taking $0$-th cohomology we see that $H^0\alpha_{\Sigma_x}$ recovers the standard identification
\[
\Hom_{A^{\heartsuit}}(C^{\ot g},x)\overset{\sim}\to \Ext^0_{A^{\heartsuit}}(C^{\ot g},x),
\]
so that all $H^0\alpha_{\Sigma_x}$ are isomorphisms and hence $H^0\alpha^{\ot}$ is a natural isomorphism of symmetric monoidal functors.
\end{proof}

\section{Derived theories with universal skeins}
\label{sect:derived_w_skeins}

We aim to produce, in these final sections, an \emph{unmarked} field theory
\begin{equation}\label{eq:4074}
\mbb{L}_{D_{fin}}:\Bord^{nc}_{3,2}\to \opn{h}\Vect
\end{equation}
from the marked theory $L_{\msc{D}_{fin}}$ constructed in Section \ref{sect:derived_lrt}. The theory $\mbb{L}_{D_{fin}}$ is obtained in two basic movements: We first localize the bordism category $\Bord_{\msc{D}_{fin}}^{nc}$ along a special class of colored merging bordisms, and produce in this way a symmetric monoidal $\infty$-category $\Bord_{\msc{D}_{fin}}^{nc}[\Theta_{\msc{D}_{fin}}^{-1}]$ of bordisms with embedded skein-like diagrams. We also argue that the derived LRT theory $L_{\msc{D}_{fin}}$ sends all such colored merging bordisms to isomorphisms in $\Vect_{\pf}$, and therefore induces a symmetric monoidal functor from this localization.  We clarify these points in the present section.  
\par

Subsequently, in Section \ref{sect:unmarked_lrt} we establish an identification between the homotopy category of localized, trivially marked bordisms $\opn{h}\Bord_{\ast}^{nc}[\Theta^{-1}_{\ast}]$ and the unmarked category $\Bord^{nc}_{3,2}$, thus yielding the promised theory \eqref{eq:4074} after restricting along a unit map $\Bord^{nc}_{\ast}\to \Bord^{nc}_{\msc{D}_{fin}}$.

\subsection{Special merging bordisms}

For any surface $\Sigma$ with markings $\zeta:D\times \mcl{I}\to \Sigma$, and index $i$ in $\mcl{I}$, we have the corresponding surface $\Sigma^i_{\zeta}$ obtained by deleting all but the $i$-th marking disk. We also consider the apparent inclusion $M\Sigma^i_{\zeta}:\Sigma^i_{\zeta}\to \Sigma_{\zeta}$ given by the $3$-manifold $\Sigma\times[0,1]$ with the embedded cylinder $D\times [0,1]\to \Sigma\times [0,1]$, $(z,t)\mapsto (\zeta(z,i),t)$.  The reader might recall at this point the notion of a merging bordism from Section \ref{sect:merge}.

\begin{definition}\label{def:special_merge}
For a connected marked surface $\zeta:D\times \mcl{I}\to \Sigma$ and positive index $i$ in $\mcl{I}$, a merging bordism $M_f:\Sigma_{\zeta}\to \Sigma^i_{\zeta}$ in $\Bord_{\ast}^{nc}$ is called special if the identity on $\Sigma^i_{\zeta}$ is recovered as a composite $id_{\Sigma^i_{\zeta}}=M_f\circ M\Sigma^i_{\zeta}$.  For a general marked surface $\Sigma_{\zeta}$, a special merging bordism $M_f:\Sigma_{\zeta}\to \Sigma_{\eta}$ is a bordism whose components are each special merging bordisms, or identities.

For a $fr\Disk$-category $q:\msc{E}^{\ot}\to fr\Disk^{\ot}$ with corresponding cocartesian fibration $q':\Bord_{\msc{E}}^{nc}\to \Bord_{\ast}^{nc}$, a special merging bordism $M_{\alpha}:\Sigma_x\to \Sigma_y$ in $\Bord_{\msc{E}}^{nc}$ is a $q'$-cocartesian edge over a special merging bordism in $\Bord_{\ast}^{nc}$.
\end{definition}

\begin{lemma}
Let $\msc{E}^{\ot}\to fr\Disk^{\ot}$ be a $fr\Disk$-category. Consider any object $\Sigma_x$ in $\Bord^{nc}_{\msc{E}}$ and positive index $i$ for the underling $\ast$-marked surface $\opn{top}(x):D\times \mcl{I}\to \Sigma$.  Let $\Sigma^c\subseteq \Sigma$ be the component of $\Sigma$ which contains the image of this $i$-th marking disk. There is a special merging bordism $M_{\alpha}:\Sigma_x\to \Sigma_y$ in which the outgoing surface $\Sigma_y$ is endowed with the single marking $\opn{top}(x)_i:D\to \Sigma^c$ on the specified component $\Sigma^c$ of $\Sigma$. 
\end{lemma}

\begin{proof}
Since the forgetful functor $\Bord_{\msc{E}}^{nc}\to \Bord_{\ast}^{nc}$ is a cocartesian fibration it suffices to prove the result in the $\ast$-marked bordism category.  We can assume further that the surface $\Sigma$ is connected. Consider $\Sigma_{\zeta}$ in $\Bord_{\ast}^{nc}$ and a positive index $i$. We take an embedding from an auxiliary disk $\tau:D\to \Sigma$ which contain all of the markings $\zeta:D\times\mcl{I}\to \Sigma$ in its image and view $\zeta$ as a map into this new disk $\zeta:D\times\mcl{I}\to D$, by intersecting with the image of $\tau$.

Let $N$ be large and consider any smooth function $c:[0,1]\to \mathbb{R}_{\geq 0}$ which is constant around $0$ and $1$, is weakly decreasing, and has values $c(0)=1$ and $c(1)=N^{-1}$.  We can take a smooth homotopy $h:[0,1]\times D\times \mcl{I}\to D$ which satisfies the following:
\begin{itemize}
\item[(a)] $h_t:D\times\mcl{I}\to D$ is a smooth embedding at all $t$.\vspace{1mm}
\item[(a)] $h_0=\zeta$ and $h(t,z,i)=\zeta(c(t)\cdot z,i)$, where $i$ is the specified index.\vspace{1mm}
\item[(b)] There exists $\epsilon>0$ for which $h(t,z,j)=h(t',z,j)$ whenever $t,t'<\epsilon$ or $t,t'>1-\epsilon$.\vspace{1mm}
\item[(c)] $h_1:D\times\mcl{I}\to D$ has image in $\zeta(D_i)$.
\end{itemize}
The corresponding bordism $M_f:\Sigma_{\zeta}\to \Sigma^i_{\zeta}$ with $M=\Sigma\times [0,1]$ and cylinders
\[
f:D\times[0,1]\times \mcl{I}\to \Sigma\times [0,1],\ \ f(z,t,i):=(h(t,z,i),t)
\]
is a special merging bordism.  (See also Proposition \ref{prop:1456}.)
\end{proof}

We have the following  analog of Lemma \ref{lem:1288}.

\begin{proposition}\label{prop:der_merge}
Let $A^{\heartsuit}$ be a modular tensor category and $\msc{D}^{\ot}\to fr\Disk^{\ot}$ be the associated $fr\Disk$-category for the derived $\infty$-category.  Any special merging bordism $M_{\alpha}:\Sigma_x\to \Sigma_y$ in $\Bord_{\msc{D}}^{nc}$ evaluates to an isomorphism
\[
L_{\msc{D}}(M_{\alpha}):L_{\msc{D}}(\Sigma_x)\overset{\sim}\to L_{\msc{D}}(\Sigma_y)
\]
in $\Vect_{\pf}$.
\end{proposition}

\begin{proof}
We can lift this special merging bordism in $\Bord_{\msc{D}}^{nc}$ along the equivalence $\Bord_{\opn{Inj}\msc{K}}^{nc}\overset{\sim}\to \Bord_{\msc{D}}^{nc}$ to a special merging bordism $M'_{\alpha'}:\Sigma_{x'}\to \Sigma_{y'}$ in $\Bord_{\opn{Inj}\msc{K}}$, at least up to an isomorphism of morphisms.  From the natural isomorphism
\[
L_{\msc{D}}|_{\Bord_{\opn{Inj}\msc{K}}^{nc}}\overset{\sim}\to L_{\opn{Inj}\msc{K}}
\]
it now suffices to show that $L_{\msc{K}}$ sends special merging bordisms to isomorphisms.  Now, since the localization functor $\Bord_{A}^{nc}\to \Bord_{\msc{K}}^{nc}$ is an essentially surjective map of cocartesian fibrations over $\Bord_{\ast}^{nc}$, we can furthermore lift special merging bordisms in $\Bord_{\msc{K}}^{nc}$ to $\Bord_{A}^{nc}$, again up to isomorphism.  We thus reduce further to proving that the functor
\[
L_{A}:\Bord_{A}^{nc}\to \Vect
\]
sends special merging bordisms to isomorphisms, and via symmetric monoidality we reduce further to considering such bordisms $M:\Sigma_x\to \Sigma_y$ for connected $\Sigma$.
\par

In the setting of $\Bord_{A}^{nc}$ we can be very explicit about things.  An object $\Sigma_x$ in $\Bord_{A}^{nc}$ is a surface $D\times \mcl{I}\to \Sigma$ with complexes $x_i$ labeling each disk $D_i$.  For each positive index, $x_i$ is a complex in $\opn{Ch}(A^{\heartsuit}_{\pf})$, and for each negative index $x_i$ is a complex in the opposite category $\opn{Ch}(A^{\heartsuit})^{op}$.  A special merging bordism $M:\Sigma_x\to \Sigma_y$ with $\Sigma$ connected is either an isomorphism or some non-trivial merging of markings on $\Sigma$.  In the former case it is immediate that $L_{A}(M)$ is an isomorphism.  So we consider the latter setting.
\par

Such a non-trivial merging bordism $M_{\alpha}:\Sigma_x\to \Sigma_y$ has the unique outgoing disk labeled by an isomorphism $\alpha_A:m_f(x)\to y$. We can write the tuple $x$ as a limit $x=\varprojlim_{\lambda} x_{\lambda}$ of objects in $A_{fin}^{\mcl{I}}$ giving the product as a limit $m_f(x)=\varprojlim_{\lambda} m_f(x_{\lambda})$, via continuity of $m_f$. Via the isomorphism $\alpha_A$ we have a corresponding expression $y=\varprojlim_{\lambda}y_{\lambda}$ with $\alpha_A$ inducing isomorphisms $\alpha_{A,\lambda}:m_f(x_{\lambda})\to y_{\lambda}$.  We therefore have corresponding merging bordisms $M_{\alpha,\lambda}:\Sigma_{x_{\lambda}}\to \Sigma_{y_{\lambda}}$ at each $\lambda$ which fit into a diagram
\[
\xymatrix{
L_{A}(\Sigma_x)\ar[rr]^{L(M_{\alpha})}\ar[d] & & L_A(\Sigma_y)\ar[d]\\
L_{A}(\Sigma_{x_{\lambda}})\ar[rr]_{L(M_{\alpha,\lambda})} & & L_A(\Sigma_{y_{\lambda}}).
}
\]

By Lemma \ref{lem:1288} the bottom maps $L_A(M_{\alpha,\lambda})$ are all isomorphisms, and by Proposition \ref{prop:dis_lrt} the vertical maps assemble into an isomorphisms
\[
L_{A}(\Sigma_x)=\varprojlim_{\lambda}L_{A}(\Sigma_{x_{\lambda}})\ \ \text{and}\ \ L_A(\Sigma_y)=\varprojlim_{\lambda}L_A(\Sigma_{y_{\lambda}}).
\]
From this we conclude that $L_A(M_{\alpha})$ is an isomorphism.
\end{proof}

\subsection{Color change bordisms}

\begin{definition}
In $\Bord_{\ast}^{nc}$ a color change bordism $M_f:\Sigma_{\zeta}\to \Sigma_{\eta}$ is one in which the following hold:
\begin{itemize}
\item The marking set $\mcl{I}'$ for $\Sigma_{\eta}$ is obtained from the marking set $\mcl{I}$ for $\Sigma_{\zeta}$ by changing some negatively colored indices in $\mcl{I}$ to positive indices.\vspace{1mm}
\item After forgetting the colors, the markings $\bar{\zeta},\bar{\eta}:D\times I\to\Sigma$ agree.\vspace{1mm}
\item The bordism $M_f$ is just $\Sigma\times [0,1]$ along with the constant cylinder embedding
\[
f:D\times[0,1]\times \mcl{I}\to \Sigma\times [0,1],\ \ f(z,t,i)=(\zeta(z,i),t).
\]
\end{itemize}
For $\msc{E}^{\ot}\to fr\Disk^{\ot}$ any $fr\Disk$-monoidal $\infty$-category, a color change bordism in $\Bord_{\msc{E}}^{nc}$ is a cocartesian lift $M_{\alpha}:\Sigma_x\to \Sigma_y$ of a color change bordism in $\Bord^{nc}_{\ast}$.
\end{definition}

\begin{lemma}\label{lem:color_change}
Let $A^{\heartsuit}$ be a finite modular tensor category.  Any color change bordism $M_{\alpha}:\Sigma_x\to \Sigma_y$ in $\Bord_{\msc{D}}^{nc}$ evaluates to an isomorphism
\[
L_{\msc{D}}(M_{\alpha}):L_{\msc{D}}(\Sigma_x)\overset{\sim}\to L_{\msc{D}}(\Sigma_y)
\]
in $\Vect_{\pf}$.
\end{lemma}

\begin{proof}
As in the proof of Proposition \ref{prop:der_merge}, it suffices to show that the map
\[
L_{A}(M_{\alpha}):L_{A}(\Sigma_x)\to L_{A}(\Sigma_y)
\]
is an isomorphism whenever $M_{\alpha}:\Sigma_x\to \Sigma_y$ is a color change bordism in $\Bord_{A}^{nc}$.  Such a bordism has outgoing boundary labeled by isomorphisms $\alpha_i:x_i\overset{\sim}\to y_i$ when the color of the $i$-th makings in $\Sigma_x$ and $\Sigma_y$ are the same, and is labeled by an isomorphism $\alpha_i:x_i^{\ast}\overset{\sim}\to y_i$ when the colors of the $i$-th markings differ.  Via the limit expressions
\[
L_{A}(\Sigma_x)=\varprojlim_{\lambda}L_{A}(\Sigma_{x_{\lambda}})\ \ \text{and}\ \ L_{A}(\Sigma_y)=\varprojlim_{\lambda} L_{A}(\Sigma_{y_{\lambda}})
\]
we reduce to the case where $x$ and $y$ are finite.
\par

We recall, in the finite setting, that the $\opn{Ch}(\opn{Vect})$-valued TQFT $\opn{L}_{A_{fin}}$ is obtained by applying a symmetric map $p:\Bord_{A_{fin}}^{nc}\to \msf{Bord}_{A_{fin}}^{nc}$ to the ``reduced" bordism category from \cite[Section 8.1]{czenkynegron}, then composing with the symmetric monoidal functor $Z^{\ast}:\msf{Bord}_{A_{fin}}^{nc}\to \opn{Ch}(\opn{Vect})$ from \cite[Theorem 9.3]{czenkynegron}.  The image $p(M_{\alpha})$ in already an isomorphism in $\msf{Bord}_{A_{fin}}^{nc}$ so that
\[
\opn{L}_{A_{fin}}(M_{\alpha})=Z^{\ast}\big(p(M_{\alpha})\big)
\]
is an isomorphism in $\opn{Ch}(\opn{Vect})$.  Consequently, after applying the localization $loc:\opn{Ch}(\opn{Vect})\to \Vect$, we have that $L_{A_{fin}}(M_{\alpha})=loc \opn{L}_{A_{fin}}(M_{\alpha})$ is an isomorphism in $\Vect$.
\end{proof}

We are interested in a class of maps which mixes both special merging bordisms and color change bordisms in $\Bord_{\msc{E}}^{nc}$.

\begin{definition}
A colored merging bordism in $\Bord_{\ast}^{nc}$ is a bordism $M_f:\Sigma_{\zeta}\to \Sigma_{\eta}$ which decomposes as a disjoint union of special merging bordisms and color change bordisms.  For any $fr\Disk$-category $\msc{E}^{\ot}\to fr\Disk^{\ot}$, a colored merging bordism in $\Bord_{\msc{E}}^{nc}$ is a cocartesian lift $M_{\alpha}:\Sigma_x\to \Sigma_y$ of a colored merging bordism in $\Bord_{\ast}^{nc}$ along the fibration $\Bord^{nc}_{\msc{E}}\to \Bord^{nc}_{\ast}$.
\end{definition}

\subsection{Inverting colored merging bordisms}

\begin{definition}
For any $fr\Disk$-monoidal $\infty$-category $\msc{E}^{\ot}\to fr\Disk^{\ot}$, let $\Theta_{\msc{E}}\subseteq \Bord^{nc}_{\msc{E}}[1]$ be the collection of colored merging bordisms.
\end{definition}

\begin{theorem}\label{thm:merge_loc}
Let $\msc{E}^{\ot}\to fr\Disk^{\ot}$ be a $fr\Disk$-category. There is a unique symmetric monoidal structure $\Bord_{\msc{E}}^{nc}[\Theta^{-1}_{\msc{E}}]^{\ot}\to \Fin_{\ast}$ on the localization $\Bord_{\msc{E}}^{nc}[\Theta^{-1}_{\msc{E}}]$ for which the following hold:
\begin{enumerate}
\item The localization functor $\Bord_{\msc{E}}^{nc}\to \Bord_{\msc{E}}^{nc}[\Theta^{-1}_{\msc{E}}]$ lifts to a symmetric monoidal functor $loc:(\Bord_{\msc{E}}^{nc})^{\ot}\to \Bord_{\msc{E}}^{nc}[\Theta^{-1}_{\msc{E}}]^{\ot}$.\vspace{1mm}
\item For any symmetric monoidal $\infty$-category $\msc{G}^{\ot}\to \Fin_{\ast}$, restriction provides a fully faithful functor
\[
loc^{\ast}:\Fun_{\Fin_{\ast}}^{cc}(\Bord_{\msc{E}}^{nc}[\Theta^{-1}_{\msc{E}}]^{\ot},\msc{G}^{\ot})\to \Fun_{\Fin_{\ast}}^{cc}((\Bord_{\msc{E}}^{nc})^{\ot},\msc{G}^{\ot})
\]
which is an equivalence onto the full subcategory spanned by those symmetric monoidal functors $F^{\ot}$ for which the underlying map $F:\Bord_{\msc{E}}^{nc}\to \msc{G}$ sends each morphism in $\Theta_{\msc{E}}$ to an isomorphism in $\msc{G}$.
\end{enumerate}
\end{theorem}

Ultimately, we obtain the symmetric monoidal $\infty$-category $\Bord_{\msc{E}}^{nc}[\Theta^{-1}_{\msc{E}}]^{\ot}$ via a vertical localization procedure (see Lemma \ref{lem:adv_loc}). We give the proof of Theorem \ref{thm:merge_loc} in Section \ref{sect:proof_merge_bord} below.  For now, let us record an immediate consequence.

\begin{corollary}\label{cor:der_lrt_skeins}
For a finite modular tensor category $A^{\heartsuit}$, with corresponding derived $\infty$-category $\msc{D}$, the derived LRT theory $L_{\msc{D}}:\Bord_{\msc{D}}^{nc}\to \Vect$ localizes to a symmetric monoidal functor
\[
L^{loc}_{\msc{D}}:\Bord_{\msc{D}}^{nc}[\Theta^{-1}_{\msc{D}}]\to \Vect_{\pf}.
\]
\end{corollary}

\begin{proof}
By Proposition \ref{prop:der_merge} and Lemma \ref{lem:color_change} the functor $L_{\msc{D}}$ sends the class $\Theta_{\msc{D}}$ into the class of isomorphisms in $\Vect_{\pf}$.  So the result follows by Theorem \ref{thm:merge_loc}.
\end{proof}

\subsection{Describing the localization via ``coarse skeins"}
\label{sect:univ_sk}

Let $\msc{E}^{\ot}\to fr\Disk^{\ot}$ be an arbitrary $fr\Disk$-category.  We note that the forgetful functor $\Bord_{\msc{E}}^{nc}\to \Bord^{nc}_{3,2}$ to the unlabeled (anomalous) bordism category sends all colored merging bordisms to the identity, and thus localizes to provide a symmetric monoidal functor $\Bord_{\msc{E}}^{nc}[\Theta^{-1}_{\msc{E}}]\to \Bord_{3,2}^{nc}$.  So the localization is some kind of decorated bordism category.
\par

The objects in the localization are still marked surfaces--or at least, they can be taken to be marked surfaces--but morphisms are now bordisms $M:\Sigma\to \Sigma'$ with a kind of ``coarse wiring diagram" which is comprised of cylinders merging and splitting in $M$. 
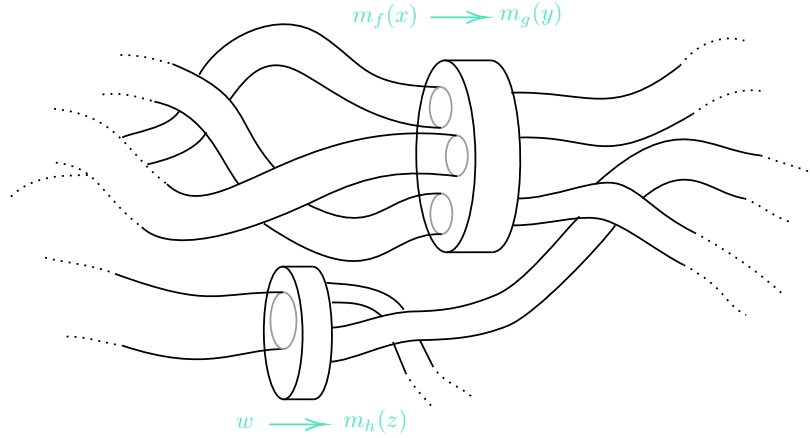
\[
\scalebox{.9}{
\tikzset{every picture/.style={line width=0.75pt}} 
\begin{tikzpicture}[x=0.75pt,y=0.75pt,yscale=-1,xscale=1]
\draw   (331.84,60.01) .. controls (340.95,59.65) and (348.62,83.3) .. (348.98,112.85) .. controls (349.34,142.39) and (342.25,166.63) .. (333.15,167) .. controls (324.04,167.36) and (316.37,143.71) .. (316.01,114.16) .. controls (315.65,84.62) and (322.73,60.38) .. (331.84,60.01) -- cycle ;
\draw    (331.84,60.01) -- (357.24,60) ;
\draw    (333.15,167) -- (358.56,166.98) ;
\draw    (357.24,60) .. controls (378.52,59.99) and (380.17,165.51) .. (358.56,166.98) ;
\draw   (241.37,175.01) .. controls (247.34,174.76) and (252.37,191.12) .. (252.61,211.55) .. controls (252.85,231.98) and (248.21,248.75) .. (242.24,249) .. controls (236.27,249.25) and (231.24,232.89) .. (231,212.46) .. controls (230.76,192.03) and (235.4,175.26) .. (241.37,175.01) -- cycle ;
\draw    (241.37,175.01) -- (258.02,175) ;
\draw    (242.24,249) -- (258.88,248.99) ;
\draw    (258.02,175) .. controls (271.96,174.99) and (273.04,247.97) .. (258.88,248.99) ;
\draw    (195,68) .. controls (206,48) and (239,30) .. (263,46) .. controls (287,62) and (286.08,74.16) .. (329.49,74.91) ;
\draw  [color={rgb, 255:red, 155; green, 155; blue, 155 }  ,draw opacity=1 ] (329.87,134.26) .. controls (333.26,134.12) and (336.06,139.08) .. (336.13,145.34) .. controls (336.21,151.6) and (333.52,156.8) .. (330.14,156.94) .. controls (326.75,157.08) and (323.95,152.12) .. (323.87,145.86) .. controls (323.8,139.59) and (326.49,134.4) .. (329.87,134.26) -- cycle ;
\draw  [color={rgb, 255:red, 155; green, 155; blue, 155 }  ,draw opacity=1 ] (338.49,101.91) .. controls (341.88,101.76) and (344.68,106.72) .. (344.75,112.99) .. controls (344.83,119.25) and (342.14,124.45) .. (338.76,124.59) .. controls (335.37,124.73) and (332.57,119.77) .. (332.49,113.5) .. controls (332.42,107.24) and (335.11,102.05) .. (338.49,101.91) -- cycle ;
\draw  [color={rgb, 255:red, 155; green, 155; blue, 155 }  ,draw opacity=1 ] (329.49,74.91) .. controls (332.88,74.76) and (335.68,79.72) .. (335.75,85.99) .. controls (335.83,92.25) and (333.14,97.45) .. (329.76,97.59) .. controls (326.37,97.73) and (323.57,92.77) .. (323.49,86.5) .. controls (323.42,80.24) and (326.11,75.05) .. (329.49,74.91) -- cycle ;
\draw    (177,129) .. controls (205,156) and (247,91) .. (338.49,101.91) ;
\draw    (163,153) .. controls (182.02,171.34) and (233,151) .. (262,137) .. controls (291,123) and (309.42,121.09) .. (338.76,124.59) ;
\draw  [dash pattern={on 0.84pt off 2.51pt}]  (114,87) .. controls (151,91) and (153,111) .. (177,129) ;
\draw  [dash pattern={on 0.84pt off 2.51pt}]  (115,117) .. controls (140,120) and (139,135) .. (163,153) ;
\draw    (212,82) .. controls (222,68) and (232,57) .. (250,66) .. controls (268,75) and (291.51,96.93) .. (329.76,97.59) ;
\draw    (255,140) .. controls (296,161) and (297,133) .. (329.87,134.26) ;
\draw    (231,151) .. controls (296,195) and (297.27,155.68) .. (330.14,156.94) ;
\draw    (180,62) .. controls (224,80) and (214,94) .. (237,120) ;
\draw    (170,83) .. controls (199,89) and (206,115) .. (217,129) ;
\draw  [dash pattern={on 0.84pt off 2.51pt}]  (141,79) .. controls (150,78) and (157,79) .. (170,83) ;
\draw  [dash pattern={on 0.84pt off 2.51pt}]  (151,57) .. controls (163,57) and (169,59) .. (180,62) ;
\draw    (166,118) .. controls (175,116) and (193,111) .. (199,100) ;
\draw    (152,103) .. controls (158,100) and (178,98) .. (184,88) ;
\draw  [dash pattern={on 0.84pt off 2.51pt}]  (90,136) .. controls (102,124) and (119,124) .. (135,124) ;
\draw    (370,152) .. controls (390,152) and (408,145) .. (421,149) .. controls (434,153) and (451,167) .. (465,174) ;
\draw    (373,137) .. controls (393,137) and (416,125) .. (429,129) .. controls (442,133) and (449,145) .. (471,156) ;
\draw    (268,228) .. controls (280.19,227.75) and (300.1,216.3) .. (318,216) .. controls (335.9,215.7) and (350.83,214.54) .. (367,208) .. controls (383.17,201.46) and (419,163) .. (427,152) ;
\draw    (269,209) .. controls (281.19,208.75) and (297.1,200.3) .. (315,200) .. controls (332.9,199.7) and (343.83,197.54) .. (360,191) .. controls (376.17,184.46) and (396,158) .. (407,147) ;
\draw    (269,195) .. controls (279,195) and (289,195) .. (297,203) ;
\draw    (266,184) .. controls (273,185) and (299,184) .. (309,199) ;
\draw    (303,217) -- (310,234) ;
\draw    (318,216) -- (325,230) ;
\draw  [dash pattern={on 0.84pt off 2.51pt}]  (310,234) .. controls (315,244) and (323,250) .. (326,255) ;
\draw  [dash pattern={on 0.84pt off 2.51pt}]  (325,230) .. controls (330,240) and (344,243) .. (347,248) ;
\draw [color={rgb, 255:red, 80; green, 227; blue, 194 }  ,draw opacity=1 ]   (324,36) -- (353,36) ;
\draw [shift={(355,36)}, rotate = 180] [color={rgb, 255:red, 80; green, 227; blue, 194 }  ,draw opacity=1 ][line width=0.75]    (10.93,-3.29) .. controls (6.95,-1.4) and (3.31,-0.3) .. (0,0) .. controls (3.31,0.3) and (6.95,1.4) .. (10.93,3.29)   ;
\draw [color={rgb, 255:red, 80; green, 227; blue, 194 }  ,draw opacity=1 ]   (234,263) -- (263,263) ;
\draw [shift={(265,263)}, rotate = 180] [color={rgb, 255:red, 80; green, 227; blue, 194 }  ,draw opacity=1 ][line width=0.75]    (10.93,-3.29) .. controls (6.95,-1.4) and (3.31,-0.3) .. (0,0) .. controls (3.31,0.3) and (6.95,1.4) .. (10.93,3.29)   ;
\draw    (152,219) .. controls (200,235) and (216,221) .. (242.12,220.95) ;
\draw  [color={rgb, 255:red, 155; green, 155; blue, 155 }  ,draw opacity=1 ] (241.75,189.31) .. controls (245.75,189.14) and (249.08,196.08) .. (249.19,204.82) .. controls (249.29,213.56) and (246.12,220.78) .. (242.12,220.95) .. controls (238.12,221.12) and (234.79,214.17) .. (234.69,205.43) .. controls (234.58,196.69) and (237.75,189.47) .. (241.75,189.31) -- cycle ;
\draw    (148,179) .. controls (203,199) and (204,189) .. (241.75,189.31) ;
\draw  [dash pattern={on 0.84pt off 2.51pt}]  (109,170) .. controls (121,173) and (136,174) .. (148,179) ;
\draw  [dash pattern={on 0.84pt off 2.51pt}]  (121,214) .. controls (134,214) and (138,217) .. (152,219) ;
\draw    (369,78) .. controls (412,77) and (425,93) .. (463,63) ;
\draw    (374,103) .. controls (417,102) and (418,121) .. (470,90) ;
\draw  [dash pattern={on 0.84pt off 2.51pt}]  (463,63) .. controls (472,51) and (486,45) .. (499,49) ;
\draw  [dash pattern={on 0.84pt off 2.51pt}]  (470,90) .. controls (482,81) and (495,75) .. (508,77) ;
\draw    (442,137) .. controls (468,115) and (480,134) .. (499,137) ;
\draw  [dash pattern={on 0.84pt off 2.51pt}]  (465,174) .. controls (491,188) and (487,196) .. (500,202) ;
\draw  [dash pattern={on 0.84pt off 2.51pt}]  (471,156) .. controls (497,170) and (511,182) .. (520,190) ;
\draw    (416,129) .. controls (462,87) and (478,110) .. (508,113) ;
\draw  [dash pattern={on 0.84pt off 2.51pt}]  (499,137) .. controls (509,141) and (513,146) .. (524,150) ;
\draw  [dash pattern={on 0.84pt off 2.51pt}]  (508,113) .. controls (522,118) and (524,118) .. (535,122) ;

\draw (214,258) node [anchor=north west][inner sep=0.75pt]  [color={rgb, 255:red, 80; green, 227; blue, 194 }  ,opacity=1 ] [align=left] {$\displaystyle w$};
\draw (274,253) node [anchor=north west][inner sep=0.75pt]  [color={rgb, 255:red, 80; green, 227; blue, 194 }  ,opacity=1 ] [align=left] {$\displaystyle m_{h}( z)$};
\draw (279,27) node [anchor=north west][inner sep=0.75pt]  [color={rgb, 255:red, 80; green, 227; blue, 194 }  ,opacity=1 ] [align=left] {$\displaystyle m_{f}( x)$};
\draw (361,26) node [anchor=north west][inner sep=0.75pt]  [color={rgb, 255:red, 80; green, 227; blue, 194 }  ,opacity=1 ] [align=left] {$\displaystyle m_{g}( y)$};

\end{tikzpicture}}
\]
The interfaces of these cylinders are labeled by morphisms from products of objects in $\msc{E}^{\pm}$, so that one can imagine them as skein diagrams.  Indeed, if one has some notion of a bordism category $\opn{SkBord}_{\msc{E}}^{nc}$ with ``skeins labeled by $\msc{E}$" then we expect a symmetric monoidal functor $\Bord_{\msc{E}}^{nc}\to \opn{SkBord}_{\msc{E}}^{nc}$ which inverts colored merging bordisms, and hence expect a canonical map from the localization
\[
\Bord_{\msc{E}}^{nc}[\Theta^{-1}_{\msc{E}}]\to \opn{SkBord}_{\msc{E}}^{nc}.
\]
We leave it to the interested reader to explore the relation between the localization $\Bord_{A_{fin}^{\heartsuit}}[\Theta^{-1}_{A^{\heartsuit}_{fin}}]$ and the skein bordism category for a ribbon category $A^{\heartsuit}$, in the discrete setting.
\par

From the above perspective, Corollary \ref{cor:der_lrt_skeins} says that the LRT theory for the derived $\infty$-category $\msc{D}(A^{\heartsuit})$ can be taken as a theory from a category of bordisms with ``coarse $\msc{D}(A^{\heartsuit})$-labeled skeins" as above (cf.\ \cite{derenzietal23,brownhaioun,costantinoetal,costantinoetalII}).

\subsection{Proving Theorem \ref{thm:merge_loc}}
\label{sect:proof_merge_bord}

\begin{lemma}\label{lem:3779}
Given a map of cocartesian fibrations
\[
\xymatrix{
\msc{C}\ar[rr]^{q}\ar[dr]_{p} & & \msc{D}\ar[dl]^{p'}\\
 & \msc{T}
}
\]
in which $q$ itself is a cocartesian fibration, and a morphism $f:s\to t$ in $\msc{T}$, the transport functor $f_!:\msc{C}_s\to \msc{C}_t$ sends $q_s$-cocartesian edges to $q_t$-cocartesian edges.
\end{lemma}

\begin{proof}
First, for any $p$-cocartesian edge $\alpha:x\to y$ in $\msc{C}$ over a map $f:s\to t$ in $\msc{T}$ we have that $q(\alpha):q(x)\to q(y)$ is $p'$-cocartesian over $f$. To compare, we take a $q$-cocartesian lift $\alpha':x\to y'$ of $q(\alpha)$, have that $\alpha'$ is $p$-cocartesian, and by uniqueness of $p$-cocartesian lifts of $f$ we obtain a $2$-simplex
\[
\xymatrix{
	& y\ar[dr]^{\sim}\\
x\ar[ur]^{\alpha}\ar[rr]_{\alpha'} & & y'.
}
\]
Since any isomorphism in $\msc{C}$ is $q$-cocartesian, the $2$-of-$3$ property \cite[\href{https://kerodon.net/tag/01TS}{01TS}]{kerodon} assures us that $\alpha$ is $q$-cocartesian.  So we see that any $p$-cocartesian edge in $\msc{C}$ is also $q$-cocartesian.

Considering now the cocartesian transformations $T_f:\Delta^1\times\msc{C}_s\to \msc{C}$ over $f$. Recall that $T_f$ evaluates to a $p$-cocartesian edge in $\msc{C}$ at each object $x$ in the fiber $\msc{C}_s$, and hence to a $p$-cocartesian edge as well.  It follows that, for each $q_s$-cocartesian edge $\beta:x\to u$, which is just a $q$-cocartesian edge which lives in the fiber, $T_f$ produces a diagram
\[
\xymatrix{
x\ar[d]_{\beta}\ar[rr]^(.45){(T_f)_x}\ar@{-->}[drr] & & f_!(x)\ar[d]^{f_!(\beta)}\\
u\ar[rr]_(.45){(T_f)_u} & & f_!(v)
}
\]
in $\msc{C}$ in which all map, save for possibly $f_!(\beta)$, are known to be $q$-cocartesian. We again reference the $2$-of-$3$ property for cocartesian edges to find that $f_!(\beta)$ is in fact $q$-cocartesian.
\end{proof}

We now provide the proof of Theorem \ref{thm:merge_loc}.

\begin{warning}
The proof of Theorem \ref{thm:merge_loc} is quite long. We note that this result is only \emph{used} in the terminal case $\msc{E}=\ast$. So a superficial reading of the arguments in the general setting will not hinder one's ability to engage with subsequent materials.
\end{warning}

\begin{proof}[Proof of Theorem \ref{thm:merge_loc}]
As the proof for general $\msc{E}$ is somewhat rambunctious, we proceed in parts. We first consider the case of $\Bord_{\ast}^{nc}$.

{\it I.\ The trivial case $\msc{E}^{\pm}=\ast$:} For each finite set $I$, any choice of ordering on $I$ provides an injective equivalence into the fiber
\begin{equation}\label{eq:3814}
(\Bord_{\ast}^{nc})^I\overset{\sim}\to (\Bord_{\ast}^{nc})^{\ot}_I.
\end{equation}
This is according to the explicit construction of the $\infty$-category $(\Bord^{nc}_{\ast})^{\ot}$ from Definitions \ref{def:e_ot} and \ref{def:bord_nc}.  Under such an identification the transport functor
\begin{equation}\label{eq:4424}
\xymatrix{
(\Bord_{\ast}^{nc})^{\ot}_I\ar[rr]^{\bar{f}_!} & & (\Bord_{\ast}^{nc})^{\ot}_J\\
(\Bord_{\ast}^{nc})^I\ar[rr]_{disj\ union}\ar[u]^{incl} & & (\Bord_{\ast}^{nc})^J\ar[u]_{incl} 
}
\end{equation}
along a map $\bar{f}:I\to J$ in $\Fin_{\ast}$ is given by the product functor on $\Bord_{\ast}^{nc}$, which we recall is the disjoint union operation (see Proposition \ref{prop:symm_infty}). To be clear, we also insert the empty surface in the $j$-th factor whenever the preimage $\bar{f}^{-1}(j)$ is empty. 
\par

Let us take $\Theta=\Theta_{\ast}$.  For each finite set $I$ we consider the collection $\Theta^I$ of tuples of colored merging bordisms in $(\Bord_{\ast}^{nc})^I$, and define the class of maps $\hat{\Theta}_I$ in $(\Bord_{\ast}^{nc})^{\ot}_I$ as those maps which are isomorphic to a map in $\Theta^I$ in $\Fun(\Delta^1,(\Bord_{\ast}^{nc})^{\ot}_I)$.  (Here we identify $\Theta^I$ with a collection in the fiber via the injective equivalence \eqref{eq:3814}.)
\par

Though the embedding \eqref{eq:3814} does depend on a choice of ordering, we note that all such embeddings are naturally isomorphic via the apparent permutation isomorphisms, and so find that the class $\hat{\Theta}_I$ does not depend on any choice of ordering.  We take now $\hat{\Theta}$ the union of all of these fibers
\[
\hat{\Theta}=\coprod_{I\text{ a fin set}}\hat{\Theta}_I\ \subseteq\ (\Bord_{\ast}^{nc})^{\ot}[1].
\]
\par

Since the class of colored merging bordisms is stable under disjoint union, by construction, we see from the diagram \eqref{eq:4424} that each transport function $\bar{f}_!:(\Bord_{\ast}^{nc})^{\ot}_I\to (\Bord_{\ast}^{nc})^{\ot}_J$ sends $\hat{\Theta}_I$ into $\hat{\Theta}_J$.  So $\hat{\Theta}$ is a vertical class in $(\Bord_{\ast}^{nc})^{\ot}$, in the sense of Defintion \ref{def:vertical}.  Furthermore, since $\hat{\Theta}_{\{0\}}$ is simply the closure of $\Theta$ up to isomorphism in $\Fun(\Delta^1,\Bord_{\ast}^{nc})$, one sees from the diagrams
\[
\xymatrix{
(\Bord_{\ast}^{nc})^{\ot}_I\ar[drr]^{transport\ equiv}\\
(\Bord_{\ast}^{nc})^I\ar[u]^{include}\ar[rr]_= & & (\Bord_{\ast}^{nc})^I
}
\]
that each $\hat{\Theta}_I$ is identified with the pullback of $\hat{\Theta}_{\{0\}}^I\subseteq (\Bord^{nc}_{\ast})^I[1]$ along the transport equivalence.  Hence, $\hat{\Theta}$ is the vertical class generated by $\hat{\Theta}_{\{0\}}$, and by Proposition \ref{prop:sym_loc} the vertical localization $(\Bord_{\ast}^{nc})^{\ot}[\hat{\Theta}^{-1}]$ is symmetric monoidal.  
\par

Since $\hat{\Theta}_{\{0\}}$ is generated by $\Theta$ under compositions with isomorphisms, the natural map $\Bord_{\ast}^{nc}[\Theta^{-1}]\to \Bord_{\ast}^{nc}[\hat{\Theta}^{-1}_{\{0\}}]$ is an equivalence.  Hence, the proposed universal property for the localization is a consequence of Lemma \ref{lem:adv_loc}.  We now turn to the general setting.
\vspace{2mm}

{\it II.\ Defining the class $\hat{\Theta}_{\msc{E}}$:} Consider an arbitrary $fr\Disk$-category $\msc{E}^{\ot}\to fr\Disk^{\ot}$ and let $q:(\Bord^{nc}_{\msc{E}})^{\ot}\to (\Bord_{\ast}^{nc})^{\ot}$ denote the corresponding cocartesian fibration provided by the forgetful functor.  Take $\hat{\Theta}_{\msc{E};I}$ in $(\Bord^{nc}_{\msc{E}})^{\ot}_I$ the collection of $q_I$-cocartesian edges which are sent to maps in $\hat{\Theta}_I$ along the pulled back fibration
\[
q_I:(\Bord_{\msc{E}}^{nc})^{\ot}_I\to (\Bord_{\ast}^{nc})^{\ot}_I.
\]
We note that each $\hat{\Theta}_{\msc{E};I}$ is stable under isomorphism in $\Fun(\Delta^1,(\Bord_{\msc{E}}^{nc})^{\ot}_I)$ since each $\hat{\Theta}_I$ is stable under isomorphism in $\Fun(\Delta^1,(\Bord_{\ast}^{nc})^{\ot}_I)$.  Take now $\hat{\Theta}_{\msc{E}}$ the union of all of the $\hat{\Theta}_{\msc{E};I}$ in $(\Bord_{\msc{E}}^{nc})^{\ot}$.
\vspace{2mm}

{\it III.\ The class $\hat{\Theta}_{\msc{E}}$ is stable under transport:} Over any map of finite pointed sets $\bar{f}:I\to J$ we have a diagram
\[
\xymatrix{
(\Bord_{\msc{E}}^{nc})^{\ot}_I\ar[rr]^{\bar{f}_!}\ar[d]_{q_I} & & (\Bord_{\msc{E}}^{nc})^{\ot}_J\ar[d]^{q_J}\\
(\Bord_{\ast}^{nc})^{\ot}_I\ar[rr]_{\bar{f}_!} & & (\Bord_{\ast}^{nc})^{\ot}_J
}
\]
in $\sCat_{\infty}$.  By our analysis in the case $\msc{E}^{\pm}=\ast$ we understand that the bottom map sends $\hat{\Theta}_I$ into $\hat{\Theta}_J$.  Since transport for $(\Bord_{\msc{E}}^{nc})^{\ot}$ sends $q_I$-cocartesian edges to $q_J$-cocartesian edges (Lemma \ref{lem:3779}) it follows that the functor $\bar{f}_!:(\Bord_{\msc{E}}^{nc})^{\ot}_I\to (\Bord_{\msc{E}}^{nc})^{\ot}_J$ sends $\hat{\Theta}_{\msc{E};I}$ into $\hat{\Theta}_{\msc{E};J}$.  Thus $\hat{\Theta}_{\msc{E}}$ is a vertical class in $(\Bord_{\msc{E}}^{nc})^{\ot}$.
\vspace{2mm}

{\it IV.\ The class $\hat{\Theta}_{\msc{E}}$ is generated by $\hat{\Theta}_{\msc{E};\{0\}}$:} Let $I$ be an arbitrary finite pointed set. Considering transport over the inert projections $\bar{\rho}_i:I\to \{0\}$, we have a diagram
\begin{equation}\label{eq:3866}
\xymatrix{
(\Bord_{\msc{E}}^{nc})^{\ot}_I\ar[rr]^{\bar{\rho}_!}_{\sim}\ar[d]_{q_I} & & (\Bord_{\msc{E}}^{nc})^I\ar[d]^{q_{\{0\}}^I}\\
(\Bord_{\ast}^{nc})^{\ot}_I\ar[rr]^{\bar{\rho}_!}_{\sim} & & (\Bord_{\ast}^{nc})^I.
}
\end{equation}
By ({\it III}) we know that $\bar{\rho}_!$ sends $\hat{\Theta}_{\msc{E};I}$ into $\hat{\Theta}_{\msc{E};\{0\}}^I$, and by our analysis in the case $\msc{E}^{\pm}=\ast$ we understand that $\hat{\Theta}_I$ is precisely the preimage of $\hat{\Theta}^I_{\{0\}}$ under the transport equivalence.
\par

Now, for any map $\alpha:x\to y$ in $(\Bord_{\msc{E}}^{nc})^{\ot}_I$ which has image in the product $\hat{\Theta}^I_{\msc{E};\{0\}}\subseteq(\Bord^{nc}_{\msc{E}})^I[1]$ we have $q_I(\alpha)$ in $\hat{\Theta}_I$.  We want to show that $\alpha$ is in $\hat{\Theta}_{\msc{E};I}$, which reduces to showing that $\alpha$ is $q_I$-cocartesian.  For this point, we can take a cocartesian lift $\alpha':x\to y'$ of $q_I(\alpha)$ then lift the simplex
\[
\xymatrix{
	& q_I(y')\ar[dr]^=\\
q_I(x)\ar[ur]\ar[rr] & & q_I(y)
}
\]
along $q_I$ to get a simplex
\[
\xymatrix{
	& y'\ar[dr]^{\tau}\\
x\ar[rr]_{\alpha}\ar[ur]^{\alpha'} & & y
}
\]
in $(\Bord_{\msc{E}}^{nc})^{\ot}_I$.
\par

Taking the image under $\bar{\rho}_!$ we obtain a $2$-simplex
\begin{equation}\label{eq:3893}
\xymatrix{
	& \bar{\rho}_!(y')\ar[dr]^{\bar{\rho}_!(\tau)}\\
\bar{\rho}_!(x)\ar[rr]_{\bar{\rho}_!(\alpha)}\ar[ur]^{\bar{\rho}_!(\alpha')} & & \bar{\rho}_!(y)
}
\end{equation}
in which the two maps from $\bar{\rho}_!(x)$ are $q^I$-cocartesian, and hence in which all of the edges are $q^I$-cocartesian by the $2$-of-$3$ property \cite[\href{https://kerodon.net/tag/01TS}{01TS}]{kerodon}.  By (non-strict) commutativity of the diagram \eqref{eq:3866} we see that the $2$-simplex \eqref{eq:3893} is sent to a $2$-simplex of the form
\[
\xymatrix{
	& b'\ar[dr]^{\sim}\\
a\ar[ur]\ar[rr] & & b
}
\]
in $(\Bord_{\ast}^{nc})^I$.  Since the lift of an isomorphism along a cocartesian fibration is a cocartesian edge if and only if it is an isomorphism, it follows that $\bar{\rho}_!(\tau)$ is an isomorphism in $(\Bord_{\msc{E}}^{nc})^I$.  Finally, since $\bar{\rho}_!$ is an equivalence we find that $\tau$ itself is an isomorphism, is therefore $q_I$-cocartesian, and by the $2$-of-$3$ property again we find that $\alpha$ is $q_I$-cocartesian.  Hence $\alpha$ is in $\hat{\Theta}_{\msc{E};I}$.  This show that the class $\hat{\Theta}_{\msc{E}}$ is in fact generated by $\hat{\Theta}_{\msc{E};\{0\}}$.

{\it V.\ The universal property for $\hat{\Theta}_{\msc{E}}$-localization:} We now understand that $\hat{\Theta}_{\msc{E}}$ is the vertical class generated by $\hat{\Theta}_{\msc{E};\{0\}}$.  Thus the localization $(\Bord^{nc}_{\msc{E}})^{\ot}[\hat{\Theta}^{-1}_{\msc{E}}]$ has the proposed universal property for symmetric monoidal functors which invert the generators $\hat{\Theta}_{\msc{E};\{0\}}$ in $\Bord_{\msc{E}}^{nc}$.  This follows by Lemma \ref{lem:adv_loc}.  So it suffices to show that maps out of $\Bord_{\msc{E}}^{nc}$ which sends the subset $\Theta_{\msc{E}}$ to isomorphisms sends all maps in $\hat{\Theta}_{\msc{E};\{0\}}$ to isomorphisms.  For this it suffices, further, to show that every map $\alpha:x\to y$ in $\hat{\Theta}_{\msc{E};\{0\}}$ fits into a diagram
\begin{equation}\label{eq:3909}
\xymatrix{
x'\ar[rr]^{\alpha'}\ar[d] & & y'\ar[d]\\
x\ar[rr]_{\alpha} & & y
}
\end{equation}
in which $\alpha':x'\to y'$ is in $\Theta_{\msc{E}}$ and the vertical maps are isomorphisms.
\par

Consider the forgetful functor $q_0=q_{\{0\}}:\Bord_{\msc{E}}^{nc}\to \Bord_{\ast}^{nc}$.  By definition, any map $\alpha:x\to y$ in $\hat{\Theta}_{\msc{E};0}$ is a $q_0$-cocartesian edge and the image $q_0(x)\to q_0(y)$ fits into a diagram
\begin{equation}\label{eq:3919}
\xymatrix{
\bar{x}'\ar[rr]^{\bar{\alpha}'}\ar[d]_{\cong}\ar[drr] & & \bar{y}'\ar[d]^{\cong}\\
q_0(x)\ar[rr]_{\alpha} & & q_0(y).
}
\end{equation}
in which $\bar{\alpha}'$ is a colored merging bordism.  Since $q_0$ is a cocartesian fibration, and hence an isofibration, we can lift the isomorphism $\bar{x}'\to q_0(x)$ to an isomorphism $x'\to x$, from some unspecified object over $\bar{x}'$.  We then take a cocartesian lift $\alpha':x'\to y'$ of $\bar{\alpha}'$ to obtain a partial diagram
\begin{equation}\label{eq:3926}
\xymatrix{
x'\ar[rr]^{\alpha'}\ar[d]_{\cong}\ar@{..>}[drr] & & y'\ar@{..>}[d]\\
x\ar[rr]_{\alpha} & & y.
}
\end{equation}
By solving some lifting problems along the cocartesian fibration $q_0:\Bord_{\msc{E}}^{nc}\to \Bord_{\ast}^{nc}$ we complete this partial diagram to a diagram $t:\Delta^1\times\Delta^1\to \Bord_{\msc{E}}^{nc}$ of the form \eqref{eq:3909}, and with $q(t):\Delta^1\times\Delta^1\to \Bord_{\ast}^{nc}$ recovering \eqref{eq:3919}.  By the $2$-of-$3$ property for cocartesian edges, we find that all edges in such a completion of \eqref{eq:3926} are $q_0$-cocartesian.  In particular, $\alpha':x'\to y'$ is in $\Theta_{\msc{E}}$ and the map $y'\to y$ is an isomorphism.  So we have the desired result, and observed the claimed universal property, i.e.\ that restriction provides an equivalence
\[
\Fun_{\Fin_{\ast}}^{cc}\left((\Bord_{\msc{E}}^{nc})^{\ot}[\hat{\Theta}_{\msc{E}}^{-1}],\msc{G}^{\ot}\right)\overset{\sim}\to \Fun_{\Fin_{\ast}}^{cc}\left((\Bord_{\msc{E}}^{nc})^{\ot},\msc{G}^{\ot}\right)_{\Theta_{\msc{E}}}
\]
at any symmetric monoidal $\infty$-category $\msc{G}^{\ot}\to \Fin_{\ast}$.
\end{proof}

\section{Emergent (unmarked) field theories}
\label{sect:unmarked_lrt}

We extract from the derived theory $L_{\msc{D}}:\Bord_{\msc{D}}^{nc}\to \Vect_{\pf}$ of Section \ref{sect:derived_lrt}--or really its $\Vect$-valued restriction $L_{\msc{D}_{fin}}$--an unmarked field theory $\mbb{L}_{D_{fin}}:\Bord_{3,2}^{nc}\to \opn{h}\Vect$. As explained in the preamble to Section \ref{sect:derived_w_skeins}, producing the theory $\mbb{L}_{D_{fin}}$ reduces to a somewhat delicate, hands-on analysis of the homotopy truncation $\opn{h}\Bord_{\ast}^{nc}[\Theta_{\ast}^{-1}]$.
\par

In human terms, we are claiming that after marking all surfaces by the unit and washing out all homotopy relations, the values of the theory $L_{\msc{D}_{fin}}$ does not actually depend on any choices of decorations on surfaces or bordisms.

\subsection{Unit labeled bordisms}

Let $q:\msc{E}^{\ot}\to fr\Disk^{\ot}$ be a $fr\Disk$-category. We consider the initial object $\emptyset$ in $fr\Disk^{\ot}$ and define the fibration of unit objects in $\msc{E}^{\ot}$ as the $\infty$-category
\[
\opn{Units}(\msc{E}^{\otimes})=\left\{\begin{array}{c}
\text{the full subcategory in the fiber product }\\
\msc{E}^{\otimes}_{\emptyset}\times_{\Fun(\{0\},\msc{E}^{\ot})}\Fun(\Delta^1,\msc{E}^{\ot})\\
\text{spanned by $q$-cocartesian edges}
\end{array}\right.
\]
equipped with the isofibration
\[
\opn{Units}(\msc{E})^{\ot}\overset{ev_1}\longrightarrow \msc{E}^{\ot}\overset{q}\to fr\Disk^{\ot}.
\]

\begin{lemma}\label{lem:4866}
The structure map $\opn{Units}(\msc{E})^{\ot}\to fr\Disk^{\ot}$ is a trivial Kan fibration.
\end{lemma}

\begin{proof}
To ease notation take $\msc{T}=fr\Disk^{\ot}$. By \cite[\href{https://kerodon.net/tag/01VK}{01VK}]{kerodon}, the cocartesian fibration $\msc{E}^{\ot}\to \msc{T}$ induces a trivial Kan fibration
\[
\opn{Units}(\msc{E})^{\ot}\to \msc{E}_{\emptyset}^{\ot}\times_{\Fun(\{0\},\msc{T})}\Fun(\Delta^1,\msc{T}),
\]
and since $\msc{E}^{\ot}_{\emptyset}$ is contractible the map $\msc{E}^{\ot}_{\emptyset}\to \{\emptyset\}$ over $\msc{T}$ is a trivial Kan fibration as well. Hence the map
\[
\msc{E}_{\emptyset}^{\ot}\times_{\Fun(\{0\},\msc{T})}\Fun(\Delta^1,\msc{T})\to \{\emptyset\}\times_{\Fun(\{0\},\msc{T})}\Fun(\Delta^1,\msc{T})
\]
is a trivial Kan fibration. Finally, since $\emptyset$ is initial in $\msc{T}$, evaluation induces a trivial Kan fibration
\[
ev_1:\{\emptyset\}\times_{\Fun(\{0\},\msc{T})}\Fun(\Delta^1,\msc{T})\to \msc{T}
\]
\cite[\href{https://kerodon.net/tag/06R2}{06R2}, \href{https://kerodon.net/tag/01KU}{01KU}]{kerodon}. As any composite of trivial Kan fibrations is a trivial Kan fibration, we obtain the desired result.
\end{proof}

By Lemma \ref{lem:4866} we observe the unique $fr\Disk$-monoidal section $fr\Disk^{\ot}\to \opn{Units}(\msc{E})^{\ot}$ of the structure map, and compose with the evaluation functor $\opn{Units}(\msc{E})^{\ot}\to \msc{E}^{\ot}$ to obtain a $fr\Disk^{\ot}$-monoidal functor
\[
\opn{unit}_q:fr\Disk^{\ot}\to \msc{E}^{\ot}.
\]
We refer to this functor as the unit map.
\par

The unit map now determines, by pullback along the outgoing boundary $\partial_{out}^{\ot}:(\Bord_{\ast}^{nc})^{\ot}\to fr\Disk^{\ot}$, a functor $(\Bord_{\ast}^{nc})^{\ot}\to (\Bord_{\msc{E}}^{nc})^{\ot}$ which is simultaneously symmetric monoidal and a map cocartesian fibrations over $(\Bord_{\ast}^{nc})^{\ot}$,
\[
\xymatrix{
(\Bord_{\ast}^{nc})^{\ot}\ar[dr]_{id}\ar[rr] & & (\Bord_{\msc{E}}^{nc})^{\ot}\ar[dl]\\
	& (\Bord_{\ast}^{nc})^{\ot}	&.
}
\]

\begin{definition}
For any $fr\Disk$-category $q:\msc{E}^{\ot}\to fr\Disk^{\ot}$ we let $\opn{unit}^{\ot}_{\msc{E}}$ denote the pullback functor
\[
\opn{unit}^{\ot}_{\msc{E}}:=(\Bord_{\ast}^{nc})^{\ot}\times_{fr\Disk^{\ot}}\opn{unit}_q:(\Bord_{\ast}^{nc})^{\ot}\to (\Bord_{\msc{E}}^{nc})^{\ot}.
\]
\end{definition}

\begin{remark}
What these unit functors do is fairly intuitive. Namely, any $fr\Disk$-category $\msc{E}^{\ot}\to fr\Disk^{\ot}$ admits unique ``unit objects" $\1_{t}$ which are the images of the transport functors $\1_{t}:\ast\cong \msc{E}^{\ot}_{\emptyset}\to \msc{E}^{\ot}_{t}$ along each (unique) morphisms $\emptyset\to t$ in $fr\Disk^{\ot}$. At the level of the underlying $\infty$-categories, the map $\opn{unit}_{\msc{E}}:\Bord_{\ast}^{nc}\to \Bord_{\msc{E}}$ labels each marked surface $\zeta:D\times \mcl{I}\to \Sigma$ by the unit object $\1_{\mcl{I}}$, and labels each bordism by the unique cocartesian map between units $\1_{\mcl{I}}$ and $\1_{\mcl{J}}$.
\end{remark}

\subsection{LRT theories along the unit}

In the case of the $fr\Disk$-category $\msc{D}_{fin}^{\ot}\to fr\Disk^{\ot}$ associated to the derived $\infty$-category of a modular tensor category $A^{\heartsuit}$, we restrict along the symmetric monoidal functor
\[
\opn{unit}_{\msc{D}_{fin}}:\Bord_{\ast}^{nc}\to\Bord_{\msc{D}_{fin}}^{nc}
\]
to obtain a field theory from the trivially labeled bordism category
\[
L_{\msc{D}_{fin}}\circ\opn{unit}_{\msc{D}_{fin}}:\Bord_{\ast}^{nc}\to \Vect.
\]

\begin{lemma}\label{lem:lrt_unit}
For any modular tensor category $A^{\heartsuit}$, the functor $L_{\msc{D}_{fin}}\circ\opn{unit}_{\msc{E}}$ sends all colored merging bordisms to isomorphisms in $\Vect$, and hence localizes to a uniquely determined symmetric monoidal functor
\[
\mbb{L}^{\opn{unit}}_{\msc{D}_{fin}}:\Bord_{\ast}^{nc}[\Theta^{-1}_{\ast}]\to \Vect.
\]
\end{lemma}

\begin{proof}
Since $\opn{unit}_{\msc{D}_{fin}}$ is a map of cocartesian fibrations over $\Bord_{\ast}^{nc}$, it sends colored merging bordisms in $\Bord_{\ast}^{nc}$ to colored merging bordisms in $\Bord_{\msc{D}_{fin}}^{nc}$.  Since $L_{\msc{D}_{fin}}$ is obtained from $L_{\msc{D}}$ via restriction, the result now follows from Corollary \ref{cor:der_lrt_skeins} and Theorem \ref{thm:merge_loc}.
\end{proof}

For this section, we are especially interested in the induced map on homotopy categories
\[
\opn{h}\mbb{L}^{\opn{unit}}_{\msc{D}_{fin}}:\opn{h}\Bord_{\ast}^{nc}[\Theta^{-1}_{\ast}]\to \opn{h}\Vect,
\]
where we note that the target $\opn{h}\Vect$ is the discrete derived category of linear cochains.  We claim, in fact, that the homotopy category $\opn{h}\Bord^{nc}_{\ast}[\Theta_{\ast}^{-1}]$ is identified with the unlabeled $3$-dimensional bordism category, and hence that $\opn{h}\mbb{L}^{\opn{unit}}_{\msc{D}_{fin}}$ provides an unlabeled topological field theory in $3$-dimensions. We argue the point below.

\subsection{Reductions in the homotopy bordism category}

Our goal is to prove that the forgetful functor $\opn{h}\Bord_{\ast}^{nc}[\Theta_{\ast}^{-1}]\to \Bord^{nc}_{3,2}$ is an equivalence. We first establish a collection of relations for morphisms in the homotopy truncation $\opn{h}\Bord_{\ast}^{nc}[\Theta^{-1}_{\ast}]$.

\begin{lemma}\label{lem:rel1}
\begin{enumerate}
\item Any colored merging bordism $M_{\tau}:\Sigma_{\zeta}\to \Sigma_{\eta}$ decomposes as a composite
\[
M_{\tau}=M'_{\tau'}\circ M''_{\tau''}
\]
where $M'_{\tau'}$ is a special merging bordism and $M''_{\tau''}$ is a color change bordism.\vspace{1mm}
\item For any special merging bordism $M_{\tau}:\Sigma_{\zeta}\to \Sigma_{\eta}$, there is a bordism $M'_{\tau'}:\Sigma_{\eta}\to \Sigma_{\zeta}$ for which $[M_{\tau}]^{-1}=[M'_{\tau'}]$ in $\opn{h}\Bord_{\ast}^{nc}[\Theta_{\ast}^{-1}]$.\vspace{1mm}
\item For any color change bordism $M_{\tau}:\Sigma_{\zeta}\to \Sigma_{\zeta'}$, and arbitrary $N_{\omega}:\Sigma_{\zeta}\to \Sigma'_{\eta}$ there is a color change bordism $M'_{\tau'}:\Sigma'_{\eta'}\to \Sigma'_{\eta}$, and some other bordism $N_{\omega'}:\Sigma_{\zeta'}\to \Sigma'_{\eta'}$ for which
\[
[N_{\omega}]\circ [M_{\tau}]^{-1}=[M'_{\tau'}]^{-1}\circ[N_{\omega'}]
\]
in $\opn{h}\Bord_{\ast}^{nc}[\Theta^{-1}_{\ast}]$.
\end{enumerate}
\end{lemma}

\begin{proof}
(1) By definition $M_{\tau}$ consists of a disjoint union of special merging bordisms on some collection of components, and color change bordisms on others.  Isolating the relevant components provides the claimed factorization.
\par

(2) By definition (Definition \ref{def:special_merge}) every special merging bordism admits a right inverse in $\Bord_{\ast}^{nc}$.   So there exists $M'_{\tau'}$ with $M_{\tau}\circ M'_{\tau'}=id_{\Sigma_{\eta}}$, and in the homotopy category we get
\[
[M_{\tau}]^{-1}=[M_{\tau}]^{-1}\circ[M_{\tau}]\circ[M'_{\tau'}]=[M'_{\tau'}].
\]
\par

(3) Let $M'_{\tau'}:\Sigma_{\zeta}\to \Sigma_{\zeta'}$ be the unique color change bordism which changes all negative markings to positive markings.  Then we have
\[
M'_{\tau'}\circ N_{\omega}=N_{\omega''}\circ M''_{\tau''}
\]
where $M''_{\tau''}:\Sigma_{\zeta}\to \Sigma_{\zeta''}$ changes all negative markings to positive markings and $N_{\omega''}$ is obtained from $N_{\omega}$ by changing all negatively colored cylinders to positively colored cylinders.  Since $M_{\tau}:\Sigma_{\zeta}\to \Sigma_{\zeta'}$ changes some select negative markings to positive markings, we have $M''_{\tau''}=M'''_{\tau'''}\circ M_{\tau}$ for the color change bordism which changes the remaining negative markings on $\Sigma_{\zeta'}$ to positives.  Hence $N_{\omega''}\circ M''_{\tau''}=N_{\omega'}\circ M_{\tau}$ where $N_{\omega'}=N_{\omega''}\circ M'''_{\tau'''}$.  We now have 
\[
[N_{\omega}]\circ [M_{\tau}]^{-1}=[M'_{\tau'}]^{-1}\circ[M'_{\tau'}]\circ[N_{\omega}]\circ[M_{\tau}]^{-1}=[M'_{\tau'}]^{-1}\circ[N_{\omega'}]
\]
in the homotopy category.
\end{proof}

\begin{lemma}\label{lem:rel2}
Let $\zeta:D\to \Sigma$ be a connected surface with a single positive (resp.\ negative) marking, and $\kappa:D'\to \Sigma$ be any auxiliary disk embedding which contains the image of $\zeta$ in its interior.  There are two smooth homotopies $h^1,h^2:[0,1]\times D\to \Sigma$ with the following properties:
\begin{itemize}
\item $h^1_0=\zeta$ and $h^2_t=\zeta$ at all times $t\in [0,1]$.\vspace{1.5mm}
\item $h^{1}$ is constant for $t<1/3$ and $t>2/3$, and for $t\in [1/3,2/3]$ it simultaneously shrinks the image of $\zeta$ by a precomposition factor $c^{-1}id_D:D\to D$, for some $c\geq 1$, and shifts the image of $D$ linearly towards $\sqrt{-1}\in D'$.\vspace{1.5mm}
\item The terminal maps $h^{1}_1,h^2_1:D\to \Sigma$ have disjoint image.\vspace{1.5mm}
\end{itemize}
Furthermore, for $\eta:D\times \{1,2\}\to \Sigma$ the resulting positively (resp.\ negatively) marked surface with $\eta|_{D_{\varepsilon}}=h^{\varepsilon}_1$, and $M^{\varepsilon}_{\tau^{\varepsilon}}:\Sigma_{\zeta}\to \Sigma_{\eta}$ the resulting bordisms with $M^{\varepsilon}=\Sigma\times[0,1]$ and
\[
\tau^{\varepsilon}:D\times[0,1]\to \Sigma\times[0,1],\ \ (z,t)\mapsto(h^{\varepsilon}_t(z),t),
\]
we have $[M^1_{\tau^1}]=[M^2_{\tau^2}]$ in the homotopy category $\opn{h}\Bord_{\ast}^{nc}[\Theta_{\ast}^{-1}]$.
\end{lemma}

\begin{proof}
The existence of such homotopies is clear. As for the latter claim, $[M^1_{\tau^1}]=[M^2_{\tau^2}]$, it suffices to find a bordism $N_{\omega}:\Sigma_{\eta}\to \Sigma_{\zeta}$ in $\Bord_{\ast}$ which is a composition of maps in $\Theta_{\ast}$ and for which $N_{\omega}\circ M^1_{\tau^1}=  N_{\omega}\circ M^2_{\tau^2}$.  Indeed, such $N_{\omega}:\Sigma_{\eta}\to \Sigma_{\zeta}$ becomes a unit in the localized homotopy category and we therefore have
\[
[M^1_{\tau^1}]=[N_{\omega}]^{-1}\circ [N_{\omega}]\circ [M^1_{\tau^1}]=[N_{\omega}]^{-1}\circ[N_{\omega}]\circ[M^2_{\tau^2}]=[M^2_{\tau^2}].
\]
\par

In the case that $\Sigma_{\zeta}$ is positively marked, we can take $N_{\omega}$ to be a special merging bordism.  Explicitly, we can construct such $N_{\omega}$ via a sequence of smooth homotopies.  First take $y^0:[0,1]\times D\times\{1,2\}\to D\times\{1,2\}$ with $y^0_0=id$ and which shrinks each disk by a factor of $1/3$, then shifts the first disk linearly towards $\sqrt{-1}$ in the ambient disk $D_1$, and shifts the second disk linearly towards $-\sqrt{-1}$ in $D_2$.  We then have the homotopy $y^1:[0,1]\times D\times\{1,2\}\to \Sigma$ which applies $h^{\varepsilon}$ on $D_{\varepsilon}$ in reverse time, and we take $y:[0,1]\times D\times \{1,2\}\to \Sigma$ to be the composite
\[
[0,1]\times D\times\{1,2\}\overset{\opn{diag}\times 1}\to [0,1]^2\times D\times\{1,2\}\overset{1\times y^0}\to [0,1]\times D\times\{1,2\}\overset{y^1}\to \Sigma.
\]
For $N=\Sigma\times [0,1]$ we then obtain $\omega$ and the function $\omega(z,t,i)=(y(t,z,i),t)$.
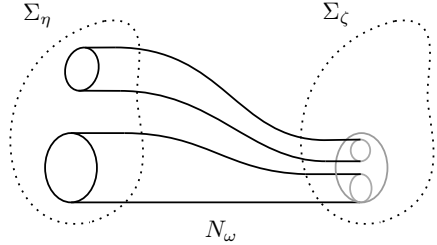
\[
\scalebox{.9}{
\tikzset{every picture/.style={line width=0.75pt}} 
\begin{tikzpicture}[x=0.75pt,y=0.75pt,yscale=-1,xscale=1]
\draw  [dash pattern={on 0.84pt off 2.51pt}]  (154.39,35.45) .. controls (184.21,31.19) and (170.51,63.91) .. (172.93,82.23) .. controls (175.35,100.55) and (185.82,135.88) .. (146.33,147.66) ;
\draw  [dash pattern={on 0.84pt off 2.51pt}]  (154.39,35.45) .. controls (81.04,53.44) and (90.72,159.43) .. (146.33,147.66) ;
\draw   (132.7,61.33) .. controls (133.74,54.69) and (138.69,49.74) .. (143.77,50.26) .. controls (148.85,50.79) and (152.13,56.59) .. (151.09,63.22) .. controls (150.06,69.85) and (145.1,74.8) .. (140.02,74.28) .. controls (134.94,73.76) and (131.67,67.96) .. (132.7,61.33) -- cycle ;
\draw   (121.37,117.28) .. controls (121.33,106.61) and (127.77,97.94) .. (135.76,97.92) .. controls (143.74,97.9) and (150.25,106.53) .. (150.29,117.21) .. controls (150.33,127.88) and (143.89,136.55) .. (135.9,136.57) .. controls (127.91,136.59) and (121.41,127.95) .. (121.37,117.28) -- cycle ;
\draw    (135.76,97.92) -- (165.67,97.93) ;
\draw    (135.9,136.57) -- (297.05,136.53) ;
\draw  [color={rgb, 255:red, 155; green, 155; blue, 155 }  ,draw opacity=0.85 ] (291.44,128.69) .. controls (291.42,124.35) and (294.1,120.82) .. (297.41,120.82) .. controls (300.72,120.81) and (303.42,124.32) .. (303.43,128.66) .. controls (303.45,133) and (300.78,136.52) .. (297.46,136.53) .. controls (294.15,136.54) and (291.46,133.03) .. (291.44,128.69) -- cycle ;
\draw    (266.42,120.83) -- (297.41,120.82) ;
\draw    (165.67,97.93) .. controls (211.62,97.93) and (222.09,119.52) .. (266.42,120.83) ;
\draw  [dash pattern={on 0.84pt off 2.51pt}]  (316.4,35.45) .. controls (346.22,31.19) and (332.52,63.91) .. (334.93,82.23) .. controls (337.35,100.55) and (347.83,135.88) .. (308.34,147.66) ;
\draw  [dash pattern={on 0.84pt off 2.51pt}]  (316.4,35.45) .. controls (243.05,53.44) and (252.72,159.43) .. (308.34,147.66) ;
\draw  [color={rgb, 255:red, 155; green, 155; blue, 155 }  ,draw opacity=1 ] (283.38,117.28) .. controls (283.34,106.61) and (289.78,97.94) .. (297.77,97.92) .. controls (305.75,97.9) and (312.26,106.53) .. (312.3,117.21) .. controls (312.33,127.88) and (305.89,136.55) .. (297.91,136.57) .. controls (289.92,136.59) and (283.42,127.95) .. (283.38,117.28) -- cycle ;
\draw    (140.02,74.28) -- (163.26,74.38) ;
\draw    (143.77,50.26) -- (160.84,50.17) ;
\draw  [color={rgb, 255:red, 155; green, 155; blue, 155 }  ,draw opacity=0.85 ] (291.8,107.6) .. controls (291.79,104.25) and (294.22,101.54) .. (297.22,101.53) .. controls (300.22,101.52) and (302.66,104.23) .. (302.67,107.57) .. controls (302.68,110.92) and (300.26,113.64) .. (297.26,113.65) .. controls (294.26,113.65) and (291.82,110.95) .. (291.8,107.6) -- cycle ;
\draw    (276.9,101.86) -- (297.22,101.53) ;
\draw    (277.71,113.63) -- (297.26,113.65) ;
\draw    (163.26,74.38) .. controls (215.65,74.38) and (235.8,113.63) .. (277.71,113.63) ;
\draw    (160.84,50.17) .. controls (221.29,50.17) and (227.74,101.2) .. (276.9,101.86) ;

\draw (107.73,24.48) node [anchor=north west][inner sep=0.75pt]   [align=left] {$\displaystyle \Sigma _{\eta }$};
\draw (208.9,145.53) node [anchor=north west][inner sep=0.75pt]   [align=left] {$\displaystyle N_{\omega }$};
\draw (274.57,23.83) node [anchor=north west][inner sep=0.75pt]   [align=left] {$\displaystyle \Sigma _{\zeta }$};

\end{tikzpicture}}
\]

It is easy to see that $N_{\omega}:\Sigma_{\eta}\to \Sigma_{\zeta}$ satisfies $N_{\omega}\circ M^2_{\tau^2}=id_{\Sigma_{\zeta}}$, so that $N_{\omega}$ is a special merging bordism, and via a mild argument about bordisms obtained from homotopies we find
\[
N_{\omega}\circ M^1_{\tau^1}=id_{\Sigma_{\zeta}}=N_{\omega}\circ M^2_{\tau^2}.
\]
This gives the claimed equality in the homotopy category.
\par

In the case where $\Sigma_{\zeta}$ is negatively marked, we first apply a color change bordism $\Sigma_{\eta}\to \Sigma_{\eta'}$ to change both negative indices to positive, then apply a special merging bordism constructed exactly as above.  This gives a bordism $N_{\omega}:\Sigma_{\eta}\to \Sigma_{\zeta}$ for which
\[
N_{\omega}\circ M^1_{\tau^1}=N_{\omega}\circ M^2_{\tau'}=\left\{\begin{array}{l}
\text{the unique color change bordism which}\\
\text{exchanges the sole negative marking}\\
\text{on $\Sigma_{\zeta}$ with a positive marking}.
\end{array}\right.
\]
\end{proof}

We also record a basic fact about localization.

\begin{lemma}\label{lem:hc_loc}
For any $\infty$-category $\msc{C}$ and class $W$ in $\msc{C}[1]$ which contains all identity maps, $\msc{C}[W^{-1}]$ is equivalent to its full subcategory which is spanned by all objects which are in the image of the the localization functor $loc:\msc{C}\to\msc{C}[W^{-1}]$.  Furthermore, any morphism $\xi:loc(x)\to loc(y)$ in the localization is a composite
\begin{equation}\label{eq:4146}
\xi=\alpha_1w_1^{-1}\alpha_2 w_2^{-1}\cdots \alpha_rw_r^{-1}
\end{equation}
of morphisms $\alpha_i$ in $\msc{C}$ and $w_i$ in $W$.
\end{lemma}

\begin{proof}
Consider $\msc{C}'\subseteq \msc{C}[W^{-1}]$ the subcategory spanned by those objects in the image of the localization functor, and all morphisms which are expressible as a composite as in \eqref{eq:4146}. It suffices to show that the map on homotopy categories $\opn{h}\msc{C}'\to \opn{h}\msc{C}[W^{-1}]$ is an equivalence.
\par

By the universal property of localization, the identity on $\msc{C}[W^{-1}]$ factors through the inclusion $\msc{C}'\to \msc{C}[W^{-1}]$.  It follows that the inclusion is essentially surjective and that the injection
\[
\pi_0\Maps_{\msc{C}'}(x,y)\to \pi_0\Maps_{\msc{C}[W^{-1}]}(x,y)
\]
is bijective. Hence the inclusion $\opn{h}\msc{C}'\to \opn{h}\msc{C}[W^{-1}]$ is fully faithful as well.
\end{proof}

\subsection{Unmarked bordisms from unit labeled bordisms}

We have the symmetric monoidal forgetful functors $\Bord_{\ast}^{nc}\to \Bord_{3,2}^{nc}$ to the usual $3$-dimensional (non-compact, anomalous) bordism category of unmarked surfaces and unmarked bordisms.  Since this map forgets all internal decorations, it localizes to a symmetric monoidal functor $\Bord_{\ast}^{nc}[\Theta_{\ast}^{-1}]\to \Bord_{3,2}^{nc}$. As the target of this localized functor is discrete, we have the induced symmetric monoidal functor from the homotopy truncation
\[
forget:\opn{h}\Bord_{\ast}^{nc}[\Theta^{-1}_{\ast}]\to \Bord_{3,2}.
\]
It is clear that this functor is essentially surjective.  We are concerned with its behavior on Hom sets.

\begin{lemma}\label{lem:4187}
For positively marked surfaces $\Sigma_{\zeta}$ and $\Sigma'_{\eta}$, the map
\[
\pi_0\Maps_{\Bord_{\ast}^{nc}}(\Sigma_{\zeta},\Sigma'_{\eta})\to \Hom_{\opn{h}\Bord_{\ast}^{nc}[\Theta_{\ast}^{-1}]}(\Sigma_{\zeta},\Sigma'_{\eta})
\]
is surjective.
\end{lemma}

\begin{proof}
By Lemmas \ref{lem:rel1} and \ref{lem:hc_loc} every map $\xi:\Sigma_{\zeta}\to \Sigma'_{\eta}$ in $\opn{h}\Bord_{\ast}^{nc}[\Theta^{-1}_{\ast}]$ is of the form
\[
\xi=[M_{\tau}]^{-1}[N_{\omega}]:\Sigma_{\zeta}\to \Sigma'_{\eta'}\to \Sigma'_{\eta}
\]
where $M_{\tau}:\Sigma'_{\eta}\to \Sigma'_{\eta'}$ is a color change bordism.  However, $\Sigma'_{\eta}$ contains no negative markings, by assumption, so that $\Sigma'_{\eta'}=\Sigma'_{\eta}$ necessarily and $M_{\tau}=id$.  Thus $\xi=[N_{\omega}]$, and we see that all morphisms in the localized category lift to $\Bord_{\ast}^{nc}$.
\end{proof}

\begin{theorem}\label{thm:unmarking}
The functor
\[
forget:\opn{h}\Bord_{\ast}^{nc}[\Theta^{-1}_{\ast}]\to \Bord_{3,2}^{nc}
\]
is a symmetric monoidal equivalence.
\end{theorem}

\begin{proof}
Since we have inverted all color change bordisms and special merging bordisms we see that every object in $\Bord_{\ast}^{nc}[\Theta^{-1}_{\ast}]$ is isomorphic to a surface $\Sigma_{\zeta}$ with a single positive marking on each component.  Hence each map on hom sets
\begin{equation}\label{eq:4213}
forget:\Hom_{\opn{h}\Bord_{\ast}^{nc}[\Theta^{-1}_{\ast}]}(\Sigma_{\zeta},\Sigma'_{\eta})\to \Hom_{\Bord^{nc}_{3,2}}(\Sigma,\Sigma')
\end{equation}
is bijective if and only if it is bijective whenever $\Sigma_{\zeta}$ and $\Sigma'_{\eta}$ carry a single positive marking on each component. We fix two such surfaces.
\par

Via our definition of non-compact bordisms, it is clear that the map \eqref{eq:4213} is surjective.  For injectivity, consider two morphisms $\xi,\xi':\Sigma_{\zeta}\to \Sigma_{\eta}$ with the same image $M:\Sigma\to \Sigma'$ in $\Bord_{3,2}$.  By Lemma \ref{lem:4187} we have $\xi=[N_{\tau}]$ and $\xi'=[N'_{\tau'}]$ for bordisms $N_{\tau},N'_{\tau'}:\Sigma_{\zeta}\to \Sigma'_{\eta}$, and so can take $N=N'=M$.  Now we have two embeddings
\[
\tau,\tau':\opn{Cyl}\to M
\]
which agree in some collars around the boundary, and we want to show $[M_{\tau}]=[M_{\tau'}]$ in $\opn{h}\Bord_{\ast}^{nc}[\Theta^{-1}_{\ast}]$.
\par

After applying homotopies in the interior of $M$, we can assume that $\tau$ and $\tau'$ have disjoint image away from a collar around the incoming boundary $\Sigma\times[0,\delta)\to M$.  By subsequently applying a homotopy in the collar we can assume further that $M_{\tau}$ and $M_{\tau'}$ decompose as composites
\[
M_{\tau}=M_{\omega}\circ M^1_{\tau^1}\ \ \text{and}\ \ M_{\tau'}=M_{\omega}\circ M^2_{\tau^2}
\]
where $\omega$ consists of two embeddings $\omega:\opn{Cyl}\times \{1,2\}\to M$ with disjoint images, and the $M^{\varepsilon}_{\tau^{\varepsilon}}$ are obtained, at each component of the incoming boundary, via homotopy as in the statement of Lemma \ref{lem:rel2}:
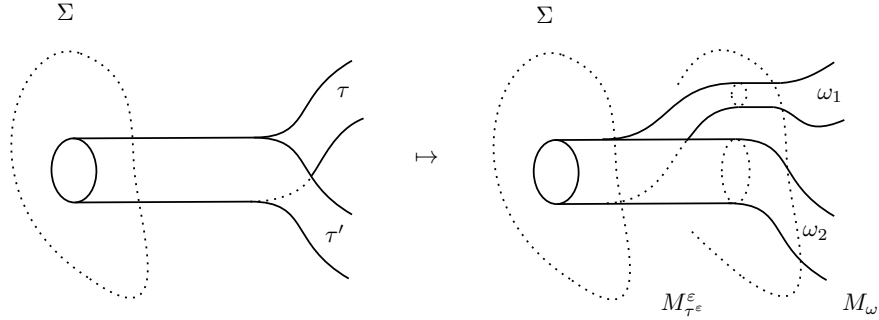
\[
\scalebox{.9}{
\tikzset{every picture/.style={line width=0.75pt}} 

\begin{tikzpicture}[x=0.75pt,y=0.75pt,yscale=-1,xscale=1]

\draw  [dash pattern={on 0.84pt off 2.51pt}]  (87.23,50.91) .. controls (116.07,62.22) and (115.61,82.92) .. (117.06,112.75) .. controls (118.5,142.59) and (144.95,215.1) .. (91.44,175.86) ;
\draw  [dash pattern={on 0.84pt off 2.51pt}]  (91.44,175.86) .. controls (19.98,129.96) and (52.23,38.78) .. (87.23,50.91) ;
\draw   (82.96,98.32) .. controls (89.87,97.99) and (95.78,105.71) .. (96.18,115.57) .. controls (96.57,125.43) and (91.29,133.7) .. (84.39,134.04) .. controls (77.48,134.38) and (71.56,126.65) .. (71.17,116.79) .. controls (70.77,106.93) and (76.05,98.66) .. (82.96,98.32) -- cycle ;
\draw [color={rgb, 255:red, 0; green, 0; blue, 0 }  ,draw opacity=1 ]   (82.96,98.32) -- (183.8,97.86) ;
\draw [color={rgb, 255:red, 0; green, 0; blue, 0 }  ,draw opacity=1 ]   (84.39,134.04) -- (182,133.57) ;
\draw    (183.8,97.86) .. controls (217.72,97.86) and (204.8,74.62) .. (238.72,54.96) ;
\draw    (216.11,119.31) .. controls (222.57,110.37) and (229.21,93.43) .. (245.18,87.14) ;
\draw    (182,133.57) .. controls (215.92,133.57) and (203,156.81) .. (236.91,176.47) ;
\draw    (183.8,97.86) .. controls (217.72,97.86) and (204.8,121.1) .. (238.72,140.76) ;
\draw  [dash pattern={on 0.84pt off 2.51pt}]  (182,133.57) .. controls (200.76,133.61) and (212.88,126.46) .. (216.11,119.31) ;
\draw  [dash pattern={on 0.84pt off 2.51pt}]  (356.13,51.8) .. controls (384.97,63.12) and (384.51,83.82) .. (385.96,113.65) .. controls (387.4,143.48) and (413.85,215.99) .. (360.34,176.75) ;
\draw  [dash pattern={on 0.84pt off 2.51pt}]  (360.34,176.75) .. controls (288.88,130.86) and (321.13,39.67) .. (356.13,51.8) ;
\draw   (351.86,99.22) .. controls (358.77,98.88) and (364.69,106.6) .. (365.08,116.46) .. controls (365.47,126.33) and (360.19,134.59) .. (353.29,134.93) .. controls (346.38,135.27) and (340.46,127.55) .. (340.07,117.69) .. controls (339.67,107.82) and (344.95,99.56) .. (351.86,99.22) -- cycle ;
\draw    (378.41,98.75) .. controls (416.37,99.65) and (412.33,66.58) .. (455.13,67.47) ;
\draw    (425.25,99.65) .. controls (431.71,92.5) and (436.55,81.77) .. (454.32,80.88) ;
\draw    (455.13,67.47) -- (477.74,67.47) ;
\draw    (454.32,80.88) -- (473.7,80.88) ;
\draw [color={rgb, 255:red, 0; green, 0; blue, 0 }  ,draw opacity=1 ]   (351.86,99.22) -- (452.71,98.75) ;
\draw [color={rgb, 255:red, 0; green, 0; blue, 0 }  ,draw opacity=1 ]   (353.29,134.93) -- (450.9,134.47) ;
\draw    (448.48,134.47) .. controls (482.39,134.47) and (469.47,157.7) .. (503.39,177.37) ;
\draw    (452.71,98.75) .. controls (486.62,98.75) and (473.7,121.99) .. (507.62,141.65) ;
\draw    (477.74,67.47) .. controls (487.43,66.58) and (493.89,64.79) .. (507.62,55.86) ;
\draw    (473.7,80.88) .. controls (489.04,82.67) and (488.24,98.75) .. (513.27,88.03) ;
\draw  [dash pattern={on 0.84pt off 2.51pt}]  (380.84,134.5) .. controls (409.1,133.61) and (415.56,111.27) .. (425.25,99.65) ;
\draw  [dash pattern={on 0.84pt off 2.51pt}] (450.44,73.66) .. controls (450.52,70.17) and (452.45,67.39) .. (454.75,67.46) .. controls (457.05,67.53) and (458.85,70.42) .. (458.77,73.91) .. controls (458.68,77.4) and (456.75,80.18) .. (454.45,80.11) .. controls (452.15,80.04) and (450.35,77.15) .. (450.44,73.66) -- cycle ;
\draw  [dash pattern={on 0.84pt off 2.51pt}] (444.94,116.63) .. controls (444.94,106.76) and (448.42,98.75) .. (452.71,98.75) .. controls (456.99,98.75) and (460.47,106.76) .. (460.47,116.63) .. controls (460.47,126.5) and (456.99,134.5) .. (452.71,134.5) .. controls (448.42,134.5) and (444.94,126.5) .. (444.94,116.63) -- cycle ;
\draw  [dash pattern={on 0.84pt off 2.51pt}]  (451.42,50.01) .. controls (480.25,61.33) and (479.8,82.03) .. (481.24,111.86) .. controls (482.69,141.69) and (509.14,214.2) .. (455.62,174.96) ;
\draw  [dash pattern={on 0.84pt off 2.51pt}]  (455.62,174.96) .. controls (444.63,165.78) and (433.32,155.95) .. (426.06,147.91) ;
\draw  [dash pattern={on 0.84pt off 2.51pt}]  (420.4,68.37) .. controls (425.25,60.32) and (433.32,45.13) .. (451.42,50.01) ;

\draw (228.68,67.94) node [anchor=north west][inner sep=0.75pt]   [align=left] {$\tau$};
\draw (222.22,144.8) node [anchor=north west][inner sep=0.75pt]   [align=left] {$\tau'$};
\draw (72.64,21.83) node [anchor=north west][inner sep=0.75pt]   [align=left] {$\Sigma $};
\draw (496.77,68.84) node [anchor=north west][inner sep=0.75pt]   [align=left] {$\omega _{1}$};
\draw (488.7,145.7) node [anchor=north west][inner sep=0.75pt]   [align=left] {$\omega _{2}$};
\draw (409.48,183.38) node [anchor=north west][inner sep=0.75pt]   [align=left] {$M_{\tau^{\varepsilon}}^{\varepsilon}$};
\draw (510.65,183.4) node [anchor=north west][inner sep=0.75pt]   [align=left] {$M_{\omega}$};
\draw (272.26,105.84) node [anchor=north west][inner sep=0.75pt]   [align=left] {$\mapsto$};
\draw (339.92,23.62) node [anchor=north west][inner sep=0.75pt]   [align=left] {$\Sigma $};
\end{tikzpicture}
}
\]
Since $[M_{\tau^1}^1]=[M_{\tau^2}^2]$ in $\opn{h}\Bord^{nc}_{\ast}[\Theta_{\ast}^{-1}]$, by Lemma \ref{lem:rel2}, it follows that
\[
[M_{\tau}]=[M_{\omega}]\circ[M_{\tau^1}^1]=[M_{\omega}]\circ[M_{\tau^2}^2]=[M_{\tau'}].
\]
This shows that the map \eqref{eq:4213} is injective, and hence bijective. So we see that the functor
\[
forget:\opn{h}\Bord_{\ast}^{nc}[\Theta^{-1}_{\ast}]\to \Bord_{3,2}
\]
is both fully faithful and essentially surjective, and thus an equivalence.
\end{proof}

\subsection{Unmarked TQFTs}

We can now take the field theory $\mbb{L}_{\msc{D}_{fin}}^{\opn{unit}}$ from Lemma \ref{lem:lrt_unit} and apply homotopy truncation to obtain a symmetric monoidal functor from the unlabeled bordism category
\[
\mbb{L}_{D_{fin}}:\Bord^{nc}_{3,2}\underset{\text{Thm \ref{thm:unmarking}}}\cong\opn{h}\Bord^{nc}_{\ast}[\Theta^{-1}_{\ast}]\overset{\opn{h}\mbb{L}^{\opn{unit}}_{\msc{D}_{fin}}}\longrightarrow \opn{h}\Vect.
\]
We recall that the homotopy category $\opn{h}\Vect$ recovers the discrete derived category of dg vector spaces $\opn{h}\Vect=D(\opn{Vect})$, so that the target for this unmarked field theory is completely familiar.

\begin{theorem}\label{thm:unmarked_lrt}
Consider a finite modular tensor category $A^{\heartsuit}$ and take $D_{fin}=D^b(A^{\heartsuit}_{fin})$. Let $C=C_{A^{\heartsuit}}$ be the canonical coend in $A^{\heartsuit}$. There is a symmetric monoidal functor from the unmarked (anomalous) $3$-dimensional bordism category
\[
\mbb{L}_{D_{fin}}:\Bord_{3,2}^{nc}\to D(\opn{Vect})
\]
which takes the value $\mbb{L}_{D_{fin}}(\Sigma)\cong \opn{RHom}_{A^{\heartsuit}}(C^{\ot g},\1)$ on any connected genus $g$ surface $\Sigma$.
\end{theorem}

\begin{remark}
When $A^{\heartsuit}$ is semisimple, the theory $\mbb{L}_{D_{fin}}$ already takes values in the full subcategory $\opn{Vect}\subseteq \Vect$.  Furthermore $\msc{K}(A^{\heartsuit})=\msc{D}(A^{\heartsuit})$ in this case, so that derived LRT is obtained directly from discrete LRT via localization. From this fact, and the fact that discrete LRT restricted to $\Bord^{nc}_{A^{\heartsuit}_{fin}}$ recovers the theory from \cite{derenzietal23}, we see that $\mbb{L}_{D_{fin}}$ is just the usual Reshetikhin-Turaev TQFT for $A^{\heartsuit}$ in the semisimple case.
\end{remark}

We conjecture that the $1$-categorical theory which we've constructed above is, in fact, a truncation of an unmarked theory which exists at the $\infty$-categorical level.

\begin{conjecture}\label{conj:0}
The truncated theory $\mbb{L}_{D_{fin}}$ from Theorem \ref{thm:unmarked_lrt} lifts to a symmetric monoidal functor $\mbb{L}_{\msc{D}_{fin}}:\Bord_{3,2}^{nc}\to \Vect$.
\end{conjecture}

We outline a proof which uses the hypothetical skein theoretic bordism category from Section \ref{sect:skein_bordisms}.

\begin{proof}[Idea of proof]
In the (conjectural) skein theoretic realization of the bordism category we'll have a field theory $L_{\msc{D}_{fin}}:\opn{SkBord}^{nc}_{\msc{D}_{fin}}\to \Vect$, by a variant of Theorem \ref{thm:der_lrt}. We then restrict along the unit $\ast\to \msc{D}_{fin}$ to obtain a field theory $L_{\msc{D}_{fin}}|_{\opn{unit}}:\opn{SkBord}^{nc}_{\ast}\to \Vect$. But now, the forgetful functor $\opn{SkBord}^{nc}_{\ast}\to \Bord^{nc}_{3,2}$ induces homotopy equivalences on mapping spaces since the $\infty$-category $\ast$ is (obviously) contractible. Hence this forgetful functor is an equivalence, and we obtain a unique symmetric monoidal functor $\mbb{L}_{\msc{D}_{fin}}:\Bord^{nc}_{3,2}\to \Vect$ which completes a diagram
\[
\xymatrixcolsep{10mm}
\xymatrix{
\opn{SkBord}^{nc}_{\ast}\ar[r]^{\opn{unit}}\ar[dr]_{forget} & \opn{SkBord}^{nc}_{\msc{D}_{fin}}\ar[r]^{L_{\msc{D}_{fin}}} & \Vect\\
	&\Bord^{nc}_{3,2}\ar[ur]_{\mbb{L}_{\msc{D}_{fin}}} & .
}
\]
\end{proof}

\subsection{Mapping class group actions in dimension $2$}
\label{sect:mcg}

If we consider any nonempty surface $\Sigma$ (with chosen Lagrangian), we have the non-anomalous bordism category $\mfk{Bord}_{3,2}$ in which the oriented diffeomorphisms of the surface generate a subgroup in the automorphism group
\[
\opn{MCG}(\Sigma)\ \subseteq\ \opn{Aut}_{\mfk{Bord}}(\Sigma),
\]
by way of mapping cylinders.  We take the pullback under the forgetful functor from $\Bord^{nc}_{3,2}$ to obtain a subgroup of anomalous automorphisms
\[
\xymatrix{
\tilde{\opn{MCG}}(\Sigma)\ar[r]\ar[d] & \opn{Aut}_{\Bord}(\Sigma)\ar[d]\\
\opn{MCG}(\Sigma)\ar[r] & \opn{Aut}_{\mfk{Bord}}(\Sigma).
}
\]
The group $\tilde{\opn{MCG}}(\Sigma)$ is a central extension of $\opn{MCG}(\Sigma)$ by $\mbb{Z}$ which we call the projective mapping class group.  The following is an immediate consequence of Theorem \ref{thm:unmarked_lrt}.

\begin{corollary}\label{cor:mcg}
For each connected genus $g$ surface $\Sigma$, the derived linear maps $\opn{RHom}_{A^{\heartsuit}}(C^{\ot g},\1)$ carry the natural structure of a representation for the projective mapping class group,
\[
\psi:\tilde{\opn{MCG}}(\Sigma)\to \opn{Aut}_{D(\opn{Vect})}\big(\opn{RHom}_{A^{\heartsuit}}(C^{\ot g},\1)\big).
\]
\end{corollary}

One can compute these mapping class group actions directly by taking resolutions of the unit. The process is exactly as outlined in Section \ref{sect:calc_3mfld} below, where we discuss similar calculations for $3$-manifold invariants.

\begin{remark}
As relayed in the introduction, mapping class group actions on cohomology were already observed in works of Lentner, Mierach, Schweigert, and Sommerh\"auser \cite{lentneretal23}, and Schweigert and Woike \cite{schweigertwoike21}. We have not checked, but we fully expect that our mapping class group actions agree with those of \cite{lentneretal23,schweigertwoike21}. Let us also note that the mapping class group actions from \cite{schweigertwoike21} occur at the level of the $\infty$-category $\Vect$, and so provide an advancement of the particular actions obtained in Corollary \ref{cor:mcg}. This $\Vect$-level action also suggests the existence of a lifted theory as in Conjecture \ref{conj:0}.
\end{remark}

\section{Conjectures}
\label{sect:conjectures}

This paper has focused exclusively on formalizations of TQFTs for derived categories of modular tensor ($1$-)categories. However, we are fundamentally interested in a space of homotopically flavored TQFTs for ``modular $\infty$-categories".  We approach this general topic through a series of conjectures.
\par

In addition to our interest in homotopical Reshetikhin-Turaev type theories, we have a special investment in \emph{couplings} of topological field theories to moduli of local systems. This kind of coupling is well-established in the physics literature, and we discuss the topic in some detail in Sections \ref{sect:locsys} and \ref{sect:locsys_formal} below.

\subsection{$3$-d TQFTs from modular $\infty$-categories}

\begin{conjecture}[cf.\ Conjecture \ref{conj:0}]\label{conj:1}
For any modular $\infty$-category $\msc{A}$, appropriately defined, there is an associated $3$-dimensional TQFT
\[
L_{\msc{A}}:\Bord^{nc}_{\msc{A}}\to \Vect_{\pf}
\]
whose state spaces, for any connected genus $g$ surface $\Sigma$, admit a natural identification $L_{\msc{A}}(\Sigma_x)\overset{\sim}\to \underline{\Maps}_{\msc{A}}(C^{\ot g},m_fx)$.  Furthermore, this TQFT induces an unmarked theory after labeling by the unit
\[
\mbb{L}_{\msc{A}_{fin}}:\Bord^{nc}_{3,2}\to \Vect.
\]
\end{conjecture}

There are two main issues here. The first is in the definition of a modular $\infty$-category. Again, we might say succinctly that such a category is a $2$-dimensional TQFT which satisfies certain non-degeneracy conditions \cite{mullerwoike23,steinebrunner}. Examples should include derived categories of finite modular tensor ($1$-)categories, and derived categories of sheaves on compact symplectic varieties, possibly quantized \cite[Section 9.2]{robertswillerton10}.

\begin{remark}
Though the $\infty$-category of sheaves $\msc{D}(X)$ on a projective symplectic variety $X$ has a strong $t$-structure, the $R$-matrix for the quantization \cite{robertswillerton10} should not respect the $t$-structure.  So we do not think of $\msc{D}(X)$ as admitting a $t$-structure when considered as a modular $\infty$-category.
\end{remark}

In the derived setting, we propose specifically that there is some way to construct our theory $L_{\msc{D}}$ for $\msc{D}=\msc{D}(A^{\heartsuit})$ in a manner which operates directly at the level of $\infty$-categories. As a related point, we conjecture that the derived theories $L_{\msc{D}}$ are natural in functors between $\infty$-categories.

\begin{conjecture}\label{conj:2}
Consider finite modular tensor categories $A^{\heartsuit}$ and $B^{\heartsuit}$ with respective derived $\infty$-categories $\msc{D}=\msc{D}(A^{\heartsuit})$ and $\msc{D}'=\msc{D}(B^{\heartsuit})$. Any map of balanced monoidal $\infty$-categories $F:\msc{D}\to \msc{D}'$, i.e.\ $fr\opn{E}_2$-monoidal $\infty$-categories, induces a symmetric monoidal transformation
\[
L_F:L_{\msc{D}}\to L_{\msc{D}'}\circ \Bord^{nc}_F
\]
which, at the level of state spaces, recovers the natural maps
\[
\underline{\Maps}_{\msc{D}}\big(C^{\ot g},m_fx\big)\to \underline{\Maps}_{\msc{D}'}\big({C'}^{\ot g},m_fF(x)\big)
\]
induced by the functor $F$ and a natural map between canonical algebras $C'\to F(C)$.  In particular, the field theory $L_{\msc{D}}$ is a derived invariant.
\end{conjecture}

For a functor $F$ which is not $t$-exact, establishing such naturality would presumably require an understanding of the theory $L_{\msc{D}}$ directly at the homotopical level.

\begin{remark}
We expect such non-$t$-exact functors to arise from Kazhdan-Lusztig type equivalences at bad levels.  See for example \cite[Theorem 7.13]{mcraeyang23}.
\end{remark}

\subsection{Extension down to $S^1$}

The TQFT $L_{\msc{D}}$ is paired with a ``conjugate" theory $L'_{\msc{D}}:\Bord^{nc}_{\msc{D}}\to \Vect$ whose state spaces are given by the continuous duals
\[
L'_{\msc{D}}(\Sigma_x)\overset{\sim}\to \underline{\Maps}_{\msc{D}}(C^{\ot g}\ot m_f(x),\1)^{\ast},
\]
where the $nc$ superscript in this case indicates a constraint on \emph{incomoing} boundaries and labels come from a $fr\Disk$-category associated to the continuous duality $-^{\ast}:\msc{D}(A^{\heartsuit}_{\pf})^{op}\overset{\sim}\to\msc{D}(A^{\heartsuit})$.  (Compare with the original theory of \cite{derenzietal23}, \cite{lagiotis}, and \cite[proof of Theorem 7.11]{czenkynegron}.) Similarly, supposing that the theories $L_{\msc{A}}$ from Conjecture \ref{conj:1} actually exist, we have the conjugate theory $L'_{\msc{A}}$.

\begin{conjecture}\label{conj:3}
Let $\msc{A}$ be a modular $\infty$-category. The conjugate theory $L'_{\msc{A}}:\Bord^{nc}_{3,2}\to \Vect$ extends down to a symmetric monoidal functor of $(\infty,2)$-categories
\[
\mfk{L}'_{\msc{A}}:\Bord^{nc}_{3,2,1}\to \msc{P}\!r_{k}
\]
which takes values in the $(\infty,2)$-category of presentable $\Vect$-module categories, and for which $\mfk{L}'_{\msc{A}}(S^1)=\msc{A}$.
\end{conjecture}

In the case of a derived $\infty$-category $\msc{A}=\msc{D}(A^{\heartsuit})$, one can sketch the construction of ``most of" an extension for $L'_{\msc{D}}$--at least in principle--by following work of De Renzi in the abelian setting \cite{derenzi21} and applying various completion processes \cite[Proposition 5.5.8.22]{lurie09} \cite[Corollary 1.4.4.5]{ha}.  (See also \cite{bartlettetal,lagiotis}.) We omit further details, for now.

In the case of a derived $\infty$-category $\msc{D}(A^{\heartsuit})$, an obvious sub-conjecture here is that $\msc{D}(A^{\heartsuit})$ appears as the circle value in a $2$-dimensional TQFT.

\begin{conjecture}\label{conj:3.5}
Any finite modular $\infty$-category $\msc{D}(A^{\heartsuit})$ appears, along with its $\opn{E}_2$-structure, as the circle value in a $2$-dimensional field theory
\[
\bar{\mfk{L}}'_{\msc{D}}:\Bord_{2,1}\to \msc{P}\!r_k.
\]
\end{conjecture}

We note that Conjecture \ref{conj:3.5} is not self-evident, though it may be accessible.

\begin{remark}
If we decide to conceptualize modular $\infty$-categories explicitly as non-degenerate $2$-dimensional field theories, then Conjectures \ref{conj:1} and \ref{conj:3} together assert that an extension \emph{up} from dimension $2$ to non-compact manifolds in dimension $3$ occurs if and only if the given $2$-dimensional field theory satisfies a finite number of checkable properties, and in particular that no additional structures are required. From the perspective of the cobordism hypothesis, say, such an assertion seems immanently reasonable.
\end{remark}

\subsection{LRT bounding a $4$-dimensional theory?}

In the semisimple setting we have the $4$-dimensional Crane-Yetter theory which is associated, modulo some details, to a given modular fusion category $A^{\heartsuit}$.  This theory is invertible, and so supports anomalous boundary theories \cite{johnsonscheimbauer17}.  It is shown in work of Ha\"ioun \cite{haioun}--following proposals of Walker, Freed and Teleman \cite{walker1,walker2,freedteleman14}--that one can obtain the (semisimple) Reshetikhin-Turaev theory for a modular fusion category $A^{\heartsuit}$ as a boundary theory to $4$-dimensional Crane-Yetter. One similarly obtains the non-semisimple abelian theory $L_{A^{\heartsuit}_{fin}}$ as a boundary to so-called non-semisimple Crane-Yetter \cite{brochieretal21,haioun}.
\par

In the derived setting, considering $\msc{D}=\msc{D}(A^{\heartsuit})$ for modular $A^{\heartsuit}$, the situation is rather unclear. We expect the $4$-dimensional theory obtained from $\msc{D}$ to be quite a wild object, and in particular non-invertible.  Specifically, the value on $S^2$ is presumably identified with the $E_3$-center for $\msc{D}$, which should be exceedingly non-trivial in general.  For example, in the case of a quantum group $A^{\heartsuit}=(\Rep_qG)_{\opn{small}}$, all expectations are that the $E_3$-center should be an intricate category of sheaves on the cone $\check{\mcl{N}}$ of nilpotent elements in the Lie algebra for the Langlands dual group. So the non-semisimple (derived) and semisimple settings seem quite dissimilar, and one should think about whether the semisimple boundary intuition is actually applicable when one works homotopically.

\begin{question}\label{quest:4}
Is there an extended $4$-dimensional theory $T_{\msc{D}}$ which recovers $L_{\msc{D}}$, or $L'_{\msc{D}}$, as a $3$-dimensional boundary?
\end{question}

As a slightly different, but maybe more fundamental set of questions, we can investigate the $4$-dimensional theory for $\msc{D}$ outright.

\begin{question}\label{quest:4.5}
To what extent can we describe the $4$-dimensional theory associated to $\msc{D}=\msc{D}(A^{\heartsuit})$, where we consider the derived $\infty$-category as an object in the $(\infty,4)$-category of $E_2$-monoidal $\infty$-categories and approach the issue--either philosophically or materially--via the cobordism hypothesis \cite{haugseng17,johnsonscheimbauer17,brochierjordansnyder21,brochieretal21}.
\par

In particular, does the fully local topological theory for $\msc{D}$ extend \emph{up} to dimension $3$ (cf.\ \cite[Theorem 1.2]{brochieretal21})? Is it, furthermore, a non-compact $4$-dimensional theory? Are its values on certain $3$-manifolds definitively non-invertible? Is there an explicit skein theory which describes its values on $3$-manifolds (cf.\ Section \ref{sect:skein_bordisms})? Can we describe factorization homology for $\msc{D}$ via a skein theory? Etc. 
\end{question}

\subsection{Coupling to local systems}
\label{sect:locsys}

In addition to the points discussed above, the other (from our perspective) tangible issue which one is required to think about is the coupling of $3$-dimensional field theories to local systems. The following conjecture is, as far as we can tell, a ``fact" in the physics literature. See for example \cite{marcus95,gaiotto19}.

\begin{conjecture}\label{conj:5}
Let $\msc{A}$ be a modular $\infty$-category, and suppose a connected algebraic group $H$ acts rationally on $\msc{A}$ via automorphisms which preserve the modular structure. Suppose also that a certain extension problem for $\msc{A}$ has a solution (see Conjecture \ref{conj:6}). Then the unmarked TQFT $\mbb{L}_{\msc{A}}:\Bord_{3,2}^{nc}\to \Vect$ deforms naturally along $H$-local systems. In particular, the state spaces $\mbb{L}_{\msc{A}}(\Sigma)$ appear as the trivial fiber of a naturally occurring dg sheaf $\mcl{L}^H_{\msc{A}}(\Sigma)$ over $\opn{LocSys}_H(\Sigma)$. 
\end{conjecture}

In the case of a derived $\infty$-category, the extent to which $\mbb{L}_{\msc{A}}$ deforms as a fully functional, mathematically formal, topological field theory is a bit subtle. So, our main ambition here is to ``look at'' Conjecture \ref{conj:5} through the lense of some examples, and to extract some more explicit claims.

Let us consider the case of quantum group representations.  Here we can think about the state spaces from the perspective of $2$-dimensional conformal field theory, via a bulk/boundary correspondence, at which point we conceptualize $2$-d CFTs through the theory of vertex tensor categories.
\par

From the perspective of the ``logarithmic KL conjecture" we anticipate an equivalence of modular tensor categories $(\Rep_qG)_{\opn{small}}\overset{\sim}\to \mcl{W}_{\kappa}(G)\text{-mod}$ where $\mcl{W}_{\kappa}(G)$ is the logarithmic $W$-algebra associated to $G$ at a prescribed level \cite{lentner21,sugimoto21,lisugimoto}. This equivalence should be $\check{G}$-equivariant for the action of the Langlands dual group $\check{G}$ on $\mcl{W}_{\kappa}(G)\text{-mod}$ induced by the natural action of $\check{G}$ on $\mcl{W}_{\kappa}(G)$ by VOA automorphisms. These $\mcl{W}_{\kappa}(G)$ should furthermore be the associated boundary algebras for a $3$-dimensional physical theory associated to $G$ \cite[Section 6.2.2]{creutzigetal24}, and it is proposed directly in \cite[Section 2.1]{gaiotto19}, for example, that this $\check{G}$-action is responsible for the coupling of the topological theory $L_{\msc{D}}$ to local systems (see also \cite{feiginlentner}).
\par

In categorical terms, one expects these kinds of deformed TQFT to be regulated by an ``$H$-crossed braided tensor $\infty$-category" $\tilde{\msc{A}}_H$ which extends $\msc{A}$, as in Turaev's HTQFTs \cite{turaev00}.

\begin{conjecture}\label{conj:6}
Given an $H$-crossed braided extension $\tilde{\msc{A}}_H$ of a modular tensor $\infty$-category $\msc{A}$, the associated topological field theory $L_{\msc{A}}$ from Conjecture \ref{conj:1} deforms along local systems to produce a mathematically formal field theory $\mcl{L}_{\msc{A}}^H$ which is coupled to local systems.
\end{conjecture}

Returning to our particular logarithmic $W$-algebra, the $\check{G}$-extension of the category of $\mcl{W}_{\kappa}(G)$-modules
\[
\mcl{W}_{\kappa}(G)\text{-mod}\to \mcl{W}_{\kappa}(G)\text{-mod}_{\check{G}}^{\opn{tw}}
\]
should be provided by the category of $\check{G}$-twists modules for $\mcl{W}(G)$. Working directly through the representation theory of the small quantum group, rather than through a Kazhdan-Lusztig equivalence, we believe that the analogous $\check{G}$-crossed extension
\[
(\Rep_q\check{G})_{\opn{small}}\to \opn{DK}_q(G)
\]
is realized by a variant of the De Concini-Kac module category. In this modified De Concini-Kac category the Frobenius center in $U^{\opn{DK}}_q(\mfk{g})$ is replaced by functons on the Langlands--rather than Frobenius--dual group (cf.\ \cite{creutzigetal24,losevtsymbaliukvu}). We place special emphasis on this case, as we expect the category $\opn{DK}_q(G)$ to be of interest in both quantum algebra and geometric representation theory.

\begin{remark}
If we are lucky, the desired category $\opn{DK}_q(G)$ can be constructed by applying some act of tensor categorical hoodoo to the recent ``corrected" De Concini-Kac quantum enveloping algebra of Losev, Tsymbaliuk, and Vu \cite{losevtsymbaliukvu}. (Compare with the corrected versions of small quantum groups at even order parameters \cite{creutziggainutdinovrunkel20,negron21,creutzigrupert22,negron26}.) Already they have produced stability for the De Concini-Kac algebra across all orders in the quantum parameter, and their $\check{G}$-integrable subalgebra $U^{fin}_q$ admits a central embedding $\msc{O}(\check{G})\to U^{fin}_q$ from the algebra of functions on the dual group.  Hence the category of $U^{fin}_q$-modules is naturally a sheaf over $\check{G}$. What we are missing is the appropriate crossed tensor structure, and the relation to small quantum groups.
\end{remark}

\begin{remark}
The $\opn{LocSys}_{\check{G}}$-coupled theory for derived quantum group representations lies, from a certain physics perspective, on the spectral side of a physical Langlands duality. The automorphic side of this duality is a corresponding theory which couples to $\opn{Bun}_{G}$. See \cite{kapustinwitten06,gaiotto18,gaiotto19}.
\end{remark}

\begin{remark}
We expect the theory from Conjecture \ref{conj:5} to accept certain types of markings, though the inclusion of such markings complicates the story.
\end{remark}

\subsection{What is a topological theory coupled to local systems, really?}
\label{sect:locsys_formal}

While we cannot make a definitive statement in this regard, let us say a few words towards the precise nature of the theory $\mcl{L}_{\msc{A}}^H$ from Conjecture \ref{conj:6}.
\par

In considering a TQFT ``coupled" to flat connections $\mcl{L}_\msc{A}^H$ we anticipate the existence of a symmetric monoidal functor
\[
\mbb{L}^H_{\msc{A}}:H\text{-}\Bord_{3,2}^{nc}\to \Vect
\]
from an $H$-bordism category consisting of surfaces $\Sigma$ equipped with a map of stacks $\lambda_{\Sigma}:\Sigma\to BH$ and bordisms $M$ with $\lambda_M:M\to BH$, as in \cite{turaev00}. Here the classifying stack is the sheaf of spaces with sections $(BH)(Y)=B(H(Y))$, over either an algebraic or analytic site, and each manifold is treated as the constant sheaf.  Equivalently, the objects and morphisms come equipped with a point on the mapping stack $\lambda_X:\opn{Spec}(k)\to \opn{LocSys}_H(X)=\Maps(X,BH)$.
\par

In top dimension the resulting functions $\mbb{L}_{\msc{A}}^H(M_{\lambda}\setminus\text{a ball}):k\to k\subseteq \mbb{L}_{\msc{A}}(S^2)$ should vary smoothly over $\opn{LocSys}_H(M)$, and in dimension $2$ the derived vector spaces $\mbb{L}^H_{\msc{A}}(\Sigma_{\lambda})$ should assemble into a dg sheaf $\mcl{L}^H_{\msc{A}}(\Sigma)$ over local systems.  As for the sheaf structure, for commutative algebras $R$ in $\Vect$ we anticipate functorially varying theories
\[
\mbb{L}^H_{\msc{A}}(R):H\text{-}\Bord_{3,2}^{nc}(R)\to \msc{M}\!od(R)
\]
for manifolds with specified $R$-families of local systems $\lambda_X:\opn{Spec}(R)\to \opn{LocSys}_H(X)$, which collectively assemble into a ``sheaf of field theories" $\mcl{L}_{\msc{A}}^H$.

To come back to the real world, we return the quantum group case $A^{\heartsuit}=(\Rep_qG)_{\opn{small}}$, $\msc{D}=\msc{D}(A^{\heartsuit})$.  If we restrict our attention only to the \emph{torus} $\check{T}\subseteq \check{G}$ we obtain a theory with \emph{abelian} connections
\[
\mbb{L}^{\check{T}}_{\msc{D}}:\check{T}\text{-}\Bord_{3,2}^{nc}\to \Vect.
\]
By all indications, this $\check{T}$-coupled theory \emph{is} the derived (and sheafified) variant of the CGP theory for the unrolled quantum group \cite{blanchetetal16,derenzigeermirand20}
\par

We note also that deformed theories over local systems were produced in \emph{$4$-dimensions} for non-semisimple (non-derived) Crane-Yetter theories in work of Kinnear \cite{kinnear}. Here, for the quantum group, one considers the $5$-dimensional theory associated to $\Rep\check{G}$ and takes $\Rep_qG$ as an invertible map in the $5$-category of $E_3$-monoidal $1$-categories. This situation enriches the $4$-dimensional theory for the fiber $\opn{Vect}\ot_{\Rep\check{G}}\Rep_qG=(\Rep_qG)_{\opn{small}}$ in sheaves over local systems.

It is possible that an analogous approach can be taken in dimension $3$ as well, though we don't currently see a clear line of argumentation in this direction.

\section{Post-credit scenes in $3$ dimension}
\label{sect:dimension_3}

From the author's perspective at least, this paper is done. We include here some post-credit materials, where we discuss sporadic outputs of our theory in dimension $3$.
\par

At the most immediate level, we obtain diffeomorphism invariants of $3$-manifolds. However, as we explain below, these $3$-manifold invariants collapse to invariants which are already calculable at the abelian level.  In order to see higher information in dimension $3$, one can performing higher genus operation to a closed $3$-manifold.  For one example, we explain how $\mbb{L}_{D_{fin}}$ can be used to produce framed knot invariants by excising tori from a $3$-sphere.

\subsection{Invariants from bored $3$-manifold}

We fix a finite modular tensor category $A^{\heartsuit}$ with finite derived category $D_{fin}=D(A^{\heartsuit}_{fin})$.  We consider the unmarked TQFT
\[
\mbb{L}_{D_{fin}}:\Bord_{3,2}^{nc}\to D(\opn{Vect})
\]
from Theorem \ref{thm:unmarked_lrt}, and fix a bordism $\mbf{\Sigma}:\emptyset\to \Sigma$ in $\Bord_{3,2}$ for a surface $\Sigma$.  So, this is just a chosen ``filling" of $\Sigma$.  Let $\bar{\mbf{\Sigma}}:\bar{\Sigma}\to \emptyset$ denote the orientation reversal.
\par

For any closed $3$-manifold with a pair of oriented embeddings
\begin{equation}\label{eq:mcheck}
\xymatrixrowsep{0mm}
\xymatrix{
& \mbf{\Sigma}\ar[ddr] & & \bar{\mbf{\Sigma}}\ar[ddl]\\
\check{M}=\\
&	& M
}
\end{equation}
we remove the interiors, and adopt signature defect $0$, to obtain a non-anomalous bordism $\Sigma\to \Sigma$.  By an abuse of notation we denote this bordism by $\check{M}:\Sigma\to \Sigma$ as well.

We note that the bordism $\check{M}:\Sigma\to \Sigma$ associated to such a configuration \eqref{eq:mcheck} only depends on the isotopy class of the embedding $\mbf{\Sigma}\amalg\bar{\mbf{\Sigma}}\to M$, up to a diffeomorphism.  Hence this bordism, considered as a map in $\Bord^{nc}_{3,2}$, is stable under isotopies for the input data \eqref{eq:mcheck}.  This follows by the isotopy extension theorem \cite{palais60,lima63}.  So, for example, when $M$ is connected and $\Sigma=S^2$ the bordism $\check{M}$ is in fact independent of the choice of embeddings completely. We record this point.

\begin{lemma}\label{lem:iet}
Given a configuration as in \eqref{eq:mcheck}, the resulting map $\check{M}:\Sigma\to \Sigma'$ in $\Bord^{nc}_{3,2}$ only depends on the isotopy class of the embedding $\mbf{\Sigma}\amalg\bar{\mbf{\Sigma}}\to M$.
\end{lemma}

Taking $\check{M}$ as above, we have the induced morphism $\mbb{L}_{D_{fin}}(\check{M}):\mbb{L}_{D_{fin}}(\Sigma)\to \mbb{L}_{D_{fin}}(\Sigma)$ and we take cohomology to get a graded endomorphism
\begin{equation}\label{eq:4601}
H^{\ast}\mbb{L}_{D_{fin}}(\check{M}):H^{\ast}\mbb{L}_{D_{fin}}(\Sigma)\to H^{\ast}\mbb{L}_{D_{fin}}(\Sigma).
\end{equation}
In the event that $\Sigma$ is connected and of genus $g$, this cohomology is identified with extensions for $A^{\heartsuit}$, as in Theorem \ref{thm:unmarked_lrt}, and we equivalently obtain an endomorphism of cohomology
\[
H^{\ast}\mbb{L}_{D_{fin}}(\check{M}):\Ext_{A^{\heartsuit}}^{\ast}(C^{\ot g},\1)\to \Ext_{A^{\heartsuit}}^{\ast}(C^{\ot g},\1).
\]

Now, if we recall how the theory $\mbb{L}_{D_{fin}}$ is produced, we understand that it is only determined up to a unique natural isomorphism.  So we should only think of the endomorphism \eqref{eq:4601} as determined up to conjugation.  Hence we consider the following.

\begin{definition}
Fix a finite modular tensor category $A^{\heartsuit}$.  For any closed $3$-manifold with embeddings $\check{M}:\mbf{\Sigma}\amalg \bar{\mbf{\Sigma}}\to M$ from fillings of a given surface $\Sigma$, we define
\[
\opn{Inv}(A^{\heartsuit}|\check{M}):=\sum_{n=0}^{\infty}\det\left(H^n\mbb{L}_{D_{fin}}(\check{M})-X\cdot id\right)t^n\ \in\ \mcl{O}(\mbb{A}^1_k)[\![t]\!],
\]
where $\mcl{O}(\mbb{A}^1_k)=k[X]$.
\end{definition}

The following is apparent.

\begin{proposition}
For each excised $3$-manifold $\check{M}$ as above, the assignment
\[
\opn{Inv}(-|\check{M}):\big\{\text{finite modular $\ot$-cats}\big\}\to \mcl{O}(\mbb{A}^1_k)[\![t]\!]
\]
is stable on equivalence classes, and hence defines an invariant of finite modular tensor categories. 
\end{proposition}

\begin{proof}
If we have an equivalence of modular tensor categories $A^{\heartsuit}\overset{\sim}\to B^{\heartsuit}$ then we obtain an equivalence between the discrete LRT theories, in particular we have a diagram
\[
\xymatrix{
\Bord_{A}^{nc}\ar[dr]_{\sim}\ar[rr]^{L_A} & & \Vect_{\pf}\\
	& \Bord_{B}^{nc}\ar[ur]_{L_B}
}
\]
in $\SM_{\infty}$, and subsequent identification between the associated topological theories $\mbb{L}_{D(A^{\heartsuit})_{fin}}\overset{\sim}\to \mbb{L}_{D(B^{\heartsuit})_{fin}}$.
\end{proof}

Of course, we conjecture that each operation $\opn{Inv}(-|\check{M})$ is in fact a \emph{derived} invariant for modular tensor categories. (See Conjecture \ref{conj:2}.)
\par

We claim that the operation $\opn{Inv}(A^{\heartsuit}|-)$ also serves as a topological invariant, at least to some extent. Below we record some basic information concerning the endomorphisms $\mbb{L}_{D_{fin}}(\check{M})$, discuss the behaviors of $\opn{Inv}(A^{\heartsuit}|-)$ in low genus, then explain how one can directly calculate this invariant via injective resolutions.

\subsection{Power series invariants and $\mbb{L}_{D_{fin}}(S^2)$-linearity}
\label{sect:3d_lin}

By basic topological considerations we understand that the sphere value $\mbb{L}_{D_{fin}}(S^2)$ is a commutative algebra and that each state space $\mbb{L}_{D_{fin}}(\Sigma)$ for a connected surface is a module over $\mbb{L}_{D_{fin}}(S^2)$ in $D(\opn{Vect})$.  The multiplication on the sphere is obtained by boring out the interiors of two small balls in a larger ball, and the action on $\mbb{L}_{D_{fin}}(\Sigma)$ is obtained by boring out the interior of a ball within the bordism $\Sigma\times [0,1]$.  Under this natural module structure, each bordism between connected surfaces $N:\Sigma\to \Sigma'$ evaluates to a map of $\mbb{L}_{D_{fin}}(S^2)$-modules under the field theory $\mbb{L}_{D_{fin}}$.
\par

Though we omit the details, we claim that under the identifications
\[
\mbb{L}_{D_{fin}}(S^2)\cong \opn{RHom}_{A^{\heartsuit}}(\1,\1)\ \ \text{and}\ \ \mbb{L}_{D_{fin}}(\Sigma)\cong \opn{RHom}_{A^{\heartsuit}}(C^{\ot g},\1)
\]
the natural algebra structure on $\mbb{L}_{D_{fin}}(S^2)$ agrees with the standard Yoneda product on extensions, and the module structure agrees with the natural tensor action of $\opn{RHom}_{A^{\heartsuit}}(\1,\1)$ on $\opn{RHom}_{A^{\heartsuit}}(C^{\ot g},\1)$.

Now, it is a well-established conjecture that, for any finite tensor category, the cohomology $\Ext^{\ast}_{A^{\heartsuit}}(x,\1)$ of an arbitrary rigid object $x$ is a finite module over $\Ext^{\ast}_{A^{\heartsuit}}(\1,\1)$ \cite{friedlandersuslin97,etingofostrik04}. So, in any setting in which one can hope to compute anything, each morphism
\[
H^{\ast}\mbb{L}_{D_{fin}}(\check{M}):H^{\ast}\mbb{L}_{D_{fin}}(\Sigma)\overset{\sim}\to H^{\ast}\mbb{L}_{D_{fin}}(\Sigma)
\]
for connected $\Sigma$ is completely determined by its values in low-degree. As a consequence, we expect the low degree coefficients
\[
\opn{Inv}(A^{\heartsuit}|\check{M})=c_0(X)+c_1(X)t+c_2(X)t^2+\cdots
\]
to carry the most important information. 

\subsection{Manifold invariants in genus $0$}
\label{sect:3d_0}

\begin{lemma}\label{lem:3d_0}
For $\Sigma=S^2$, any choice of a connected closed $3$-manifold $M$, and arbitrary embedding $\check{M}:B^3\amalg\bar{B}^3\to M$, the polynomial
\[
\opn{Inv}(A^{\heartsuit}|M):=\opn{Inv}(A^{\heartsuit}|\check{M})
\]
is a diffeomorphism invariant of $M$.
\end{lemma}

\begin{proof}
This follows from the fact that the resulting bordism $\check{M}:S^2\to S^2$ only depends on the isotopy class of the given embedding (Lemma \ref{lem:iet}).
\end{proof}

Though this invariant is a natural by-product of the unmarked field theory $\mbb{L}_{D_{fin}}$, and so is worth considering in the abstract, it does not actually carry any information from the higher cohomologies $H^{>0}\mbb{L}_{D_{fin}}(S^2)$.

\begin{proposition}\label{prop:inv_zero}
For $\Sigma=S^2$, the scalar term $\lambda_M=\opn{Inv}(A^{\heartsuit}|M)|_{X=t=0}$ determines the entire invariant $\opn{Inv}(A^{\heartsuit}|M)$.  In particular, we have
\[
\opn{Inv}(A^{\heartsuit}|\check{M}):=\sum_{n\geq 0}(\lambda_M-X)^{d_n} t^n
\]
where $d_n=\opn{dim}\Ext^n_{A^{\heartsuit}}(\1,\1)$.
\end{proposition}

\begin{proof}
Since $H^{\ast}\mbb{L}_{D_{fin}}(\check{M})$ is a $H^{\ast}\mbb{L}_{D_{fin}}(S^2)$-module endomorphism of $H^{\ast}\mbb{L}_{D_{fin}}(S^2)$ in the category of graded vector spaces, it is determined by its restriction to degree zero
\[
H^0\mbb{L}_{D_{fin}}(\check{M}):H^0\mbb{L}_{D_{fin}}(S^2)=k\to H^0\mbb{L}_{D_{fin}}(S^2)=k.
\]
This restriction just returns a scalar $\lambda_M$, and we have $H^{\ast}\mbb{L}_{D_{fin}}(\check{M})=\lambda_M\cdot id_{\mbb{L}(S^2)}$.
\end{proof}

We note that the scalar $\lambda_M$ recovers the abelian (renormalized) Lyubashenko invariant from \cite{derenzietal23}.  In particular, we have
\[
\lambda_M=L_{A^{\heartsuit}_{fin}}(M):L_{A^{\heartsuit}_{fin}}(S^2)=k\to L_{A^{\heartsuit}_{fin}}(S^2)=k.
\]
See Lemma \ref{lem:calc} and Section \ref{sect:calc_justif} below.

\subsection{Framed knot invariants in genus $1$}
\label{sect:3d_1}

Take a knot in $3$-space $\nu:S^1\to \mbb{R}^3$, which we may consider as a knot in $S^3$ which avoids the point $\infty$. We suppose the differential $\partial\nu$ is non-vanishing, and assume $\nu$ comes equipped with an orthonormal framing in which the first vector is the differential $\partial\nu/\partial t$. From this framing we deduce an embedding from the solid torus $\nu':\mbf{T}^2\to \mbb{R}^3\subseteq S^3$ where the center of $\mbf{T}^2$ traverses $\nu$. We take another embedding $\mu:\bar{\mbf{T}}^2\to \mbb{R}^3\subseteq S^3$ which wraps around the exterior of $\nu'$ at an arbitrary time $t\in S^1$.  In this way we obtain an assignment $\nu\mapsto \check{S}^3_\nu$ from a framed knot $S^2\to \mbb{R}^3$ to a bordism $\check{S}^3:T^2\to T^2$.

\begin{proposition}\label{prop:inv_knots}
For any finite modular tensor category $A^{\heartsuit}$, the assignment
\[
\opn{Embed}_{fr}(S^1,\mbb{R}^3)\to \mcl{O}(\mbb{A}^1_k)[\![t]\!],\ \ \nu\mapsto \opn{Inv}\left(A^{\heartsuit}|\check{S}^3_\nu\right),
\]
from the space of smooth framed embeddings is continuous.  In particular, $\opn{Inv}\left(A^{\heartsuit}|\check{S}^3_\nu\right)$ is a framed knot invariant.
\end{proposition}

We can consider, as a fairly explicit example, the category $A^{\heartsuit}=\Rep(u_q)$ of small quantum group representations for $\opn{SL}_2$ at an odd order root of unity. We assume, for the sake of simplicity, that $\opn{ord}(q)=3$. Here the Hochschild cohomology $H^{\ast}\mbb{L}_{D_{fin}}(T^2)\cong \Ext^{\ast}_{u_q}(C,\1)$ is generated over the extensions $\Ext^{\ast}_{u_q}(\1,\1)$ in degrees $0$ and $1$, and the self-extensions of the unit are identified with algebraic functions $\Ext^{\ast}_{u_q}(\1,\1)=\mcl{O}(\mcl{N})$ on the cone of nilpotent matrices in $\mfk{sl}_2$.  We have, specifically, four module generators for $\Ext^{\ast}_{u_q}(C,\1)$ in degree $0$ and four generators in degree $1$.
\par

Both cohomologies come equipped with natural $\opn{PSL}_2$-actions and, as a $\opn{PSL}_2$-representation, we have
\begin{equation}\label{eq:5248}
\Ext^{0}_{u_q}(C,\1)=L(0)\otimes_{\mbb{C}}(\mbb{C}^2\oplus \mbb{C}^2)\ \ \text{and}\ \ \Ext^{1}_{u_q}(C,\1)=(L(0)\oplus L(2))\otimes_{\mbb{C}}\mbb{C}^2
\end{equation}
where $\opn{PSL}_2$ acts on the left hand factors \cite[Section 5]{lachowskaqi21}. We claim that the map $H^{\ast}\mbb{L}_{D_{fin}}(\check{S}^3_{\nu})$ is $\opn{PSL}_2$-linear \cite{sajan}, so that $H^1\mbb{L}_{D_{fin}}(\check{S}^3_{\nu})$ is specified by two $2\times 2$ matrices $\Lambda_0$ and $\Lambda_2$ with eigenvalues $\lambda_{0i}$ and $\lambda_{2i}$. Hence we obtain
\[
\opn{Inv}(u_q|\check{S}^3_{\nu})=c_0(X)-\prod_{i=1}^2(X-\lambda_{0i})(X-\lambda_{2i})^3\cdot t+O(t^2).
\] 
The coefficient $c_0(X)$ can be calculated via the abelian theory $L_{A^{\heartsuit}_{fin}}$.

\begin{remark}
Obviously, the (unordered) pairs of scalars $(\lambda_{01},\lambda_{02})$ and $(\lambda_{01},\lambda_{02})$ are themselves framed knot invariants.
\end{remark}

\begin{remark}
The dependence of the invariant $\opn{Inv}(A^{\heartsuit}|\check{S}^3_{\nu})$ on the framing is precisely captured by the $\opn{SL}_2(\mbb{Z})$-action on Hochschild cohomology. Namely, precomposing $H^1\mbb{L}_{D_{fin}}(\check{S}^3_{\nu})$ with the action of the $T$ matrix in $\opn{SL}_2(\mbb{Z})$ produces the endomorphisms
\[
H^1\mbb{L}_{D_{fin}}(\check{S}^3_{\nu})\circ T^n=H^1\mbb{L}_{D_{fin}}(\check{S}^3_{\nu|n})
\]
associated to the new framing on $\nu$ obtained by introducing an $n$-th order twisting around the image of $\nu$. So the framed knot produces a tuple of power series
\[
\opn{Inv}_{unfr}(A^{\heartsuit}|\check{S}^3_{\nu})=\opn{Inv}(A^{\heartsuit}|\check{S}^3_{\nu|-}):\mbb{Z}\to \mcl{O}(\mbb{A}^1_k)[\![t]\!]
\]
which is uniquely determined up to a translation in $\mbb{Z}$. In particular, any asymptotic behaviors of this function produce unframed invariants.
\end{remark}

\begin{remark}
The decomposition of \eqref{eq:5248} should be one of $\opn{PSL}_2(\mbb{C})\times \opn{SL}_2(\mbb{Z})$-representations, with $\opn{PSL}_2$ acting on the left factor and $\opn{SL}_2(\mbb{Z})$ acting on the right. The $\mbb{C}^2$ factor in the first Hochschild cohomology is a twist of the standard representation for $\opn{SL}_2(\mbb{Z})$ by a nontrivial character, as is one of the $\mbb{C}^2$ summands in degree $0$. See \cite{lachowskaqi21}.
\end{remark}

Let us explain, finally, how one can calculate these invariants in practice.

\subsection{Calculating $\opn{Inv}$ via resolutions}
\label{sect:calc_3mfld}

Though we won't calculate any of the above invariants explicitly, let us at least explain \emph{how} one can calculate them.  The process is rather simple.  Let us fix modular $A^{\heartsuit}$ and closed connected $M$ with embeddings $\mbf{\Sigma}\to M\leftarrow \bar{\mbf{\Sigma}}$.
\par

Firt, take your favorite injective resolution $\1\to Q$ of the unit in $A^{\heartsuit}$ and endow the surface $\Sigma$ with a positive disk labeled by $Q$.  Then produce a bordism
\[
\check{M}_{tube}:\Sigma_Q\to \Sigma_Q
\]
in $\Bord_{A}^{nc}$ by connecting the markings via a single tube which traverses the interior of $\check{M}$ arbitrarily.  We then evaluate under discrete LRT to get an endomorphism of cochains
\[
L_{A}(\check{M}_{tube}):\Hom^{\ast}_{A^{\heartsuit}}(C^{\ot g},Q)\to \Hom^{\ast}_{A^{\heartsuit}}(C^{\ot g},Q).
\]
This endomorphism can be calculated, degree-by-degree, via the original abelian LRT theory from \cite{derenzietal23}.  Take cohomology to get an endomorphism
\[
H^{\ast}L_A(\check{M}_{tube}):\Ext^{\ast}_{A^{\heartsuit}}(C^{\ot g},\1)\to \Ext^{\ast}_{A^{\heartsuit}}(C^{\ot g},\1).
\]

\begin{lemma}\label{lem:calc}
The endomorphisms $H^{\ast}L_A(\check{M}_{tube})$ and $H^{\ast}\mbb{L}_{D_{fin}}(\check{M})$ are conjugate.  Hence the coefficient of $t^n$ in $\opn{Inv}(A^{\heartsuit}|\check{M})$ is calculated as the determinant of the matrix $H^nL_A(\check{M}_{tube})-X\cdot id$.
\end{lemma}

For the proof of this result we provide a more leisurely explanation below.

\subsection{Justification for Lemma \ref{lem:calc}}
\label{sect:calc_justif}

Let us explain why this claim holds.  First note that we have a diagram
\[
\xymatrix{
	& \opn{h}\Bord_{\msc{D}_{fin}}^{nc}\ar[dr]^{\opn{h}L_{\msc{D}_{fin}}}\\
\opn{h}\Bord_{\ast}^{nc}\ar[dr]_{forget}\ar[rr]^{\opn{h}L_{\msc{D}_{fin}}unit}\ar[ur]^{unit} & & D(\opn{Vect})\\
	& \Bord_{3,2}^{nc}\ar[ur]_{\mbb{L}_{D_{fin}}}
}
\]
so that we can calculate $\mbb{L}_{D_{fin}}(\check{M})$ by lifting to a marked bordism $\check{M}_{tube}'$ in $\Bord_{\ast}^{nc}$ along the forgetful functor.  This just involves introducing arbitrary markings on the bounding surfaces, though they should be \emph{the same} markings so that $\check{M}_{tube}'$ remains an endomorphism, and arbitrary cylinder embeddings connecting these markings through $\check{M}$.  To simplify the situation we should assume that intersecting with the outgoing boundary reproduces the identity map $\partial_{out}(\check{M}_{tube}')=id_{\mcl{I}}:D\times\mcl{I}\to D\times\mcl{I}$ in $fr\Disk$.
\par

As for the map $\Bord_{\ast}^{nc}\to \Bord_{\msc{D}_{fin}}$, we have the diagram
\[
\xymatrix{
	& \Bord_{\opn{Units}(\msc{D}_{fin})}^{nc}\ar[dr]\ar[dl]_{triv\ kan\ fib}\\
\Bord_{\ast}^{nc}\ar[rr]_{unit} & & \Bord_{\msc{D}_{fin}}^{nc}
}
\]
so that, we can calculate $\opn{unit}(\check{M}_{tube}'):\Delta^1\to \Bord_{\ast}$, up to isomorphism, by first lifting along the forgetful functor from $\Bord_{\opn{Units}(\msc{D}_{fin})}^{nc}$ then projecting down to $\Bord_{\msc{D}_{fin}}^{nc}$.  This can be done explicitly by marking both bounding surfaces by \emph{the same} choice of a unit object $\1_{\msc{D}}$ in $\msc{D}_{fin}$, and all outgoing boundaries by the identity morphism.  This gives us a $\msc{D}_{fin}$-marked bordism $M''_{tube}=\opn{unit}(\check{M}_{tube}')$.  We note that the particular choice of unit does not matter here, but the fact that all markings are the same is important.
\par

Now to calculate the value $L_{\msc{D}_{fin}}(\check{M}''_{tube})$ we have the identification
\[
L_{\msc{D}_{fin}}(\check{M}''_{tube})=L_{\msc{K}}(R\check{M}''_{tube})
\]
where $R$ is the (restricted) right adjoint $R:\Bord_{\msc{D}_{fin}}^{nc}\to \Bord_{\msc{K}}^{nc}$.  The bordism $R\check{M}''_{tube}$ is determined, up to conjugation by an isomorphism, as the $\ast$-marked bordism $\check{M}_{tube}'$ now labeled by an arbitrary injective resolution $Q$ of the unit $\1_{\msc{K}}$, along with identity morphisms at the boundary.  To calculate this morphism we can lift finally along the map $\Bord_{A}^{nc}\to \Bord_{\msc{K}}^{nc}$ to get
\[
L_{\msc{K}}(R\check{M}''_{tube})=L_A(\check{M}_{tube})\ \text{up to conjugation},
\]
where $\check{M}_{tube}$ is obtained from $\check{M}'_{tube}$ by labeling by $Q$ at each marking, and the identity at the outgoing boundaries.  In total we find
\[
H^{\ast}\mbb{L}_{D_{fin}}(\check{M})=H^{\ast}L_A(\check{M}_{tube})\ \text{up to conjugation}.
\]
We can therefore calculate the proposed determinants in the definition of the invariant $\opn{Inv}(A^{\heartsuit}|\check{M})$ via $H^{\ast} L_A(\check{M}_{tube})$.

\appendix

\section{About $fr\E_2$}
\label{sect:more_fre2}

We discuss equivalences between maps in $fr\E_2$, while paying special attention to issues of smoothness.  We then prove Proposition \ref{prop:a_pm}, which claims that the inclusion $Linfr\E_2\to fr\E_2$ is an equivalence, and discuss related issues for bordism categories. Finally, we provide the proofs of Lemmas \ref{lem:full_to_symm} and \ref{lem:structure_transp} from the body of the text. The reader might recall at this point our explicit construction of the $\infty$-operad $fr\E_2$ from Section \ref{sect:fr_e2}.

\subsection{The homotopy category $\opn{h}fr\E_2$}

The following results are well-known to experts. Corollary \ref{cor:h_e2}, in particular, follow by an equivalence between the space of smooth embeddings $\opn{Embed}(D\times I,D)\subseteq \Maps(D\times I,D)$ and the space of framed configurations of points $\opn{Conf}^{fr}_I(D)$.\footnote{See for example this mathoverflow \href{https://mathoverflow.net/questions/490076/embedding-space-of-disjoint-disks-and-framed-configuration-spaces}{answer} form Najib Idrissi.} We give some simple-minded proofs in any case.

\begin{lemma}\label{lem:1430}
Any two oriented smooth embeddings $f_0,f_1:D\to D$ are homotopic via a homotopy $h:|\Delta^1|\times D\to D$ whose evaluation $h_t$ at each time $t\in |\Delta^1|$ is a smooth embedding.  Furthermore, the homotopy $h$ itself can be chosen to be smooth.
\end{lemma}

\begin{proof}
Consider an arbitrary smooth, oriented embedding $f:D\to D$.  Via the homotopy $f(\frac{z}{1+t}):[0,1]\times D\to D$ we may assume that $f$ has image in the interior to $D$. By applying the diffeomorphism $\frac{z}{1-|z|}:D^o\to \mbb{R}^2$ we reduce to considering embeddings $f:D\to \mbb{R}^2$.  We claim that any such $f$ is homotopic to the standard inclusion, and in particular that all oriented embeddings are homotopic.
\par

By translating we assume $f(0)=0$, and consider the differential $T_0f:\mbb{R}^2\to \mbb{R}^2$.  Since $T_0f$ is a matrix with positive determinant, we can find a (smooth) path $A(t):[0,1]\to \opn{GL}_2(\mbb{R})$ with $A(0)=I_2$ and $A(1)=(T_0f)^{-1}$.  Then by composing $A(t)f:[0,1]\times D\to \mbb{R}^2$ we reduce to the case where $T_0f=I_2$.
\par

Now consider the open path
\[
h^o:[-1,0)\times D\to D,\ \ h^o_t(z)=\frac{f(-tz)}{-t}.
\]
We have
\[
\lim_{t\to 0} h^o_t(z)=\frac{\partial f(tz)}{\partial t}|_{t=0}=T_0f(z)=I_2(z)=z.
\]
Hence $h^o(t)$ extends continuously to a homotopy $h:[-1,0]\times D\to D$ with $h_0=id_D$.  To see that $h$ is smooth, we have the alternate expression
\[
f(-tz)=-\int_0^tx\frac{\partial f}{\partial x}(-uz)+y\frac{\partial f}{\partial y}(-uz)du
\]
\[
=-t\int^1_0x\frac{\partial f}{\partial x}(-vtz)+y\frac{\partial f}{\partial y}(-vtz)dv
\]
\[
\Rightarrow\ h(t,z)=-\int_0^1x\frac{\partial f}{\partial x}(-vtz)+y\frac{\partial f}{\partial y}(-vtz)dv.
\]

We now have a sequence of composable smooth homotopies $h^0,\dots,h^m:[0,1]\times D\to D$, for any choice of a smooth embedding $f:D\to D$, with $h^0_0=f$ and $h^m_1=f_{\infty}$.  Here $f_{\infty}$ is fixed to be the preimage of the identity $D\to D\subseteq \mbb{R}^2$ under the diffeomorphism $\frac{z}{1-|z|}:D^o\to \mbb{R}^2$.  After choosing a smooth weakly increasing function $b(t):[0,1]\to [0,1]$ with $b(t)=0$ in a neighborhood around $0$ and $b(t)=1$ in a neighborhood around $1$, we obtain a sequence of smooth homotopies $H^l(t,z)=h^l(b(t),z)$ who compose to a smooth homotopy
\[
H_{f}=H^m\circ\dots\circ H^0:[0,m+1]\times D\to D
\]
between $f$ and $f_{\infty}$.  For an arbitrary pair of smooth embeddings $f_0,f_1:D\to D$ we now construct a smooth homotopy between $f_0$ and $f_1$ via the compositie $H_{f_1}^{-1}\circ H_{f_0}$.
\end{proof}

\begin{proposition}\label{prop:1456}
Suppose two maps $f_0,f_1:D\times\mcl{I}_{f}\to D\times \mcl{J}$ in $fr\E_2$ have the same image $f:\mcl{I}\to \mcl{J}$ in $\mbf{Fin}_{\ast}$. There is a smooth homotopy $h:|\Delta^1|\times D\times \mcl{I}\to D\times \mcl{J}$ with $h_0=f_0$ and $h_1=f_1$, and for which the evaluation $h_t:D\times \mcl{I}_{f}\to D\times\mcl{J}$ at each $t\in |\Delta^1|$ is a smooth embedding.
\end{proposition}

\begin{proof}
It suffices to prove the result when $\mcl{J}$ is a singleton and $f:\mcl{I}\to \{\pm 0\}$ is active, since we can independently resolve the problem for the preimages $D\times \mcl{I}_j$ of the disks in the codomain one-by-one.  Furthermore, after shrinking the images of the disks, then translating if necessary, we may assume each of the images $f_{0}(D_i)$ and $f_1(D_{i'})$ are disjoint in $D$ for each $i$ and $i'$ in $\mcl{I}_{f}$.  Then, by isolating each of the pairs $f_0(D_i)\cup f_1(D_i)$ in disjoint neighborhoods $\tau_i:D\hookrightarrow D$ we reduce to the case where $\mcl{I}$ is a singleton.  This reduces us to the situation considered in Lemma \ref{lem:1430}, where the problem was already resolved.
\end{proof}

\begin{corollary}\label{cor:h_e2}
The forgetful functor $fr\E_2\to \mbf{Fin}_{\ast}$ induces an equivalence from the homotopy truncation $\opn{h}fr\E_2\overset{\sim}\to \mbf{Fin}_{\ast}$.  
\end{corollary}

\begin{proof}
Morphisms in the homotopy category are given by the connected components in the mapping space $\Maps_{fr\E_2}(\mcl{I},\mcl{J})$.  We have the forgetful functor $\Maps_{fr\E_2}(\mcl{I},\mcl{J})\to \Hom_{\mbf{Fin}_{\ast}}(\mcl{I},\mcl{J})$ which separates this mapping space into a disjoint union
\[
\Maps_{fr\E_2}(\mcl{I},\mcl{J})=\amalg_{f}\Maps_{fr\E_2}(\mcl{I},\mcl{J})_{f}
\]
over $f$ in $\Hom_{\mbf{Fin}_{\ast}}(\mcl{I},\mcl{J})$.  So it suffices to show that each space $\Maps_{fr\E_2}(\mcl{I},\mcl{J})_{f}$ is connected.  As two maps in this space $f_0,f_1:D\times\mcl{I}_{f}\to D\times \mcl{J}$ are connected by a path if and only if there is a homotopy $h:|\Delta^1|\times D\times \mcl{I}_{f}\to D\times \mcl{J}$ with $h_0=f_0$ and $h_1=f_1$, the result follows by Proposition \ref{prop:1456}.
\end{proof}

\subsection{Linearly arrangements of disks}

\begin{lemma}\label{lem:linfr_fr}
The inclusion $Linfr\E_2\to fr\E_2$ is an equivalence.
\end{lemma}

\begin{proof}
The inclusion is essentially surjective, as $Linfr\E_2$ contains all objects. Given any pair of finite colored sets $\mcl{I}$ and $\mcl{J}$, we need only show that the inclusion of mapping spaces
\[
\Maps_{Linfr\E_2}(\mcl{I},\mcl{J})\to \Maps_{fr\E_2}(\mcl{I},\mcl{J})
\]
is a homotopy equivalence.  For this we can consider the corresponding simplicial category $Linfr\uE_2$ in $fr\uE_2$ and show that the inclusion
\[
\uHom_{Linfr\uE_2}(\mcl{I},\mcl{J})\to \uHom_{fr\uE_2}(\mcl{I},\mcl{J})
\]
is an equivalence (see \cite[\href{https://kerodon.net/tag/02LN}{02LN}]{kerodon}).  Since the subspace $\uHom_{Linfr\uE_2}(\mcl{I},\mcl{J})$ is full in $\uHom_{fr\uE_2}(\mcl{I},\mcl{J})$ we need only show that the inclusion is surjective on connected components.  However, this is a consequence of Proposition \ref{prop:1456}, which shows that there is only one map in $\uHom_{fr\uE_2}(\mcl{I},\mcl{J})$, up to homotopy, over each map of finite colored sets $\mcl{I}\to \mcl{J}$.
\end{proof}

We similarly have the subcategory $Linfr\Disk^{\ot}$ in $fr\Disk^{\ot}$ which fits into a pullback diagram
\[
\xymatrix{
Linfr\Disk^{\ot}\ar[r]\ar[d] & fr\Disk^{\ot}\ar[d]\\
Linfr\E_2\ar[r] & fr\E_2.
}
\]
Since the disjoint union functor $Linfr\Disk^{\ot}\to fr\E_2$ is an isofibration, the existence of such a pullback diagram tells us that the inclusion $Linfr\Disk^{\ot}\to fr\Disk^{\ot}$ is an equivalence.

\begin{lemma}\label{lem:linfr_disk}
The inclusion $Linfr\Disk^{\ot}\to fr\Disk^{\ot}$ is an equivalence.
\end{lemma}

\subsection{Linear arrangements for bordisms}

\begin{lemma}\label{lem:1939}
Fix two smooth oriented embeddings $f_0,f_1:D\times \mcl{I}\to D$.  There is a smooth homotopy $h:|\Delta^1|\times \opn{Cyl}\times \mcl{I}\to D$, and subsequent homotopy
\[
H:|\Delta^1|\times \opn{Cyl}\times \mcl{I}\to \opn{Cyl},\ \ H(s,z,t,i)\mapsto (h(s,z,t,i),t),
\]
which satisfies the following:
\begin{itemize}
\item There exists $\varepsilon>0$ for which $h(s,z,t,i)=h(s,z,t',i)$ whenever $t$ and $t'$ are simultaneously in the neighborhood $[0,\varepsilon)$, or in $(1-\varepsilon,1]$.\vspace{1mm}
\item $h(0,z,t,i)=h(s,z,0,i)=f_0(z,i)$ for all $s$ and $t$.\vspace{1mm}
\item $h(1,z,1,i)=f_1(z,i)$.
\end{itemize}
\end{lemma}

\begin{proof}
Take $\mu:|\Delta^1|\times D\times \mcl{I}\to D$ a smooth homotopy which has $\mu_0=f_0$, $\mu_1=f_1$, and for which $\mu_t$ is a smooth embedding at all $t$ in $|\Delta^1|$.  Such a homotopy $\mu$ exists by Proposition \ref{prop:1456}. Let $b:[0,1]\to \mathbb{R}$ be a smooth, weakly increasing function for which $b(t)$ is of constant value $0$ in a neighborhood of $0$ and of constant value $1$ in a neighborhood of $1$.  Let $c:[0,1]\to \mathbb{R}$ be the identity $c(s)=s$. Take now
\[
h:|\Delta^1|\times \opn{Cyl}\times \mcl{I}\to D\ \ \text{and}\ \ H:|\Delta^1|\times \opn{Cyl}\times \mcl{I}\to \opn{Cyl}
\]
defined by $h(s,z,t,i)=\mu(c(s)b(t),z,i)$ and $H(s,z,t,i)=\Big(\mu(c(s)b(t),z,i),t\Big)$.
\end{proof}

\begin{proposition}\label{prop:whatever}
Consider any $\ast$-marked surfaces $\Sigma_{\zeta}$ and $\Sigma'_{\eta}$, and bordism $M_f:\Sigma_{\zeta}\to \Sigma'_{\eta}$ in $\Bord_{\ast}$.  Let $f_0=\partial_{out}(M_f):D\times \mcl{I}\to D\times \mcl{J}$ be the corresponding map in $fr\E_2$, and take $t:\mcl{I}\to \mcl{J}$ its image in $\mbf{Fin}_{\ast}$.  For any other map $f_1:D\times\mcl{I}\to \mcl{J}$ in $fr\E_2$ over $t$, there is an isomorphic bordism $M_f\overset{\sim}\to M_{f'}$ in $\Maps_{\Bord_{\ast}}(\Sigma_{\zeta},\Sigma_{\eta})$ with $\partial_{out}(M_{f'})=f_1$.
\end{proposition}

\begin{proof}
Consider the case where $\Sigma'=\Sigma$, $\eta:D\times \mcl{J}\to \Sigma$ are disk which collectively contain the image of $\zeta:D\times \mcl{I}\to \Sigma$, $M=\Sigma\times [0,1]$, and
\[
f:\opn{Cyl}\times \mcl{I}\to M
\]
is constant in the linear factor of $\opn{Cyl}=D\times[0,1]$.  So, $f$ is explicitly the map $f(z,t,i)=f_0(z,i)$.
\par

For each $j$ in $\mcl{J}$ we have the corresponding embedded cylinder $\hat{\eta}_j:\opn{Cyl}_j\to M$, $\hat{\eta}_j(t,z)=(t,\eta(t,z))$, and we intersect with this cylinder to obtain a map
\[
\hat{f}^j:\opn{Cyl}\times \mcl{I}_j\to \opn{Cyl}_j
\]
which is constant in the linear factor, $\hat{f}^j(z,t,i)=f_0(z,i)$, where we now restrict to $z$ in $D\times \mcl{I}_j\subseteq D\times\mcl{I}$.  By Lemma \ref{lem:1939} and Proposition \ref{prop:whatever} we can find, for each $j$, a smooth homotopy
\[
H^j:|\Delta^1|\times \opn{Cyl}\times\mcl{I}_j\to \opn{Cyl}_j\to M
\]
whose terminal map $H^j_1:\opn{Cyl}\times \mcl{I}_j\to M$ intersects with the outgoing boundary disk to reproduce $f_1|_{D\times \mcl{I}_j}$, and which together provide an isomorphism
\[
H:|\Delta^1|\times \opn{Cyl}\times\mcl{I}\to M
\]
between $M_f$ and $M_{f'}$ in $\Maps_{\Bord_{\ast}}(\Sigma_{\zeta},\Sigma_{\eta})$.  (Here $f'=H_1$.)  This homotopy provides, in particular, a diagram of the form
\[
\xymatrix{
 & \Sigma_{\eta}\ar[dr]^{id}\\
\Sigma_{\zeta}\ar[ur]^{M_f}\ar[rr]_{M_{f'}} & & \Sigma_{\eta},
}
\]
which determines such an isomorphism in the mapping space.
\par

In the case of a general bordism $M_f:\Sigma_{\zeta}\to \Sigma'_{\eta}$, we can cut $M_f$ along its outgoing collar to decompose it as a composite $M_f=N_g\circ M_e$ where
\[
N_g:\Sigma'_{\xi}\to \Sigma'_{\eta}
\]
is of the type considered above.  We then have an isomorphism $H:N_g\overset{\sim}\to N_{g'}$ of the prescribed form which composes to provide the desired isomorphism $(M_e)^{\ast}H:M_f\overset{\sim}\to M_{f'}$.
\end{proof}

\begin{corollary}\label{cor:whatever}
For any $fr\Disk$-category $\msc{E}^{\ot}\to fr\Disk^{\ot}$ the projection
\[
Linfr\Disk^{\ot}\times_{fr\Disk^{\ot}}\Bord_{\msc{E}}^{\ot}\to \Bord_{\msc{E}}^{\ot}
\]
is an equivalence of symmetric monoidal $\infty$-categories.
\end{corollary}

\begin{proof}
It suffices to show that the underlying map on $\infty$-categories is an equivalence. Consider the case of the trivial famed pair
\begin{equation}\label{eq:2015}
Linfr\Disk\times_{fr\Disk}\Bord_{\ast}\to \Bord_{\ast}.
\end{equation}
The left-hand $\infty$-category is the subcategory in $\Bord_{\ast}$ consisting of all marked surfaces, and bordisms $M_f:\Sigma_{\zeta}\to \Sigma'_{\eta}$ whose outgoing map $\partial_{out}(M_f)$ lies in $Linfr\Disk$.  Equivalently, the fiber product is the nerve of the simplicial subcategory $Lin\underline{\Bord}_{\ast}$ in $\underline{\Bord}_{\ast}$ containing all objects, and whose mapping complexes
\[
\uHom_{Lin\underline{\Bord}_{\ast}}(\Sigma_{\zeta},\Sigma_{\eta})\ \subseteq\ \uHom_{\underline{\Bord}_{\ast}}(\Sigma_{\zeta},\Sigma_{\eta})
\]
are the full subspaces spanned by all bordisms $M_f$ with $\partial_{out}(M_f)$ in $Linfr\uE_2$.  The corresponding inclusions of mapping complexes are clearly fully faithful, and by (the proof of) Proposition \ref{prop:whatever} they are essentially surjective.  Hence each inclusion is a homotopy equivalence, and the corresponding functor on homotopy coherent nerves is an equivalence \cite[\href{https://kerodon.net/tag/01Q1}{01Q1}]{kerodon}.

For a general $fr\Disk$-category, we have a pullback diagram
\[
\xymatrix{
Linfr\Disk\times_{fr\Disk}\Bord_{\msc{E}}\ar[r]\ar[d] & \Bord_{\msc{E}}\ar[d]\\
Linfr\Disk\times_{fr\Disk}\Bord_{\ast}\ar[r] & \Bord_{\ast}
}
\]
in which the map $\Bord_{\msc{E}}\to \Bord_{\ast}$ is a cocartesian fibration, and hence an isofibration.  Since the bottom projection is an equivalence, it follows that the top projection is an equivalence as well.
\end{proof}

\subsection{Proof of Lemma \ref{lem:structure_transp}}
\label{sect:str_proof}

\begin{proof}[Proof of Lemma \ref{lem:structure_transp}]
By Theorem \ref{thm:hk}, we have an equivalence of $fr\Disk$-categories $\msc{E}^{\ot}\overset{\sim}\to \opn{Env}(\msc{E}^{\odot})$ for some $fr\E_2$-monoidal $\infty$-category $\msc{E}^{\odot}$. So it suffices to prove the result in the case of such an envelope $\msc{E}^{\ot}=\opn{Env}(\msc{E}^{\odot})$.
\par

For $\pi:fr\Disk^{\ot}\to fr\E_2$ the disjoint union functor, and any map $\tilde{f}:s\to t$ in $fr\Disk^{\ot}$ with image $f=\pi(\tilde{f}):\mcl{I}\to \mcl{J}$, we have the diagram
\[
\xymatrixrowsep{5mm}
\xymatrix{
\msc{E}^{\ot}_s\ar[rr]^{\tilde{f}_!}\ar[dd]_(.4){=}\ar[dr]^{\rho_!} & & \msc{E}^{\ot}_t\ar[dd]|{\hole}_(.4)=\ar[dr]^{\rho_!}\\
	& \msc{E}^{\mcl{I}}\ar[rr]_{\opn{induced}} & & \msc{E}^{\mcl{J}}\\
\msc{E}^{\odot}_{\mcl{I}}\ar[rr]_{f_!}\ar[ur]_{\rho_!} & & \msc{E}^{\odot}_{\mcl{J}}\ar[ur]_{\rho_!} & .
}
\]
So it suffices to show that the transport functor $f_!$ along an arbitrary map $f:\mcl{I}\to \mcl{J}$ in $fr\E_2$ is expressible as a composite of projections onto select factors, transition maps, products, and unit maps.
\par

Write out $f$ explicitly as a disk embedding $f:D\times \mcl{I}_f\to D\times \mcl{J}$. We have that $f$ factors as an inert projection $\rho_f:\mcl{I}\to \mcl{I}_f$ composed with a coloring of the identity $f_0:\mcl{I}_f\to \mcl{I}_0$, then composed with an active map $f_1:\mcl{I}_0\to \mcl{J}$ which preserves colors. Transport along the inert projection simply recovers the projection $\msc{E}^{\mcl{I}}\to \msc{E}^{\mcl{I}_f}$, and transport along the colored identity just applies the transition functor $d:\msc{E}^-\to \msc{E}^+$ in the relevant factors. So that we may assume $f$ itself is active and color preserving, i.e.\ is defined by a disk embedding $f:D\times\mcl{I}\to D\times\mcl{J}$ which sends positive disks to positive disks, and negative disks to negative disks.
\par

Now, by Corollary \ref{cor:h_e2}, $f$ is isomotopic to a composite sequence $f\simeq \iota f_n\dots f_0$ where each $f_r:D\times \mcl{I}^r\to D\times \mcl{I}^{r+1}$ is just the identity $D\to D$ on all disks save for two components $D_{i_r}\amalg D_{i'_r}$, and on these two components we apply the linear disk embedding $D\amalg D\to D$, $z\mapsto \frac{1}{3}z\pm \frac{1}{2}$. The map $\iota:D\times \mcl{I}^{n}\to D\times \mcl{J}$ is defined by an inclusion of colored sets $\mcl{I}^n\to \mcl{J}$. 
\par

Composing $\iota$ with the inert projection $\rho_j:\mcl{J}\to \{0\}$ onto an index $j$ in the complement $\mcl{J}\setminus \opn{im}(\mcl{I}^n)$ recovers the zero map $\mcl{I}^ln\to \emptyset\to \{0\}$. From this one sees that the transport functor $\iota_!:\msc{E}^{\mcl{I}^n}\to \msc{E}^{\mcl{J}}$ just inserts the unit $\ast\to \msc{E}^{\pm}$ at all excess factors. Transport along each $f_i$ applies the product $\msc{E}^{\mcl{I}^r}\to \msc{E}^{\mcl{I}^{r+1}}$ on the associated factors. Taking $m_f=\iota_!(f_n)_!\cdots (f_0)_!$, we obtain the desired diagram
\[
\xymatrix{
\msc{E}^{\otimes}_{\mcl{I}}\ar[r]^{f_!}\ar[d] & \msc{E}^{\ot}_{\mcl{J}}\ar[d]\\
\msc{E}^{\mcl{I}}\ar[r]_{m_f} & \msc{E}^{\mcl{J}}.
}
\]
\end{proof}

\subsection{Proof of Lemma \ref{lem:full_to_symm}}
\label{sect:full_to_symm}

\begin{proof}
We first show that the inherited map $\msc{L}^{\ot}\to \Fin_{\ast}$ is a cocartesian fibration. Let $z$ be an object in $\msc{L}^{\ot}$ with images $t=(\mcl{I},\mu)$ and $K$ in $fr\Disk^{\ot}$ and $\opn{Fin}_{\ast}$, respectively. Consider a map $\bar{f}:K\to L$ in $\opn{Fin}_{\ast}$. We claim, in particular, that any cocartesian lift $\alpha:z\to z'$ of $\bar{f}$ to $\msc{E}^{\ot}$ has $z'$ in $\msc{L}^{\ot}$, and hence that $\alpha$ lines in $\msc{L}^{\ot}$ as well.
\par

First, lift $\bar{f}$ to the cocartesian edge $f:t=(\mcl{I},\mu)\to t'=(\mcl{I}_{\bar{f}},\bar{f}\mu)$ in $fr\Disk^{\ot}$. Here $\mcl{I}_{\bar{f}}$ is the primage of the specified subset $K_{\bar{f}}$ along the map $\mu:\mcl{I}\to K$ and $f$ is explicitly given by the diagram
\[
\xymatrix{
D\times \mcl{I}_{\bar{f}}\ar[r]^{id}\ar[d]_{\mu|_{\mcl{I}_f}} & D\times \mcl{I}_{\bar{f}}\ar[d]^{\bar{f}\mu|_{\mcl{I}_{\bar{f}}}}\\
K_{\bar{f}}\ar[r]_{\bar{f}} & L.
}
\]
For the separations $\tilde{t}=(\mcl{I},id_I)$ and $\tilde{t}'=(\mcl{I}_{\bar{f}},id_{I_{\bar{f}}})$, and any index $i$ in $\mcl{I}_{\bar{f}}$, we have the diagram
\[
\xymatrix{
t_i& \tilde{t}\ar[r]\ar[d]\ar[l]_{\rho_i} & t\ar[d]\\
 & \tilde{t}'\ar[r]\ar[ul]^{\rho_i} & t'.
}
\]
Here $t_i$ is just the singleton $\{\pm 0\}$ with color $\opn{color}(i)$ and $\rho_i$ is the inert projection. These diagrams imply diagrams for the transport functors
\[
\xymatrix{
\msc{E}^{\ot}_t\ar[rr]^{f_!}\ar[dr]_{\pi_i\rho_!} & & \msc{E}^{\ot}_{t'}\ar[dl]^{\pi_i\rho_!}\\
	& \msc{E}^{\pm}
}
\]
which, taken together, provide a diagram
\[
\xymatrix{
\msc{E}^{\ot}_t\ar[rr]^{f_!}\ar[dr]_{\pi_{\bar{f}}\rho_!} & & \msc{E}^{\ot}_{t'}\ar[dl]^{\rho_!}\\
	& \msc{E}^{\mcl{I}_{\bar{f}}}	& .
}
\]
(Here $\pi_{\bar{f}}$ is the projection $\msc{E}^{\mcl{I}}\to \msc{E}^{\mcl{I}_{\bar{f}}}$.) From this we conclude that $f_!$ sends $\msc{L}^{\ot}_t$ into $\msc{L}^{\ot}_{t'}$. Evaluating at our original object $z$ in $\msc{L}^{\ot}_t$, we observe a cocartesian edge $z\to f_!(z)$ in $\msc{L}^{\ot}$ over $f:t\to t'$ in $fr\Disk^{\ot}$, and hence over $\bar{f}:K\to L$ as well.
\par

We've now shown that the induced map $\msc{L}^{\ot}\to \Fin_{\ast}$ is a cocartesian fibration and that the inclusion $\msc{L}^{\ot}\to \msc{E}^{\ot}$ is a map of cocartesian fibrations over $\Fin_{\ast}$. To conclude symmetric monoidality, we need to show that that transport along the inert projections $\bar{\rho}_k:K\to \{0\}$ provide an equivalence $\msc{L}^{\ot}_K\to \msc{L}_{\{0\}}^K$.
\par

Since $\msc{L}^{\ot}_K$ and $\msc{L}^K_{\{0\}}$ are full in $\msc{E}^{\ot}_K$ and $\prod_k\msc{E}^{\ot}_{\{0\}}$, this functor inherits fullness from fullness of the corresponding functor for the ambient category $\msc{E}^{\ot}$. So we need only verify essential surjectivity. For this choose any tuple of objects $(z_k:k\in K)$ in $\msc{L}^K_{\{0\}}$ and consider its image $(t_k:k\in K)$ in $fr\Disk^K$.
\par

Each $t_k$ is just a tuple of colored disks, and taking the disjoint union we obtain an object $t$ in the fiber $fr\Disk_K^{\ot}$ with inert projections $\rho_i:t\to t_k$. Writing $t=(\mcl{I},\mu)$, we can argue as above to observe a diagram
\[
\xymatrix{
\msc{L}^{\ot}_t\ar[rr]\ar[dr]_{\rho_!} & & \prod_k\msc{L}^{\ot}_{t_k}\ar[dl]^{\prod_k\rho_!}\\
 & \msc{L}^{\mcl{I}}
 }
\]
where the maps to $\msc{L}^{\mcl{I}}$ are both equivalences. Hence transport provides an equivalence $\msc{L}^{\ot}_t\overset{\sim}\to \prod_k\msc{L}^{\ot}_{t_k}$, and there is an object $z$ in $\msc{L}^{\ot}_t$ which admits cocartesian edges $\varrho_k:z\to z_k$ over the projections $\rho_k:t\to t_k$. Composing to $\Fin_{\ast}$, we have cocartesian edges $\varrho_k$ over the inert projections $\bar{\rho}_k:K\to \{0\}$. By uniqueness of cocartesian edges it follows that the tuple $(z_k:k\in K)$ is in the essential image of the functor $\msc{L}^{\ot}_K\to \prod_k\msc{L}^{\ot}_{\{0\}}$, as required.
\par

Finally, as for stability of $\msc{L}^{\ot}$ under isomorphism, choose an isomorphism $\alpha:z\to w$ in $\msc{E}^{\ot}$. Let $f:t=(\mcl{I},\mu)\to s=(\mcl{J},\nu)$ be the image of $\alpha$ in $fr\Disk^{\ot}$ and note that $f$ is an isomorphism as well. Hence $f$ is defined by a collection of disk embeddings $D_i\to D_j$, each of which is homotopic to the identity (see Lemma \ref{lem:1430}), and we have a diagram
\[
\xymatrix{
\msc{E}^{\ot}_t\ar[rr]^{f_!}\ar[d]_{\sim} & & \msc{E}^{\ot}_s\ar[d]^{\sim}\\
\msc{E}^{\mcl{I}}\ar[rr] & & \msc{E}^{\mcl{J}}
}
\]
where the bottom map sends the $i$-th factor in $\msc{E}^{\mcl{I}}$ identically to the $\pi_0(f)(i)$-th factor in $\msc{E}^{\mcl{J}}$. Since $\alpha$ is a cocartesian lift of $f$, we have $f_!(t)\cong w$ in the fiber $\msc{E}^{\ot}_s$, and from the above diagram it follows that $z$ is in $\msc{L}^{\ot}_t$ if and only if $w$ is in $\msc{L}^{\ot}_s$.
\end{proof}

\section{Framed pairs from balanced categories}
\label{sect:framed_pairs}

We construct a canonical $fr\E_2$-monoidal $\infty$-category $q:E^{\odot}\to fr\E_2$ associated to a map of balanced monoidal categories $d:E^{-}\to E^+$. As a preliminary point, we recall how one constructs an (uncolored) $fr\opn{E}_2$-monoidal $\infty$-category from any balanced monoidal category. Pulling back along the disjoint union functor, we obtain a corresponding $fr\Disk$-category $E^{\ot}=\opn{Env}(E^{\odot})\to fr\Disk^{\ot}$.

\subsection{Framed $\opn{E}_2$-algebras from balanced categories}

We consider the standard operad $fr\opn{E}_2$ of (uncolored) disk embeddings.  This operad is defined via a corresponding simplicial operad $fr\underline{\opn{E}}_2$, exactly as in Section \ref{sect:fr_e2}, but without colors.  So, objects are finite sets, morphisms $f:I\to J$ are oriented disk embeddings $f:D\times I_{f}\to D\times J$, $2$-simplices are generally non-commuting diagrams
\begin{equation}\label{eq:5787}
\xymatrix{
	& J\ar[dr]^g\\
I\ar[rr]_h\ar[ur]^f& & K
}
\end{equation}
whose underlying maps of finite pointed sets commute $\bar{g}\bar{f}=\bar{h}$, and which are equipped with a homotopy $\sigma:|\Delta^1|\times D\times I_{\bar{h}}\to D\times K$ between $gf$ and $h$.  As in the colored setting, we have the non-full subcategory
\[
Linfr\opn{E}_2\to fr\opn{E}_2
\]
spanned by all objects and linear disk arrangements $f:D\times I_{\bar{f}}\to D\times J$.  The same arguments employed in the colored case tell us that this inclusion is an equivalence.
\par

From any $2$-simplex $\sigma:\Delta^2\to Linfr\opn{E}_2$ as in \eqref{eq:5787} the underlying homotopy traces out a framed braid $[\sigma]=(\beta,n)$ in $\opn{Br}_{fr}(gf)$ which transaltes between the orderings on the subsets $I_k\subseteq I$ specified by the two maps $gf$ and $h$.  This framed braid is defined exactly as in Section \ref{sect:framed_ribbon}.  The following is well-known, and can be obtained from a calculation of the $2$-categorical truncation $\pi_1 Linfr\E_2$ via framed braids.  See for example \cite[Definition 7.1, Propositions 7.4 and 7.6]{salvatorewahl03} (cf.\ \cite[Proposition 6.2.2]{fresse17}).

\begin{proposition}\label{prop:balanced_fre2}
Any balanced monoidal $(1-)$category $B$ has an associated functor $F_B:Linfr\opn{E}_2\to \sCat$ with the following values:
\begin{itemize}
\item $F_B(I)=B^I$.\vspace{1mm}
\item For $f:I\to J$ in $Linfr\opn{E}_2$, $F_B(f)=m_f:B^I\to B^J$.\vspace{1mm}
\item For a $2$-simplex $\sigma$ as in \eqref{eq:5787}, with associated framed braid $(\beta,n)$,
\[
F_B(\sigma)=\zeta_{\beta,n}:m_gm_f\cong m_{gf}\overset{\sim}\to m_h.
\]
\end{itemize}
\end{proposition}

Here $m_f:B^I\to B^J$ applies of the product functors $B^{I_j}\to B$ in each factor, where each $I_j$ is ordered by the given map $f$ as in \eqref{eq:m_f}, and $\zeta_{\beta,n}$ is defined via the braid group action on the multiplication functors $m_l:B^I\to B^K$ provided by the braiding on $B$ and the twist, as in Section \ref{sect:braids},
\[
\zeta_{\beta,n}=\beta\circ m_{gf}(\theta^n)\circ\opn{assoc}.
\]
We sketch a proof.

\begin{proof}[Sketch proof]
We have the simplicial operad $Lin\underline{\opn{Conf}}^{fr}$ whose objects are finite pointed sets, and maps are linear embeddings of disks $f:D\times I_{\bar{f}}\to D\times J$.  An $n$-simplex in the mapping complex $\uHom(I,J)$ consists of the choice of a map $\bar{f}:I\to J$ in $\Fin_{\ast}$, choices of linear disk embeddings $f_r$ over $\bar{f}$ for each $0\leq r\leq n$, and the choice of an $n$-simplex $\bar{\sigma}:\Delta^n\to \prod_{j\in J}\opn{Conf}^{fr}_{I_j}(D)$ in the space of framed embeddings of points $I_j\to D$ for which $\bar{\sigma}_r=f_r|_{\{0\}\times I_j}$ at each index $r$.  Composition is defined as in $fr\underline{\opn{E}}_2$.
\par

We apply the homotopy category functor on the mapping complexes to obtain a strict $2$-category $\opn{h}_2Lin\underline{\opn{Conf}}^{fr}$.  The objects in this $2$-category are finite pointed sets, maps are linear disk embeddings as above, and $2$-morphisms $(\beta,n):f\to f'$ are choices of framed braids which translate between the orderings on the preimages $I_j\to \{j\}$ provided by the initial disk embedding $f:D\times I_{f}\to D\times J$ and the target embedding $f':D\times I_{f}\to D\times J$.
\par

Remembering only the centers of disks and their inherited framings gives a forgetful functor $Linfr\underline{\opn{E}}_2\to Lin\underline{\opn{Conf}}^{fr}$ and taking the homotopy truncation provides a map of simplicial operads $Linfr\underline{\opn{E}}_2\to \opn{h}_2Lin\underline{\opn{Conf}}^{fr}$.  We apply the nerve to obtain a map of $\infty$-operads $Linfr\opn{E}_2\to \opn{h}_2Lin\opn{Conf}^{fr}$, where the latter simplicial set is described entirely by is $0$, $1$ and $2$-simplices as in Lemma \ref{lem:2cat}.
\par

We have the well defined assignment of $0$, $1$ and $2$-simplices $\bar{F}_B[d]:\opn{h}_2Lin\opn{Conf}^{fr}[d]\to \sCat[d]$ as prescribed in the statement, and it is a practical matter to check that these assignments define a map of simplicial sets $\bar{F}_B$, via Lemma \ref{lem:2cat}.  We compose with the projection from $Linfr\opn{E}_2$ to obtain the promised map $F_B$.
\end{proof}

We note that the functor $F_B$ extends uniquely to a functor from the ambient category $fr\opn{E}_2$.

\begin{lemma}\label{lem:6017}
For any balanced monoidal category $B$, there is a unique functor $F'_B:fr\opn{E}_2\to \sCat$ which completes a strictly commuting diagram
\[
\xymatrix{
Linfr\opn{E}_2\ar[rr]^{F_B}\ar[dr] & & \sCat\\
	& fr\opn{E}_2\ar[ur]_{F'_B}
}
\] 
\end{lemma}

\begin{proof}
Since the map $Linfr\opn{E}_2\to fr\opn{E}_2$ is an injective equivalence the restriction functor
\[
\Fun(fr\opn{E}_2,\sCat_{\infty})\to \Fun(Linfr\opn{E}_2,\sCat)
\]
is a equivalence which is also an isofibration \cite[\href{https://kerodon.net/tag/01F3}{01F3}]{kerodon}. It follows that the projection from the fiber
\[
\{F_B\}\times_{\Fun(Linfr\opn{E}_2,\sCat)}\Fun(fr\opn{E}_2,\sCat)\to \ast
\]
is also an equivalence.  In particular the above fiber, which classifies such completing functors, is a contractible Kan complex.
\end{proof}

\begin{definition}\label{def:balanced_ot}
For any balanced monoidal category we let
\[
B^{\odot}=\int_{fr\opn{E}_2}F'_B\to fr\opn{E}_2
\]
denote the cocartesian fibration determined by the functor $F'_B:fr\opn{E}_2\to \sCat\subseteq \sCat_{\infty}$ from Proposition \ref{prop:balanced_fre2} and Lemma \ref{lem:6017}.
\end{definition}

Explicitly, $B^{\odot}$ is defined by pulling back the universal fibration $\msc{P}\!\sCat_{\infty}\to \sCat_{\infty}$ along the map $F'_B$.  In particular, from the diagram in Lemma \ref{lem:6017}, the pullback along the inclusion $Linfr\opn{E}_2\to fr\opn{E}_2$ recovers the fibration
\[
B^{\odot}|_{Linfr\opn{E}_2}=B^{\odot}\times_{fr\opn{E}_2}Linfr\opn{E}_2=\int_{Linfr\opn{E}_2}F_B.
\]
Also the transport functors along maps $f:I\to J$ in $Linfr\opn{E}_2$ recover the product functors $m_f:B^I\to B^J$.  It is therefore clear that the given fibration realizes $B^{\odot}$ as a $fr\opn{E}_2$-monoidal $\infty$-category \cite[Definition 2.1.2.13]{ha}.  We describe $B^{\odot}$ explicitly in Section \ref{sect:bal_description} below.

\subsection{Describing the fibration $B^{\odot}\to fr\opn{E}_2$}
\label{sect:bal_description}

Consider a balanced monoidal category $B$ and let $B^{\odot}\to fr\opn{E}_2$ be the associated fibration from Definition \ref{def:balanced_ot}.  By definition, $B^{\odot}$ is the pullback of the universal fibration $\msc{P}\!\sCat_{\infty}\to \sCat_{\infty}$ along the composite of $F'_B$ with the inclusion $\sCat\to \sCat_{\infty}$.  Taking $\msc{P}\!\sCat=\sCat\times_{\sCat_{\infty}}\msc{P}\!\sCat_{\infty}$ we now obtain $B^{\odot}$ as the pullback
\[
\xymatrix{
B^{\odot}=fr\opn{E}_2\times_{\sCat}\msc{P}\!\sCat\ar[r]\ar[d] & \msc{P}\!\sCat\ar[d]\\
fr\opn{E}_2\ar[r] & \sCat.
}
\]
So, to describe $B^{\odot}$ it suffices to describe the $\infty$-category $\msc{P}\!\sCat$.
\par

Explicitly, we have the homotopy coherent nerve $\opn{N}(\underline{\opn{Cat}})$ and the undercategory $\opn{N}(\underline{\opn{Cat}})_{\ast/}\to \opn{N}(\underline{\opn{Cat}})$.  (The simplicial set $\opn{N}(\underline{\opn{Cat}})$ is an $(\infty,2)$-category, not an $\infty$-category, though the point is irrelevant for us.)  The $\infty$-category $\msc{P}\!\sCat$ can be identified with the fiber product
\[
\msc{P}\!\sCat=\sCat\times_{\opn{N}(\underline{\opn{Cat}})}\opn{N}(\underline{\opn{Cat}})_{\ast/}
\]
\cite[\href{https://kerodon.net/tag/0213}{0213}]{kerodon}.

Now, since the simplicial category $\underline{\opn{Cat}}$ is a strict $2$-category, its homotopy coherent nerve can be described explicitly in accordance with Lemma \ref{lem:2cat}. This explicit description yields an explicit description of the simplices of the slice category $\opn{N}(\underline{\opn{Cat}})_{\ast/}$, and hence of $\msc{P}\!\sCat$.
\par

Precisely, objects in $\msc{P}\!\sCat$ are categories equipped with a chosen object $x:\ast\to C$.  Morphisms are functors $F:C\to C'$ equipped with a map between the chosen objects $\alpha:F(x)\to x'$, and $2$-simplices are triples of functors $F_{ij}:C_i\to C_j$ with a choice of natural isomorphism $\xi:F_{12}F_{01}\to F_{02}$ which produces a diagram
\begin{equation}\label{eq:5993}
\xymatrix{
F_{12}F_{01}(x_0)\ar[d]_{\xi}\ar[rr]^{F_{12}\alpha_{01}} & & F_{12}(x_1)\ar[d]^{\alpha_{12}}\\
F_{02}(x_0)\ar[rr]_{\alpha_{02}} & & x_2
}
\end{equation}
for the chosen objects $x_i$ and morphisms $\alpha_{ij}$. A general $n$-simplex consists of a choice of pointed categories $x_i:\ast\to C_i$, functors $F_{ij}:C_i\to C_j$, morphisms $\alpha_{ij}:F_{ij}(x_i)\to x_j$, and natural isomorphisms $\xi_{ijk}:F_{jk}F_{ij}\to F_{ik}$ for all $0\leq i<j<k\leq n$ which complete diagrams \eqref{eq:5993} at all such triples, and which also satisfy the compatibility
\[
\xi_{jkl}(id_{F_{kl}}\circ \xi_{ijk})=\xi_{ijl}(\xi_{jkl}\circ id_{F_{ij}}):F_{kl}F_{jk}F_{ij}\to F_{il}
\]
at all quadruples $0\leq i<j<k<l\leq n$.  The cocartesian edges in $\msc{P}\!\sCat$, along the universal fibration $\msc{P}\sCat\to\sCat$, are precisely those pairings of a functor $F:\msc{C}\to \msc{C}'$ with an isomorphism $\alpha:F(x)\overset{\sim}\to x'$ between the chosen objects.
\par

Pulling back along the functor $F_B:Linfr\opn{E}_2\to \sCat$, which is precisely the portion of the functor $F'_B$ which we understand, we obtain an explicit description of the restriction $B^{\odot}|_{Linfr\opn{E}_2}$. Objects in $B^{\odot}|_{Linfr\E_2}$ are choices of a finite colored set $I$ and an objects $x$ in the exponent $B^{I}$, and morphisms $\alpha:\{x,I\}\to \{y,J\}$ consist of the choice of a map $f:I\to J$ in $Linfr\opn{E}_2$ and a map $\alpha:m_f(x)\to y$ in $B^{J}$.  A general $n$-simplex $\sigma:\Delta^n\to B^{\odot}|_{Linfr\opn{E}_2}$ consists of an $n$-simplex $\bar{\sigma}:\Delta^n\to fr\opn{E}_2$, which specifies a linear disk arrangement $f_{ij}:D\times I^i_{f_{ij}}\to D\times I^j$ and homotopies $\bar{\sigma}_{ijk}:f_{jk}f_{ij}\to f_{ik}$, and a choice of the following data:
\begin{itemize}
\item An object $x_i$ in $B^{I^i}$ for each index $0\leq i\leq n$.\vspace{1mm}
\item A choice of morphism $\alpha_{ij}:x_i\to x_j$ for each $i\leq j$.
\end{itemize}
These data are required to produce a diagram
\[
\xymatrix{
m_{jk}m_{ij}(x_i)\ar[d]_{\zeta_{ijk}}\ar[rr]^{m_{jk}\alpha_{ij}} & & m_{jk}(x_j)\ar[d]^{\alpha_{jk}}\\
m_{ik}(x_i)\ar[rr]_{\alpha_{ik}} & & x_k
}
\]
whenever $i<j<k$ and to produce equalities
\[
\zeta_{jkl}(id_{m_{kl}}\circ \zeta_{ijk})=\zeta_{ijl}(\zeta_{jkl}\circ id_{m_{ij}}):m_{kl}m_{jk}m_{ij}\to m_{il}
\]
whenever $i<j<k<l$, where $m_{ij}=m_{f_{ij}}$ and $\zeta_{ijk}:m_{jk}m_{ij}\to m_{ik}$ is the transformation specified by the framed braid extracted from the homotopy $\bar{\sigma}_{ijk}$.

It is a general fact that the transport functors for the cocartesian fibration $\int_{\msc{K}}F\to \msc{K}$ associated to a functor $F:\msc{K}\to \sCat_{\infty}$ recover the functors $F(f):F(s)\to F(t)$ determined by $F$ \cite[\href{https://kerodon.net/tag/027K}{027}]{kerodon}. In particular, the product functors $m_f:B^I\to B^J$ are recovered as transport for the fibration $B^{\odot}\to fr\opn{E}_2$, along any map in $Linfr\opn{E}_2$.  Taking $\rho_i:I\to \{0\}$ the inert maps in $fr\opn{E}_2$, i.e.\ the tautological identifications $\rho_i:D\times\{i\}\to D$, we have that the transport functor along $\rho_i$ is simply the projection $m_{\rho_i}=p_i:B^I\to B$. Hence transport along the inert projections recovers the identity map $\rho_!:B^{\odot}_I=B^I\to B^I$. We therefore see that $B^{\odot}\to fr\opn{E}_2$ is in fact a $fr\opn{E}_2$-monoidal $\infty$-category.

\begin{proposition}\label{prop:bal_fre2_monidal}
For any balanced monoidal category $B$, the cocartesian fibration $B^{\odot}\to fr\opn{E}_2$ from Definition \ref{def:balanced_ot} is a $fr\opn{E}_2$-monoidal $\infty$-category. Furthermore, for each map $f:I\to J$ in $Linfr\opn{E}_2$, the transport functor along $f$ recovers the associated product functor $m_f:B^I\to B^J$.
\end{proposition}

\subsection{$fr\E_2$-monoidal $\infty$-categories from balanced functors}
\label{sect:e_ot_const}

Given a balanced monoidal functor $d:E^-\to E^+$ we have the associated fibration $p:E_d\to \Delta^1$ obtained from the weighted nerve of $d$.  Specifically, $E_d$ is the ($1$-)category whose fibers over $0$ and $1$ are $E^-$ and $E^+$ respectively, and whose maps from the negative to positive fiber are
\[
\Hom_{E_d}(x,y)=\Hom_{E^+}(d(x),y).
\]
Note that there are no maps from the positive to negative fiber.
\par

The balanced structure on $d$ endows $E_d$ with the unique structure of a balanced monoidal category so that the two inclusions $E^{\pm}\to E_d$ are maps of (non-unital) balanced monoidal categories, with trivial tensor compatibility, and with the product between the fibers given by
\[
x\ot y=d(x)\ot y\ \ \text{and}\ \ y\ot x=y\ot d(x)\ \ \text{whenever $x$ is over $0$ and $y$ is over }1.
\]
The unit for this category is $\1= \1_{E^-}$, and the associator between the fibers is provided by the monoidal structure on $d$.

\begin{remark}
The inclusion $E^-\to E_d$ is unital while the inclusion from $E^+$ is non-unital. Indeed, the unit for $E^+$ is sent to the object $\1_{E^+}$ in $(E_d)_{\{1\}}$, and this object admits a distinguished $p$-cocartesian morphism $\1\to \1_{E^+}$ provided by the unit structure $d(\1_{E^-})\overset{\sim}\to \1_{E^+}$ on $d$.
\end{remark}

We have the associated fibration $q_d:E_d^{\odot}\to fr\opn{E}_2$ from Definition \ref{def:balanced_ot} which admits two apparent maps of operads
\[
\xymatrix{
(E^-)^{\odot}\ar[r]\ar[dr]_{q^-} & E_d^{\odot}\ar[d]^{q_d} & (E^+)^{\odot}\ar[l]\ar[dl]^{q^+}\\
	& fr\opn{E}_2 & .
}
\]
We pull back along the forgetful functor $fr\E_2\to fr\opn{E}_2$ to get a cocartesian fibration $fr\E_2\times_{fr\opn{E}_2}E_d^{\odot}\to fr\E_2$.

\begin{definition}
Given a map $d:E^-\to E^+$ between balanced monoidal ($1$-)categories and associated fibration $q_d:E_d^{\odot}\to fr\opn{E}_2$, we take
\[
E^{\odot}\ \subseteq\ fr\E_2\times_{fr\opn{E}_2}E_d^{\odot}
\]
the full subcategory whose fibers are the colored exponents $E^{\odot}_{\mcl{I}}=E^{\mcl{I}}\ \subseteq\ E_d^I$. We consider $E^{\odot}$ along with the inner fibration $q:E^{\odot}\to fr\E_2$ given by restricting the cocartesian fibration $fr\E_2\times_{fr\opn{E}_2}E_d^{\odot}\to fr\E_2$ to $E^{\odot}$.
\end{definition}

We note that $E^{\odot}$ is stable under isomorphism in $fr\E_2\times_{fr\opn{E}_2}E_d^{\odot}$, so that the map $q:E^{\odot}\to fr\E_2$ is in fact an isofibration.  
\par

Recall that each $f:\mcl{I}\to\mcl{J}$ in $Linfr\E_2$ determines an associated colored product $m_f:E^{\mcl{I}}\to E^{\mcl{J}}$, and the corresponding map $\bar{f}:I\to J$ in $Linfr\opn{E}_2$ determines an uncolored product $m_{\bar{f}}:E_d^I\to E_d^J$. Given such $f$ in $Linfr\E_2$ and $x$ in the fiber $E^{\odot}_{\mcl{I}}=E^{\mcl{I}}$, we define the map $1_f:x\to m_f(x)$ in $E^{\odot}$ over $f$ which we recall is defined by an underlying map $1_f:m_{\bar{f}}(x)\to m_f(x)$ in $E_d^J$ as follows:
\begin{itemize}
\item For each negatively colored index $j\in\mcl{J}$, we have that the $j$-th components agree $m_{\bar{f}}(x)_j=m_f(x)_j$ and we take $(1_f)_j:m_{\bar{f}}(x)_j\to m_f(x)_j$ to be the identity $(1_f)_j=id_{m_f(x)_j}$.\vspace{1mm}
\item For each positively colored index $j\in\mcl{J}$ whose preimage $\mcl{I}_j$ also contains a positive index, we again have $m_{\bar{f}}(x)_j=m_{f}(x)_j$ and take $(1_f)_j=id_{m_f(x)_j}$.\vspace{1mm}
\item For each positively colored index $j\in\mcl{J}$ whose preimage $\mcl{I}_j$ contains \emph{no} positive indices, we have $m_{f}(x)_j=m_{\bar{f}}(d^{I_j}x)_j$, where $d^{I_j}$ applies $d$ at each index in $I_j\subseteq I$, and we let $(1_f)_j:m_{\bar{f}}(x)_j\to m_f(x)_j$ be the map in $E_d$ given by the natural isomorphism $d(m_{\bar{f}}(x)_j)\overset{\sim}\to m_{\bar{f}}(d^{I_j}x)_j=m_f(x)$ given by the monoidal structure on $d$.
\end{itemize}

\begin{proposition}\label{prop:balanced1}
For any balanced monoidal functor $d:E^-\to E^+$, the associated isofibration $q:E^{\odot}\to fr\E_2$ from Definition \ref{def:balanced_ot} is a cocartesian fibration which gives $E^{\odot}$ the structure of a $fr\E_2$-monoidal $\infty$-category.  Furthermore, for each $f:\mcl{I}\to\mcl{J}$ in $Linfr\E_2$ and $x$ in $E^{\odot}$, the map $1_f:x\to m_f(x)$ is a $q$-cocartesian edge over $f$.
\end{proposition}

\begin{proof}
Call a map $f:\mcl{I}\to \mcl{J}$ in $fr\E_2$ consistently colored if the preimage $\mcl{I}_j$ of each positive index $j$ in $\mcl{J}$ is nonempty and contains a positive index $i$. Since this property is determined by the underlying morphism of finite pointed colored sets, we see that such consistency is stable under isomorphisms of maps in $fr\E_2$.
\par 

For consistently colored $f:\mcl{I}\to \mcl{J}$ in $Linfr\E_2$ we note that the uncolored product $m_{\bar{f}}:E_d^I\to E_d^J$ for $E_d$ sends $E^{\mcl{I}}$ into $E^{\mcl{J}}$, and in fact recovers the colored product $m_f:E^{\mcl{I}}\to E^{\mcl{J}}$ under restriction.  Hence for any $x$ in $E^{\mcl{I}}\subseteq E^I_d$ the distinguished cocartesian lift $1_{f}:x\to m_f(x)$ in $E_d^{\odot}$ provided by the identity on $m_f(x)$ is already in $E^{\odot}$.  For general consistently colored $f':\mcl{I}\to \mcl{J}$ in $fr\E_2$ we have $m_{f'}\cong m_f$ for some consistently colored $f$ in $Linfr\E_2$, so that again the cocartesian lift $x\to m_{f'}(x)$ in $E_d^{\odot}$ lies in $E^{\odot}$, by fullness of $E^{\odot}$ and stability under isomorphism.
\par

Now, by a color change map in $fr\E_2$ we mean a map $\omega:\mcl{I}\to \mcl{I}'$ in $fr\E_2$ which maps to $id_I:I\to I'=I$ in $fr\opn{E}_2$.  We note that every map $f'$ in $fr\E_2$ is a composite $f'=f\omega u$ of consistently colored maps $f$ and $u$, and a color change map $\omega$.  Here $u$ is just an inclusion $D\times \mcl{I}^0\to D\times \mcl{I}^1$ induced by an inclusion of sets $\mcl{I}^0\to \mcl{I}^1$ in which all added indices in $\mcl{I}^1$ are negatively colored. So, via the $2$-of-$3$ property for cocartesian edges \cite[\href{https://kerodon.net/tag/01TS}{01TS}]{kerodon}, it suffices to show that each color change map admits a $q$-cocartesian lift.
\par

For a color change map $\omega:\mcl{I}\to\mcl{I}'$ let $d_{\omega}:E^{\mcl{I}}\to E^{\mcl{I}'}$ be the functor which is the identity at each index $i$ with $\opn{color}(i)=\opn{color}(\omega i)$, and which applies $d$ at each index $i$ with $\opn{color}(i)=-\opn{color}(\omega i)$. (This map $d_{\omega}$ is equal to the colored product $m_{\omega}$ over $\omega$.)  For $x$ in $E^{\mcl{I}}$ take $1_{\omega}:x\to d_{\omega}(x)$ the map over $\omega$ in $E_d^{\odot}$ which is just the identity on each factor.  We claim that $1_{\omega}$ is $q$-cocartesian.
\par

Consider an $n$-simplex $\tau:\Delta^n\to fr\E_2$, with specified objects $\mcl{I}^r$, maps $f'_{rs}:\mcl{I}^r\to \mcl{I}^s$, and homotopies $\tau_{rst}:f'_{st}f'_{rs}\to f'_{rt}$.  Suppose $f'_{01}$ is a color change map $f'_{01}=\omega:\mcl{I}^0\to \mcl{I}^1$, so that each $f'_{0s}$ admits a strict factorization $f'_{0s}=f_{0s}\omega$ for a uniquely specified map $f_{0s}$.  In fact, taking $f_{rs}=f'_{rs}$ when $r>0$, the $n$-simplex $\sigma$ specifies a corresponding $n$-simplex $\nu$ with maps $f_{rs}$ and the same higher data, as this higher data is purely topological, i.e.\ has nothing to do with colorings.
\par

For $\tau$ as above we consider a lifting problem
\begin{equation}\label{eq:6122}
\xymatrix{
\Lambda^n_0\ar[r]^{\bar{\sigma}}\ar[d] & E^{\odot}\ar[d]\\
\Delta^n\ar[r]_{\tau}\ar@{..>}[ur] & fr\E_2
}
\end{equation}
in which $\bar{\sigma}|_{\Delta^{\{0,1\}}}=1_{\omega}:x_0\to d_{\omega}(x_0)=x_1$.  The horn $\bar{\sigma}$ is specified by maps $\alpha_{rs}:m'_{rs}(x_r)\to x_s$ which complete diagrams
\[
\xymatrix{
m'_{st}m'_{rs}(x_r)\ar[rr]^{m'_{st}\alpha_{rs}}\ar[d]_{\zeta_{rst}} & & m'_{st}(x_s)\ar[d]^{\alpha_{st}}\\
m'_{rt}(x_r)\ar[rr]_{\alpha_{rt}} & & x_t
}
\]
where $m'_{rs}=m_{f'_{rs}}$ and $\zeta_{rst}$ is the transformation specified by the framed braid $(\beta_{rst},n_{rst})$ associated to the homotopy $\tau_{rst}$.  Each $m'_{0s}$ factors as $m_{0s}d_{\omega}$ for $m_{0s}=m_{f_{0s}}$, and taking $m_{rs}=m'_{rs}$ for $r>s$ we obtain a corresponding horn $\bar{\mu}:\Lambda^n_0\to E^{\odot}$ specified by the objects $y_r$ with $y_0=d_{\omega}(x_0)=x_1$ and $y_r=x_r$ for all $r>0$, the same maps $\alpha_{rs}$, and the same natural transformations $\zeta_{rst}$.  Note that $\alpha_{01}=id_{x_1}$ here.
\par

The horn $\bar{\mu}$ fits into a lifting problem
\[
\xymatrix{
\Lambda^n_0\ar[r]^{\bar{\mu}}\ar[d] & E^{\odot}\ar[d]\\
\Delta^n\ar[r]_{\nu}\ar@{..>}[ur] & fr\E_2
}
\]
and since $id_{x_1}$ is $q$-cocartesian we can solve this lifting problem to produce an $n$-simplex $\mu:\Delta^n\to E^{\odot}$. The two layers of data, and third order compatibilities which specify this $n$-simplex specify an $n$-simplex $\sigma:\Delta^n\to E^{\odot}$ which solves the original lifting problem \eqref{eq:6122}.  The generic existence of such solutions verifies that $1_{\omega}$ is $q$-cocartesian.
\par

As for the claim that $1_{f}:x\to m_f(x)$ is always a cocartesian lift of $f$ in $Linfr\E_2$, this is validated immediately in the case where $f$ is consistently colored.  For $f':\mcl{I}'\to \mcl{J}$ not consistently colored we have a factoring $f'=f\omega u$ where $u$ is induced by an inclusion $\mcl{I}'\to \mcl{I}$ in which all added indices are negative, $\omega$ is a color change map, and $f$ is consistently colored.   We can specifically take
\[
\mcl{I}''=\{j\in \mcl{J}:j\text{ positive with empty preimage in }\mcl{I}'\},
\]
recolor this set negatively $-\mcl{I}''$, and take $\mcl{I}=\mcl{I}'\amalg -\mcl{I}''$.  We then have the cocartesian lift of $f'$ given as the composite map
\[
1_f\circ 1_{\omega}\circ 1_u:x\to z\to d_{\omega}(z)\to m_fd_{\omega}(z)
\]
where $z$ in $E^{\mcl{I}}$ is obtained from $x$ in $E^{\mcl{I}'}$ by inserting a copy of $\1_{E^-}$ at each index in $-\mcl{I}''$.  The object $m_fd_{\omega}(z)$ differs from $m_{f'}(x)$ only at the indices in $\mcl{I}''\subseteq \mcl{J}$, and we have the isomorphism $\xi:m_fd_{\omega}(z)\overset{\sim}\to m_{f'}(x)$ in $E^{\mcl{J}}$ which just applies the isomorphism $d(\1_{E^-})\cong \1_{E^+}$ at all such indices.  The isomorphism $\xi$ is a cocartesian lift of its image $id_{\mcl{J}}$ in $fr\E_2$ so that we obtain the new cocartesian lift
\[
1_{f'}=\xi\circ 1_f\circ 1_{\omega}\circ 1_u:x\to m_{f'}(x)
\]
of $f'$ along $q$, as desired. 
\end{proof}

\begin{proposition}\label{prop:balanced2}
Let $d:E^-\to E^+$ be a balanced monoidal functor.  For any map $f:\mcl{I}\to \mcl{J}$ in $Linfr\E_2$ the associated colored product functor
\[
m_f:E^{\mcl{I}}\to E^{\mcl{J}}
\]
is the transport functor for $q:E^{\odot}\to fr\E_2$ along $f$.
\end{proposition}

\begin{proof}
When $f$ is consistently colored, in the sense that $\mcl{I}_j$ contains a positive index whenever $j$ is positive in $\mcl{J}$, then the colored product is obtained by restricting the corresponding uncolored product on $E_d$.  Since the uncolored product is reproduced via transport along $f$ for $E_d^{\odot}\to fr\opn{E}_2$, it follows that the colored product is also recovered via transport along $f$.  (Here we are using the fact that the map $1_{f}:x\to m_f(x)$ is a cocartesian lift of $f$ in both $E^{\odot}$ and $E_d^{\odot}$.)
\par

For a color change map $\omega:\mcl{I}\to \mcl{I}'$, i.e.\ a map over the identity in $fr\opn{E}_2$, recall that $m_{\omega}:E^{\mcl{I}}\to E^{\mcl{J}}$ just applies $d$ at all indices where the color changes from negative to positive.  We have the uniquely specified cocartesian transformation
\[
T_{\omega}:\Delta^1\times E^{\mcl{I}}\to E^{\odot}
\]
with $T_{\omega}|_{\{0\}}=id$, $T_{\omega}|_{\{1\}}=m_{\omega}$, $T_{\omega}(0<1,\alpha:x\to y)=m_{\omega}(\alpha):x\to m_{\omega}(y)$.  Note that for each $x$ in $E^{\mcl{I}}$ this transformation has $T_{\omega}|_{\Delta^1\times\{x\}}:\Delta^1\to E^{\odot}$ equal to the $q$-cocartesian edge $1_{\omega}$.  Hence, restricting along $\{1\}\to \Delta^1$, we see that $m_{\omega}$ provides the transport functor for $E^{\odot}$ along $\omega$.
\par

For general $f':\mcl{I}_0\to \mcl{J}$ which is not consistently colored we have a factorization $f'=f\omega u$ for color change $\omega$, $f$ consistently colored, and $u$ induced by an inclusion $\mcl{I}_0\to \mcl{I}_1$.  Here we can take $\mcl{I}_1$ obtained from $\mcl{I}_0$ by adding in those indices in $\mcl{J}$ which are positively colored but have vanishing preimage, and negating the colors on these indices.  Note that $u$ is consistently colored.  Since composites of transport functors are transport functors, we find that the composite $m_fm_{\omega}m_u:E^{\mcl{I}}\to E^{\mcl{J}}$ is a transport functor along $f'$.
\par

The functor $m_fm_{\omega}m_u$ inserts a copy of $d(\1_{E^{-}})$ at each index $j$ in $\mcl{I}$, and the structural isomorphism $d(\1_{E^-})\cong \1_{E^+}$ induces a natural isomorphism $\xi:m_fm_{\omega}m_u\overset{\sim}\to m_{f'}$ which composes with the cocartesian transformations to produce a cocartesian lift
\[
\xi\circ T_f\circ T_{\omega}\circ T_u:\Delta^1\times E^{\mcl{I}}\to E^{\odot}
\]
of the map $\Delta^1\times E^{\mcl{I}}\to \Delta^1\overset{f}\to fr\E_2$. This transformation realizes $m_{f'}$ as another transport functor along $f'$ in $E^{\odot}$.
\end{proof}

Proposition \ref{prop:balanced2} tells us, in particular, that transport along the inert projections $\rho_i:\mcl{I}\to \{\pm 0\}$ recovers the identity map on the colored exponent $\rho_!=id:E^{\odot}_{\mcl{I}}=E^{\mcl{I}}\to E^{\mcl{I}}$.

\begin{corollary}\label{cor:balanced3}
For any balanced monoidal functor $d:E^-\to E^+$, the associated cocartesian fibration $q:E^{\odot}\to fr\E_2$ gives $E^{\odot}$ the structure of a $fr\E_2$-monoidal $\infty$-category.
\end{corollary}

We now record the proof of Proposition \ref{prop:a_pm}, for the sake of completeness.

\begin{proof}[Proof of Proposition \ref{prop:a_pm}]
Define the $fr\Disk$-monoidal $\infty$-category $E^{\ot}\to fr\Disk^{\ot}$ as the envelope of the $fr\E_2$-monoidal $\infty$-category $E^{\odot}\to fr\E_2$ constructed above.
\end{proof}

\section{Module categories and inner-Homs}
\label{sect:module_cats}

We cover basics information regarding module categories over symmetric monoidal $\infty$-categories.  As a primary point of interest, we calculate inner-Homs for the action of $\Vect_{\pf}$ on $\msc{D}_{\pf}$. For this section one does need to have some familiarity with $\infty$-operads, and we refer the reader to \cite[Sections 2.1.1 and 2.1.2]{ha} for the details.

\subsection{$\infty$-Categories of cocartesian fibrations}
\label{sect:cocart_sm}

For any $\infty$-category $\msc{T}$ and cocartesian fibrations $q_i:\msc{C}_i\to \msc{T}$ we have the complexes
\[
\Fun_{\msc{T}}(\msc{C}_0,\msc{C}_1)=\{q_0\}\times_{\Fun(\msc{C}_0,\msc{T})}\Fun(\msc{C}_0,\msc{C}_1)
\]
of functors over $\msc{T}$, and consider the full subcategory $\Fun_{\msc{T}}^{cc}(\msc{C}_0,\msc{C}_1)$ spanned by those maps which preserve cocartesian edges.  Since the map $(q_1)_{\ast}:\Fun(\msc{C}_0,\msc{C}_1)\to \Fun(\msc{C}_0,\msc{T})$ is an inner fibration \cite[\href{https://kerodon.net/tag/01BV}{01BV}]{kerodon}, the fiber $\Fun_{\msc{T}}(\msc{C}_0,\msc{C}_1)$ is an $\infty$-category, as is the full subcategory $\Fun_{\msc{T}}^{cc}(\msc{C}_0,\msc{C}_1)$.
\par

The composition operations on the ambient complexes $\Fun(\star,\star)$ restrict to provide composition operations for $\Fun_{\msc{T}}^{cc}$, so that we obtain a simplicial category $\underline{\opn{Cocart}}(\msc{T})$ which is enriched in $\infty$-categories.  Applying the associated Kan complex functor on mapping complexes, then applying the homotopy coherent nerve, we obtain an $\infty$-category $\msc{C}\!ocart(\msc{T})=\opn{N}(\underline{\opn{Cocart}}(\msc{T})^+)$.
\par

We note that pulling back along any map $w:\msc{S}\to \msc{T}$ provides a functor $\underline{\opn{Cocart}}(\msc{T})\to \underline{\opn{Cocart}}(\msc{S})$ between $\opn{Cat}_{\infty}$-enriched categories, so that we obtain a pullback functor
\[
w^{\ast}:\msc{C}\!ocart(\msc{T})\to \msc{C}\!ocart(\msc{S}).
\]

\begin{definition}\label{def:sm_infty}
Given an $\infty$-operad $\msc{O}^{\ot}\to \Fin_{\ast}$, the $\infty$-category $\msc{M}\!on_{\infty}^{\msc{O}^{\ot}}$ is the full subcategory in $\msc{C}\!ocart(\msc{O}^{\ot})$ spanned by $\msc{O}$-monoidal $\infty$-categories.  In the particular case $\msc{O}^{\ot}=\Fin_{\ast}$ we take
\[
\SM_{\infty}:=\msc{M}\!on_{\infty}^{\Fin_{\ast}}.
\]
\end{definition}

\subsection{Module categories over symmetric monoidal $\infty$-categories}

\begin{definition}\label{def:mod}
Let $\opn{Mod}$ denote the subcategory of $\mbf{Fin}_{\ast}$ (see Section \ref{sect:col_sets}) consisting of all finite colored sets $\mcl{I}$ and all morphisms $f:\mcl{I}\to \mcl{J}$ which satisfy the following:
\begin{enumerate}
\item[(a)] Every positively colored index in the specified subset $\mcl{I}_{f}\subseteq \mcl{I}$ is sent to a positively colored index in $\mcl{J}$. 
\item[(b)] For each positively colored $j$ in $\mcl{J}$, the preimage $\mcl{I}_j=f^{-1}(j)$ is nonempty and contains precisely one positively colored index.
\end{enumerate}
We consider $\opn{Mod}$ as a category over $\Fin_{\ast}$ via the forgetful functor $\opn{Mod}\to \Fin_{\ast}$.
\end{definition}

\begin{remark}
Condition (a) in Definition \ref{def:mod} is redundant, in the sense that it is already included in the notion of morphisms for $\mbf{Fin}_{\ast}$.
\end{remark}

\begin{lemma}
The forgetful functor $\opn{Mod}\to \Fin_{\ast}$ gives $\opn{Mod}$ the structure of an $\infty$-operad. 
\end{lemma}

\begin{proof}
Any inert map $\bar{\rho}:I\to J$ in $\Fin_{\ast}$ specifies a bijection from the specified subset $I'=I_{\bar{\rho}}\overset{\cong}\to J$.  For any coloring $\mcl{I}$ of $I$ we have the corresponding colored subset $\mcl{I}'$ and coloring $\mcl{J}$ of $J$ induced by the given bijection.  We thus obtain the colored lift $\rho^{\mcl{I}}:\mcl{I}\to \mcl{J}$ of $\bar{\rho}$ which is observed directly to be cocartesian over $\bar{\rho}$.  As the coloring of $\mcl{J}$ is canonically determined by the coloring of $\mcl{I}$ and the map $\bar{\rho}$ in $\Fin_{\ast}$, we have the associated transport functor
\[
\bar{\rho}_!:\opn{Mod}_I\to \opn{Mod}_J,\ \ \mcl{I}\mapsto \mcl{J},
\]
with cocartesian transformation $T_{\bar{\rho}}:\Delta^1\times \opn{Mod}_I\to \opn{Mod}$ provided by the $\rho^{\mcl{I}}$.  The remaining requirements from \cite[Definition 2.1.1.10]{ha}l are now also observed directly.
\end{proof}

We note that there is a unique section $i:\Fin_{\ast}\to \opn{Mod}$ of the forgetful functor which sends each finite pointed set $I$ to $I$ along with its complete negative coloring.  This section is a map of $\infty$-operads, in the sense that it preserves inert edges.  In particular, any $\opn{Mod}$-monoidal $\infty$-category pulls back along $i$ to produce a symmetric monoidal $\infty$-category.

\begin{definition}
Let $\msc{F}^{\ot}\to \Fin_{\ast}$ be a symmetric monoidal $\infty$-category.  A module category over $\msc{F}$ is a $\opn{Mod}$-monoidal $\infty$-category $\msc{M}^{\ot}\to \opn{Mod}$ which is equipped with an equivalence of symmetric monoidal $\infty$-categories
\[
\xymatrix{
\msc{F}^{\ot}\ar[rr]^(.4){\xi}\ar[dr] & & \Fin_{\ast}\times_{\opn{Mod}}\msc{M}\ar[dl]\\
 & \Fin_{\ast} & .
}
\]
The underlying $\infty$-category for an $\msc{F}$-module category is the fiber over the positive singleton $\msc{M}=\msc{M}^{\ot}_{\{+0\}}$, and a map of $\msc{F}$-module categories $F:\msc{M}^{\ot}\to \msc{N}^{\ot}$ is a map of $\opn{Mod}$-monoidal categories equipped with a choice of a diagram
\[
\xymatrix{
\Fin_{\ast}\times_{\opn{Mod}}\msc{M}^{\ot}\ar[rr]^{F|_{\Fin_{\ast}}} & & \Fin_{\ast}\times_{\opn{Mod}}\msc{N}^{\ot}\\
	&\msc{F}^{\ot}\ar[ul]^{\sim}\ar[ur]_{\sim}
}
\]
in $\SM_{\infty}$.
\end{definition}

\begin{remark}
We are abusing notation, as we should really be speaking of $\msc{F}^{\ot}$-module categories, or even $q$-module categories for the given fibration $q:\msc{F}^{\ot}\to \Fin_{\ast}$.
\end{remark}

For any $\msc{F}$-module $\infty$-category $\msc{M}$, we have the action functor
\[
\opn{act}_{\msc{M}}:\msc{F}\times \msc{M}\to \msc{M}
\]
which is produced via transport along the unique active map $\{-1,+0\}\to \{+0\}$ in $\opn{Mod}$.

\begin{remark}\label{rem:3259}
Any equivalence of symmetric monoidal $\infty$-categories $\msc{F}^{\ot}\overset{\sim}\to \opn{Fin}_{\ast}\times_{\opn{Mod}}\msc{M}^{\ot}$ realizes the fiber product as an initial object in the undercategory $(\SM_{\infty})_{\msc{F}^{\ot}/}$ \cite[\href{https://kerodon.net/tag/02J2}{02J2}]{kerodon}.  For $(\SM_{\infty})_{\msc{F}^{\ot}/}^{\opn{init}}$ the full subcategory spanned by initial objects, we can define the $\infty$-category of $\msc{F}$-module categories formally as the fiber product
\[
\msc{F}^{\ot}\text{-}\opn{mod}=(\SM_{\infty})^{\opn{init}}_{\msc{F}^{\ot}/}\times_{\SM_{\infty}}\msc{M}\!od^{\opn{Mod}}_{\infty}.
\]
Since the subcategory of initial objects is contractible \cite[\href{https://kerodon.net/tag/02HM}{02HM}]{kerodon} and since the forgetful functor $(\SM_{\infty})_{\msc{F}^{\ot}/}\to \SM_{\infty}$ is a left fibration \cite{joyal02} \cite[\href{https://kerodon.net/tag/018F}{018F}]{kerodon}, the above fiber product is equivalent to the homotopy fiber product $\{\msc{F}^{\ot}\}\times^{\opn{htop}}_{\SM_{\infty}}\msc{M}\!od^{\opn{Mod}}_{\infty}$.
\end{remark}

\subsection{Module $\infty$-categories from discrete module categories}

Let $S$ be a discrete symmetric monoidal category and $M$ be a discrete module category over $S$, in the usual $1$-categorical sense of the term.  (See for example \cite[Definition 7.1.1]{egno15}.)  We construct the corresponding cocartesian fibration $M^{\ot}\to \opn{Mod}$ whose objects consist of a choice of finite colored set $\mcl{I}$, a linear ordering on the negatively colored indices in $\mcl{I}$, and a choice of an object $x$ in $M^{\mcl{I}}=S^{\opn{color}^{-1}(-)}\times M^{\opn{color}^{-1}(+)}$.  Morphisms $\alpha:x\to x'$ consist of a map $f:\mcl{I}\to \mcl{J}$ in $\opn{Mod}$ and, for each $j$ in $\mcl{J}$, a choice of map $x_{\mcl{I}_j}\to x'_j$ in $S$ or $M$ where
\[
x_{\mcl{I}_j}=\bigotimes_{i\in \mcl{I}_j}x_i.
\]
Here the product is ordered via the chosen ordering on $\mcl{I}_j$ when $j$ is negative, and otherwise is taken in the unique ordering which extends the given ordering on the negative indices and takes the positive index in $\mcl{I}_j$ to be maximal.  In this way any module category in the $1$-categorical sense becomes a module category in the $\infty$-categorical sense.

\begin{definition}\label{def:lin_ab}
By a finitely generated linear abelian category $A^{\heartsuit}$ we mean a presentable, compactly generated, abelian category whose compact objects are precisely the finite length objects, and for which
\[
\opn{dim}\Hom_{A^{\heartsuit}}(x,y)<\infty\ \ \text{whenever $x$ and $y$ are compact}.
\]
\end{definition}

We note that any finitely generated linear abelian category is a module category over $\opn{Vect}$.  Indeed, we might take this as our definition of linearity.  For our purposes, however, we want to consider actions of the pro-finite category $\opn{Vect}_{\pf}$.
\par

The linear action on $A^{\heartsuit}$ completes to an action of $\opn{Vect}_{\pf}$ on $A^{\heartsuit}_{\pf}$ which is defined by the expected formula
\[
v^0\ot x^0=\varprojlim_{\lambda,\mu}v^0_{\lambda}\ot_k x^0_{\mu}.
\]
(See Sections \ref{sect:pf_obj} and \ref{sect:pf_complexes} for basic information on pro-finite completion.) The action of $\opn{Vect}_{\pf}$ on $A^{\heartsuit}_{\pf}$ now extends to an action of unbounded cochains $\opn{Ch}(\opn{Vect}_{\pf})=\opn{Ch}(\opn{Vect})_{\pf}$ on $A_{\pf}=\opn{Ch}(A^{\heartsuit})_{\pf}$, which is again given by the expected formula
\[
v\ot x=\text{the complex with degree $n$ subobject }\prod_{n_1+n_2=n}v^{n_1}\ot x^{n_2}.
\]
Hence we obtain a $\opn{Ch}(\opn{Vect})_{\pf}$-module category of unbounded cochains over $A^{\heartsuit}_{\pf}$.

\begin{proposition}
For any finitely generated linear abelian category $A^{\heartsuit}$, the corresponding category $A_{\pf}$ of pro-finite cochains admits the natural structure of a module $\infty$-category ${_kA}^{\ot}_{\pf}\to \opn{Mod}$ over the symmetric monoidal $\infty$-category $\opn{Ch}(\opn{Vect})_{\pf}^{\ot}\to \Fin_{\ast}$.
\end{proposition}

\subsection{Inner-Homs}

For symmetric $\msc{F}$ and an $\msc{F}$-module category $\msc{M}^{\ot}\to \opn{Mod}$ we have the action map $\opn{act}_{\msc{M}}:\msc{F}\times \msc{M}\to \msc{M}$.  Hence, for any Hom-functor $H_{\msc{M}}:\msc{M}^{op}\times\msc{M}\to \sKan$ we have the corresponding sequence
\[
\msc{F}^{op}\times \msc{M}^{op}\times \msc{M}\to \msc{M}^{op}\times \msc{M}\to \sKan
\]
given by composing the action with $H_{\msc{M}}$.  Via adjunction this specifies a map
\[
\opn{act}^{\vee}:\msc{M}^{op}\times \msc{M}\to \Fun(\msc{F}^{op},\sKan),\ \ (m,m')\mapsto \{x\mapsto H_{\msc{M}}(x\ot m,m')\}.
\]

\begin{remark}
Recall that a Hom functor $H_{\msc{M}}$ is simply a ``bi-functorialization" of the mapping spaces
\[
\Maps_{\msc{M}}(x,x')=\{x,x'\}\times_{\Fun(\partial\Delta^1,\msc{M})}\Fun(\Delta^1,\msc{M}).
\]
Such bifunctors are obtained by straightening the twisted arrow fibration $\opn{TW}(\msc{M})\to \msc{M}^{op}\times \msc{M}$.  See \cite[\href{https://kerodon.net/tag/03K1}{03K1}, \href{https://kerodon.net/tag/03K6}{03K6}]{kerodon}.
\end{remark}

Recall the Yoneda embedding $\msc{F}\to \Fun(\msc{F}^{op},\sKan)$, which is an equivalence onto the full subcategory spanned by representable functors \cite[\href{https://kerodon.net/tag/03NF}{03NF}, \href{https://kerodon.net/tag/03M5}{03M5}]{kerodon}.

\begin{definition}
An $\msc{F}$-module category $\msc{M}$ is said to admit inner-Homs if there is a bifunctor $\underline{\Maps}:\msc{M}^{op}\times \msc{M}\to \msc{F}$ which completes a diagram
\[
\xymatrix{
	& \msc{F}\ar[dr]^{Yon}\\
\msc{M}^{op}\times \msc{M}\ar[rr]_{\opn{act}^{\vee}}\ar[ur]^{\underline{\Maps}} & & \Fun(\msc{F}^{op},\sKan)
}
\]
in $\sCat_{\infty}$.
\end{definition}

Clearly inner-Homs exist if and only if each functor $\opn{act}^{\vee}(x,x')=H_{\msc{M}}(-\ot x,x')$ is representable.  Furthermore, in this case the representing functor $\underline{\opn{Maps}}$ is uniquely determined up to a contractible space of choices.  For us, existence of inner-Homs over all of $\msc{M}$ is a bit too much to ask for.  So we consider the following.

\begin{definition}
An $\msc{F}$-module category $\msc{M}$ is said to admit restrictive inner-Homs, relative to a specified full subcategory $\msc{M}_{fin}$ in $\msc{M}$, if there is a bifunctor $\underline{\Maps}:\msc{M}^{op}_{fin}\times \msc{M}\to \msc{F}$ which completes a diagram
\[
\xymatrix{
	& \msc{F}\ar[dr]^{Yon}\\
\msc{M}^{op}_{fin}\times \msc{M}\ar[rr]_{\opn{act}^{\vee}|_{\msc{M}_{fin}^{op}\times \msc{M}}}\ar[ur]^{\underline{\Maps}} & & \Fun(\msc{F}^{op},\sKan)
}
\]
in $\sCat_{\infty}$.
\end{definition}

We note that the full subcategory $\msc{M}_{fin}$ is \emph{not} assumed to be an $\msc{F}$-module subcategory, so that the restriction of $\opn{act}^{\vee}$ is obtained from the sequence
\[
\msc{F}^{op}\times\msc{M}_{fin}^{op}\times \msc{M}\overset{\opn{act}\times 1}\to \msc{M}^{op}\times \msc{M}\overset{H_{\msc{M}}}\to \sKan.
\]

\begin{lemma}\label{lem:a_innhom}
For any finitely generated linear abelian category $A^{\heartsuit}$, the $\opn{Ch}(\opn{Vect})_{\pf}$-module category $A_{\pf}$ admits restrictive inner-Homs, relative to the subcategory $A_{fin}$ of finite length bounded cochains.  Furthermore, these inner-Homs are given by the standard cochain Homs
\[
\underline{\Maps}(x,y)=\Hom^{\ast}_{A}(x,y).
\]
\end{lemma}

In the above expression $x$ is finite and bounded, and $y$ is expressed as a limit $y=\varprojlim_{\alpha}y_{\alpha}$ of finite $y_{\alpha}$. Hence the Homs are naturally topologized via the limit expression
\[
\Hom^{\ast}_{A}(x,y)=\varprojlim_{\alpha}\Hom^{\ast}_{A_{fin}}(x,y_{\alpha}).
\]

\begin{proof}
All categories under consideration are discrete, so that the standard set-valued Hom functors provide $\infty$-categorical Hom functors.  Also, for $x$ in $A_{fin}$ the topological product $v\ot x=v\hat{\ot}_kx$ agrees with the usual product $v\ot_k x$.  Hence the result follows by the standard adjunction
\[
\begin{array}{c}
\Hom_{A_{\pf}}(v\ot x,y)\overset{\sim}\to \Hom_{\opn{Ch}(\opn{Vect})_{\pf}}(v,\Hom^{\ast}_{A}(x,y))\vspace{2mm}\\
f\mapsto (\nu\mapsto f(\nu,-)).
\end{array}
\]
\end{proof}

\subsection{Localizations of module $\infty$-categories}

One should recall at this point the general theory of vertical localization outlined in Section \ref{sect:vert_loc}.

\begin{proposition}\label{prop:mod_loc}
Let $\msc{M}^{\ot}\to \opn{Mod}$ be an $\msc{F}$-module category and $W^-\subseteq \msc{F}$, $W^+\subseteq \msc{M}$, be classes of maps which are stable under isomorphism and which are preserved by the monoidal product and action functors
\[
\msc{F}\times\msc{F}\to \msc{F}\ \ \text{and}\ \ \msc{F}\times\msc{M}\to \msc{M}.
\]
For each finite colored set $\mcl{I}$, let $W_{\mcl{I}}$ be the preimage of the product class $\prod_{i\in \mcl{I}}W^{\opn{color}(i)}$ under the transport equivalence
\begin{equation}\label{eq:3381}
\rho_!:\msc{M}^{\ot}_{\mcl{I}}\overset{\sim}\to \msc{M}^\mcl{I}
\end{equation}
and take $W=\cup_{\mcl{I}}W_{\mcl{I}}\subseteq \msc{M}^{\ot}[1]$.  The collection $W$ forms a vertical class over $\opn{Mod}$ and, for $W_{\Fin_{\ast}}$ the class in $\msc{F}^{\ot}$ generated by $W^{-}$, the localization $\msc{M}^{\ot}[W^{-1}]\to \opn{Mod}$ is a module category over the symmetric monoidal $\infty$-category $\msc{F}^{\ot}[W_{\Fin_{\ast}}^{-1}]\to \Fin_{\ast}$.
\end{proposition}

Before offering a proof, let us try to make sense of the statement.  We have the structural equivalence of symmetric monoidal $\infty$-categories $\msc{F}^{\ot}\overset{\sim}\to \Fin_{\ast}\times_{\opn{Mod}}\msc{M}^{\ot}$ so that, in particular, we have an equivalence $\msc{F}\overset{\sim}\to \msc{M}^{\ot}_{\{-0\}}$. Hence the colored exponents admit uniquely determined equivalences
\[
\msc{M}^{\mcl{I}}\overset{\sim}\to \msc{F}^{I_-}\times \msc{M}^{I_+},
\]
where $I_{\pm}$ are the subsets of $\pm$-colored indices in $\mcl{I}$, and trasport therefore specifies an equivalence
\[
\msc{M}^{\ot}_{\mcl{I}}\overset{\sim}\to \msc{F}^{I_-}\times\msc{M}^{I_+}.
\]
So our definitions of $W_{\mcl{I}}$ and $W$ make sense.
\par

Now, supposing that $W$ is in fact vertical over $\opn{Mod}$, in particular that the fibers $W_{\mcl{I}}$ are stable under transport, we have the pullback class $W_{\Fin_{\ast}}$ in the fiber product $\Fin_{\ast}\times_{\opn{Mod}}\msc{M}^{\ot}$ which pulls back to a class of maps in $\msc{F}^{\ot}$ along the structural equivalence
\[
\msc{F}^{\ot}\overset{\sim}\to \Fin_{\ast}\times_{\opn{Mod}}\msc{M}^{\ot}.
\]
We also call this class $W_{\Fin_{\ast}}$, by an abuse of notation.

By Proposition \ref{prop:sym_loc}, the localization $\msc{F}^{\ot}[W_{\Fin_{\ast}}^{-1}]\to \Fin_{\ast}$ is in fact symmetric monoidal and hence we obtain an equivalence of symmetric monoidal $\infty$-categories
\[
\msc{F}^{\ot}[W_{\opn{Fin}_{\ast}}^{-1}]\overset{\sim}\to (\Fin_{\ast}\times_{\opn{Mod}}\msc{M}^{\ot})[W_{\Fin_{\ast}}^{-1}]\overset{\sim}\to \Fin_{\ast}\times_{\opn{Mod}}(\msc{M}^{\ot}[W^{-1}]),
\]
where the second equivalence follows by Proposition \ref{prop:vert_pullback}.  So we observe that $\msc{M}^{\ot}[W^{-1}]$ is in fact a module category over $\msc{F}^{\ot}[W^{-1}_{\Fin_{\ast}}]$, provided $W$ is vertical over $\opn{Mod}$.

\begin{proof}[Proof of Proposition \ref{prop:mod_loc}]
One sees directly that all maps in $\opn{Mod}$ are obtained as (products of) composites of inert maps and the generating active maps $\{-1,-0\}\to \{-0\}$ and $\{-1,+0\}\to \{0\}$.  Hence stability of all $W_{\mcl{I}}$ under transport follows by stability of $W^-$ under the product on $\msc{F}$ and stability of $W^-\times W^+$ under the action functor for $\msc{M}$.  According to the requirements of Definition \ref{def:vertical}--or more directly the hypotheses of Proposition \ref{prop:fib_localize}--we have that $W$ is vertical along the fibration $\msc{M}^{\ot}\to \opn{Mod}$.  The fact that the fibers $(\msc{M}^{\ot}[W^{-1}])_{\mcl{I}}$ are identified with the colored products $(\msc{M}[W_{\{-0\}}^{-1}])^{\mcl{I}}$ via transport, and hence that $\msc{M}^{\ot}[W^{-1}]$ is realized as a $\opn{Mod}$-monoidal category, follows by the fiber calculation of Proposition \ref{prop:fib_localize} (iv).
\par

There is a structural equivalence
\[
\xi_W:\msc{F}^{\ot}[W_{\Fin_{\ast}}^{-1}]\overset{\sim}\to \Fin_{\ast}\times_{\opn{Mod}}(\msc{M}^{\ot}[W^{-1}])
\]
induced by the analogous equivalence for $\msc{M}^{\ot}$, which we obtained precisely as outlined above.  So we see that the localization $\msc{M}^{\ot}[W^{-1}]$ carries the natural structure of an $\msc{F}[W_{\Fin_{\ast}}^{-1}]$-module category in this setting.
\end{proof}

\begin{lemma}\label{lem:3421}
Let $A^{\heartsuit}$ be a finitely generated linear abelian category and $A_{\pf}$ be the category of pro-finite cochains.  The product and action functors,
\[
\opn{Ch}(\opn{Vect})_{\pf}\times \opn{Ch}(\opn{Vect})_{\pf}\to \opn{Ch}(\opn{Vect})_{\pf}
\]
and
\[
\opn{Ch}(\opn{Vect})_{\pf}\times A_{\pf}\to A_{\pf},
\]
preserve both homotopy equivalences and quasi-isomorphisms, in each factor.  Furthermore, the classes of homotopy equivalences and quasi-isomorphisms agree in $\opn{Ch}(\opn{Vect})_{\pf}$.
\end{lemma}

\begin{proof}
We can proceed abstractly, but let us make things very concrete.  We can identify $A^{\heartsuit}$ with the category $\opn{Corep}(\Lambda)$ of comodules over a coalgebra $\Lambda$ \cite[Theorem 5.1]{takeuchi77}. For $\Lambda'$ the co-opposite coalgebra, we have the equivalence $-^{\ast}:\opn{Corep}(\Lambda')^{op}\overset{\sim}\to \opn{Corep}(\Lambda)_{\pf}$ provided by linear duality, which then extends to a duality equivalence at the level of cochains
\begin{equation}\label{eq:3442}
\opn{Ch}(\Lambda')^{op}\overset{\sim}\to \opn{Ch}(\Lambda)_{\pf},
\end{equation}
just as in Sections \ref{sect:pf_obj} and \ref{sect:pf_complexes}.  The inverse is provided by the continuous dual.
\par

This equivalence intertwines the $\opn{Ch}(\opn{Vect})$ action on the left with the $\opn{Ch}(\opn{Vect})_{\pf}$ action on the right, in the sense that we have a diagram
\[
\xymatrix{
\opn{Ch}(\opn{Vect})^{op}\times\opn{Ch}(\Lambda')^{op}\ar[r]^(.45){\sim}\ar[d] & \opn{Ch}(\opn{Vect})_{\pf}\times\opn{Ch}(\Lambda)_{\pf}\ar[d]\\
\opn{Ch}(\Lambda')^{op}\ar[r]^(.45){\sim} & \opn{Ch}(\Lambda)_{\pf}.
}
\]

As with any equivalence of abelian categories, duality is additive and exact, and in particular commutes with cohomology.  Hence the equivalence \eqref{eq:3442} identifies homotopy equivalence in the source category with homotopy equivalence in the target category, and similarly preserves quasi-isomorphisms.  So, the fact that the usual action
\[
\opn{Ch}(\opn{Vect})\times \opn{Ch}(\Lambda')\to \opn{Ch}(\Lambda')
\]
preserves homotopy equivalence and quasi-isomorphisms implies that the pro-finite action also preserves homotopy equivalence and quasi-isomorphisms (see also Proposition \ref{prop:tensor_qiso}). Considering the case $\Lambda=k$, we see that the product on $\opn{Ch}(\opn{Vect})_{\pf}$ similarly preserves these classes of maps, and that homotopy equivalence agrees with quasi-isomorphisms in $\opn{Ch}(\opn{Vect})_{\pf}$ since they agree in $\opn{Ch}(\opn{Vect})$.
\end{proof}

\begin{corollary}\label{cor:htop_lin}
Let $A^{\heartsuit}$ be a finitely generated linear abelian category.  The pro-finite homotopy and derived $\infty$-categories $\msc{K}_{\pf}=\msc{K}(A^{\heartsuit}_{\pf})$ and $\msc{D}_{\pf}=\msc{D}(A^{\heartsuit}_{\pf})$ are naturally module categories over $\Vect_{\pf}$.  Furthermore, the action maps for $\msc{K}_{\pf}$ and $\msc{D}_{\pf}$ fit into commuting diagrams
\[
\xymatrixcolsep{0mm}
\xymatrix{
	& \opn{Ch}(\opn{Vect})_{\pf}\times A_{\pf}\ar[dl]\ar[dr]\ar[d]\\
\Vect_{\pf}\times \msc{K}_{\pf}\ar[d] & A_{\pf}\ar[dl]\ar[dr] & \Vect_{\pf}\times \msc{D}_{\pf}\ar[d]\\
\msc{K}_{\pf} & & \msc{D}_{\pf}
}
\]
in $\sCat_{\infty}$.
\end{corollary}

\begin{proof}
Immediate from Proposition \ref{prop:mod_loc} and Lemma \ref{lem:3421}.
\end{proof}

We let $_k\msc{K}^{\ot}_{\pf}\to \opn{Mod}$ and $_k\msc{D}^{\ot}_{\pf}\to \opn{Mod}$ denote the fibrations for the $\Vect_{\pf}$-module categories $\msc{K}_{\pf}$ and $\msc{D}_{\pf}$, respectively. By the universal property of (vertical) localization we have a diagram
\[
\xymatrix{
	& {_kA}_{\pf}^{\ot}\ar[dl]\ar[dr] \\
{_k\msc{K}}^{\ot}_{\pf}\ar[rr] & & {_k\msc{D}}^{\ot}_{\pf}
}
\]
in $\msc{M}\!od^{\opn{Mod}}_{\infty}$.  This diagram pulls back to a diagram
\[
\xymatrix{
 & \opn{Ch}(\opn{Vect})_{\pf}^{\ot}\ar[dl]\ar[dr]\\
\Fin_{\ast}\times_{\opn{Mod}}{_k\msc{K}}^{\ot}_{\pf}\ar[rr] & &\Fin_{\ast}\times_{\opn{Mod}}{_k\msc{D}}^{\ot}_{\pf}
}
\]
in $\SM_{\infty}$, which then localizes to a diagram
\[
\xymatrix{
	& \Vect_{\pf}^{\ot}\ar[dl]_{\sim}\ar[dr]^{\sim}\\
\Fin_{\ast}\times_{\opn{Mod}}{_k\msc{K}}^{\ot}_{\pf}\ar[rr] & & \Fin_{\ast}\times_{\opn{Mod}}{_k\msc{D}}^{\ot}_{\pf}.
}
\]
This is all to say, the localization functor $\msc{K}_{\pf}\to \msc{D}_{\pf}$ is $\Vect_{\pf}$-linear.

\subsection{Calculating inner-Homs under advantageous localization}

\begin{proposition}\label{prop:loc_innerhoms}
Suppose that an $\msc{F}$-module $\infty$-category $\msc{M}$ admits restrictive inner-Homs, relative to a full subcategory $\msc{M}_{fin}\subseteq \msc{M}$.  Let $W^-$ and $W^+$ be classes of maps in $\msc{F}$ and $\msc{M}$, respectively, which satisfy the hypotheses of Proposition \ref{prop:mod_loc}, and suppose that for each $x$ in $\msc{M}_{fin}$ the functor
\[
\underline{\Maps}_{\msc{M}}(x,-):\msc{M}\to \msc{F}
\]
sends maps in $W^+$ to maps in $W^-$.  Then the localized module category $\msc{M}[(W^+)^{-1}]$ admits restrictive inner-Homs over $\msc{F}[(W^-)^{-1}]$, relative to the full subcategory spanned by the image of $\msc{M}_{fin}$ under localization, and the inner-Hom functors are calculated via localization.  More precisely, for each $x$ in $\msc{M}_{fin}$ the functors $\underline{\Maps}_{\msc{M}[(W^+)^{-1}]}(loc(x),-)$ is the unique one which completes a diagram
\[
\xymatrix{
	& \msc{M}[(W^+)^{-1}]\ar[dr]^(.55){\underline{\Maps}(loc(x),-)}\\
\msc{M}\ar[ur]^{loc}\ar[rr]_{\underline{\Maps}(x,-)} & & \msc{F}.
}
\]
\end{proposition}

\begin{proof}
Take $S=W^-$ and $T=W^+$.  The action maps for $\msc{M}$ and its localization fit into a diagram
\[
\xymatrix{
\msc{F}\times\msc{M}\ar[d]\ar[rr]^(.4){loc\times loc} & & \msc{F}[S^{-1}]\times \msc{M}[T^{-1}]\ar[d]\\
\msc{M}\ar[rr]_(.4){loc} & & \msc{M}[T^{-1}]
}
\]
so that, for any $x$ in $\msc{M}$, the localization of the action map
\[
loc\circ (-\ot x)[S^{-1}]:\msc{F}[S^{-1}]\to \msc{M}[T^{-1}]
\]
is isomorphic to the action $-\ot loc(x)$.  Hence, to say that we have restrictive inner-Homs relative to the full subcategory $\msc{M}[T^{-1}]_{fin}$ spanned by the image of $\msc{M}_{fin}$ is to say that, for each $x$ in $\msc{M}_{fin}$, the localization of the functor $-\ot x:\msc{F}\to \msc{M}$ admits a right adjoint.  (See \cite[\href{https://kerodon.net/tag/045S}{045S}]{kerodon}.)
\par

For $x$ in $\msc{M}_{fin}$, the functor $-\ot x$ specifies a cocartesian fibration
\[
q:\msc{Z}\to \Delta^1
\]
with fibers $\msc{Z}_0= \msc{F}$ and $\msc{Z}_1= \msc{M}$ and transport recovering $-\ot x$.  The specification of $\underline{\Maps}_{\msc{M}}(x,-)$ as right adjoint to $-\ot x$ tells us that the fibration $q$ is also cartesian with cotransport recovering $\underline{\Maps}_{\msc{M}}(x,-)$.  (The fibration $q:\msc{Z}\to \Delta^1$ can be explicitly constructed by the weighted nerve \cite[\href{https://kerodon.net/tag/025X}{025X}, \href{https://kerodon.net/tag/02FN}{02FN}]{kerodon}.)
\par

Take $U\subseteq \msc{Z}[1]$ the union of the collections $U_0=S$ and $U_1=T$ and consider the weighted nerve $q':\msc{Z}'\to \Delta^1$ of the localized action functor.  Since $U$ is vertical for $q$, considered both as a cocartesian fibration and a cartesian fibration, we have the vertical localization $\msc{Z}[U^{-1}]\to \Delta^1$ and the diagram
\begin{equation}\label{eq:3510}
\xymatrix{
\msc{Z}\ar[rr]^{\opn{loc}}\ar[dr]^{i}\ar[ddr]\ar[rr] & & \msc{Z}'\ar[ddl]\\
	& \msc{Z}[U^{-1}]\ar[ur]^{i'}\ar[d]\\
	&\Delta^1
}
\end{equation}
in which $i'$ is an equivalence of cocartesian fibrations, by Proposition \ref{prop:fib_localize} (iv) and \cite[\href{https://kerodon.net/tag/028B}{028B}]{kerodon}.
\par

By Proposition \ref{prop:fib_localize}, the cocartesian, or dually cartesian localization promised in Proposition \ref{prop:fib_localize} is obtained by taking any localization $F:\msc{Z}\to\msc{Z}[U^{-1}]$ which is equipped with an isofibration to $\Delta^1$ which fits into a strictly commuting diagram
\[
\xymatrix{
\msc{Z}\ar[rr]^F\ar[dr]_q & & \msc{Z}[U^{-1}]\ar[dl]^{q_U}\\
	& \Delta^1 & .
}
\]
In particular, the vertical localization $\msc{Z}[U^{-1}]\to \Delta^1$ is simultaneously cartesian and cocartesian.  From the equivalence $i':\msc{Z}[U^{-1}]\to \msc{Z}'$ in \eqref{eq:3510} we therefore deduce that $q':\msc{Z}'\to \Delta^1$ is both cocartesian and cartesian as well, and that the map $\msc{Z}\to \msc{Z}'$ induced by the diagram
\[
\xymatrix{
\msc{F}\ar[rr]^{-\ot x}\ar[d]_{loc} & & \msc{M}\ar[d]^{loc}\\
\msc{F}[S^{-1}]\ar[rr]_{-\ot loc(x)} & & \msc{M}[T^{-1}]
}
\]
in $\sCat_{\infty}$ is a map of both cocartesian and cartesian fibrations.  It follows that the functor $-\ot loc(x)$ admits a right adjoint $\opn{M}^{loc(x)}$ and that the right adjoint fits into a diagram
\begin{equation}\label{eq:3535}
\xymatrix{
\msc{M}\ar[rr]^{\underline{\Maps}(x,-)}\ar[d]_{loc} & & \msc{F}\ar[d]^{loc}\\
\msc{M}[T^{-1}]\ar[rr]_{\opn{M}^{loc(x)}} & & \msc{F}[S^{-1}]
}
\end{equation}
\cite[\href{https://kerodon.net/tag/02FN}{02FN}]{kerodon}. Thus $\opn{M}^{loc(x)}$ is obtained from the functor $\underline{\Maps}(x,-)$ via localization.
\par

We now have that the functor $-\ot loc(x)$ admits a right adjoint whenever $x$ is in $\msc{M}_{fin}$, and we conclude that the functor $-\ot x':\msc{F}[S^{-1}]\to \msc{M}[T^{-1}]$ admits a right adjoint whenever $x'$ is in the full subcategory $\msc{M}[T^{-1}]_{fin}$ spanned by the image of $\msc{M}_{fin}$.  So the functor
\[
\msc{M}[T^{-1}]_{fin}^{op}\times \msc{M}[T^{-1}]\to \Fun(\msc{F}[S^{-1}]^{op},\sKan)
\]
\[
(x',x'')\mapsto \{z\mapsto H_{\msc{M}}(z\ot x',x'')\}
\]
has image in the full subcategory of representable functors, and so the $\msc{F}[S^{-1}]$-module category $\msc{M}[T^{-1}]$ admits restrictive inner-Homs.
\par

For $x'=loc(x)$ in $\msc{M}[T^{-1}]_{fin}$, the inner-Homs
\[
\underline{\Maps}_{\msc{M}[T^{-1}]}(loc(x),-):\msc{M}[T^{-1}]\to \msc{F}[S^{-1}] 
\]
are right adjoint to the action functor $-\ot loc(x)$.  We therefore find
\[
\underline{\Maps}_{\msc{M}[T^{-1}]}(loc(m),-)=\opn{M}^{loc(x}=\left\{\begin{array}{l}
\text{the functor obtained from}\\
\underline{\Maps}_{\msc{M}}(x,-)\ \text{via localization},
\end{array}\right.
\]
by a consideration of the diagram \eqref{eq:3535}.
\end{proof}

\begin{corollary}\label{cor:K_innerhoms}
For any finitely generated linear abelian category $A^{\heartsuit}$, the category $\msc{K}_{\pf}=\msc{K}(A^{\heartsuit}_{\pf})$ admits restrictive inner-Homs, relative to the subcategory $\msc{K}_{fin}$ of finite length bounded complexes.  Furthermore, for each finite complex $x$, the inner-Homs
\[
\underline{\Maps}_{\msc{K}_{\pf}}(x,-):\msc{K}\to \Vect_{\pf}
\]
are obtained from Hom-complexes $\Hom^{\ast}_{A}(x,-)$ via localization.
\end{corollary}

\subsection{Inner-Homs in the derived setting}

Our aim in this subsection is to prove the following.

\begin{proposition}\label{prop:D_innerhom}
For a finitely generated linear abelian category $A^{\heartsuit}$ which has enough compact injectives, the $\Vect_{\pf}$-module $\infty$-category $\msc{D}_{\pf}=\msc{D}(A^{\heartsuit}_{\pf})$ admits restrictive inner-Homs, relative to the subcategory $\msc{D}_{fin}$ of complexes with finite length, bounded cohomology.
\end{proposition}

Recall that a complex $x$ in $A_{\pf}=\opn{Ch}(A^{\heartsuit}_{\pf})$ is called $K$-injective if the Hom-complex $\Hom^{\ast}_{A_{\pf}}(-,x)$ preserves quasi-isomorphisms.  One similarly defines $K$-projective complexes.  Note that $K$-injectives and $K$-projectives are stable under homotopy equivalence.

\begin{lemma}\label{lem:inj_linear}
The full subcategory $\opn{Inj}\msc{K}_{\pf}$ of $K$-injectives in $\msc{K}_{\pf}$ is stable under the action of $\Vect_{\pf}$, and hence forms a $\Vect_{\pf}$-module subcategory in $\msc{K}_{\pf}$.
\end{lemma}

\begin{proof}
Via duality and localization, it suffices to show that the subcategory $\opn{Proj}B$ of $K$-projectives in $B=\opn{Ch}(B^{\heartsuit})$ is stable under the action of $\opn{Ch}(\opn{Vect})$, for any finitely generated abelian category $B$ with enough compact projectives. However, this is clear as any linear cochain complex $v$ is homotopy equivalent to its cohomology $H^\ast(v)$.  Hence the product with a $K$-projective $v\ot P$ is homotopy equivalent to a possibly-infinite sum of shifts of $P$. Since the class of $K$-projectives is stable under arbitrary sums and shifting, we have that $v\ot P$ is in fact $K$-projective.
\par

We consider the fibration ${_k\msc{K}}^{\ot}_{\pf}\to \opn{Mod}$ and the full subcategory ${_k\opn{Inj}\msc{K}}_{\pf}^{\ot}$ in ${_k\msc{K}}^{\ot}_{\pf}$ spanned by those objects in the fibers $({_k\msc{K}}^{\ot}_{\pf})_{\mcl{I}}$ which map to objects in $\Vect^{I_-}\times \opn{Inj}\msc{K}^{I_+}$ under the transport equivalence
\[
({_k\msc{K}}^{\ot}_{\pf})_{\mcl{I}}\overset{\sim}\to ({_k\msc{K}}_{\pf})^{\mcl{I}}=\Vect^{I_-}_{\pf}\times \msc{K}^{I_+}_{\pf}.
\]
Stability of $\opn{Inj}\msc{K}_{\pf}$ under the action of $\Vect_{\pf}$ tells us that the fibers $({_k\opn{Inj}\msc{K}}^{\ot}_{\pf})_{\mcl{I}}$ are stable under transport, and hence that the sequence ${_k\opn{Inj}\msc{K}}^{\ot}_{\pf}\to {_k\msc{K}}^{\ot}_{\pf}\to \opn{Mod}$ is a cocartesian fibration which gives ${_k\opn{Inj}\msc{K}}^{\ot}_{\pf}$ the structure of a $\opn{Mod}$-monoidal $\infty$-category.  Furthermore we have the equivalence
\[
\Vect^{\ot}_{\pf}\overset{\sim}\to \Fin_{\ast}\times_{\opn{Mod}}({_k\opn{Inj}\msc{K}}^{\ot}_{\pf})=\Fin_{\ast}\times_{\opn{Mod}}({_k\msc{K}}^{\ot}_{\pf}),
\]
so that the inclusion ${_k\opn{Inj}\msc{K}}^{\ot}_{\pf}\to {_k\msc{K}}^{\ot}_{\pf}$ is an inclusion of $\Vect_{\pf}$-module $\infty$-categories.
\end{proof}

\begin{lemma}\label{lem:enough_inj}
If $A^{\heartsuit}$ has enough compact injectives, then the category $A_{\pf}=\opn{Ch}(A^{\heartsuit})_{\pf}$ has enough $K$-injectives.
\end{lemma}

Of course, by enough $K$-injectives we mean that every object $x$ admits a quasi-isomorphism $x\to I(x)$ into a $K$-injective complex.

\begin{proof}
As in the proof of Lemma \ref{lem:3421}, we can find a finitely generated linear abelian category $B^{\heartsuit}$ with $(B^{\heartsuit})^{op}\cong A^{\heartsuit}_{\pf}$ and $B^{op}\cong A_{\pf}$, for $B=\opn{Ch}(B^{\heartsuit})$.  By our hypotheses $B^{\heartsuit}$ has enough compact projectives, and therefore enough projectives.  It follows that $B$ has enough $K$-projectives (cf.\ \cite[Lemma 13.3]{drinfeld04}), and by duality $A_{\pf}$ has enough $K$-injectives.
\end{proof}

Existence of enough $K$-injectives tells us that the full subcategory $\opn{Inj}\msc{K}_{\pf}$ spanned by $K$-injectives in $\msc{K}_{\pf}$ is a reflective \cite[\href{https://kerodon.net/tag/02F5}{02F5}]{kerodon}.  We therefore apply standard results for reflective subcategories--as in Section \ref{sect:kan_htop_der}--to see that the derived $\infty$-category $\msc{D}_{\pf}$ is identified with the subcategory of $K$-injective complexes in $\msc{K}_{\pf}$, and that the localization functor $\msc{K}_{\pf}\to \msc{D}_{\pf}$ is defined via an application of functorial $K$-injective resolutions.

\begin{corollary}\label{cor:dercat_loc}
Suppose $A^{\heartsuit}$ has enough compact injectives.  The inclusion $\opn{Inj}\msc{K}_{\pf}\to \msc{K}_{\pf}$ admits a left adjoint $L:\msc{K}_{\pf}\to \opn{Inj}\msc{K}_{\pf}$ which realizes $\opn{Inj}\msc{K}_{\pf}$ as a localization of $\msc{K}_{\pf}$ along the class of quasi-isomorphisms.  In particular, the composite functor
\[
\xymatrixrowsep{3mm}
\xymatrix{
	&\msc{K}_{\pf}\ar@{-->}[dr]\\
\opn{Inj}\msc{K}_{\pf}\ar[rr]\ar@{-->}[ur] & & \msc{D}_{\pf}
}
\]
is an equivalence, and the localization functor $\msc{K}_{\pf}\to\msc{D}_{\pf}$ admits a fully faithful right adjoint $R:\msc{D}_{\pf}\to \msc{K}_{\pf}$ which is an equivalence onto the full subcategory of $K$-injectives.
\end{corollary}

\begin{proof}
This is a special case of \cite[\href{https://kerodon.net/tag/02FD}{02FD}, \href{https://kerodon.net/tag/02FE}{02FE}, \href{https://kerodon.net/tag/04JL}{04JL}]{kerodon}.
\end{proof}

\begin{remark}
In general, the category $A_{\pf}$ will have enough $K$-\text{projectives}, and one can construct the derived $\infty$-category $\msc{D}_{\pf}$ via $K$-projective resolutions.  For us the injective setting is more relevant.
\end{remark}

We can now prove our main result.

\begin{proof}[Proof of Proposition \ref{prop:D_innerhom}]
We have the linear functor ${_k\msc{K}}^{\ot}_{\pf}\to {_k\msc{D}}^{\ot}_{\pf}$ which restricts to a linear functor ${_k\opn{Inj}\msc{K}}^{\ot}_{\pf}\to {_k\msc{D}}^{\ot}_{\pf}$ whose fibers are equivalence, and hence which is an equivalence of $\Vect_{\pf}$-module categories \cite[\href{https://kerodon.net/tag/028C}{028C}]{kerodon}.  So it suffices to show that $\opn{Inj}\msc{K}_{\pf}$ admits restrictive inner-Homs relative to the subcategory of complexes with finite cohomology.
\par

Any complex $x$ with finite cohomology admits a pair of quasi-isomorphisms $x\to x'\leftarrow x''$ in which $x''$ is a bounded complex of finite length objects.  (To see this one can check the analogous claim in the dual category $\opn{Ch}(B^{\heartsuit})$ for $B^{\heartsuit}=(A^{\heartsuit}_{\pf})^{op}$, where the result is standard.)  Since Hom-complexes with $K$-injective target preserve quasi-isomorphisms, and since the action of $\Vect_{fin}$ preserves quasi-isomorphisms by Lemma \ref{lem:3421}, we have a natural isomorphisms of functors
\[
H_{\msc{K}_{\pf}}(-\ot x,y)\cong H_{\msc{K}_{\pf}}(-\ot x'',y):\Vect_{\pf}\to \sKan
\]
whenever $y$ is in $\opn{Inj}\msc{K}_{\pf}$.  (Here $H_{\msc{K}_{\pf}}$ is a Hom functor for the homotopy $\infty$-category.)  By Corollary \ref{cor:K_innerhoms} the functor $H_{\msc{K}_{\pf}}(-\ot x'',y)$ is representable, and hence $H_{\msc{K}_{\pf}}(-\ot x,y)$ is representable.  Thus $\opn{Inj}\msc{K}$ admits the proposed inner-Homs.
\end{proof}

From the proof we also deduce a formula for the inner-Homs.

\begin{proposition}\label{prop:D_innerhom_formula}
Suppose $A^{\heartsuit}$ admits enough compact injectives, and let $R:\msc{D}_{\pf}\to \msc{K}_{\pf}$ be right adjoint to the localization functor $loc:\msc{K}_{\pf}\to \msc{D}_{\pf}$.  The $\Vect_{\pf}$-module category $\msc{D}_{\pf}$ admits restrictive inner-Homs, relative to the full subcategory $\msc{D}_{fin}$ spanned by the image of $\msc{K}_{fin}$ under localization, and for any $x$ in $\msc{K}_{fin}$ there is a natural isomorphism
\[
\underline{\Maps}_{\msc{D}_{\pf}}(loc(x),-)\cong \underline{\Maps}_{\msc{K}_{\pf}}(x,R-)
\]
in $\Fun(\Vect_{\pf},\sKan)$.
\end{proposition}

\bibliographystyle{abbrv}

\end{document}